\documentclass[12pt]{article}
\usepackage{mathrsfs}
\usepackage{amsmath}
\usepackage{amssymb}
\usepackage{amsfonts}
\usepackage{mathtools}
\usepackage[top=2cm,bottom=2cm,left=2.5cm,right=2cm]{geometry}
\usepackage{graphicx}
\usepackage[colorinlistoftodos]{todonotes}
\presetkeys{todonotes}{inline,color=blue!15,bordercolor=blue}{}
\usepackage{amsthm}
\usepackage{bbm}
\usepackage{enumerate}
\usepackage{enumitem}
\usepackage{verbatim}
\usepackage[
  backend=biber,
  style=alphabetic,
  sorting=nyt,
  maxnames=99,
  giveninits=true,
  doi=true,
  url=false,
  isbn=false,
  eprint=true
]{biblatex}
\usepackage{float}
\usepackage{tikz-cd}
\usepackage[all]{xy}
\usepackage[T1]{fontenc}
\usepackage{xcolor}
\usepackage{tcolorbox}
\usepackage{adjustbox}
\usepackage[bottom]{footmisc}

\newtcolorbox{commentbox}[1][]{colback=yellow!10!white, colframe=red!50!black, title=#1}

\usepackage[affil-it]{authblk}
\usepackage{commands}
\usepackage[colorlinks=true, allcolors=blue]{hyperref}

\title{The Unipotent Chabauty--Kim--Kantor Method for Relative Completions}
\author[1]{David Corwin\thanks{corwind@bgu.ac.il}}
\author[2]{Sa'ar Zehavi\thanks{saarzehavi@gmail.com}}
\affil[1]{Ben-Gurion University of the Negev}
\affil[2]{Ben-Gurion University of the Negev}

\begin{document}

\maketitle

\begin{abstract}
We develop a new $p$-adic analytic method for studying integral points on hyperbolic curves,
building on Kantor's relative-completion approach to unifying the
Chabauty--Kim and Lawrence--Venkatesh methods. 
Its central mechanism converts dimension inequalities between algebraic global
and local Galois cohomology stacks into nonzero rigid-analytic functions
vanishing on the integral points within a residue disk. Assuming the relevant Bloch--Kato
dimension formulas, we verify these inequalities for curves of genus at
least $2$, as well as for certain modular curves, 
and thereby obtain a new conditional proof of the theorems of Faltings and Siegel. Our method produces genuinely new functions, different from those produced by Chabauty--Kim, pointing to an improved effective approach.

Our key technical contributions include the following. We prove that
Kantor's $p$-adic period map is analytic and algebraically Zariski
dense in the relevant period domain. We also bypass the conjectural
representability of Kantor's Selmer stacks by decomposing the
integral-point locus into finitely many Kummer strata. On each stratum,
the relative Kummer map factors through a neutral-fibre Selmer stack,
which we prove is an algebraic Artin stack of finite type.
We also develop a theory of finite-type motivic quotients of the relative completion: these
quotients are cofinal among all finite-type quotients and recover the
relative completion as their inverse limit. For Kodaira--Parshin
families over curves of genus at least $2$, we construct a canonical
adjoint-supported tower with fixed finite-dimensional abelianization
and controlled graded pieces. At sufficiently deep levels, the
Bloch--Kato formulas imply the required dimension inequality and
produce genuinely relative rigid-analytic functions that do not arise
from the ordinary unipotent completion.
\end{abstract}

\tableofcontents

﻿\section{Introduction}
This paper develops a relative Chabauty--Kim theory for studying integral points on hyperbolic
curves, building on the work of Kantor~\cite{KANTORTHESIS}. 
Here, \emph{relative} refers to the functorial completion of a group
relative to a monodromy representation with reductive image. Replacing the
unipotent completion of the classical Chabauty--Kim method by this relative
completion yields a more flexible framework for studying the arithmetic
of hyperbolic curves, providing a bridge between the Chabauty--Kim and Lawrence--Venkatesh approaches; see~\cite{KIM,LV}.

We begin with our motivation in developing the theory; our main results are summarized in~\S\ref{subsection:intro_main_results}, 
see especially Theorems~\ref{theorem:main_intro},
\ref{theorem:intro_finite_type_quotients},
\ref{theorem:intro_adjoint_supported_tower},
and~\ref{theorem:main_main}.

\subsection{Motivation: Effective Faltings}

Let \(X/K\) be a smooth quasi-projective curve over a number field
\(K\). Write \(\mathcal O_K\) for the ring of integers of \(K\), and
let \(S\) be a finite set of finite places such that \(X\) admits a
smooth model
\[
\mathcal X
\longrightarrow
\Spec\mathcal O_{K,S},
\qquad
\mathcal O_{K,S}\coloneqq \mathcal O_K[1/S].
\]
The theorems of Siegel and Faltings imply that
\[
\mathcal X(\mathcal O_{K,S})
\]
is finite whenever \(X(\mathbb C)\) has negative Euler
characteristic. The general problem of determining this finite set
effectively remains open; we refer to it, broadly, as the
\emph{Effective Faltings problem}.

Kim's non-abelian Chabauty method
\cite{KIM2,KIM}
approaches this problem by replacing the Jacobian of
Chabauty--Coleman with successively deeper finite-dimensional
quotients \(U_n\) of the unipotent fundamental group. At level \(n\),
one compares a global Selmer variety with its local de Rham
realization. If the global Selmer image has smaller dimension than
the local target, an equation for that image pulls back along the Zariski dense
\(p\)-adic unipotent Albanese map to a nonzero locally analytic
function on a residue disk that vanishes on the \(S\)-integral
points. Such Chabauty--Kim functions have been computed in a
substantial range of cases; see, for example,
\cite{BALAKRISHNAN_DOGRA,BALAKRISHNAN_QUADRATIC}
for quadratic Chabauty and
\cite{CKone,MTMUE,CDC2,ZEHAVI}
for motivic Chabauty--Kim. In general, however, the existence of a
level satisfying the required dimension inequality is conditional 
on conjectural control of global Selmer groups; in the applications of this paper, 
the input is Conjecture~\ref{conj:BCL_BK}/Bloch--Kato.

Lawrence and Venkatesh
\cite{LV}
prove Faltings' theorem unconditionally by a different \(p\)-adic
dimension argument. Given a Kodaira--Parshin family
\[
f:Y\longrightarrow X
\]
and a place \(v\mid p\) of good reduction, fix
\(b\in\mathcal X(\mathcal O_{K_v})\). Write
\(\underline b\in\mathcal X(k_v)\) for its reduction and
\(]\underline b[\) for the corresponding rigid-analytic residue tube,
so that \(]\underline b[(K_v)\) is the \(v\)-adic residue disk of
\(b\). The Gauss--Manin connection identifies the de Rham cohomologies
of the fibres over this tube with
\[
H^1_{\dR}(Y_b/K_v).
\]
The resulting rigid-analytic period map records the varying Hodge
filtration:
\[
\Phi_v:
]\underline b[
\longrightarrow
\operatorname{Flag}
\left(
H^1_{\dR}(Y_b/K_v)
\right).
\]
For an \(S\)-integral point \(x\), the representation
\[
H^1_{\et}
\left(
Y_{\overline x},\mathbb Q_p
\right)
\]
arises by localization from a global \(p\)-adic Galois
representation of fixed dimension and weight, unramified outside a
fixed set of places and crystalline at \(v\). Faltings' finiteness
theorem leaves only finitely many possible semisimple global
representations of this kind.

On a fixed residue disk, the underlying $\varphi$-module
\((V,\varphi)\) is constant and only its Hodge filtration varies.
Two filtrations define isomorphic filtered $\varphi$-modules exactly
when they lie in the same orbit of the Frobenius centralizer
\[
Z(\varphi)
\coloneqq 
\operatorname{Aut}(V,\varphi)
\]
on the relevant flag variety. The period image of the integral
points is therefore contained in a finite union of
\(Z(\varphi)\)-orbits. The Lawrence--Venkatesh argument combines the
Zariski density of the period map with the positive codimension of
these orbits.

The difficulty is that the Frobenius centralizer can be too large
relative to the original period domain. Lawrence and Venkatesh
overcome this by an auxiliary extension construction: they arrange
that the cohomology of the relevant fibre decomposes over a
sufficiently large unramified extension and exploit the resulting
semilinearity of Frobenius. The associated Weil-restricted period
domain grows with the extension degree, while the centralizer remains
comparatively controlled. This produces the required codimension,
but the same auxiliary extension substantially enlarges the
arithmetic data that an effective implementation must handle; see
\cite[\S\S2--3]{BALAKRISHNAN_SURVEY}.

Kantor's thesis
\cite{KANTORTHESIS}
proposes a common generalization of the Chabauty--Kim and
Lawrence--Venkatesh methods. Let \(\mathcal R\) be the reductive
monodromy group of the Kodaira--Parshin Gauss--Manin system. Its relative completion is a \emph{pro-unipotent extension} of $\Rc$, meaning an extension
\[
1
\longrightarrow
\mathcal U
\longrightarrow
\mathcal G
\longrightarrow
\mathcal R
\longrightarrow
1,
\]
where \(\mathcal U\) is pro-unipotent. Taking
\(\mathcal R=1\) recovers the unipotent completion underlying
Chabauty--Kim, whereas passing to the reductive quotient recovers the
period geometry of Lawrence--Venkatesh. In some sense,
Kantor's
method is to Lawrence--Venkatesh what Kim's method is to
Chabauty--Coleman.

This analogy suggests a different response to the Frobenius
centralizer problem. Lawrence--Venkatesh enlarge an auxiliary
extension so that the reductive period domain grows relative to its
centralizer orbits. Kantor instead climbs the relatively unipotent
tower
\[
\cdots
\longrightarrow
\mathcal G_n
\longrightarrow
\mathcal G_{n-1}
\longrightarrow
\cdots
\longrightarrow
\mathcal R.
\]
The de Rham targets acquire new relatively unipotent directions,
while the Frobenius-fixed contribution has constant dimension:
\[
\dim\mathcal G_n^{\varphi=1}
=
\dim\mathcal R^{\varphi=1};
\]
see
\cite[Proposition~6.4.7 and Remark~27]{KANTORTHESIS}.
Thus the decisive dimension inequality may be sought by enlarging the
relative completion while keeping both the base field and the
Frobenius-centralizer penalty fixed. Kantor's framework also suggests
a possibility of organizing the global representation theory by residual
pseudorepresentations, potentially replacing the enumeration of
characteristic-zero representations by a coarser residual
stratification. This residual stratification is not used in the present paper; it
belongs to the broader program discussed in
\S\ref{subsection:kantor_program}.

The aim of this paper is to turn
this effectivity-motivated relative unipotent strategy into a cohesive
Chabauty--Kim theory.

\begin{remark}[Relative Chabauty--Skolem]
The first relative unipotent level is already of independent
interest. The quotient
\[
\mathcal G/\overline{[\mathcal U,\mathcal U]}
\]
is an extension of \(\mathcal R\) by
\(\mathcal U^{\mathrm{ab}}\). We call the specialization of the
relative method to this quotient \emph{relative
Chabauty--Skolem}. It combines the reductive variation exploited by
Lawrence--Venkatesh with the abelian level of Chabauty--Kim, without
passing beyond the abelianization of the unipotent radical. The
modular Eisenstein application of~\S\ref{section:eis} occurs
precisely at this level, since its finite-type Eisenstein quotients
factor through \(\mathcal G/[\mathcal U,\mathcal U]\). One might hope that this method can be made geometric in a similar way to Chabauty--Coleman, and we hope to explore this in future work.
\end{remark}

\subsection{Main Results}
\label{subsection:intro_main_results}

\subsubsection{The Relative Cutter}
\label{subsubsection:intro_relative_box_cutter}

Retain the preceding notation, and fix a Kodaira--Parshin family
\[
f:Y\longrightarrow X;
\]
thus \(f\) is smooth and proper of positive relative dimension and
has symplectic monodromy; see
\S\ref{section:rmc_realizations}. After enlarging \(S\) if necessary,
assume that the fixed Kodaira--Parshin data extend over
\(\mathcal O_{K,S}\) with the smooth proper or logarithmic models, as
appropriate, fixed in
Convention~\ref{convention:fixed_kp_family}. Choose a rational prime
\(p\) such that no place above \(p\) lies in \(S\), a place \(v\mid p\),
and a basepoint
\[
b\in\mathcal X(\mathcal O_{K,S}).
\]
Set
\[
S_p\coloneqq \{\,w:w\mid p\,\},
\qquad
T\coloneqq S\cup S_p,
\]
and let
\[
K_0\coloneqq \operatorname{Frac}W(k_v),
\]
where \(k_v\) is the residue field of \(K_v\). Write \(G_T\) for the Galois group of the maximal extension of \(K\)
unramified outside \(T\), and \(G_v\subseteq G_T\) for the chosen
decomposition group at \(v\).

The Gauss--Manin system of \(f\) defines relative completions in the
étale, crystalline, and de Rham realizations:
\[
1
\longrightarrow
\mathcal U^\bullet
\longrightarrow
\mathcal G^\bullet
\longrightarrow
\mathcal R^\bullet
\longrightarrow
1,
\qquad
\bullet\in\{\et,\cri,\dR\}.
\]
Here
\[
\mathcal G^{\et}/\mathbb Q_p,
\qquad
\mathcal G^{\cri}/K_0,
\qquad
\mathcal G^{\dR}/K_v
\]
are affine pro-algebraic group schemes, with pro-unipotent radicals
\(\mathcal U^\bullet\) and reductive quotients
\(\mathcal R^\bullet\). They carry the structures native to their
realizations: a Galois action in the étale realization, a
\(\sigma\)-semilinear Frobenius in the crystalline realization, and
a Hodge filtration in the de Rham realization. The non-abelian
comparison isomorphisms relate these structures and the corresponding
path torsors; see
\S\ref{section:rmc_realizations} and
Appendix~\ref{appendix:rmc}.

A \emph{finite-type motivic quotient} of the relative completion is
a compatible system
\begin{equation}
\label{eq:intro_finite_level_quotient}
1
\longrightarrow
\mathcal U_n^\bullet
\longrightarrow
\mathcal G_n^\bullet
\longrightarrow
\mathcal R^\bullet
\longrightarrow
1,
\qquad
\bullet\in\{\et,\cri,\dR\},
\end{equation}
in which the étale quotient is Galois equivariant and the
crystalline and de Rham quotients correspond to it under the
comparison isomorphisms, with their Frobenius and Hodge
structures. Fix such a quotient for the rest of \S\ref{subsubsection:intro_relative_box_cutter}.

The crystalline--de Rham comparison gives a canonical isomorphism
\[
\iota_n:
\mathcal G_n^{\cri}\otimes_{K_0}K_v
\xrightarrow{\sim}
\mathcal G_n^{\dR}.
\]
The Frobenius-fixed functor of \(\mathcal G_n^{\cri}\) is represented
by the affine \(\mathbb Q_p\)-group
\[
\mathcal G_n^{\varphi=1}
\coloneqq 
\left(
\operatorname{Res}_{K_0/\mathbb Q_p}
\mathcal G_n^{\cri}
\right)^{\varphi=1}.
\]
After extension to \(K_v\), evaluation followed by \(\iota_n\)
defines a canonical closed immersion
\[
\mathcal G_{n,K_v}^{\varphi=1}
\coloneqq 
\mathcal G_n^{\varphi=1}
\otimes_{\mathbb Q_p}K_v
\hookrightarrow
\mathcal G_n^{\dR}.
\]
Together with the Hodge subgroup
\[
F^0\mathcal G_n^{\dR}
\subseteq
\mathcal G_n^{\dR},
\]
this defines the Hodge filtration space and the de Rham period stack
\[
\operatorname{Fil}_n^{\dR}
\coloneqq 
\mathcal G_n^{\dR}/F^0\mathcal G_n^{\dR},
\qquad
\mathscr M_n^{\dR}
\coloneqq 
\left[
\mathcal G_{n,K_v}^{\varphi=1}
\backslash
\operatorname{Fil}_n^{\dR}
\right].
\]
The stack \(\mathscr M_n^{\dR}\) classifies
\(\mathcal G_n^{\dR}\)-torsors equipped with Hodge and Frobenius
reductions; see
\S\ref{section:dr_moduli}.

Fix
\[
S_p\subseteq T'\subseteq T.
\]
We construct two moduli stacks over \(K_v\),
\[
\mathscr S_{n,T'},
\qquad
\mathscr S_{n,v},
\]
which we call the global and local \emph{neutral Selmer stacks}, respectively. 
The global stack classifies \(G_T\)-equivariant \(\mathcal G_n^{\et}\)-torsors that
are Bloch--Kato finite at every place above \(p\), unramified at the
places in \(T\setminus T'\), and whose reductive pushout (a.k.a. \emph{reductive (Galois) type}; cf. Definition ~\ref{definition:reductive_galois_type}) is neutral. 

The local stack is defined analogously for \(G_v\), with the
Bloch--Kato condition at \(v\). By
Theorem~\ref{theorem:neutral_selmer_algebraicity} and Corollary~\ref{corollary:neutral_bk_algebraicity}, both are
finite-type Artin stacks with affine inertia. Their trivialized
neutral fibres
\[
\widetilde{\mathscr S}_{n,T'}
\longrightarrow
\mathscr S_{n,T'},
\qquad
\widetilde{\mathscr S}_{n,v}
\longrightarrow
\mathscr S_{n,v},
\]
which retain an equivariant trivialization of the reductive pushout,
are algebraic spaces of finite type.

For a global or local integral point \(x\), let
\[
\left[{}^{\vphantom{\et}}_{x}P_{b,n}^{\et}\right]
\]
denote the \(K_v\)-point of the corresponding non-abelian cohomology
stack defined by its relative étale path torsor. On the residue tube
of \(b\), set
\[
\Sigma_b^{\mathrm{glob}}
\coloneqq 
\left\{
x\in\mathcal X(\mathcal O_{K,S})
\ \middle|\
\begin{array}{l}
\operatorname{loc}_v(x)\in]\underline b[(K_v),\\
\left[{}^{\vphantom{\et}}_{x}P_{b,n}^{\et}\right]
\in\mathscr S_{n,T'}(K_v)
\end{array}
\right\}
\]
and
\[
\Sigma_b^{\mathrm{loc}}
\coloneqq 
\left\{
x\in]\underline b[(K_v)
\ \middle|\
\left[{}^{\vphantom{\et}}_{x}P_{b,n}^{\et}\right]
\in\mathscr S_{n,v}(K_v)
\right\}.
\]
Note that, in particular, elements of $\Sigma_b^{\mathrm{glob}}$ (resp. $\Sigma_b^{\mathrm{loc}}$) have the \emph{same reductive global (resp. local) Galois type} as $b$.

Localization sends
\(\Sigma_b^{\mathrm{glob}}\) to
\(\Sigma_b^{\mathrm{loc}}\), and the path torsors define Kummer maps
\[
\kappa_{b,n,T'}^{\mathrm{neu}}:
\Sigma_b^{\mathrm{glob}}
\longrightarrow
\mathscr S_{n,T'}(K_v),
\qquad
\kappa_{b,n,v}^{\mathrm{neu}}:
\Sigma_b^{\mathrm{loc}}
\longrightarrow
\mathscr S_{n,v}(K_v).
\]

Applying our relative \(p\)-adic Hodge realization for families of
torsors gives algebraic neutral Bloch--Kato logarithms
(Proposition~\ref{proposition:local_neutral_bk_logarithm} and
Definition~\ref{definition:neutral_bk_logarithm}):
\[
\BK_{n,v}:
\mathscr S_{n,v}
\longrightarrow
\mathscr M_n^{\dR},
\qquad
\BK_{n,T'}
\coloneqq 
\BK_{n,v}\circ\operatorname{loc}_v:
\mathscr S_{n,T'}
\longrightarrow
\mathscr M_n^{\dR}.
\]

The universal relative de Rham path torsor defines a rigid-analytic
period map
\[
[\Phi_{n,v}]:
]\underline b[
\longrightarrow
\left(
\mathscr M_n^{\dR}
\right)^{\an},
\]
together with its canonical lift
\[
\Phi_{n,v}:
]\underline b[
\longrightarrow
\left(
\operatorname{Fil}_n^{\dR}
\right)^{\an};
\]
see Definition~\ref{definition:v_adic_period_map} and
Proposition~\ref{proposition:an_rep}.

Our principal geometric result is the following:

\begin{theorem}
[The relative cutter;
Theorem~\ref{theorem:bk_period_compatibility},
Theorem~\ref{theorem:dim_ineq},
and Corollary~\ref{corollary:dim_ineq}]
\label{theorem:main_intro}
With the preceding notation, the following statements hold.

\begin{enumerate}
\item
The Kummer maps, localization, the local Bloch--Kato logarithm, and
the \(v\)-adic period map form a \(2\)-commutative diagram
\[
\begin{tikzcd}[column sep=large,row sep=large]
\Sigma_b^{\mathrm{glob}}
  \arrow[r, "\operatorname{loc}_v"]
  \arrow[d, "\kappa_{b,n,T'}^{\mathrm{neu}}"']
&
\Sigma_b^{\mathrm{loc}}
  \arrow[r, hook]
  \arrow[d, "\kappa_{b,n,v}^{\mathrm{neu}}"']
&
]\underline b[(K_v)
  \arrow[d, "{[\Phi_{n,v}]}"]
\\
\mathscr S_{n,T'}(K_v)
  \arrow[r, "\operatorname{loc}_v"']
&
\mathscr S_{n,v}(K_v)
  \arrow[r, "\BK_{n,v}"']
&
\mathscr M_n^{\dR}(K_v).
\end{tikzcd}
\]
In particular, the bottom composite is the global logarithm
\(\BK_{n,T'}\).

\item
Suppose that
\begin{equation}
\label{eq:intro_dimension_inequality}
\dim\widetilde{\mathscr S}_{n,T'}
+
\dim\mathcal G_{n,K_v}^{\varphi=1}
+
\dim F^0\mathcal G_n^{\dR}
<
\dim\mathcal G_n^{\dR}.
\end{equation}
Then
\[
\#\Sigma_b^{\mathrm{glob}}<\infty.
\]

More precisely, there exist a closed affinoid disk
\[
D_b^+
\subseteq
]\underline b[
\]
containing \(]\underline b[(K_v)\) and a nonzero analytic function
\[
F_{n,b}\in\mathcal O(D_b^+)
\]
such that
\[
F_{n,b}
\left(
\operatorname{loc}_v(x)
\right)
=
0
\qquad
\text{for all }
x\in\Sigma_b^{\mathrm{glob}}.
\]
\end{enumerate}
\end{theorem}

This is the relative analogue of Kim's
cutter~\cite[p.5 and Conjecture~1]{KIM}, so named both because it cuts
out rational points and because the diagram resembles a box cutter. Pulling back the Hodge filtration space along the logarithm gives the
Selmer pullback
\[
\widetilde{\mathscr C}_n
\coloneqq 
\widetilde{\mathscr S}_{n,T'}
\times_{\mathscr M_n^{\dR}}
\operatorname{Fil}_n^{\dR}.
\]
It is a
\(\mathcal G_{n,K_v}^{\varphi=1}\)-torsor over
\(\widetilde{\mathscr S}_{n,T'}\), and hence
\[
\dim\widetilde{\mathscr C}_n
=
\dim\widetilde{\mathscr S}_{n,T'}
+
\dim\mathcal G_{n,K_v}^{\varphi=1}.
\]
The commutativity of the box cutter diagram places the lifted period image of
\(\Sigma_{b}^{\mathrm{glob}}\) in the scheme-theoretic image of
\(\widetilde{\mathscr C}_n\) inside the connected component of
\(\operatorname{Fil}_n^{\dR}\) containing the identity coset.
Under
\eqref{eq:intro_dimension_inequality}, this image is a proper closed
subscheme. The lifted period map is algebraically Zariski dense in
that component by
Proposition~\ref{proposition:p-adic-dense}; pulling back the defining
ideal therefore produces \(F_{n,b}\).

\subsubsection{Finite-Type Motivic Quotients}
\label{subsubsection:intro_finite_type_quotients}
The inequality (\ref{eq:intro_dimension_inequality}) requires finite-dimension, so its application requires finite-type motivic quotients. For a pro-unipotent group
\(\mathcal U\), write
\[
\mathcal F_{\DCS}^1\mathcal U\coloneqq \mathcal U,
\qquad
\mathcal F_{\DCS}^{n+1}\mathcal U
\coloneqq 
\overline{[\mathcal U,\mathcal F_{\DCS}^n\mathcal U]},
\qquad
\mathcal{U}_n \coloneqq  \mathcal{U}/F_{\DCS}^{n+1}\mathcal U,
\qquad
\Gr^{\DCS}_n \mathcal{U} \coloneqq  F_{\DCS}^{n}\mathcal{U}/F_{\DCS}^{n+1}\mathcal U
\]
for its descending central series and associated (sub)quotients. For the ordinary unipotent
completion, this series provides a canonical tower of finite-type
quotients. For a relative completion, however, the abelianization of
its pro-unipotent radical, and hence its descending-central-series (DCS)
graded pieces, may already be infinite-dimensional. Consequently,
the quotients
\[
\mathcal G^\bullet/
\mathcal F_{\DCS}^{n+1}\mathcal U^\bullet
\]
need not be of finite type.

We prove the cofinality of finite-type motivic
quotients and construct specific families needed in our
arithmetic applications. The basic finiteness input is
Theorem~\ref{thm:ftg_basic}: under the finite-type reductive and
finite-dimensional \(\operatorname{Ext}^1\) hypotheses satisfied here,
bounded iterated self-extensions of a fixed semisimple object generate
an algebraic Tannakian category. Applied
to the relative completion, this yields the following approximation
theorem.

\begin{theorem}
[Finite-type motivic approximation;
Theorem~\ref{thm:inv_limit_of_ft_revised}]
\label{theorem:intro_finite_type_quotients}
Finite-type motivic quotients form a cofiltered inverse system, cofinal 
among the finite-type quotients of the relative completion. In each
realization, the canonical morphism
\[
\mathcal G^\bullet
\longrightarrow
\varprojlim_\alpha
\mathcal G_\alpha^\bullet
\]
is an isomorphism. Thus
\[
\mathcal G^{\et}
\simeq
\varprojlim_\alpha\mathcal G_\alpha^{\et},
\qquad
\mathcal G^{\cri}
\simeq
\varprojlim_\alpha\mathcal G_\alpha^{\cri},
\qquad
\mathcal G^{\dR}
\simeq
\varprojlim_\alpha\mathcal G_\alpha^{\dR},
\]
where \(\alpha\) ranges over the finite-type motivic quotients.
\end{theorem}

The arithmetic applications over \(\mathbb Q\) use the following
dimension consequence of the Bloch--Kato conjecture; see
\cite[Lemma~2.10]{BETTS_CORWIN}. For every pure geometric
\(p\)-adic representation \(W\) of weight at most \(-3\),
Conjecture~\ref{conj:BCL_BK}, a consequence of the Bloch--Kato Conjectures, predicts
\[
\dim H^1_f(\mathbb Q,W)
=
\dim H^1_f(\mathbb Q_p,W)
-
\dim W^+,
\]
where \(W^+\) denotes the subspace fixed by complex conjugation.
Thus every nonzero \(W^+\) contributes a strict deficit between the
global and local Selmer dimensions. The formula applies only in
weights at most \(-3\), so it does not control the first two DCS graded pieces.

Consequently, cofinality alone is not enough for the arithmetic
applications. We need a tower whose low-weight contribution remains
uniformly controlled, while nonzero conjugation-invariant higher-weight graded pieces occur
in arbitrarily large degrees.

\subsubsection{Arithmetic Applications}
\label{subsubsection:intro_arithmetic_applications}

We provide different approaches to this in two different cases, as described in Convention~\ref{convention:fixed_kp_family}: the \emph{modular case}, in which we take $X$ to be a modular curve and $f$ to be the universal elliptic curve, and the \emph{genus $\ge 2$ case}, where we take $X$ to be a smooth curve of genus $\ge 2$ and $f$ to be the restriction of a Kodaira--Parshin family on its compactification as in Lawrence--Venkatesh. Note that for $X$ a modular curve of genus $\ge 2$, we may apply either the \emph{modular case} to $X$ or the \emph{genus $\ge 2$ case} to its smooth compactification to get the same finiteness result.


The relative cutter is applied to the intersection of each Kummer
stratum with each residue tube at \(v\). In the modular case,
Theorem~\ref{theorem:faltings_finiteness} and
Proposition~\ref{proposition:kummer_unip} give a finite partition into
neutral Kummer strata with \(T'=T\). In the genus \(\geq2\) case,
Theorem~\ref{thm:kummer_stratification} gives such a finite partition
with \(T'=S_p\). Refining either partition by reduction modulo \(v\)
produces only finitely many further pieces, indexed by
\(\underline b\in\mathcal X(k_v)\). After choosing an integral base
point \(b\) in each nonempty piece, the cutter applies to the
corresponding locus \(\Sigma_b^{\mathrm{glob}}\). Consequently, global
finiteness reduces to verifying
\eqref{eq:intro_dimension_inequality}
on finitely many such loci.

For the genus $\ge 2$ case, we construct a distinguished tower of
adjoint-supported quotients. Here \emph{adjoint-supported} means that
every simple \(\mathcal R\)-constituent of every
DCS graded piece of the Lie algebra of the
unipotent radical occurs in
\[
\Lie\left((\mathcal R^\circ)^{\mathrm{der}}\right),
\]
the adjoint representation of the connected semisimple monodromy
group; see Definition~\ref{definition:adjoint_support}. The tower is
obtained by truncating the maximal quotient with this property.

\begin{theorem}
[The adjoint-supported motivic tower;
Theorem~\ref{theorem:special_motivic_exists} and
Lemma~\ref{lemma:adjoint_supported_complex_conjugation}]
\label{theorem:intro_adjoint_supported_tower}
Let
\[
Y\longrightarrow X'\longrightarrow X
\]
be a Kodaira--Parshin family over a smooth projective curve \(X/K\)
of genus at least \(2\). There is a canonical inverse system
\[
\mathcal G^\bullet
\twoheadrightarrow
\cdots
\twoheadrightarrow
\mathcal G_{\mathrm{ad},n+1}^\bullet
\twoheadrightarrow
\mathcal G_{\mathrm{ad},n}^\bullet
\twoheadrightarrow
\cdots
\twoheadrightarrow
\mathcal G_{\mathrm{ad},1}^\bullet
\]
of adjoint-supported finite-type motivic quotients, obtained as the
DCS truncations of the maximal
adjoint-supported quotient. Every transition morphism is faithfully
flat and is not an isomorphism.

For every \(k\geq1\), the graded piece
\[
W_k^\bullet
\coloneqq 
\Gr_k^{\DCS}
\Lie\left(\mathcal U_{\mathrm{ad},n}^\bullet\right),
\qquad n\geq k,
\]
is independent of \(n\), finite-dimensional, and nonzero. Its étale
realization is pure of weight \(-k\). If \(K=\mathbb Q\), then
\[
\dim\left(W_k^{\et}\right)^+\geq1
\qquad(k\geq1).
\]
\end{theorem}

The fixed finite-dimensional abelianization controls the
contribution outside the range of the Bloch--Kato dimension formula,
while the nonzero complex-conjugation invariants in the higher graded
pieces provide the dimension growth required for (\ref{eq:intro_dimension_inequality}) to hold.

We now specialize to \(K=\mathbb Q\). For the modular application, let
\[
f:\mathcal E\longrightarrow X
\]
be the universal elliptic curve. Recall that \(X\) has
\emph{enough motivic Eisenstein classes} if, for infinitely many
integers \(m\geq0\), the compatible realization system
\[
H^1_\bullet
\left(
X,\Sym^m R^1_\bullet f_*\mathbf 1
\right),
\qquad
\bullet\in\{\et,\cri,\dR\},
\]
contains a compatible rank-one direct summand whose étale realization
is the pure Tate representation \(\mathbb Q_p(-m-1)\); see
Definitions~\ref{definition:motivic_eisenstein_class}
and~\ref{definition:eis_quotient}.

\begin{theorem}
[Arithmetic finiteness;
Corollaries~\ref{corollary:modular_finiteness}
and~\ref{corollary:hyperbolic_main}]
\label{theorem:main_main}
Let
\[
\mathcal X
\longrightarrow
\Spec\mathbb Z[1/S]
\]
be a smooth model, and let \(p\notin S\) satisfy the standing
good-reduction hypotheses. Assume
Conjecture~\ref{conj:BCL_BK}. Suppose that either
\begin{enumerate}
\item
\(X\) is a modular curve with enough motivic Eisenstein classes, or

\item
the smooth proper compactification of \(X\) has genus at least \(2\).
\end{enumerate}
Then there exists a finite type motivic quotient for which (\ref{eq:intro_dimension_inequality}) holds; i.e.,
\[
\dim\widetilde{\mathscr S}_{n,T'}
+
\dim\mathcal G_{n,K_v}^{\varphi=1}
+
\dim F^0\mathcal G_n^{\dR}
<
\dim\mathcal G_n^{\dR}.
\]

In particular, it follows from our Theorem~\ref{theorem:main_intro}(2) and Theorem~\ref{theorem:faltings_finiteness} that
\[
\#\mathcal X(\mathbb Z[1/S])<\infty.
\]
\end{theorem}

\begin{remark}
We emphasize that the theorem gives a new route from the Bloch--Kato conjectures to
Faltings' theorem. In the ordinary Chabauty--Kim setting, the requisite
dimension estimates are known unconditionally in several special
motivic situations: for
\(\mathbb P^1\setminus\{0,1,\infty\}\) through its mixed-Tate
fundamental group~\cite{KIM2}, for curves whose Jacobians are
geometrically of CM type~\cite{COATES_KIM}, and for superelliptic curves over
\(\mathbb Q\)~\cite{ELLENBERG}. While such control is not
available for a general curve, the additional freedom afforded by
relative completion through the choice of reductive variation and
a finite-type quotient may open a route to unconditional, and
ultimately effective, results for broader classes of curves.
\end{remark}

In the modular case, a motivic Eisenstein class of degree \(m\)
determines the compatible \emph{Eisenstein quotient}
\(\mathcal Q_m^\bullet\) of the abelianized unipotent radical and
hence a finite-type motivic quotient \(\mathcal G_m^\bullet\) of the
relative completion; see \eqref{eq:eis_quotient}.
Theorem~\ref{theorem:modular_main} shows that
\(\mathcal G_m^\bullet\) satisfies
\eqref{eq:intro_dimension_inequality} for $m$ sufficiently large; no higher unipotent
truncation is required. Thus the modular application is an instance of relative
Chabauty--Skolem. By
Fact~\ref{fact:what_has_enough_eis}, the hypothesis on Eisenstein
classes holds for \(Y_1(N)\) whenever \(N\geq4\), and in particular
for
\[
Y_1(4)_{\mathbb Z[1/2]}
\simeq
\mathbb P^1_{\mathbb Z[1/2]}\setminus\{0,1,\infty\}.
\]
In the genus \(\geq2\) case, we show that every sufficiently deep truncation of
the adjoint-supported motivic tower of Theorem \ref{theorem:intro_adjoint_supported_tower} satisfies the dimension
inequality under Conjecture \ref{conj:BCL_BK}. We emphasize that this tower is disjoint from the purely unipotent tower of Chabauty--Kim, which implies by our Zariski density result that we get genuinely new analytic functions vanishing on the global points. The construction of the Selmer spaces, motivic
quotients, period maps, Bloch--Kato logarithms, and the relative box
cutter is unconditional; Conjecture \ref{conj:BCL_BK} is used only to
verify the dimension inequalities in these applications.

\subsection{Kantor's Program}
\label{subsection:kantor_program}

Kantor's thesis proposes a broader program than the finite-level
implementation developed here. Two directions for future work
suggested by this broader program are especially relevant to
effectivity: a substantially coarser arithmetic stratification of
the Selmer stack and a period map varying over the entire
\(p\)-adic integral locus.

\paragraph{Residual strata and algebraicity.}

Lawrence--Venkatesh and the present work partition the integral locus
according to reductive Galois types with $\Qp$ (hence characteristic-zero) coefficients.

Kantor instead proposes fixed-type Selmer strata indexed, in part, by
residual pseudorepresentations; see
\cite[Definition~5.1.2, Definition~7.0.2, Proposition~7.0.5]{KANTORTHESIS}. Over an algebraically closed residue field, a pseudorepresentation records the semisimplification of the residual representation. Consequently, a single residual stratum may contain many distinct characteristic-zero reductive types. Replacing the enumeration of the latter by a residual stratification could therefore yield a substantial computational advantage over both the
Lawrence--Venkatesh stratification and the finer Kummer
stratification used in this paper.

The principal obstacle is algebraicity. Kantor proves that the
ambient continuous non-abelian cohomology stacks are
strat-representable by rigid-analytic stacks, but the corresponding
Bloch--Kato loci are not presently known, in this generality, to be
algebraizable; see
\cite[
Theorem~4.0.1,
Chapters~4--5,
and Remark~14
]{KANTORTHESIS}.
This distinction is structural rather than merely technical. A
one-dimensional rigid-analytic map may have algebraically Zariski
dense image in a higher-dimensional algebraic target. This
flexibility is precisely what makes density of the period map
possible, but it also means that analytic dimension alone cannot
supply the algebraic dimension cutter: one needs algebraizable
strata on which scheme-theoretic images and algebraic dimensions can
be controlled.

At the reductive level, the required algebraization of the
potentially crystalline fixed-type loci would follow from the
FGAMS conjecture and the results of Wang--Erickson; see
\cite[
\S3.3]{CARL_WE}
and
\cite[Proposition~7.0.5]{KANTORTHESIS}.
Propagating this algebraicity through the relatively unipotent tower
requires coherence results for Galois cohomology in rigid-analytic
families with the relevant \(p\)-adic Hodge-theoretic local
conditions. Pottharst proves the corresponding finiteness and
coherence statement without these local conditions
\cite[Theorem on p.~1]{POTTHARST_SELMER}, while Shah establishes the
interpolation of Hodge--Tate and de Rham periods
\cite{SHAH_PERIODS}. The crystalline analogue required in Kantor's
inductive construction remains open; see
\cite[Remark~14]{KANTORTHESIS}.

Our neutral-fibre construction takes a different route. After fixing
an integral basepoint \(b\), we restrict to the fibre over the neutral
characteristic-zero reductive Galois type determined by \(b\). The remaining
directions are governed by finite-level unipotent cohomology, which
allows us to prove algebraicity unconditionally; see
Theorem~\ref{theorem:neutral_selmer_algebraicity}. The price is a
finer arithmetic stratification. An algebraization of Kantor's
residual strata would make it possible to retain the coarser
stratification while applying the same dimension-cutting philosophy.

\begin{remark}[General fixed reductive fibres]
The neutral-fibre construction suffices for qualitative finiteness,
since one may choose a basepoint in each reductive Galois
type in the image of the Kummer map and work relative to it. For effectivity, it would be preferable to construct Selmer stacks over arbitrary fixed reductive types. This would allow each type to be analyzed directly, without first selecting an integral point representing it. This should not present significant technical hurdles, as it should be representable by the same kind of non-abelian Galois cohomology of a pro-unipotent group; we leave this twisted
fixed-fibre theory to future work.
\end{remark}

\paragraph{Global period maps and relative Coleman integration.}

A second ambitious feature of Kantor's program is its intended
globality. His relative unipotent Albanese map is formulated
pointwise on the full \(p\)-adic integral locus and is locally
analytic around every point; see
\cite[\S\S6.2 and~6.5]{KANTORTHESIS}. By contrast, both the
Lawrence--Venkatesh period map and the period map constructed in this
paper are defined after fixing a residue tube.

In the ordinary unipotent setting, the passage between residue disks
is supplied by Besser--Coleman integration. Its Tannakian
construction uses Frobenius-invariant crystalline paths
\cite{BESSER_COLEMAN}. For a torsor under a negatively weighted
unipotent group, the Frobenius Lang map is an isomorphism, and hence
the Frobenius-invariant path exists uniquely; see also
\cite[Theorem~6.2.5]{KANTORTHESIS}. This uniqueness gives canonical
parallel transport and permits Kim's unipotent Albanese map to be
continued across residue disks.

For a full relative path torsor, neither conclusion remains valid.
Even when Frobenius-invariant paths exist, they form a torsor under
the generally nontrivial group
\(\mathcal G^{\varphi=1}\), and hence need not determine a canonical
transport. More fundamentally, the Frobenius reduction may be empty.

The Legendre family illustrates the latter obstruction. Let
\[
f:
\mathcal E
\longrightarrow
X
=
\mathbb P^1\setminus\{0,1,\infty\},
\]
and work over a \(p\)-adic field whose residue field contains smooth
ordinary and supersingular fibres of \(f\). Choose lifts
\[
b,x\in X(\mathcal O_{K_v})
\]
of an ordinary point \(\underline b\) and a supersingular point
\(\underline x\), respectively. If the reductive pushout of the
crystalline path torsor
\[
{}^{\vphantom{\cri}}_{\underline x}P_{\underline b}^{\cri}
\]
had a Frobenius-fixed point, it would induce, after extension of
scalars, an isomorphism of \(F\)-isocrystals
\[
H^1_{\cris}
\left(
\mathcal E_{\underline b}/K_0
\right)
\xrightarrow{\sim}
H^1_{\cris}
\left(
\mathcal E_{\underline x}/K_0
\right).
\]
This is impossible: the ordinary fibre has Frobenius slopes
\[
\{0,1\},
\]
whereas the supersingular fibre has slopes
\[
\left\{
\frac12,\frac12
\right\}.
\]
By the Dieudonn\'e--Manin classification after base change to an
algebraic closure of the residue field, the slope multiset is invariant
under Frobenius-equivariant isomorphism and is preserved under
extension of the residue field; see
\cite{KEDLAYA_ISOC}. Thus the Frobenius reduction of the relative de
Rham path torsor
\[
{}^{\vphantom{\dR}}_{x}P_b^{\dR}
\]
is empty, and this torsor is not admissible. In particular, a
period map valued in the moduli of admissible relative torsors cannot
be extended across arbitrary residue disks.

The present paper avoids this obstruction by working on a fixed
residue tube. For
\[
x\in]\underline b[(K_v),
\]
one has
\[
\underline x=\underline b,
\]
so the corresponding crystalline path torsor is the trivial torsor at
\(\underline b\), whose identity supplies the required canonical
Frobenius-fixed point. We globalize this construction over the tube
using the universal crystalline and de Rham path torsors and their
family-level comparison. This yields a rigorous rigid-analytic
period map and proves its analyticity and algebraic Zariski density;
see
Proposition~\ref{proposition:an_rep}
and
Proposition~\ref{proposition:p-adic-dense}.

Thus the present work realizes a finite-level, algebraic, and
residue-local part of Kantor's program. It trades the coarseness of
the conjectural residual stratification and the globality of a
relative Coleman theory for unconditional algebraicity and canonical
Frobenius transport. A relative Besser--Coleman theory should, more
generally, identify the locus on which cross-disk Frobenius transport
exists and control its nonuniqueness. We intend to
return to this problem in future work.

\subsection{Plan of Paper}
The paper is organized as follows.
\begin{itemize}
    \item
    \S\ref{chapter:relative_malcev_completion} reviews relative
    completion from the Tannakian viewpoint and specializes the
    construction to the realizations associated with Kodaira--Parshin
    families.

    \item
    \S\ref{section:finite_type_motivic_quotients} develops the
    finiteness theory of motivic quotients and constructs the
    adjoint-supported finite-type tower used in the applications.

    \item
    \S\ref{chapter:selmer_stacks} constructs neutral-fibre Selmer
    stacks and relative Kummer maps, proves their finite-level
    algebraicity, and establishes the required finite Kummer
    stratifications.

    \item
    \S\ref{section:dr} constructs the $v$-adic period map and the
    neutral Bloch--Kato logarithm, proves their compatibility and
    density properties, and derives the dimension-cutter criterion.

    \item
    \S\ref{chapter:applications} applies this criterion to modular
    curves with enough Eisenstein classes and to curves of genus at
    least $2$, conditional on the relevant Bloch--Kato formulas.

    \item
    Appendix~\ref{appendix:rmc} establishes the comparison and
    semisimplicity results for the Kodaira--Parshin realizations used
    throughout the paper.

    \item
    Appendix~\ref{appendix:filtered_de_rham_path_groupoid} constructs
    the Hodge filtration on the universal relative de Rham path torsor
    and descends it to the finite-level motivic quotients.

    \item
    Appendix~\ref{appendix:tannakian_finiteness} collects the auxiliary
    Tannakian finiteness results used in the proof of
    Theorem~\ref{thm:ftg_basic} and gives an alternative proof of that
    theorem using universal extensions.

    \item
    Appendix~\ref{appendix:adjoint_supported_quotients} develops the
    formalism of supported quotients of pro-unipotent extensions used
    to construct the adjoint-supported tower.
\end{itemize}

\section*{Acknowledgements}
The authors are deeply grateful to Ishai Dan-Cohen and Noam Kantor for their generosity
with their time and ideas, which substantially shaped this work. They also
thank Amnon Besser, Nadav Groper, John Pridham, and Mark Shusterman for helpful discussions. 
D.C.~thanks Will Sawin for comments on the proof of
Proposition~\ref{prop:embed_hereditary}. S.Z.~thanks Ishai Dan-Cohen
for many conversations on the theoretical foundations of the
subject.

S.Z. acknowledges support from the Israel Science Foundation
(grant No.~621/21), and the European Union
(ERC, Function Fields, 101161909).
﻿\section{The Relative Completion}
\label{chapter:relative_malcev_completion}
This chapter is dedicated to the theory of relative completions.
\begin{itemize}
    \item \S\ref{section:abstract_theory} focuses on the abstract aspects of the theory from a Tannakian point of view.
    \item \S\ref{section:rmc_realizations} specializes our discussion to relative completions 
    in the Kodaira--Parshin setting.
\end{itemize}

\subsection{Abstract Theory}
\label{section:abstract_theory}

\subsubsection{Tannakian Categories}
References for Tannakian categories may be found
in~\cite{DelMil_TC22,DelMilOgShi82}. Let \(k\) be a Hausdorff
topological field, and let \((\mathcal T,\omega)\) be a neutral
Tannakian category over \(k\).

\begin{definition}
\label{definition:tannakian_fundamental_group}
If \(\omega'\) is another fibre functor, we denote by
\[
\operatorname{Isom}^{\otimes}(\omega,\omega')
\]
the affine \(k\)-scheme of tensor isomorphisms between \(\omega\) and
\(\omega'\). When \(\omega=\omega'\), the affine \(k\)-group scheme
\[
\operatorname{Aut}^{\otimes}(\mathcal T,\omega)
\coloneqq 
\operatorname{Isom}^{\otimes}(\omega,\omega)
\]
is called the Tannakian fundamental group of \(\mathcal T\) based at
\(\omega\).
\end{definition}

\begin{remark}
\label{remark:torsors}
The scheme \(\operatorname{Isom}^{\otimes}(\omega,\omega')\) is a right
\(\operatorname{Aut}^{\otimes}(\mathcal T,\omega)\)-torsor and a left
\(\operatorname{Aut}^{\otimes}(\mathcal T,\omega')\)-torsor.
\end{remark}

Let \(G\) be a topological (resp. pro-algebraic) group. We write
\[
\operatorname{Rep}_k(G)
\]
for the neutral Tannakian category of finite-dimensional continuous (resp. rational)
\(k\)-linear representations of \(G\), equipped with its
forgetful fibre functor \(\omega\). Finite-dimensional \(k\)-vector
spaces carry their standard topology; when \(G\) is discrete,
this is the usual category of finite-dimensional representations. When $G$ is pro-algebraic, all such representations factor through
a finite-type algebraic quotient.

If $\Gamma$ is a topological group and $H$ pro-algebraic, we equip $H(k)$ with the inverse-limit topology induced by its algebraic quotients, and all
homomorphisms
\[
\Gamma\longrightarrow H(k)
\]
are understood to be continuous. The Tannakian fundamental group
\[
\operatorname{Aut}^{\otimes}
\bigl(\operatorname{Rep}_k(\Gamma),\omega\bigr)
\]
is the algebraic envelope of \(\Gamma\) with respect to its continuous
finite-dimensional representations. The canonical homomorphism from
\(\Gamma\) to its \(k\)-points is continuous and Zariski dense.

\begin{definition}
\label{def:algebraic_tannakian}
A neutral Tannakian category $(\mathcal{Q}, \omega)$ is called \emph{algebraic} 
if it is equivalent to the category of representations of an affine algebraic group scheme. 
Equivalently (see~\cite[Proposition~2.20]{DelMil_TC22}), $\mathcal{Q}$ is algebraic if it is 
tensor generated by a single object $\mathcal{M}\in \mathcal{Q}$, i.e., $\mathcal{Q} = \langle \mathcal{M} \rangle_{\otimes}$.
\end{definition}

\subsubsection{The Relative Completion}
Fix \(\mathcal L\in\operatorname{Rep}_k(\Gamma)\), and write
\[
\mathcal D\coloneqq \langle\mathcal L\rangle_{\otimes}
\]
for the smallest full Tannakian subcategory containing \(\mathcal L\).
Thus \(\mathcal D\) is closed under finite direct sums, tensor products,
duals, and subquotients. Let
\[
\mathcal R
\coloneqq 
\operatorname{Aut}^{\otimes}
\bigl(\mathcal D,\omega|_{\mathcal D}\bigr).
\]

Since \(\mathcal D\) is generated by a single object, \(\mathcal R\) is
an affine algebraic \(k\)-group scheme. We shall assume that
\(\mathcal R\) is reductive. When \(\operatorname{char}(k)=0\), this is
equivalent to \(\mathcal L\) being semisimple, and also to every object
of \(\mathcal D\) being semisimple
(see~\cite[Proposition~2.23]{DelMil_TC22}).

\begin{definition}
\label{definition:thick}
The thick tensor category generated by \(\mathcal L\), denoted
\[
\mathcal C
\coloneqq 
\langle\mathcal L\rangle_{\otimes}^{\mathrm{ext}},
\]
is the smallest full Tannakian subcategory of
\(\operatorname{Rep}_k(\Gamma)\) containing \(\mathcal D\) and closed
under extensions. Its objects are precisely the representations
admitting a finite filtration
\[
0=V_0\subset V_1\subset\cdots\subset V_n=V
\]
such that
\[
V_i/V_{i-1}\in\mathcal D
\qquad
(1\leq i\leq n).
\]
\end{definition}

The inclusion \(\mathcal D\subseteq\mathcal C\) is full, and
\(\mathcal D\) is closed under subquotients taken in \(\mathcal C\).

Let
\[
\rho:\Gamma\longrightarrow\mathcal R(k)
\]
be the continuous Zariski-dense monodromy homomorphism induced by
\(\mathcal L\).

\begin{definition}
The relative Malcev completion of \((\Gamma,\rho)\) is the pair
\[
(\mathcal G,\widetilde\rho),
\]
where
\[
\mathcal G
\coloneqq 
\operatorname{Aut}^{\otimes}
\bigl(\mathcal C,\omega|_{\mathcal C}\bigr)
\]
and
\[
\widetilde\rho:\Gamma\longrightarrow\mathcal G(k)
\]
is the canonical continuous Zariski-dense homomorphism induced by the
action of \(\Gamma\) on the fibre functor.
\end{definition}

We also refer to \(\mathcal G\) as the relative completion of
\(\Gamma\). The inclusion \(\mathcal D\subseteq\mathcal C\) induces a
morphism
\[
\pi:\mathcal G\longrightarrow\mathcal R
\]
satisfying
\[
\pi(k)\circ\widetilde\rho=\rho.
\]
Since \(\mathcal D\) is full and closed under subquotients in
\(\mathcal C\), the morphism \(\pi\) is faithfully flat by
\cite[Proposition~2.21(a)]{DelMil_TC22}. Writing
\(\mathcal U\coloneqq \ker(\pi)\), we obtain an exact sequence
\begin{equation}
\label{ses:rel_comp}
1\longrightarrow\mathcal U
\longrightarrow\mathcal G
\xrightarrow{\pi}\mathcal R
\longrightarrow 1.
\end{equation}
The group scheme \(\mathcal U\) is pro-unipotent
(see~\cite[Proposition~2.1.2]{KANTORTHESIS}). If
\(\operatorname{char}(k)=0\), the sequence admits a noncanonical Levi
splitting, and hence a noncanonical isomorphism
\begin{equation}
\label{rel_completion_semidirect}
\mathcal G\simeq\mathcal U\rtimes\mathcal R
\end{equation}
by~\cite[Proposition~3.1]{Hain_HodgeDeRham_Chapter2016}.

\begin{proposition}[{Universal property of the relative completion;
cf.~\cite[Proposition~2.1.4]{KANTORTHESIS}}]
\label{prop:rmc_univ_property}
Let
\[
1\longrightarrow U'
\longrightarrow G'
\xrightarrow{p}\mathcal R
\longrightarrow 1
\]
be an exact sequence of affine pro-algebraic \(k\)-group schemes, with
\(U'\) pro-unipotent. For every continuous homomorphism
\[
\varphi:\Gamma\longrightarrow G'(k)
\]
such that \(p(k)\circ\varphi=\rho\), there exists a unique morphism
\[
F_\varphi:\mathcal G\longrightarrow G'
\]
such that
\[
p\circ F_\varphi=\pi,
\qquad
F_\varphi(k)\circ\widetilde\rho=\varphi.
\]
Equivalently, composition with \(\widetilde\rho\) induces a bijection
\[
\operatorname{Hom}_{\mathcal R}(\mathcal G,G')
\xrightarrow{\sim}
\left\{
\varphi\in
\operatorname{Hom}_{\mathrm{cts}}\bigl(\Gamma,G'(k)\bigr)
\ \middle|\
p(k)\circ\varphi=\rho
\right\}.
\]
\end{proposition}

\begin{proof}
The displayed map is well defined because
\[
p\circ F=\pi
\qquad\Longrightarrow\qquad
p(k)\circ F(k)\circ\widetilde\rho=\rho.
\]
We prove that it is surjective and injective.

\smallskip
\noindent\textbf{Surjectivity.}
Let
\[
\varphi:\Gamma\longrightarrow G'(k)
\]
be a continuous lift of \(\rho\). Every
\(V\in\operatorname{Rep}_k(G')\) factors through an algebraic quotient
of \(G'\). Since the image of \(U'\) in this quotient is unipotent and
normal, the usual invariant-vector argument gives a finite
\(G'\)-stable filtration of \(V\) whose graded pieces have trivial
\(U'\)-action. As \(G'/U'\simeq\mathcal R\), these graded pieces are
pulled back from \(\operatorname{Rep}_k(\mathcal R)\). After restriction
along \(\varphi\), they therefore belong to \(\mathcal D\), because
\(p(k)\circ\varphi=\rho\). Hence restriction along \(\varphi\) defines
an exact tensor functor
\[
\Phi_\varphi:
\operatorname{Rep}_k(G')
\longrightarrow
\mathcal C
\simeq
\operatorname{Rep}_k(\mathcal G)
\]
preserving the fibre functors. By Tannakian
duality~\cite[Corollary~2.9]{DelMil_TC22}, this functor is induced by a
morphism
\[
F_\varphi:\mathcal G\longrightarrow G'.
\]
By construction,
\[
F_\varphi(k)\circ\widetilde\rho=\varphi.
\]
Consequently,
\[
p(k)\circ F_\varphi(k)\circ\widetilde\rho
=
p(k)\circ\varphi
=
\rho
=
\pi(k)\circ\widetilde\rho.
\]
Since \(\widetilde\rho(\Gamma)\) is Zariski dense in \(\mathcal G\), it
follows that \(p\circ F_\varphi=\pi\).

\smallskip
\noindent\textbf{Injectivity.}
Suppose that
\[
F_1,F_2\in\operatorname{Hom}_{\mathcal R}(\mathcal G,G')
\]
satisfy
\[
F_1(k)\circ\widetilde\rho
=
F_2(k)\circ\widetilde\rho.
\]
Then \(F_1\) and \(F_2\) agree on the Zariski-dense subset
\(\widetilde\rho(\Gamma)\subseteq\mathcal G(k)\), and hence
\(F_1=F_2\). Thus the displayed map is bijective.
\end{proof}

\begin{corollary}[{Functoriality of relative Malcev completion;
cf.~\cite[p.~12]{KANTORTHESIS}}]
\label{corollary:rmc_functoriality}
Let
\[
\rho:\Gamma\longrightarrow\mathcal R(k),
\qquad
\rho':\Gamma\longrightarrow\mathcal R'(k)
\]
be continuous Zariski-dense homomorphisms into reductive affine
algebraic \(k\)-group schemes, and let
\[
f:\mathcal R\longrightarrow\mathcal R'
\]
satisfy \(\rho'=f(k)\circ\rho\). If
\((\mathcal G,\widetilde\rho)\) and
\((\mathcal G',\widetilde\rho')\) are the corresponding relative
completions, then there exists a unique morphism
\[
\widetilde f:\mathcal G\longrightarrow\mathcal G'
\]
making the diagram
\[
\begin{tikzcd}
\Gamma
  \arrow[r, "\widetilde\rho"]
  \arrow[rd, "\widetilde\rho'"']
&
\mathcal G
  \arrow[r, "\pi"]
  \arrow[d, "\widetilde f"']
&
\mathcal R
  \arrow[d, "f"]
\\
&
\mathcal G'
  \arrow[r, "\pi'"]
&
\mathcal R'
\end{tikzcd}
\]
commute. These morphisms are compatible with identities and
composition.
\end{corollary}

\begin{proof}
Form the pullback
\[
\mathcal H
\coloneqq 
\mathcal R\times_{\mathcal R'}\mathcal G'.
\]
Then \(\mathcal H\to\mathcal R\) is a pro-unipotent extension, and
\[
\Gamma\longrightarrow\mathcal H(k),
\qquad
\gamma\longmapsto
\bigl(\rho(\gamma),\widetilde\rho'(\gamma)\bigr),
\]
is a continuous lift of \(\rho\). Proposition
\ref{prop:rmc_univ_property} gives a unique morphism
\(\mathcal G\to\mathcal H\); composing with
\(\mathcal H\to\mathcal G'\) gives \(\widetilde f\). Any other
morphism with the same properties induces another morphism
\(\mathcal G\to\mathcal H\) with the prescribed lift, and is therefore
equal to \(\widetilde f\). Compatibility with identities and
composition follows from the same uniqueness.
\end{proof}

In the following examples, \(\Gamma\) is given the discrete topology,
\(\mathcal R=\mathbf 1\), and
\[
\rho:\Gamma\longrightarrow\mathcal R(k)=\{1\}
\]
is the trivial representation.

\begin{example}[{\cite[\S2]{KNUDSON}}]
If \(k=\mathbb Q\), then \(\mathcal G\) is canonically isomorphic to
the Malcev, or pro-unipotent, completion of \(\Gamma\).
\end{example}

\begin{example}[{\cite[\S3.3]{KNUDSON}}]
Let \(k=\mathbb F_p\), and suppose that
\[
\dim_{\mathbb F_p}H_1(\Gamma,\mathbb F_p)<\infty;
\]
this holds, for example, when \(\Gamma\) is finitely generated. Then
\(\mathcal G(\mathbb F_p)\), equipped with its inverse-limit topology,
is canonically isomorphic to the pro-\(p\) completion of \(\Gamma\):
\[
\mathcal G(\mathbb F_p)
\simeq
\widehat\Gamma^{(p)}
\simeq
\varprojlim_r\Gamma/\Gamma_p^r\Gamma,
\]
where \(\Gamma_p^\bullet\Gamma\) denotes the lower \(p\)-central
series.
\end{example}

\subsection{Relative Completions Associated with Kodaira--Parshin Families}
\label{section:rmc_realizations}

\paragraph{Kodaira--Parshin families.}
Let \(F\) be a field of characteristic
zero equipped with an embedding \(F\hookrightarrow\mathbb C\), and
let \(C/F\) be a smooth quasi-projective curve.

\begin{definition}
[{\cite[Definition~5.1 and the discussion following it]{LV};
cf.~\cite[Definitions~7.1.1 and~7.2.1]{JI_LV}}]
\label{definition:kp_family}
An \emph{abelian-by-finite family} over \(C\) is a factorization
\[
Y\xrightarrow{g}C'\xrightarrow{\pi}C,
\]
where \(g\) is a polarized abelian scheme and \(\pi\) is finite
\'{e}tale. For \(x_0\in C(\mathbb C)\), there is a decomposition
\[
H^1_{\Betti}(Y_{x_0},\mathbb Q)
\simeq
\bigoplus_{x'\in\pi^{-1}(x_0)}
H^1_{\Betti}(Y_{x'},\mathbb Q).
\]
The family has \emph{full monodromy} if, for one (equivalently, every)
\(x_0\), the Zariski closure of the image of
\[
\pi_1\bigl(C(\mathbb C),x_0\bigr)
\longrightarrow
\operatorname{GL}
\bigl(H^1_{\Betti}(Y_{x_0},\mathbb Q)\bigr)
\]
contains the block-diagonal subgroup
\[
\prod_{x'\in\pi^{-1}(x_0)}
\operatorname{Sp}
\bigl(H^1_{\Betti}(Y_{x'},\mathbb Q),\omega_{x'}\bigr),
\]
where \(\omega_{x'}\) is induced by the polarization. We call an abelian-by-finite family with full monodromy a
\emph{Kodaira--Parshin family}; cf.~\cite[\S\S7--8]{LV} and \cite[\S7.3]{JI_LV}.
\end{definition}

\begin{fact}[{\cite[Theorem~8.1]{LV}}]
\label{fact:kp_monodromy}
Let \(C\) be a smooth projective curve of genus \(g\geq 2\) over a
number field \(F\). Then \(C\) admits a Kodaira--Parshin family
\[
Y\longrightarrow C'\longrightarrow C
\]
for which \(Y\to C'\) has relative dimension \(2g-1\).
\end{fact}

We now return to the standing arithmetic setup of
\S\ref{subsubsection:intro_relative_box_cutter}.

\begin{convention}[The fixed Kodaira--Parshin family]
\label{convention:fixed_kp_family}
We fix the following geometric data.

In the \emph{modular case} (cf. \S\ref{subsubsection:intro_arithmetic_applications}), we take
\[
Y=\mathcal E,
\qquad
X'=X,
\qquad
f:Y\longrightarrow X
\]
to be the universal elliptic curve.

In the \emph{genus $\ge 2$ case}, we assume that the smooth projective compactification
\(\widehat X\) has genus at least \(2\), and fix, by
Fact~\ref{fact:kp_monodromy}, a Kodaira--Parshin family
\[
\widehat Y
\xrightarrow{\widehat g}
\widehat X'
\xrightarrow{\widehat\pi}
\widehat X,
\qquad
\widehat f\coloneqq \widehat\pi\circ\widehat g.
\]
We denote its restriction to \(X\) by
\[
Y\longrightarrow X'\longrightarrow X,
\qquad
f:Y\longrightarrow X.
\]

After enlarging \(S\) away from the places above \(p\), we fix a
smooth proper model
\[
\widehat{\mathcal X}/\mathcal O_{K,S}
\]
of \(\widehat X\) such that
\[
\mathcal D
\coloneqq 
\widehat{\mathcal X}\setminus\mathcal X
\]
is a relative normal-crossings divisor.

In the modular case, the universal elliptic curve extends to a
polarized abelian scheme over \(\mathcal X\), and we fix its standard
semistable generalized-elliptic-curve compactification, proper and
log smooth over
\[
(\widehat{\mathcal X},\mathcal D).
\]
In the genus $\ge 2$ case, we fix an extension
\[
\widehat{\mathcal Y}
\longrightarrow
\widehat{\mathcal X}'
\longrightarrow
\widehat{\mathcal X}
\]
of the chosen Kodaira--Parshin family, whose first morphism is a
polarized abelian scheme and whose second is finite étale.

At every \(v\mid p\), we use the models obtained by base change to
\(\mathcal O_{K_v}\). For the comparison theorems, we use the
cuspidal log pair in the modular case and the empty boundary on
\(\widehat{\mathcal X}\) in the genus $\ge 2$ case.
\end{convention}

\paragraph{Realization categories and Gauss--Manin objects.}

To treat the two cases uniformly, set
\[
(X^\dagger,Y^\dagger,f^\dagger)
\coloneqq 
\begin{cases}
(X,Y,f),
& \text{in the modular case},\\[2mm]
(\widehat X,\widehat Y,\widehat f),
& \text{in the genus $\ge 2$ case}.
\end{cases}
\]
Thus
\[
f^\dagger:Y^\dagger\longrightarrow X^\dagger
\]
is smooth and proper. The points considered below lie in the open
subscheme
\[
X\subseteq X^\dagger,
\]
where this inclusion is the identity in the modular case. All
Tannakian categories, relative completions, and path torsors below
are formed on \(X^\dagger\) and evaluated at points of \(X\).

Fix \(v\mid p\), let \(k\) be the residue field of \(K_v\), and set
\[
K_0\coloneqq \operatorname{Frac}W(k).
\]
Write
\[
\overline f^\dagger:
Y^\dagger_{\overline K}
\longrightarrow
X^\dagger_{\overline K},
\qquad
f_v^\dagger:
Y^\dagger_{K_v}
\longrightarrow
X^\dagger_{K_v},
\qquad
\underline f^\dagger:
Y^\dagger_k
\longrightarrow
X^\dagger_k
\]
for the geometric generic, \(K_v\)-generic, and special fibres of
the fixed models. The basepoint \(b\) similarly determines
\[
\overline b\in X^\dagger_{\overline K}(\overline K),
\qquad
b_v\in X^\dagger_{K_v}(K_v),
\qquad
\underline b\in X^\dagger_k(k).
\]

We use the neutral Tannakian categories
\[
\begin{aligned}
\mathcal T^{\et}(X^\dagger)
&\coloneqq 
\operatorname{Lisse}_{\mathbb Q_p}
\left(
X^\dagger_{\overline K,\et}
\right),\\
\mathcal T^{\dR}(X^\dagger)
&\coloneqq 
\operatorname{VIC}
\left(
X^\dagger_{K_v}
\right),\\
\mathcal T^{\cri}(X^\dagger)
&\coloneqq 
\operatorname{Isoc}
\left(
X^\dagger_k/K_0
\right).
\end{aligned}
\]
Thus \(\mathcal T^{\et}(X^\dagger)\) is the category of
finite-dimensional lisse
\(\mathbb{Q}_p\)-sheaves on $X^\dagger$;
\(\mathcal T^{\dR}(X^\dagger)\) is the category of vector bundles on
\(X^\dagger_{K_v}\) equipped with an integrable connection; and
\(\mathcal T^{\cri}(X^\dagger)\) is the category of convergent
\(K_0\)-isocrystals on \(X^\dagger_k\). Their morphisms are, respectively, morphisms of lisse
\(\mathbb Q_p\)-sheaves, horizontal morphisms of vector bundles with
integrable connection, and morphisms of isocrystals. Their fibre functors at
\(b\) are
\[
\begin{aligned}
\omega_b^{\et}\coloneqq \overline b^{\,*}
&:
\mathcal T^{\et}(X^\dagger)
\longrightarrow
\operatorname{Vec}_{\mathbb Q_p},\\
\omega_b^{\dR}\coloneqq b_v^*
&:
\mathcal T^{\dR}(X^\dagger)
\longrightarrow
\operatorname{Vec}_{K_v},\\
\omega_b^{\cri}\coloneqq \underline b^{\,*}
&:
\mathcal T^{\cri}(X^\dagger)
\longrightarrow
\operatorname{Vec}_{K_0}.
\end{aligned}
\]
When the realization is clear, we write simply \(\omega_b\).

\begin{remark}
The geometric basepoint \(\overline b\) gives a tensor equivalence
\[
\mathcal T^{\et}(X^\dagger)
\xrightarrow{\sim}
\operatorname{Rep}^{\mathrm{cts}}_{\mathbb Q_p}
\left(
\pi_1^{\et}
\left(
X^\dagger_{\overline K},
\overline b
\right)
\right),
\]
under which \(\omega_b^{\et}\) corresponds to the forgetful fibre
functor. After fixing an embedding
\[
\overline K\hookrightarrow\overline{K_v}
\]
and a compatible geometric basepoint \(\overline b_v\), base change
gives an isomorphism
\[
\pi_1^{\et}
\left(
X^\dagger_{\overline K},
\overline b
\right)
\xrightarrow{\sim}
\pi_1^{\et}
\left(
X^\dagger_{\overline{K_v}},
\overline b_v
\right);
\]
see~\cite[Theorem~1.1]{LANDESMAN}. Thus the same geometric
\'{e}tale category is used in the global and local settings; only
its ambient arithmetic action differs, through \(G_K\) globally and
\(G_v\) locally.
\end{remark}

Define the geometric Gauss--Manin objects by
\[
\begin{aligned}
\mathcal L^{\et}
&\coloneqq 
R^1\overline f^\dagger_{\et,*}\mathbb Q_p,\\
\mathcal L^{\dR}
&\coloneqq 
R^1_{\dR}(f_v^\dagger)_*\mathbf 1,\\
\mathcal L^{\cri}
&\coloneqq 
R^1_{\cri}\underline f^\dagger_*\mathbf 1.
\end{aligned}
\]
These are objects of
\[
\mathcal T^{\et}(X^\dagger),
\qquad
\mathcal T^{\dR}(X^\dagger),
\qquad
\mathcal T^{\cri}(X^\dagger),
\]
respectively.

\begin{proposition}
\label{prop:semisimplicity_kp}
For every
\[
\bullet\in\{\et,\dR,\cri\},
\]
the object \(\mathcal L^\bullet\) is semisimple in
\(\mathcal T^\bullet(X^\dagger)\).
\end{proposition}

\begin{proof}
See Appendix~\ref{section:semisimplicity}.
\end{proof}

\paragraph{Relative completions and path torsors.}
Retain the Kodaira--Parshin family \(f^\dagger:Y^\dagger\to X^\dagger\) as in Convention~\ref{convention:fixed_kp_family}.

\begin{definition}
\label{definition:rel_comp}
For \(\bullet\in\{\et,\dR,\cri\}\), set
\[
\begin{aligned}
\mathcal R^\bullet
&\coloneqq 
\operatorname{Aut}^{\otimes}
\bigl(
\langle\mathcal L^\bullet\rangle_\otimes,
\omega_b^\bullet
\bigr),
\\
\mathcal G^\bullet
&\coloneqq 
\operatorname{Aut}^{\otimes}
\bigl(
\langle\mathcal L^\bullet\rangle_\otimes^{\mathrm{ext}},
\omega_b^\bullet
\bigr),
\\
\mathcal U^\bullet
&\coloneqq 
\ker\bigl(\mathcal G^\bullet\to\mathcal R^\bullet\bigr).
\end{aligned}
\]
Thus \(\mathcal G^\bullet\) is the relative completion associated with
\(\mathcal L^\bullet\), and there is an exact sequence
\[
1\longrightarrow\mathcal U^\bullet
\longrightarrow\mathcal G^\bullet
\longrightarrow\mathcal R^\bullet
\longrightarrow 1,
\]
with \(\mathcal R^\bullet\) reductive and
\(\mathcal U^\bullet\) pro-unipotent.

If \(x\) defines a fibre functor \(\omega_x^\bullet\) on
\(\langle\mathcal L^\bullet\rangle_\otimes^{\mathrm{ext}}\), define
the corresponding path torsor by
\[
{}^{\vphantom{\bullet}}_{x}P_b^\bullet
\coloneqq 
\operatorname{Isom}^{\otimes}
\bigl(\omega_b^\bullet,\omega_x^\bullet\bigr).
\]
It is a right \(\mathcal G^\bullet\)-torsor, and
\[
{}^{\vphantom{\bullet}}_{b}P_b^\bullet=\mathcal G^\bullet.
\]
\end{definition}

\paragraph{Realization structures and comparison.}

The relative completions and their path torsors carry structures
native to their respective realizations.

\begin{proposition}
[The \'{e}tale Galois action;
cf.~{\cite[Example~2.1.5 and Example~5.1.6]{KANTORTHESIS}}]
\label{proposition:etale_torsor_galois}
Let $F\in\{K,K_v\}$. For every $x,b\in X(\overline F)$ and every 
\(\sigma\in G_F\), there is a natural isomorphism
\[
\sigma:
{}^{\vphantom{\et}}_{x}P_b^{\et}
\xrightarrow{\sim}
{}^{\vphantom{\et}}_{\sigma x}P_{\sigma b}^{\et}.
\]
Consequently, \({}^{\vphantom{\et}}_{x}P_b^{\et}\) carries a continuous \(G_F\)-action
when \(b,x\in X(F)\). This action is compatible with composition of paths. 
If, moreover,
\[
x,b\in\mathcal X(\mathcal O_{K,S}),
\]
then, under the fixed good-model hypotheses, the \(G_K\)-action
factors through \(G_T\); see~\cite[Proposition~3.3.13]{KANTORTHESIS}.
\end{proposition}

\begin{proposition}[{The de Rham Hodge filtration;
cf.~\cite[\S6.4]{KANTORTHESIS}}]
\label{proposition:de_rham_torsor_filtration}
For \(b,x\in X(K_v)\), the coordinate algebra of the de Rham path
torsor carries a functorial decreasing, exhaustive, and multiplicative
Hodge filtration
\[
\cdots
\supseteq
F^p\mathcal O({}^{\vphantom{\dR}}_{x}P_b^{\dR})
\supseteq
F^{p+1}\mathcal O({}^{\vphantom{\dR}}_{x}P_b^{\dR})
\supseteq
\cdots.
\]
The coordinate Hopf algebra \(\mathcal O(\mathcal G^{\dR})\) carries
the analogous filtration, and the torsor coaction
\[
\mathcal O({}^{\vphantom{\dR}}_{x}P_b^{\dR})
\longrightarrow
\mathcal O({}^{\vphantom{\dR}}_{x}P_b^{\dR})
\otimes_{K_v}
\mathcal O(\mathcal G^{\dR})
\]
is a morphism of filtered algebras.

Define the Hodge locus by
\[
F^0{}^{\vphantom{\dR}}_{x}P_b^{\dR}
\coloneqq 
\operatorname{Spec}
\left(
\frac{
\mathcal O({}^{\vphantom{\dR}}_{x}P_b^{\dR})
}{
\left\langle
F^1\mathcal O({}^{\vphantom{\dR}}_{x}P_b^{\dR})
\right\rangle
}
\right),
\]
where
\[
\left\langle
F^1\mathcal O({}^{\vphantom{\dR}}_{x}P_b^{\dR})
\right\rangle
\]
denotes the ideal generated by the positive Hodge filtration.
Equivalently, for a \(K_v\)-algebra \(A\),
\[
F^0{}^{\vphantom{\dR}}_{x}P_b^{\dR}(A)
=
\left\{
p\in{}^{\vphantom{\dR}}_{x}P_b^{\dR}(A)
\ \middle|\
\operatorname{ev}_p:
\mathcal O({}^{\vphantom{\dR}}_{x}P_b^{\dR})
\longrightarrow A
\text{ is filtered}
\right\},
\]
where \(A\) carries the trivial filtration.
The filtrations are compatible with composition of paths and, when
\(x=b\), with the Hopf-algebra structure on
\(\mathcal O(\mathcal G^{\dR})\).
\end{proposition}
For the finite-level quotients used below, this construction and its
compatibility with the torsor structure are established in
Definition~\ref{definition:finite_level_hodge_loci} and
Corollary~\ref{corollary:finite_level_universal_hodge_reduction} of
Appendix~\ref{appendix:filtered_de_rham_path_groupoid}.

\begin{proposition}
[Crystalline Frobenius;
cf.~{\cite[\S6.2, especially the discussion following
Definition~6.2.2 and preceding Theorem~6.2.5, and
Example~6.4.16]{KANTORTHESIS}}]
\label{proposition:crystalline_torsor_frobenius}
Let
\[
\underline b,\underline x\in X_k(k),
\]
and let \(\sigma\) denote the Witt-vector Frobenius on \(W(k)\),
extended to
\[
K_0=\operatorname{Frac}W(k).
\]
The canonical \(F\)-structure on \(\mathcal L^{\cri}\), together with
the Frobenius compatibility of the fibre functors
\(\omega_{\underline b}^{\cri}\) and
\(\omega_{\underline x}^{\cri}\), induces isomorphisms
\[
\varphi:
\sigma^*\mathcal G^{\cri}
\xrightarrow{\sim}
\mathcal G^{\cri},
\qquad
\varphi:
\sigma^*{}^{\vphantom{\cri}}_{\underline x}P_{\underline b}^{\cri}
\xrightarrow{\sim}
{}^{\vphantom{\cri}}_{\underline x}P_{\underline b}^{\cri}.
\]
Equivalently, \(\mathcal G^{\cri}\) and
\({}^{\vphantom{\cri}}_{\underline x}P_{\underline b}^{\cri}\) carry
\(\sigma\)-semilinear Frobenius automorphisms. These are compatible
with composition of paths.
\end{proposition}

\begin{remark}
When \(\underline b\) and \(\underline x\) are the reductions of
chosen \(\mathcal O_{K_v}\)-valued lifts \(b\) and \(x\), we may,
where no confusion can arise, write
\[
{}^{\vphantom{\cri}}_{x}P_b^{\cri}
\quad\text{for}\quad
{}^{\vphantom{\cri}}_{\underline x}P_{\underline b}^{\cri}.
\]
\end{remark}

Appendix~\ref{appendix:rmc} proves that the \'{e}tale, crystalline and de Rham
Gauss--Manin objects are associated in the sense of $p$-adic Hodge theory 
(see Definition~\ref{definition:associated_etale_crystalline}), and deduces the
specialized comparison isomorphisms stated below.

\begin{theorem}
[The \'{e}tale--crystalline comparison]
\label{theorem:etale_crystalline_comparison}
Let \(B_{\cris}\) denote Fontaine's crystalline period ring, equipped
with its commuting \(G_v\)-action and Frobenius \(\varphi\). For
endpoints \(b,x \in \mathcal{X}(\mathcal{O}_{K_v})\), there is a natural comparison isomorphism
\[
\iota_{x,b}^{\cris}:
{}^{\vphantom{\et}}_{\overline x}P_{\overline b}^{\et}
 \otimes_{\mathbb Q_p} B_{\cris}
\xrightarrow{\sim}
{}^{\vphantom{\cri}}_{\underline x}P_{\underline b}^{\cri}
 \otimes_{K_0} B_{\cris}.
\]

This isomorphism is \(G_v\)-equivariant, where \(G_v\) acts
diagonally on the left and through its action on \(B_{\cris}\) on the
right. It is also Frobenius-equivariant, where Frobenius acts on the
left through its action on \(B_{\cris}\), and diagonally on the right.

The comparison isomorphisms are compatible with composition of paths.
\end{theorem}

\begin{proof}
See
\S\ref{section:etale_crystalline_comparison_appendix}.
\end{proof}

\begin{theorem}
[The crystalline--de Rham comparison]
\label{theorem:crystalline_de_rham_comparison}
Evaluation on the rigid generic fibre of the $p$-adic formal completion of the fixed
$\mathcal O_{K_v}$-model, after extension of coefficients from
$K_0$ to $K_v$, induces tensor equivalences
\[
\left(
\langle\mathcal L^{\cri}\rangle_\otimes
\right)_{K_v}
\xrightarrow{\sim}
\langle\mathcal L^{\dR}\rangle_\otimes
\]
and
\[
\left(
\langle\mathcal L^{\cri}\rangle_\otimes^{\mathrm{ext}}
\right)_{K_v}
\xrightarrow{\sim}
\langle\mathcal L^{\dR}\rangle_\otimes^{\mathrm{ext}},
\]
compatible with the corresponding fibre functors. Consequently,
there are canonical comparison isomorphisms
\[
\mathcal R^{\cri}\otimes_{K_0}K_v
\xrightarrow{\sim}
\mathcal R^{\dR},
\qquad
\mathcal G^{\cri}\otimes_{K_0}K_v
\xrightarrow{\sim}
\mathcal G^{\dR}.
\]
For $b,x\in\mathcal X(\mathcal O_{K_v})$, there is likewise a
canonical comparison isomorphism
\[
{}^{\vphantom{\cri}}_{\underline x}P_{\underline b}^{\cri}
\otimes_{K_0}K_v
\xrightarrow{\sim}
{}^{\vphantom{\dR}}_{x}P_b^{\dR}.
\]

Under these identifications, the crystalline realization retains
its $K_0$-semilinear Frobenius, while its de Rham realization
carries the Hodge filtration; together they form the corresponding
filtered Frobenius realization.
\end{theorem}

\begin{proof}
See
\S\ref{section:crystalline_de_rham_comparison_appendix}.
\end{proof}
\section{Finite-Type Motivic Quotients}
\label{section:finite_type_motivic_quotients}

The relative completion is generally pro-algebraic, whereas our
applications require compatible finite-type quotients in the
\'{e}tale, crystalline, and de Rham realizations. Unlike in the ordinary unipotent
setting, DCS truncations need not be of finite
type, since the abelianization of the relative completion may be
infinite-dimensional. This section develops a general finiteness mechanism 
and then constructs a distinguished tower with controlled graded pieces.

\begin{itemize}
\item
\S\ref{section:ft_quotients} proves that, for a fixed semisimple
object and a bounded extension length, the corresponding iterated
extensions generate an algebraic Tannakian category. 
An alternative proof using universal extensions is given in
Appendix~\ref{appendix:tannakian_finiteness}.
\item
\S\ref{sec:finite-type_quotients} applies this result to show that
finite-type motivic quotients are cofinal among the finite-type (that is, algebraic)
quotients of the relative completion, which is therefore recovered as
their inverse limit.

\item
\S\ref{sec:param_families} specializes to the Kodaira--Parshin
setting and constructs a canonical tower supported on the adjoint
representation of the connected semisimple monodromy group. This is
the tower used in the applications.
\end{itemize}

\subsection{Finiteness for Bounded Iterated Extensions}
\label{section:ft_quotients}

The key finiteness input for the construction of finite-type motivic
quotients is that bounded iterated self-extensions of a fixed
semisimple object generate an algebraic Tannakian category; see Theorem~\ref{thm:ftg_basic}.

Let \(F\) be a field of characteristic zero, and let
\((\Tc,\omega)\) be a neutral Tannakian category over \(F\). Write
\[
\Gc\coloneqq \pi_1(\Tc,\omega),
\]
let \(\Uc \coloneqq  \operatorname{Rad}_u(\Gc)\) be its pro-unipotent radical, and set
\[
\Rc\coloneqq \Gc/\Uc.
\]
Let \(\Tc_0\subseteq\Tc\) denote the full Tannakian subcategory of
objects on which \(\Uc\) acts trivially; thus
\[
\Tc_0\simeq\operatorname{Rep}_F(\Rc).
\]
All tensor products are taken over \(F\).

\begin{definition}
\label{definition:iterated_exts}
Fix \(\Vc\in\Tc_0\). For \(k\geq 0\), let
\(\IExt^k(\Vc)\) denote the collection of objects \(E\in\Tc\)
admitting a filtration
\[
0=E_{-1}\subset E_0\subset\cdots\subset E_k=E
\]
such that
\[
E_i/E_{i-1}\simeq\Vc
\qquad
(0\leq i\leq k).
\]
We call such an object a \(k\)-fold iterated extension of \(\Vc\) by
itself.
\end{definition}

\begin{definition}\label{def:tvk}
Set
\[
\Tc_{\Vc,k}
\coloneqq 
\left\langle
\IExt^k(\Vc)\cup\Tc_0
\right\rangle_\otimes
\subseteq\Tc,
\]
and write
\[
\Gc_{\Vc,k}\coloneqq \pi_1(\Tc_{\Vc,k},\omega),
\qquad
\Uc_{\Vc,k}\coloneqq \operatorname{Rad}_u(\Gc_{\Vc,k}),
\qquad
\uf_{\Vc,k}\coloneqq \Lie(\Uc_{\Vc,k}).
\]
\end{definition}

If \(k\geq 1\) and \(E_k\in\IExt^k(\Vc)\), a choice of such a
filtration gives exact sequences
\begin{equation}
\label{ses:e_k}
0\longrightarrow E_{k-1}
\longrightarrow E_k
\longrightarrow\Vc
\longrightarrow 0
\end{equation}
and
\begin{equation}
\label{ses:e_k'}
0\longrightarrow\Vc
\longrightarrow E_k
\longrightarrow E'_{k-1}
\longrightarrow 0,
\end{equation}
where \(E_{k-1},E'_{k-1}\in\IExt^{k-1}(\Vc)\). The copy of
\(\Vc\) in \eqref{ses:e_k'} is the first nonzero term of the chosen
filtration and is therefore contained in the subobject \(E_{k-1}\)
appearing in \eqref{ses:e_k}.

\begin{theorem}
\label{thm:ftg_basic}
Assume that \(\Rc\) is of finite type over \(F\) and that
\[
\dim_F\Ext^1_{\Tc}(M,N)<\infty
\]
for all \(M,N\in\Tc\). Then, for every \(k\geq 0\), the Tannakian
category \(\Tc_{\Vc,k}\) is algebraic. Equivalently,
\(\Gc_{\Vc,k}\) is an affine algebraic \(F\)-group scheme.
\end{theorem}

We prove Theorem~\ref{thm:ftg_basic} by controlling the
abelianization and nilpotence class of \(\Uc_{\Vc,k}\). An alternative
proof using universal extensions is given in
\S\ref{sec:proof_using_univ_ext}.

\subsubsection{Proof via Abelianizations}
\label{sec:proof_using_abel}

We first bound the nilpotence class of the unipotent radical.

\begin{lemma}
\label{lemma:T_k,V_kth_DCS}
For every \(k\geq 0\), the subgroup
\[
\Fc_{\DCS}^{k+1}\Uc
\]
acts trivially on every object of \(\IExt^k(\Vc)\). Consequently,
\[
\Fc_{\DCS}^{k+1}\Uc_{\Vc,k}=1.
\]
\end{lemma}

\begin{proof}
We argue by induction on \(k\). The case \(k=0\) follows from
\(\Vc\in\Tc_0\). Let \(k\geq 1\), and suppose the assertion is known
in degree \(k-1\). Fix \(E_k\in\IExt^k(\Vc)\), together with the exact
sequences
\[
0\longrightarrow E_{k-1}\longrightarrow E_k
\longrightarrow\Vc\longrightarrow0
\]
and
\[
0\longrightarrow\Vc\longrightarrow E_k
\longrightarrow E'_{k-1}\longrightarrow0
\]
from \eqref{ses:e_k} and \eqref{ses:e_k'}, respectively.

Let \(A\) be an \(F\)-algebra, \(a\in\Fc_{\DCS}^{k}\Uc(A)\), and
\(b\in\Uc(A)\). Since \(\Uc\) acts trivially on \(\Vc\), the operator
\(b-1\) sends \(E_{k,A}\) into \(E_{k-1,A}\). By induction, \(a-1\)
vanishes on \(E_{k-1,A}\), and hence
\[
(a-1)(b-1)=0
\qquad\text{on }E_{k,A}.
\]
Similarly, \(a-1\) sends \(E_{k,A}\) into \(\Vc_A\), because it acts
trivially on \(E'_{k-1,A}\), while \(b-1\) vanishes on \(\Vc_A\).
Thus
\[
(b-1)(a-1)=0
\qquad\text{on }E_{k,A}.
\]
Since $ab = a + b - 1 = ba$ on $E_{k,A}$, \(ab\) and \(ba\) act identically on \(E_{k,A}\), so
\([a,b]\) acts trivially. Since this holds functorially in \(A\),
\[
[\Fc_{\DCS}^{k}\Uc,\Uc]
=
\Fc_{\DCS}^{k+1}\Uc
\]
acts trivially on \(E_k\). The final assertion follows because
\(\Tc_{\Vc,k}\) is generated by \(\Tc_0\) and the objects of
\(\IExt^k(\Vc)\).
\end{proof}

Write
\[
\uf\coloneqq \Lie(\Uc).
\]
We use the following standard description of its abelianization.

\begin{fact}
[{\cite[(3.4)]{Hain_HodgeDeRham_Chapter2016};
cf.~\cite[Proposition~3.1.1]{ESKANDARI_B},
\cite[Proposition~6.1]{BrownNotes17}}]
\label{fact:u_ab_revised}
For every \(M\in\Tc_0\), there is a natural isomorphism
\[
\Hom_{\Rc}\bigl(\uf^{\ab},\omega(M)\bigr)
\simeq
\Ext^1_{\Tc}(\mathbf 1,M).
\]
In particular, when the simple objects of \(\Tc_0\) are absolutely
simple\footnote{This is always the case when the reductive quotient is connected split of finite type over $F$.}, there is a canonical \(\Rc\)-equivariant isomorphism
\[
\uf^{\ab}
\simeq
\prod_{M\ \mathrm{simple}}
\omega(M)\otimes_F
\Ext^1_{\Tc}(\mathbf 1,M)^\vee.
\]
\end{fact}

\begin{definition}
\label{defn:n(V)}
Let \(\Sigma_{\Vc}\) be the finite set of isomorphism classes of
simple constituents of
\[
\Vc\otimes\Vc^\vee.
\]
Let
\[
A_{\Vc}
\]
be the maximal quotient of the pro-\(\Rc\)-representation
\(\uf^{\ab}\) whose simple constituents belong to
\(\Sigma_{\Vc}\), and set
\[
n(\Vc)\coloneqq \dim_F A_{\Vc}.
\]
\end{definition}

The vector space \(A_{\Vc}\) is finite-dimensional: the set
\(\Sigma_{\Vc}\) is finite, and Fact~\ref{fact:u_ab_revised}, together
with the hypotheses of Theorem~\ref{thm:ftg_basic}, shows that each
of its isotypic multiplicity spaces is finite-dimensional. When the
relevant simple objects are absolutely simple,
\[
n(\Vc)
=
\sum_{S\in\Sigma_{\Vc}}
\dim_F\omega(S)\,
\dim_F\Ext^1_{\Tc}(\mathbf 1,S).
\]

For a Tannakian subcategory
\(\Tc'\subseteq\Tc\) containing \(\Tc_0\), write
\[
\Gc(\Tc')\coloneqq \pi_1(\Tc',\omega|_{\Tc'}),
\qquad
\Uc(\Tc')\coloneqq \operatorname{Rad}_u\Gc(\Tc'),
\qquad
\uf(\Tc')\coloneqq \Lie\Uc(\Tc').
\]

\begin{definition}
\label{defn:property_V}
We say that \(\Tc'\) satisfies \emph{Property~\(\mathbf{V}\)} if the natural
surjection
\[
\uf^{\ab}\longrightarrow\uf(\Tc')^{\ab}
\]
factors through \(A_{\Vc}\). Equivalently, every simple constituent
of \(\uf(\Tc')^{\ab}\) belongs to \(\Sigma_{\Vc}\).
\end{definition}

In particular, if \(\Tc'\) satisfies Property~\(\mathbf{V}\), then
\begin{equation}
\label{eq:property_V_dimension}
\dim_F\uf(\Tc')^{\ab}\leq n(\Vc).
\end{equation}

We next show that Property~\(\mathbf{V}\) is preserved when one adjoins a single
iterated extension. We use the following result on the fundamental
group of an extension.

\begin{fact}[{The Tannakian kernel of an extension;
cf.~\cite[\S3, especially Theorem~3.3.1]{ESKANDARI_A}}]
\label{fact:eskandari_murty}
Let
\[
0\longrightarrow L\longrightarrow M
\longrightarrow\mathbf 1\longrightarrow0
\]
be exact, and suppose that
\[
\Tc'\coloneqq \langle L\rangle_\otimes
\]
contains \(\Tc_0\). Set
\[
\Tc''\coloneqq \langle M\rangle_\otimes.
\]
The canonical morphism
\[
\Gc(\Tc'')\longrightarrow\Gc(\Tc')
\]
has commutative unipotent kernel \(\Uc_\eta\) and induces an exact
sequence
\[
1\longrightarrow\Uc_\eta
\longrightarrow\Uc(\Tc'')
\longrightarrow\Uc(\Tc')
\longrightarrow1.
\]
Writing
\[
\uf_\eta\coloneqq \Lie(\Uc_\eta),
\]
conjugation makes \(\uf_\eta\) an object
\(\underline{\uf_\eta}\in\Tc'\). Its natural action on \(M\)
identifies \(\underline{\uf_\eta}\) with the smallest subobject
\(L_0\subseteq L\) such that
\[
M/L_0\in\Tc'.
\]
\end{fact}

\begin{lemma}
\label{lemma:induction_ftg}
Let \(k\geq 1\), and suppose that \(\Tc'\subseteq\Tc\) is an
algebraic Tannakian subcategory containing \(\Tc_{\Vc,k-1}\) and
satisfying Property~\(\mathbf{V}\). For every
\(E_k\in\IExt^k(\Vc)\), the category
\[
\Tc''
\coloneqq 
\left\langle\Tc'\cup\{E_k\}\right\rangle_\otimes
\]
is algebraic and satisfies Property~\(\mathbf{V}\).
\end{lemma}

\begin{proof}
Choose a tensor generator \(E'\) of \(\Tc'\). Tensor
\eqref{ses:e_k} with \(\Vc^\vee\) and pull back along the coevaluation
map
\[
\mathbf 1\longrightarrow\Vc\otimes\Vc^\vee.
\]
This gives an exact sequence
\begin{equation}
\label{eq:ext_k_tilde}
0\longrightarrow E_{k-1}\otimes\Vc^\vee
\longrightarrow\widetilde E_k
\longrightarrow\mathbf 1
\longrightarrow0.
\end{equation}

The categories generated by \(E_k\) and \(\widetilde E_k\) over
\(\Tc'\) coincide. Indeed,
\(\widetilde E_k\) is a subobject of \(E_k\otimes\Vc^\vee\), so
\[
\widetilde E_k\in\langle\Tc',E_k\rangle_\otimes.
\]
Conversely, tensoring \eqref{eq:ext_k_tilde} with \(\Vc\) and pushing
out along the evaluation map
\[
E_{k-1}\otimes\Vc^\vee\otimes\Vc
\longrightarrow E_{k-1}
\]
recovers the extension \eqref{ses:e_k}. Hence
\[
E_k\in\langle\Tc',\widetilde E_k\rangle_\otimes.
\]

Set
\[
L\coloneqq (E_{k-1}\otimes\Vc^\vee)\oplus E',
\qquad
M\coloneqq \widetilde E_k\oplus E'.
\]
Then
\[
0\longrightarrow L\longrightarrow M
\longrightarrow\mathbf 1\longrightarrow0
\]
is exact, and
\[
\langle L\rangle_\otimes=\Tc',
\qquad
\langle M\rangle_\otimes=\Tc''.
\]
In particular, \(\Tc''\) is algebraic.

Let
\[
\Uc_\eta
=
\ker\bigl(\Gc(\Tc'')\longrightarrow\Gc(\Tc')\bigr),
\qquad
\uf_\eta\coloneqq \Lie(\Uc_\eta).
\]
Both \(\Tc'\) and \(\Tc''\) have reductive quotient \(\Rc\), and hence
there is an exact sequence
\[
1\longrightarrow\Uc_\eta
\longrightarrow\Uc(\Tc'')
\longrightarrow\Uc(\Tc')
\longrightarrow1.
\]
Passing to Lie algebras and then to abelianizations gives a right
exact sequence of \(\Rc\)-representations
\begin{equation}
\label{eq:abelianization_induction}
\uf_\eta
\longrightarrow
\uf(\Tc'')^{\ab}
\longrightarrow
\uf(\Tc')^{\ab}
\longrightarrow0.
\end{equation}

It remains to control the simple constituents of \(\uf_\eta\).
The chosen filtration of \(E_k\) gives an inclusion
\[
\Vc\hookrightarrow E_{k-1},
\]
and hence
\[
\Vc\otimes\Vc^\vee
\hookrightarrow
E_{k-1}\otimes\Vc^\vee
\hookrightarrow L.
\]
Moreover,
\[
\widetilde E_k/(\Vc\otimes\Vc^\vee)
\hookrightarrow
(E_k\otimes\Vc^\vee)/(\Vc\otimes\Vc^\vee)
\simeq
E'_{k-1}\otimes\Vc^\vee.
\]
The object on the right belongs to \(\Tc'\), since
\[
E'_{k-1},\Vc\in\Tc_{\Vc,k-1}\subseteq\Tc'.
\]
It follows that
\[
M/(\Vc\otimes\Vc^\vee)\in\Tc'.
\]
Fact~\ref{fact:eskandari_murty} therefore gives
\[
\underline{\uf_\eta}
\subseteq
\Vc\otimes\Vc^\vee.
\]
Thus every simple constituent of \(\uf_\eta\) belongs to
\(\Sigma_{\Vc}\).

By hypothesis, the same is true of
\(\uf(\Tc')^{\ab}\). Since \(\Rep_F(\Rc)\) is semisimple,
\eqref{eq:abelianization_induction} implies that every simple
constituent of \(\uf(\Tc'')^{\ab}\) belongs to \(\Sigma_{\Vc}\).
Hence \(\Tc''\) satisfies Property~\(\mathbf{V}\).
\end{proof}

\begin{proof}[Proof of Theorem~\ref{thm:ftg_basic}]
We prove by induction on \(k\) the stronger assertion that
\(\Tc_{\Vc,k}\) is algebraic and satisfies Property~\(\mathbf{V}\).

For \(k=0\), we have
\[
\Tc_{\Vc,0}=\Tc_0,
\]
which is algebraic because \(\Rc\) is of finite type, and its
unipotent radical is trivial.

Suppose \(k\geq1\), and assume that
\(\Tc_{\Vc,k-1}\) is algebraic and satisfies Property~\(\mathbf{V}\). For a finite
subset
\[
T\subseteq\IExt^k(\Vc),
\]
set
\[
\Tc_T
\coloneqq 
\left\langle
\Tc_{\Vc,k-1}\cup T
\right\rangle_\otimes,
\qquad
\Uc_T\coloneqq \Uc(\Tc_T).
\]
Repeated application of Lemma~\ref{lemma:induction_ftg} shows that
\(\Tc_T\) is algebraic and satisfies Property~\(\mathbf{V}\). Therefore
\[
\dim_F\Lie(\Uc_T)^{\ab}
\leq n(\Vc).
\]
By Lemma~\ref{lemma:T_k,V_kth_DCS},
\[
\Fc_{\DCS}^{k+1}\Uc_T=1.
\]
Corollary~\ref{cor:prounip_dim_bound} now gives a uniform finite bound
\begin{equation}
\label{eq:finite_T_dimension_bound}
\dim\Uc_T\leq d_{n(\Vc)}(k).
\end{equation}

As \(T\) ranges over the finite subsets of \(\IExt^k(\Vc)\), the
groups \(\Uc_T\) form a cofiltered inverse system with faithfully flat
transition maps, and
\[
\Uc_{\Vc,k}\simeq\varprojlim_T\Uc_T.
\]
The dimensions of the terms are uniformly bounded by
\eqref{eq:finite_T_dimension_bound}. Hence
Lemma~\ref{lemma:unipotent_bounded_dim} implies that
\(\Uc_{\Vc,k}\) is algebraic and
\[
\dim\Uc_{\Vc,k}\leq d_{n(\Vc)}(k).
\]
Since \(\Gc_{\Vc,k}\) is an extension of the algebraic group \(\Rc\)
by \(\Uc_{\Vc,k}\), it is algebraic. Equivalently,
\(\Tc_{\Vc,k}\) is algebraic.

Finally, \(\Tc_{\Vc,k}\) is generated by the directed union of the
subcategories \(\Tc_T\). By
Lemma~\ref{lemma:tann_noetherian}, there is a finite subset
\(T_0\subseteq\IExt^k(\Vc)\) such that
\[
\Tc_{\Vc,k}=\Tc_{T_0}.
\]
Since \(\Tc_{T_0}\) satisfies Property~\(\mathbf{V}\), so does
\(\Tc_{\Vc,k}\). This completes the induction and proves the theorem.
\end{proof}

\subsection{Cofinality of Finite-Type Motivic Quotients}
\label{sec:finite-type_quotients}

We now show that the relative completion is recovered from its
finite-type motivic quotients. Throughout this subsection, we work in the global geometric
\(p\)-adic \'{e}tale realization and suppress the superscript
\(\et\). Set
\[
\Tc
\coloneqq 
\langle\mathcal L\rangle_\otimes^{\mathrm{ext}},
\qquad
\Tc_0
\coloneqq 
\langle\mathcal L\rangle_\otimes,
\qquad
\Gc
\coloneqq 
\pi_1(\Tc,\omega_b).
\]
Thus
\[
1\longrightarrow\Uc
\longrightarrow\Gc
\longrightarrow\Rc
\longrightarrow1
\]
is the relative completion associated with \(\mathcal L\), which is a semisimple object, \(\mathcal L\in \Tc_0\)
equipped with a \(G_T\)-action. The
natural \(G_T\)-action on \(\Tc\) is by exact tensor
autoequivalences and preserves \(\Tc_0\) (as it stabilizes $\mathcal L$).

An algebraic quotient of \(\Gc\) is \(G_T\)-equivariant precisely
when its corresponding Tannakian subcategory of \(\Tc\) is
\(G_T\)-stable. By the functorial comparison results of
Appendix~\ref{appendix:rmc}, such quotients admit compatible
crystalline and de Rham realizations, equipped respectively with
Frobenius and Hodge filtration. We refer to these compatible
quotients as \emph{motivic}.

\begin{theorem}[Cofinality of finite-type motivic quotients]
\label{thm:inv_limit_of_ft_revised}
Let \(\mathscr M\) denote the collection of \(G_T\)-stable algebraic
Tannakian subcategories
\[
\Tc_0\subseteq\mathsf{S}\subseteq\Tc.
\]
Then \(\mathscr M\) is filtered and cofinal among the algebraic
Tannakian subcategories of \(\Tc\). Consequently,
\[
\Gc
\xrightarrow{\sim}
\varprojlim_{\mathsf{S}\in\mathscr M}
\Gc_{\mathsf{S}},
\qquad
\Gc_{\mathsf{S}}
\coloneqq 
\pi_1(\mathsf{S},\omega_b).
\]
\end{theorem}

\begin{proof}
We first verify the finiteness hypothesis of
Theorem~\ref{thm:ftg_basic}. Since \(\Tc\) is a full subcategory of
the category of \(p\)-adic \'{e}tale local systems on \(X^\dagger_{\overline{K}}\)
and is closed under extensions, for \(M,N\in\Tc\) there is a natural
isomorphism
\[
\Ext^1_{\Tc}(M,N)
\simeq
H^1_{\et}
\bigl(
X^\dagger_{\overline{K}},
M^\vee\otimes N
\bigr).
\]
The group on the right is finite-dimensional by standard theorems in \'etale cohomology. Hence
Theorem~\ref{thm:ftg_basic} applies to every object of \(\Tc_0\).

Let \(\Ac\subseteq\Tc\) be an algebraic Tannakian subcategory. After
replacing \(\Ac\) by
\[
\langle\Ac\cup\Tc_0\rangle_\otimes,
\]
which is again algebraic, we may assume that
\(\Tc_0\subseteq\Ac\). Choose a tensor generator \(E\) of \(\Ac\).
Since \(E\in\Tc\), it admits a finite filtration
\[
0=F_{-1}E\subset F_0E\subset\cdots\subset F_kE=E
\]
whose graded pieces
\[
W_i\coloneqq F_iE/F_{i-1}E
\]
belong to \(\Tc_0\).

For each \(i\), choose a finite direct sum \(P_i\) of tensor
products of copies of \(\mathcal L\) and \(\mathcal L^\vee\) that
contains \(W_i\) as a direct summand. Such a choice is possible
because \(\Tc_0\) is semisimple and tensor-generated by
\(\mathcal L\); moreover, \(P_i\) is automatically \(G_T\)-stable.
Set
\[
\Vc\coloneqq \bigoplus_{i=0}^k P_i.
\]
Then \(\Vc\in\Tc_0\) is \(G_T\)-stable, and every \(W_i\) is a
direct summand of \(\Vc\). For each $i$, choose \(C_i\in\Tc_0\) such that
\[
\Vc\simeq W_i\oplus C_i.
\]

Define
\[
E'
\coloneqq 
E\oplus\bigoplus_{i=0}^k C_i
\]
and equip it with the filtration
\[
F_iE'
\coloneqq 
F_iE\oplus\bigoplus_{j=0}^i C_j.
\]
Its \(i\)-th graded piece is
\[
\Gr_i^F E'
\simeq
W_i\oplus C_i
\simeq
\Vc.
\]
Thus
\[
E'\in\IExt^k(\Vc).
\]
Moreover, \(E'\in\Ac\), while \(E\) is a direct summand of \(E'\);
hence
\[
\Ac=\langle E'\rangle_\otimes.
\]

Set
\[
\Ac^{\mathrm{mot}}
\coloneqq 
\Tc_{\Vc,k}.
\]
Then
\[
\Ac\subseteq\Ac^{\mathrm{mot}}\subseteq\Tc.
\]
The category \(\Ac^{\mathrm{mot}}\) is algebraic by
Theorem~\ref{thm:ftg_basic}. It is also \(G_T\)-stable: the
\(G_T\)-action preserves \(\Tc_0\), and, since \(\Vc\) is
\(G_T\)-stable, it preserves the collection
\(\IExt^k(\Vc)\). Therefore
\[
\Ac^{\mathrm{mot}}\in\mathscr M.
\]
This proves that \(\mathscr M\) is cofinal among all algebraic
Tannakian subcategories of \(\Tc\).

Finally, if \(\mathsf{S}_1,\mathsf{S}_2\in\mathscr M\), then
\[
\langle\mathsf{S}_1\cup\mathsf{S}_2\rangle_\otimes
\]
is \(G_T\)-stable and algebraic. It is \(G_T\)-stable because both
\(\mathsf S_1\) and \(\mathsf S_2\) are \(G_T\)-stable. It is
algebraic: if \(E_i\) is any tensor generator of \(\mathsf S_i\),
then \(E_1\oplus E_2\) tensor-generates
\(\langle\mathsf S_1\cup\mathsf S_2\rangle_\otimes\).

Thus \(\mathscr M\) is filtered. Since every object of \(\Tc\) lies
in an algebraic Tannakian subcategory,
\[
\Tc=\bigcup_{\mathsf{S}\in\mathscr M}\mathsf{S}.
\]
Tannakian duality therefore gives
\[
\Gc
\simeq
\varprojlim_{\mathsf{S}\in\mathscr M}\Gc_{\mathsf{S}},
\]
as claimed.
\end{proof}

\subsection{An Adjoint-Supported Tower of Finite-Type Motivic Quotients}
\label{sec:param_families}
Theorem~\ref{thm:inv_limit_of_ft_revised} produces a cofinal family
of finite-type motivic quotients of the relative completion, but gives
little control over the \(\mathcal R\)-representation structure of their graded pieces. Such
control is needed in our applications, notably in the analysis of
complex conjugation on the DCS graded pieces.
One could use instead the finite-level truncations of the ordinary
unipotent completion, viewed as motivic quotients of the relative
completion, but this would recover Kim's unipotent Chabauty theory.
Our aim is to exploit genuinely relative quotients that retain
nontrivial reductive monodromy and thereby give rise to new
rigid-analytic functions cutting out \(S\)-integral points. 
We therefore construct a canonical tower supported on the adjoint
representation of the connected monodromy group. 
Remark~\ref{rem:why_ad_supp} below motivates our choice to work with
the adjoint representation.

We work in the geometric \'{e}tale realization fixed in
\S\ref{section:rmc_realizations}, equipped with its global
\(G_T\)-action from
Proposition~\ref{proposition:etale_torsor_galois}. Set
\[
\Tc_0\coloneqq \langle\mathcal L^{\et}\rangle_\otimes,
\qquad
\Tc\coloneqq \langle\mathcal L^{\et}\rangle_\otimes^{\mathrm{ext}},
\]
and write
\[
1\longrightarrow\Uc^{\et}
\longrightarrow\Gc^{\et}
\longrightarrow\Rc^{\et}
\longrightarrow1
\]
for the associated relative completion. We suppress the superscript
\(\et\) throughout this subsection.

Full monodromy identifies \(\Rc^\circ\) with a product of symplectic
groups. In particular,
\[
\Rc^\circ=(\Rc^\circ)^{\mathrm{der}}
\]
is semisimple. Set
\[
\uf\coloneqq \Lie(\Uc),
\qquad
\mathfrak a\coloneqq \Lie(\Rc^\circ).
\]
Since \(\Rc^\circ\) is characteristic in \(\Rc\), conjugation makes
\(\mathfrak a\) canonically an \(\Rc\)-representation.

\begin{definition}
\label{definition:adjoint_support}
Let \(\Irr(\Rc)\) denote the set of isomorphism classes of simple
objects of \(\operatorname{Rep}_{\mathbb Q_p}(\Rc)\), and let
\[
\Sigma_{\mathrm{ad}}\subseteq\Irr(\Rc)
\]
be the set of simple constituents of \(\mathfrak a\). An
\(\Rc\)-representation is \emph{adjoint-supported} if all its simple
constituents belong to \(\Sigma_{\mathrm{ad}}\). 
For an \(\Rc\)-representation or pro-representation \(W\), write
\[
Q_{\mathrm{ad}}(W)
\]
for its maximal adjoint-supported quotient; see
\S\ref{appendix:support_projectors} and
Proposition~\ref{proposition:support_projector}.

A quotient extension
\[
\Gc\twoheadrightarrow\Hc\longrightarrow\Rc
\]
is \emph{adjoint-supported} if, writing
\[
\vf\coloneqq \Lie\ker(\Hc\to\Rc),
\]
each \(\Gr_n^{\DCS}\vf\) is adjoint-supported.
\end{definition}
This is
Definition~\ref{def:supported_quotient_extension} specialized to
\(\Sigma=\Sigma_{\mathrm{ad}}\).

\begin{remark}\label{rem:why_ad_supp}
The choice of adjoint support is distinguished by a nonvanishing mechanism. The $\Rc$-equivariant Lie algebra $\af$ is nonzero and perfect, and the adjoint-supported abelianization constructed below surjects onto it. Perfectness then forces every bracket degree of the maximal adjoint-supported quotient of the free Lie algebra on this abelianization to surject onto $\af$, hence be nonzero; see Proposition~\ref{proposition:perfect_target_all_degrees}.
\end{remark}

\begin{proposition}
\label{proposition:maximal_adjoint_supported_quotient}
There is a canonical maximal adjoint-supported quotient extension
\[
1\longrightarrow\Uc_{\mathrm{ad}}
\longrightarrow\Gc_{\mathrm{ad}}
\longrightarrow\Rc
\longrightarrow1
\]
of \(\Gc\), through which every adjoint-supported quotient extension
of \(\Gc\) over \(\Rc\) factors uniquely.

Set
\[
\qf_{\mathrm{ad}}\coloneqq \Lie(\Uc_{\mathrm{ad}}),
\qquad
A_{\mathrm{ad}}\coloneqq Q_{\mathrm{ad}}(\uf^{\ab}).
\]
If \(A_{\mathrm{ad}}\) is finite-dimensional, then every truncation
\[
\Gc_{\mathrm{ad},n}
\coloneqq 
\Gc_{\mathrm{ad}}/
\Fc_{\DCS}^{n+1}\Uc_{\mathrm{ad}}
\]
is of finite type. The construction is functorial under automorphisms
of the extension \(\Gc\to\Rc\) that preserve
\(\Sigma_{\mathrm{ad}}\).
\end{proposition}
The proof is given in
Appendix~\ref{appendix:adjoint_supported_quotients}.
Proposition~\ref{proposition:maximal_supported_pronilpotent_quotient}
constructs the maximal supported pronilpotent quotient and proves its
compatibility with truncation;
Proposition~\ref{proposition:supported_abelianization_and_graded} and
Corollary~\ref{corollary:finite_supported_truncations} give the
finite-type assertion; and
Proposition~\ref{proposition:maximal_supported_extension_appendix},
together with
Remark~\ref{remark:adjoint_support_is_canonical}, gives functoriality,
including \(G_T\)-equivariance.

We shall use the following presentation of the relative-completion
Lie algebra in the projective case. For a vector space or pro-vector space \(V\), write
\(\widehat{\mathbb L}(V)\) for the free pronilpotent Lie algebra on
\(V\), equivalently the DCS completion of the
free Lie algebra on $V$.

\begin{fact}
[Quadratic presentation;
{\cite[Theorem~13.14]{HAIN_HDR_RMC};
cf.~{\cite[Proposition~2.2 and Remark~2.3]
{Hain_HodgeDeRham_Chapter2016}}}]
\label{fact:quadratic_presentation_adjoint}
Assume that \(X\) is smooth and projective. Let
\(\mathcal O(\Rc)\) denote the local system associated with the
regular representation of \(\Rc\). Under the canonical identification
\[
\uf^{\ab}
\simeq
H^1_{\et}\bigl(X_{\overline{K}},\mathcal O(\Rc)\bigr)^\vee,
\]
the dual of the cup product induces an \(\Rc\)-equivariant map
\[
\check\cup:
H^2_{\et}\bigl(X_{\overline{K}},\mathcal O(\Rc)\bigr)^\vee
\longrightarrow
\Gr_2^{\DCS}\widehat{\mathbb L}(\uf^{\ab}).
\]
There is an \(\Rc\)-equivariant minimal presentation
\[
p:
\widehat{\mathbb L}(\uf^{\ab})
\twoheadrightarrow
\uf
\]
such that
\[
\ker(p)
=
\overline{
\left\langle
\operatorname{im}(\check\cup)
\right\rangle_{\mathrm{Lie}}
}.
\]
For a smooth projective curve, \(\operatorname{im}(\check\cup)\) is a
one-dimensional \(\Rc\)-trivial quadratic relation space.
\end{fact}

\begin{theorem}[The adjoint-supported Kodaira--Parshin tower]
\label{theorem:special_motivic_exists}
Let
\[
Y\longrightarrow X'\longrightarrow X
\]
be a Kodaira--Parshin family over a smooth projective curve \(X/K\)
of genus \(g\geq2\). Then the truncations
\[
\Gc_{\mathrm{ad},n}
=
\Gc_{\mathrm{ad}}/
\Fc_{\DCS}^{n+1}\Uc_{\mathrm{ad}},
\qquad n\geq1,
\]
form an inverse system
\[
\Gc
\twoheadrightarrow\cdots\twoheadrightarrow
\Gc_{\mathrm{ad},n+1}
\twoheadrightarrow
\Gc_{\mathrm{ad},n}
\twoheadrightarrow\cdots\twoheadrightarrow
\Gc_{\mathrm{ad},1}
\]
of finite-type \(G_T\)-equivariant quotients. For every \(n\geq1\), there is an \(\Rc\)-equivariant surjection
\[
\Gr_n^{\DCS}\qf_{\mathrm{ad}}
\twoheadrightarrow
\mathfrak a.
\]
In particular, every DCS graded piece is
nonzero, and every transition map in the displayed inverse system, while faithfully flat, is not an isomorphism.

Under the comparison isomorphisms of
Theorems~\ref{theorem:etale_crystalline_comparison}
and~\ref{theorem:crystalline_de_rham_comparison}, these quotients
admit compatible crystalline and de Rham realizations. They therefore
form an infinite tower of finite-type motivic quotients.
\end{theorem}

\begin{proof}
Full monodromy identifies the connected semisimple monodromy group
with a product of symplectic groups. Consequently, its Lie algebra
\[
\mathfrak a
=
\Lie\bigl((\Rc^\circ)^{\mathrm{der}}\bigr)
\]
is a finite multiplicity-free sum
\[
\mathfrak a
=
\bigoplus_{\lambda\in\Lambda}S_\lambda
\]
of nontrivial simple \(\Rc\)-representations. Let
\(\mathbb S_\lambda\in\Tc_0\) denote the local system corresponding
to \(S_\lambda\).

The standard description of the abelianization of the relative
completion gives (see Fact~\ref{fact:u_ab_revised})
\[
\Hom_{\Rc}\bigl(\uf^{\ab},S_\lambda\bigr)
\simeq
H^1_{\et}\bigl(X_{\overline{K}},\mathbb S_\lambda\bigr).
\]
Since \(\mathbb S_\lambda\) is nontrivial and irreducible,
Poincar\'e duality and the Euler characteristic formula give
\[
H^0_{\et}(X_{\overline{K}},\mathbb S_\lambda)
=
H^2_{\et}(X_{\overline{K}},\mathbb S_\lambda)
=
0
\]
and
\[
\dim_{\mathbb Q_p}
H^1_{\et}\bigl(X_{\overline{K}},\mathbb S_\lambda\bigr)
=
(2g-2)\dim_{\mathbb Q_p}S_\lambda
>0.
\]
It follows that the adjoint-supported quotient
\[
A_{\mathrm{ad}}
\coloneqq 
Q_{\mathrm{ad}}(\uf^{\ab})
\]
is finite-dimensional and admits an \(\Rc\)-equivariant surjection
\[
A_{\mathrm{ad}}\twoheadrightarrow\mathfrak a.
\]

Set
\[
\mathfrak f_{\mathrm{ad}}
\coloneqq 
\mathfrak q_{\Sigma_{\mathrm{ad}}}
\bigl(\widehat{\mathbb L}(A_{\mathrm{ad}})\bigr).
\]
Since \(\mathfrak a\) is perfect,
Proposition~\ref{proposition:perfect_target_all_degrees} gives
surjections
\[
\Gr_n^{\DCS}\mathfrak f_{\mathrm{ad}}
\twoheadrightarrow
\mathfrak a
\qquad(n\geq1).
\]

By Fact~\ref{fact:quadratic_presentation_adjoint}, the kernel of the
minimal presentation
\[
\widehat{\mathbb L}(\uf^{\ab})
\twoheadrightarrow
\uf
\]
is generated by an \(\Rc\)-trivial quadratic subrepresentation.
Since the tensor unit is not adjoint-supported,
Lemma~\ref{lemma:unsupported_relations_vanish} shows that the
canonical map
\[
\widehat{\mathbb L}(\uf^{\ab})
\longrightarrow
\mathfrak f_{\mathrm{ad}}
\]
factors through a surjection
\[
\uf\twoheadrightarrow\mathfrak f_{\mathrm{ad}}.
\]
The universal property of the maximal adjoint-supported quotient
therefore gives a surjection
\[
\qf_{\mathrm{ad}}
\twoheadrightarrow
\mathfrak f_{\mathrm{ad}},
\]
and hence
\[
\Gr_n^{\DCS}\qf_{\mathrm{ad}}
\twoheadrightarrow
\mathfrak a
\qquad(n\geq1).
\]
In particular, every graded piece is nonzero, so each transition map
\[
\Gc_{\mathrm{ad},n+1}
\longrightarrow
\Gc_{\mathrm{ad},n}
\]
is nontrivial.

Since \(A_{\mathrm{ad}}\) is finite-dimensional,
Proposition~\ref{proposition:maximal_adjoint_supported_quotient}
implies that every \(\Gc_{\mathrm{ad},n}\) is of finite type. The
\(G_T\)-action preserves
\[
(\Rc^\circ)^{\mathrm{der}},
\qquad
\mathfrak a,
\qquad
\Sigma_{\mathrm{ad}},
\]
so the same proposition shows that the quotient
\(\Gc\twoheadrightarrow\Gc_{\mathrm{ad}}\), and hence every
truncation, is \(G_T\)-equivariant.
\end{proof}

\begin{lemma}[Purity of the adjoint-supported graded pieces]
\label{lemma:purity_adjoint_supported_gradeds}
For every \(k\geq1\), the \(G_T\)-representation
\[
\Gr_k^{\DCS}\qf_{\mathrm{ad}}
\]
is pure of weight \(-k\). Under the comparison isomorphisms, its
\'{e}tale, crystalline, and de Rham realizations form a compatible
system of weight \(-k\).
\end{lemma}

\begin{proof}
The adjoint local systems are pure of weight \(0\). Proper purity,
Fact~\ref{fact:u_ab_revised}, and
Proposition~\ref{proposition:maximal_adjoint_supported_quotient}
therefore show that
\[
\Gr_1^{\DCS}\qf_{\mathrm{ad}}
=
\qf_{\mathrm{ad}}^{\ab}
=
A_{\mathrm{ad}}
\]
is pure of weight \(-1\). The iterated Lie bracket induces a
surjection
\[
A_{\mathrm{ad}}^{\otimes k}
\twoheadrightarrow
\Gr_k^{\DCS}\qf_{\mathrm{ad}},
\]
so the latter is pure of weight \(-k\). Compatibility of the
realizations follows from the comparison isomorphisms.
\end{proof}
﻿\section{Neutral-Fibre Selmer Stacks and Relative Kummer Maps}
\label{chapter:selmer_stacks}
This section constructs the neutral-fibre Selmer stacks and relative
Kummer maps used in our applications.

\begin{itemize}
\item
\S\ref{section:galois_coh_stacks} recalls continuous non-abelian
cohomology stacks and defines the Bloch--Kato, unramified, and
restricted-ramification local conditions on a general Grothendieck
site.

\item
\S\ref{section:neutral_fibre_selmer} specializes these constructions
to relative completions. It defines the neutral and trivialized
neutral fibres, establishes their quotient presentations, and
constructs the relative Kummer maps and reductive Galois strata.

\item
\S\ref{section:finite_descent_kummer} uses finite descent and
inertia-invariant paths to construct, in the genus $\ge 2$ setting, a finite Kummer stratification on which the
global relative Kummer map factors through the restricted
neutral-fibre Selmer stack.

\item
\S\ref{section:neutral_selmer_representability} proves the
finite-level algebraicity of these targets. The trivialized neutral
Selmer stacks are algebraic spaces of finite type, while the neutral
targets are finite-type Artin stacks with affine and explicitly
described inertia.
\end{itemize}

\subsection{Continuous Non-Abelian Cohomology and Selmer Conditions}
\label{section:galois_coh_stacks}

We recall the construction of continuous non-abelian cohomology
stacks and formulate the local conditions used below. Our conventions
follow \cite[Chapter~3]{KANTORTHESIS}.

Let \((\mathcal C,\tau)\) be a Grothendieck site, let \(\Gamma\) be a
profinite group, and let
\[
\mathcal A:
\mathcal C^{\op}\longrightarrow\mathbf{TopGrp}
\]
be a group-valued \(\tau\)-sheaf. A \emph{continuous
\(\Gamma\)-action} on \(\mathcal A\) is a natural family of actions by
group automorphisms such that
\[
\Gamma\times\mathcal A(S)\longrightarrow\mathcal A(S)
\]
is continuous for every \(S\in\mathcal C\).

In our applications, \((\mathcal C,\tau)\) is either the category of
\(\Qp\)-schemes with the \'{e}tale topology or the category of rigid
analytic spaces over \(\Spm(\Qp)\) with the Tate topology. For an
affine algebraic test object \(S=\Spec(A)\), we give \(A\) the
inductive-limit topology over its finite-dimensional
\(\Qp\)-subspaces. For an affinoid test object
\(S=\Spm(\Lambda)\), we use the canonical affinoid topology on
\(\Lambda\). If \(H/\Qp\) is affine of finite type, we equip \(H(A)\)
or \(H(\Lambda)\) with the topology induced by the corresponding
coefficient algebra. More generally, if
\[
\mathcal A=\varprojlim_i\mathcal A_i
\]
is affine pro-algebraic, we equip its groups of points with the
inverse-limit topology; see
\cite[Example~3.0.5]{KANTORTHESIS} and
\cite[\S2.2]{BETTS_WEIGHT}. The functors below are thus defined on affinoid or affine test objects and then stackified with respect to $\tau$.

\begin{definition}
[Continuous non-abelian cohomology;
{\cite[Definitions~3.0.2 and~3.0.8]{KANTORTHESIS}}]
\label{definition:nonabelian_h1}
For \(S\in\mathcal C\), let
\[
Z^1\bigl(\Gamma,\mathcal A(S)\bigr)
\]
be the set of continuous maps
\(c:\Gamma\to\mathcal A(S)\) satisfying
\[
c(\gamma\delta)
=
c(\gamma)\,{}^\gamma c(\delta).
\]
The group \(\mathcal A(S)\) acts by coboundaries:
\[
(a\cdot c)(\gamma)
=
a^{-1}c(\gamma)\,{}^\gamma a.
\]

The quotient groupoids
\[
S\longmapsto
\left[
Z^1\bigl(\Gamma,\mathcal A(S)\bigr)/
\mathcal A(S)
\right]
\]
form a pseudofunctor on \(\mathcal C\). Its
\(\tau\)-stackification is denoted
\[
\mathscr H^1(\Gamma,\mathcal A).
\]
It is pointed by the class of the trivial cocycle.
\end{definition}

We reserve the notation \(H^1(\Gamma,U)\), without the script
letter, for the finite-level scheme of isomorphism classes introduced
later.

Let
\[
\mathcal A\text{-}\mathrm{Tors}_\Gamma
\]
denote the stack of \(\tau\)-locally trivial right
\(\mathcal A\)-torsors \(P\) equipped with a compatible continuous
\(\Gamma\)-action, where compatibility means
\[
\gamma(p\cdot a)
=
\gamma(p)\cdot{}^\gamma a.
\]

\begin{proposition}
[Cocycle--torsor correspondence;
{\cite[Proposition~3.1.2]{KANTORTHESIS}}]
\label{proposition:cocycle_torsor_equiv}
There is a natural equivalence of stacks
\[
\mathscr H^1(\Gamma,\mathcal A)
\xrightarrow{\sim}
\mathcal A\text{-}\mathrm{Tors}_\Gamma.
\]
\end{proposition}

For later use, let
\[
\mathscr B_\Gamma(\mathcal A)
\subseteq
\mathscr H^1(\Gamma,\mathcal A)
\]
denote the neutral gerbe of \(\Gamma\)-equivariant
\(\mathcal A\)-torsors that are \(\tau\)-locally isomorphic to the
trivial torsor. Its distinguished object is the trivial torsor, whose
automorphism sheaf is \(\mathcal A^\Gamma\). Hence there is a
canonical equivalence
\[
\mathscr B_\Gamma(\mathcal A)
\simeq
B(\mathcal A^\Gamma),
\]
where \(B(\mathcal A^\Gamma)\) denotes the classifying stack of
\(\mathcal A^\Gamma\)-torsors. If \(\mathcal A^\Gamma\) is represented
on the algebraic site by an affine group scheme \(H/\Qp\), then
\[
B(\mathcal A^\Gamma)
=
BH
=
[\Spec(\Qp)/H].
\]
On the Tate site, we use the corresponding rigid-analytic
classifying stack
\[
(BH)^{\an}
=
[\Spm(\Qp)/H^{\an}].
\]

\paragraph{Local conditions.}
Let \(w\mid p\), and assume that
\[
\mathcal A=\varprojlim_i\mathcal A_i
\]
is an affine pro-algebraic \(\Qp\)-group functor equipped with a
continuous \(G_w\)-action. For a topological \(\Qp\)-algebra \(B\),
write
\[
\mathcal A(B)\coloneqq \varprojlim_i\mathcal A_i(B),
\]
equipped with the inverse-limit topology. For an affine or affinoid
test object \(S\), set
\[
\mathcal A_{B_{\cris}}(S)
\coloneqq 
\begin{cases}
\mathcal A\bigl(A\otimes_{\Qp}B_{\cris}\bigr),
&
S=\Spec(A),
\\[2mm]
\mathcal A\bigl(\Lambda\widehat\otimes_{\Qp}B_{\cris}\bigr),
&
S=\Spm(\Lambda).
\end{cases}
\]
where the first tensor product carries its natural inductive-limit
tensor-product topology and the second its completed tensor-product
topology. In both cases, \(\mathcal A_{B_{\cris}}(S)\) is equipped
with the diagonal \(G_w\)-action.

\begin{definition}
[Local Bloch--Kato Selmer stack;
{\cite[Definition~3.3.10]{KANTORTHESIS}}]
\label{definition:kantor_bk_selmer}
The \emph{local Bloch--Kato Selmer stack}
\(\mathscr H^1_f(G_w,\mathcal A)\) is the
\(\tau\)-stackification of the pseudofunctor
\[
S\longmapsto
\ker\left(
\left[
Z^1\bigl(G_w,\mathcal A(S)\bigr)/\mathcal A(S)
\right]
\longrightarrow
\left[
Z^1\bigl(G_w,\mathcal A_{B_{\cris}}(S)\bigr)/
\mathcal A_{B_{\cris}}(S)
\right]
\right),
\]
where the kernel is the full subgroupoid of objects whose image is
isomorphic to the trivial object. Equivalently,
\(\mathscr H^1_f(G_w,\mathcal A)\) classifies equivariant
\(\mathcal A\)-torsors that become $\tau$-locally trivial after extension to
\(B_{\cris}\); no trivialization is retained.
\end{definition}

Recall from \S\ref{subsubsection:intro_relative_box_cutter} that
\[
T=S\cup S_p,
\qquad
S_p\coloneqq \{w:w\mid p\}.
\]
Let \(\mathcal A\) be equipped with a continuous \(G_T\)-action, and
fix compatible decomposition-group embeddings
\[
G_w\hookrightarrow G_T
\qquad(w\in T).
\]

\begin{definition}[Global Bloch--Kato Selmer stack]
\label{definition:kantor_global_selmer}
The \emph{global Bloch--Kato Selmer stack} is
\[
\mathscr H^1_f(G_T,\mathcal A)
\coloneqq 
\mathscr H^1(G_T,\mathcal A)
\times_{\prod_{w\mid p}\mathscr H^1(G_w,\mathcal A)}
\prod_{w\mid p}\mathscr H^1_f(G_w,\mathcal A).
\]
It is the strictly full substack of global equivariant
\(\mathcal A\)-torsors whose localization at every place
\(w\mid p\) is Bloch--Kato finite.
\end{definition}

For a finite place \(w\nmid p\), restriction to inertia gives
\[
\operatorname{res}_{I_w}:
\mathscr H^1(G_w,\mathcal A)
\longrightarrow
\mathscr H^1(I_w,\mathcal A).
\]

\begin{definition}[Unramified local condition]
\label{definition:unramified_local_condition}
The \emph{unramified local cohomology stack} is
\[
\mathscr H^1_{\mathrm{ur}}(G_w,\mathcal A)
\coloneqq 
\mathscr H^1(G_w,\mathcal A)
\times_{\mathscr H^1(I_w,\mathcal A)}
\mathscr B_{I_w}(\mathcal A).
\]
It classifies \(G_w\)-equivariant \(\mathcal A\)-torsors whose
restriction to \(I_w\) is $\tau$-locally isomorphic to the trivial inertia torsor.
No inertial trivialization is retained.
\end{definition}

Fix a set of places \(T'\) satisfying
\[
S_p\subseteq T'\subseteq T.
\]

\begin{definition}[Restricted global Selmer stack]
\label{definition:restricted_ramification_stack}
Define
\[
\mathscr H^1_{f,T'}(G_T,\mathcal A)
\coloneqq 
\mathscr H^1_f(G_T,\mathcal A)
\times_{
\prod_{w\in T\setminus T'}
\mathscr H^1(G_w,\mathcal A)
}
\prod_{w\in T\setminus T'}
\mathscr H^1_{\mathrm{ur}}(G_w,\mathcal A).
\]
Thus it classifies global \(\mathcal A\)-torsors that are
Bloch--Kato finite at every place above \(p\) and unramified at every
place in \(T\setminus T'\).
\end{definition}

\subsection{Relative Kummer Maps and Neutral Fibres}
\label{section:neutral_fibre_selmer}

We now specialize the preceding constructions to the geometric
\'{e}tale relative completion and discuss the relative Kummer maps.

\paragraph{Relative Kummer maps.}

Fix an integral basepoint
\[
b\in\mathcal X(\mathcal O_{K,S}),
\qquad
\Pi\coloneqq \pi_1^{\et}(X_{\overline{K}},\overline b).
\]
Let
\[
\Pi^{\mathrm{alg}}
\coloneqq 
\operatorname{Aut}^{\otimes}\left(\omega_{\overline b}^{\et}\right)
\]
be the pro-algebraic completion of \(\Pi\) over \(\Qp\), where the
fibre functor is taken on the Tannakian category of all lisse
\(\Qp\)-sheaves on \(X_{\overline{K}}\).

For \(F\in\{K,K_w\}\) and \(x\in X(F)\), let
\[
{}^{\vphantom{\prof}}_{x}P_b^{\prof}
\]
denote the profinite \'{e}tale path torsor from
\(\overline b\) to \(X_{\overline{K}}\), and set
\[
{}^{\vphantom{\mathrm{alg}}}_{x}P_b^{\mathrm{alg}}
\coloneqq 
\operatorname{Isom}^{\otimes}
\left(
\omega_{\overline b}^{\et},
\omega_{\overline x}^{\et}
\right),
\]
with the fibre functors again taken on the category of all lisse
\(\Qp\)-sheaves on \(X_{\overline{K}}\). Then
\({}^{\vphantom{\mathrm{alg}}}_{x}P_b^{\mathrm{alg}}\) is a right
\(\Pi^{\mathrm{alg}}\)-torsor equipped with a compatible
\(G_F\)-action. Every profinite \'{e}tale path induces a tensor isomorphism, giving
a canonical \(G_F\)-equivariant map
\[
{}^{\vphantom{\prof}}_{x}P_b^{\prof}
\coloneqq 
\pi_1^{\et}
\left(
X_{\overline{K}};\overline b,\overline x
\right)
\longrightarrow
{}^{\vphantom{\mathrm{alg}}}_{x}P_b^{\mathrm{alg}}(\Qp).
\]

\begin{definition}[\(\mathcal A\)-valued Kummer maps]
\label{definition:abs_kummer_maps}
Let \(\mathcal A\) be an affine pro-algebraic \(\Qp\)-group functor
with continuous \(G_F\)-action, and let
\[
\mu:\Pi\longrightarrow\mathcal A(\Qp)
\]
be a continuous \(G_F\)-equivariant homomorphism. By the universal
property of pro-algebraic completion, \(\mu\) extends uniquely to a
\(G_F\)-equivariant morphism
\[
\mu^{\mathrm{alg}}:
\Pi^{\mathrm{alg}}
\longrightarrow
\mathcal A.
\]
The associated \(\mathcal A\)-valued Kummer map is
\[
\kappa_F^\mu:
X(F)\longrightarrow
\mathscr H^1(G_F,\mathcal A)(\Qp),
\qquad
x\longmapsto
\left[
{}^{\vphantom{\mathrm{alg}}}_{x}P_b^{\mathrm{alg}}
\times^{\Pi^{\mathrm{alg}},\mu^{\mathrm{alg}}}
\mathcal A
\right].
\]
These maps are compatible with localization under the fixed
embeddings of decomposition groups.
\end{definition}

Let
\[
\widetilde\rho:
\Pi\longrightarrow\mathcal G^{\et}(\Qp)
\]
be the relative-completion homomorphism. It extends uniquely to a
morphism
\[
\widetilde\rho^{\mathrm{alg}}:
\Pi^{\mathrm{alg}}
\longrightarrow
\mathcal G^{\et}.
\]
The Galois action of
Proposition~\ref{proposition:etale_torsor_galois} makes these
morphisms equivariant. By the Tannakian definition of the relative
path torsor, restriction of the two endpoint fibre functors gives a
canonical identification
\[
{}^{\vphantom{\mathrm{alg}}}_{x}P_b^{\mathrm{alg}}
\times^{\Pi^{\mathrm{alg}},\widetilde\rho^{\mathrm{alg}}}
\mathcal G^{\et}
\xrightarrow{\sim}
{}^{\vphantom{\et}}_{x}P_b^{\et}.
\]
The resulting maps
\[
\kappa_F^{\mathrm{rel}}:
X(F)\longrightarrow
\mathscr H^1(G_F,\mathcal G^{\et})(\Qp),
\qquad
x\longmapsto[{}^{\vphantom{\et}}_{x}P_b^{\et}],
\]
are the \emph{relative Malcev Kummer maps}.

Let
\[
\Delta^{\et}\coloneqq \pi_0(\mathcal R^{\et}),
\qquad
{}^{\vphantom{\Delta}}_{x}P_b^\Delta
\coloneqq 
{}^{\vphantom{\et}}_{x}P_b^{\et}
\times^{\mathcal G^{\et}}
\Delta^{\et}.
\]
Define
\[
\mathcal X(\mathcal O_{K,S})_{\Delta\text{-triv}}
\coloneqq 
\left\{
x\in\mathcal X(\mathcal O_{K,S})
\ \middle|\
\left[
\left.{}^{\vphantom{\Delta}}_{x}P_b^\Delta\right|_{G_w}
\right]
=1
\text{ in }
H^1(G_w,\Delta^{\et})
\text{ for every }w\mid p
\right\}
\]
and
\[
\mathcal X_v(\mathcal O_{K_v})_{\Delta\text{-triv}}
\coloneqq 
\left\{
x\in\mathcal X_v(\mathcal O_{K_v})
\ \middle|\
[{}^{\vphantom{\Delta}}_{x}P_b^\Delta]
=1
\text{ in }
H^1(G_v,\Delta^{\et})
\right\}.
\]

\begin{proposition}
[Component-trivial relative Kummer maps;
cf.~{\cite[Proposition~3.3.13 and its proof]{KANTORTHESIS}}]
\label{proposition:kantor_kummer}
The relative Malcev Kummer maps restrict to factorizations
\[
\kappa_T^{\mathrm{rel}}:
\mathcal X(\mathcal O_{K,S})_{\Delta\text{-triv}}
\longrightarrow
\mathscr H^1_f(G_T,\mathcal G^{\et})(\Qp)
\]
and
\[
\kappa_v^{\mathrm{rel}}:
\mathcal X_v(\mathcal O_{K_v})_{\Delta\text{-triv}}
\longrightarrow
\mathscr H^1_f(G_v,\mathcal G^{\et})(\Qp).
\]
These factorizations are compatible with localization and with
pushout along every $G_T$-equivariant motivic quotient of
$\mathcal G^{\et}$.
\end{proposition}

\begin{remark}
Kantor assumes that the crystalline reductive monodromy group is
semisimple and simply connected, in which case its component group
is trivial. In the Kodaira--Parshin setting, its identity component
is semisimple and simply connected, but the full group may have a
finite permutation component. On the loci above, the induced
component-group torsor is trivial; after choosing a point of this
torsor, the remaining reductive torsor is a torsor under the identity
component, and Kantor's argument applies.
\end{remark}

Fix henceforth a finite-type \(G_T\)-equivariant motivic quotient
\begin{equation}
\label{eq:finite_level_relative_extension}
1\longrightarrow U^{\et}
\longrightarrow G^{\et}
\xrightarrow{\pi}R^{\et}
\longrightarrow1,
\end{equation}
and write
\[
{}^{\vphantom{G}}_{x}P_b^G
\coloneqq 
{}^{\vphantom{\et}}_{x}P_b^{\et}\times^{\mathcal G^{\et}}G^{\et}.
\]
We suppress the superscript \(\et\) when no confusion can arise.

\paragraph{Neutral reductive fibres.}

For a profinite group \(\Gamma\), recall the neutral gerbe
\[
\mathscr B_\Gamma(R)
\subseteq
\mathscr H^1(\Gamma,R)
\]
of equivariant \(R\)-torsors locally isomorphic to the trivial
torsor. Its distinguished object is the trivial torsor, whose
automorphism group is \(R^\Gamma\).

\begin{definition}[Neutral and trivialized neutral fibres]
\label{definition:neutral_reductive_fibre}
Let
\[
\{\mathbf1_R\}
\longrightarrow
\mathscr H^1(\Gamma,R)
\]
denote the point classifying the trivial equivariant \(R\)-torsor.
Define
\[
\mathscr N_\Gamma(G/R)
\coloneqq 
\mathscr H^1(\Gamma,G)
\times_{\mathscr H^1(\Gamma,R)}
\mathscr B_\Gamma(R)
\]
and
\[
\widetilde{\mathscr N}_\Gamma(G/R)
\coloneqq 
\mathscr H^1(\Gamma,G)
\times_{\mathscr H^1(\Gamma,R)}
\{\mathbf1_R\}.
\]

Thus \(\mathscr N_\Gamma(G/R)\) classifies equivariant \(G\)-torsors
whose reductive pushout is locally isomorphic to the trivial
\(R\)-torsor, without a chosen trivialization. The trivialized fibre
\(\widetilde{\mathscr N}_\Gamma(G/R)\) classifies pairs
\[
(P,\alpha),
\qquad
\alpha:
P\times^G R\xrightarrow{\sim}\mathbf1_R,
\]
consisting of an equivariant \(G\)-torsor and a chosen equivariant
trivialization of its reductive pushout.
\end{definition}

\begin{proposition}
[Neutral-fibre quotient]
\label{proposition:neutral_fibre_quotient}
Extension of structure group induces a canonical equivalence
\[
\mathscr H^1(\Gamma,U)
\xrightarrow{\sim}
\widetilde{\mathscr N}_\Gamma(G/R).
\]
Moreover, changing the chosen reductive trivialization defines an
action of \(R^\Gamma\), and there is a canonical equivalence
\[
\mathscr N_\Gamma(G/R)
\simeq
\left[
\widetilde{\mathscr N}_\Gamma(G/R)/R^\Gamma
\right]
\simeq
\left[
\mathscr H^1(\Gamma,U)/R^\Gamma
\right].
\]
\end{proposition}

\begin{proof}
An object of the rigidified fibre is a pair \((P,\alpha)\), where
\(P\) is an equivariant \(G\)-torsor and
\[
\alpha:
P\times^G R\xrightarrow{\sim}\mathbf1_R
\]
is an equivariant trivialization. The inverse image of the identity
section under \(\alpha\) is an equivariant \(U\)-torsor.

Conversely, extension of structure group sends an equivariant
\(U\)-torsor \(Q\) to
\[
Q\times^U G,
\]
whose reductive pushout is canonically trivialized. These
constructions are mutually inverse. Finally,
\[
R^\Gamma=\operatorname{Aut}_\Gamma(\mathbf1_R)
\]
acts by changing the chosen trivialization, and passage from the
point \(\{\mathbf1_R\}\) to its residual gerbe gives the asserted
quotient. Equivalently, this is the fibre description 
of~\cite[Proposition~3.1.4]{KANTORTHESIS}, applied to
\[
1\longrightarrow U\longrightarrow G\longrightarrow R
\longrightarrow1
\]
over the trivial class in \(\mathscr H^1(\Gamma,R)\).
\end{proof}

\paragraph{Neutral-fibre Selmer stacks.}

\begin{definition}[Neutral and trivialized neutral Selmer fibres]
\label{definition:neutral_fibre_selmer}
For \(w\mid p\), define
\[
\mathscr N_{f,w}(G/R)
\coloneqq 
\mathscr H^1_f(G_w,G)
\times_{\mathscr H^1(G_w,R)}
\mathscr B_{G_w}(R)
\]
and
\[
\widetilde{\mathscr N}_{f,w}(G/R)
\coloneqq 
\mathscr H^1_f(G_w,G)
\times_{\mathscr H^1(G_w,R)}
\{\mathbf1_R\}.
\]

More generally, for
\[
S_p\subseteq T'\subseteq T,
\]
define
\[
\mathscr N_{f,T'}(G/R)
\coloneqq 
\mathscr H^1_{f,T'}(G_T,G)
\times_{\mathscr H^1(G_T,R)}
\mathscr B_{G_T}(R)
\]
and
\[
\widetilde{\mathscr N}_{f,T'}(G/R)
\coloneqq 
\mathscr H^1_{f,T'}(G_T,G)
\times_{\mathscr H^1(G_T,R)}
\{\mathbf1_R\}.
\]

The stacks \(\mathscr N_{f,\star}(G/R)\) classify relative
\(G\)-torsors satisfying the indicated local conditions whose
reductive pushout is neutral, whereas
\(\widetilde{\mathscr N}_{f,\star}(G/R)\) also includes a chosen
trivialization of that pushout.
\end{definition}

Since
\[
\mathscr B_\Gamma(R)
\simeq
[\{\mathbf1_R\}/R^\Gamma],
\]
there are canonical equivalences
\[
\mathscr N_{f,w}(G/R)
\simeq
\left[
\widetilde{\mathscr N}_{f,w}(G/R)/R^{G_w}
\right]
\]
and
\[
\mathscr N_{f,T'}(G/R)
\simeq
\left[
\widetilde{\mathscr N}_{f,T'}(G/R)/R^{G_T}
\right].
\]

\begin{proposition}[Bloch--Kato neutral-fibre comparison]
\label{proposition:bk_neutral_fibre_comparison}
Extension of structure group induces canonical equivalences
\[
\mathscr H^1_f(G_w,U)
\xrightarrow{\sim}
\widetilde{\mathscr N}_{f,w}(G/R)
\qquad(w\mid p)
\]
and
\[
\mathscr H^1_f(G_T,U)
\xrightarrow{\sim}
\widetilde{\mathscr N}_{f,T}(G/R).
\]
Consequently,
\[
\mathscr N_{f,w}(G/R)
\simeq
\left[
\mathscr H^1_f(G_w,U)/R^{G_w}
\right]
\]
and
\[
\mathscr N_{f,T}(G/R)
\simeq
\left[
\mathscr H^1_f(G_T,U)/R^{G_T}
\right].
\]
\end{proposition}

\begin{proof}
We first prove the local equivalence for a fixed \(w\mid p\).
The ordinary neutral-fibre equivalence of
Proposition~\ref{proposition:neutral_fibre_quotient} identifies the
trivialized neutral fibre with \(\mathscr H^1(G_w,U)\). We verify that
it identifies the corresponding Bloch--Kato substacks.

Let \(S\) be an affine or affinoid test object, and set
\[
B_{\cris,S}
\coloneqq 
\begin{cases}
A\otimes_{\Qp}B_{\cris},
&
S=\Spec(A),
\\[2mm]
\Lambda\widehat\otimes_{\Qp}B_{\cris},
&
S=\Spm(\Lambda).
\end{cases}
\]
Let \(Q\) be a \(U\)-torsor over \(S\), and put
\[
P\coloneqq Q\times^U G.
\]
If \(Q\) becomes trivial over \(B_{\cris,S}\), then so does \(P\).

Conversely, suppose that \(P\) becomes trivial over \(B_{\cris,S}\)
and that its reductive pushout is equipped with an equivariant
trivialization
\[
\alpha:
P\times^G R\xrightarrow{\sim}\mathbf1_R.
\]
The corresponding \(U\)-reduction is \(Q=\alpha^{-1}(1)\).

Let
\[
1\longrightarrow U_w^{\cri}
\longrightarrow G_w^{\cri}
\longrightarrow R_w^{\cri}
\longrightarrow1
\]
denote the crystalline realization of
\eqref{eq:finite_level_relative_extension}. It admits a Levi section
\[
s^{\cri}:R_w^{\cri}\longrightarrow G_w^{\cri};
\]
cf.~\cite[Proposition~3.1]{Hain_HodgeDeRham_Chapter2016}. After
extension to \(B_{\cris,S}\) and transport through the comparison
isomorphism, it induces a \(G_w\)-equivariant section
\[
s_{B_{\cris,S}}:
R_{B_{\cris,S}}\longrightarrow G_{B_{\cris,S}}.
\] 
Working \(\tau\)-locally on \(S\), choose
\[
p\in P(B_{\cris,S})^{G_w},
\qquad
r\coloneqq \alpha(\overline p)\in R(B_{\cris,S})^{G_w}.
\]
Then
\[
q\coloneqq p\cdot s_{B_{\cris,S}}(r^{-1})
\]
belongs to \(Q(B_{\cris,S})^{G_w}\), since
\(\alpha(\overline q)=1\). Thus \(Q\) becomes
\(\tau\)-locally trivial over \(B_{\cris,S}\). This construction is
functorial in \(S\), and hence proves the local equivalence of stacks.

The global equivalence follows by imposing the local Bloch--Kato
condition at every \(w\mid p\). The quotient presentations then follow
from the preceding neutral-fibre equivalences.
\end{proof}

\begin{definition}[Relative neutral inertia locus]
\label{definition:relative_neutral_inertia_locus}
Let \(w\nmid p\). Define
\[
\mathscr Z_w(G/U)
\coloneqq 
\mathscr H^1(I_w,U)
\times_{\mathscr H^1(I_w,G)}
\mathscr B_{I_w}(G).
\]
Thus \(\mathscr Z_w(G/U)\) classifies equivariant \(U\)-torsors
whose extension of structure group to \(G\) is locally isomorphic to
the trivial equivariant \(G\)-torsor.
\end{definition}

\begin{proposition}[Restricted neutral-fibre presentation]
\label{proposition:restricted_neutral_fibre_presentation}
For every \(w\nmid p\), there is a canonical equivalence
\[
\mathscr Z_w(G/U)
\simeq
\left[
R^{I_w}/G^{I_w}
\right],
\]
where \(G^{I_w}\) acts on \(R^{I_w}\) through
\(G^{I_w}\to R^{I_w}\). Moreover, for every
\(S_p\subseteq T'\subseteq T\), there is a canonical equivalence
\[
\widetilde{\mathscr N}_{f,T'}(G/R)
\simeq
\mathscr H^1_f(G_T,U)
\times_{
\prod_{w\in T\setminus T'}
\mathscr H^1(I_w,U)
}
\prod_{w\in T\setminus T'}
\mathscr Z_w(G/U).
\]
\end{proposition}

\begin{proof}
The stack \(\mathscr Z_w(G/U)\) classifies \(U\)-reductions of a
locally trivial \(I_w\)-equivariant \(G\)-torsor. After choosing an
equivariant trivialization of the \(G\)-torsor, such reductions are
parametrized by
\[
(G/U)^{I_w}=R^{I_w}.
\]
Changing the trivialization acts through
\[
G^{I_w}\longrightarrow R^{I_w}.
\]
Forgetting the trivialization therefore gives
\[
\mathscr Z_w(G/U)\simeq[R^{I_w}/G^{I_w}],
\]
which proves the first equivalence.

By Proposition~\ref{proposition:bk_neutral_fibre_comparison}, the
trivialized global Bloch--Kato neutral fibre is
\(\mathscr H^1_f(G_T,U)\). Under this identification, requiring the
associated \(G\)-torsor to be unramified at \(w\) is precisely the
condition that its restricted \(U\)-torsor belong to
\(\mathscr Z_w(G/U)\). Imposing these conditions for all
\(w\in T\setminus T'\) gives the second equivalence.
\end{proof}

\paragraph{Reductive strata.}

Let \(\Gamma\) be \(G_T\) or \(G_v\), and let \(x,y\) be
\(\mathcal O_{K,S}\)-integral points of \(\mathcal X\) or
\(\mathcal O_{K_v}\)-integral points of \(\mathcal X_v\), respectively.
Basing the relative
completion at \(y\), let
\[
{}^{\vphantom{R}}_{x}P_y^R
\coloneqq 
{}^{\vphantom{G}}_{x}P_y^G\times^G R
\]
denote the reductive pushout of the relative path torsor. It is a
\(\Gamma\)-equivariant right \(R\)-torsor.

\begin{definition}[Reductive Galois type]
\label{definition:reductive_galois_type}
We say that \(x\) and \(y\) have the same \emph{reductive Galois
type}, and write
\[
x\sim_R y,
\]
if
\[
[{}^{\vphantom{R}}_{x}P_y^R]=1
\qquad\text{in}\qquad
H^1(\Gamma,R).
\]
Equivalently, their fibre functors on
\(\langle\mathcal L^{\et}\rangle_\otimes\) are isomorphic as
\(\Gamma\)-equivariant tensor fibre functors.
\end{definition}

Composition and inversion of path torsors show that
\(\sim_R\) is an equivalence relation.

\begin{theorem}
[{Finiteness of reductive Galois types;
cf.~\cite[Proof of Theorem~5]{FALTINGS_FINITENESS}
and~\cite[Lemmas~2.3 and~2.6]{LV}}]
\label{theorem:faltings_finiteness}
The relation \(\sim_R\) has only finitely many equivalence classes on
\[
\mathcal X(\mathcal O_{K,S}).
\]
\end{theorem}

\begin{proof}
The finite \'{e}tale part of the Kodaira--Parshin construction has
fixed degree and is unramified outside \(T\), so only finitely many
possibilities occur by Hermite--Minkowski. For each such possibility,
the polarized abelian factors have fixed dimension and polarization
type and have good reduction outside \(T\). Faltings's finiteness
theorem therefore gives only finitely many isomorphism classes of the
resulting polarized abelian-by-finite Galois representations.

By the full-monodromy condition, the connected reductive monodromy
group is the product of the symplectic groups of the cohomological
blocks, while its component group records the finite permutation
action on those blocks. Consequently, the block decomposition,
polarizations, and Galois action determine the associated
\(G_T\)-equivariant tensor fibre functor on
\[
\langle\mathcal L^{\et}\rangle_\otimes.
\]
There are therefore only finitely many isomorphism classes of these
fibre functors, which, by
Definition~\ref{definition:reductive_galois_type}, are precisely the
equivalence classes for \(\sim_R\).
\end{proof}

\begin{definition}[Reductive strata]
\label{definition:global_partition}
Choose representatives
\[
b_1,\ldots,b_s
\in
\mathcal X(\mathcal O_{K,S})
\]
for the finitely many nonempty reductive Galois types of
Theorem~\ref{theorem:faltings_finiteness}, and set
\[
\Sigma_i^{\mathrm{red}}
\coloneqq 
\left\{
x\in\mathcal X(\mathcal O_{K,S})
\ \middle|\
x\sim_R b_i
\right\}.
\]
Then
\[
\mathcal X(\mathcal O_{K,S})
=
\bigsqcup_{i=1}^s
\Sigma_i^{\mathrm{red}}.
\]

At the distinguished place \(v\mid p\), define
\[
\Sigma_{i,v}^{\mathrm{red}}
\coloneqq 
\left\{
x\in\mathcal X_v(\mathcal O_{K_v})
\ \middle|\
x\sim_R b_i
\right\},
\]
where \(b_i\) is viewed as an \(\mathcal O_{K_v}\)-point. These loci
contain the localizations of the corresponding global strata, but
need not be distinct or exhaust
\(\mathcal X_v(\mathcal O_{K_v})\).
\end{definition}

For each \(i\), rebase the finite-level relative completion at
\(b_i\), and continue to denote the resulting extension by
\[
1\longrightarrow U
\longrightarrow G
\longrightarrow R
\longrightarrow1.
\]

\begin{proposition}[Canonical neutral-fibre Kummer maps]
\label{proposition:kummer_unip}
For every nonempty reductive stratum, the relative Kummer maps
restrict canonically to
\[
\kappa^{\mathrm{neu}}_{i,T}:
\Sigma_i^{\mathrm{red}}
\longrightarrow
\mathscr N_{f,T}(G/R)(\Qp),
\qquad
x\longmapsto[{}^{\vphantom{G}}_{x}P_{b_i}^G],
\]
and
\[
\kappa^{\mathrm{neu}}_{i,v}:
\Sigma_{i,v}^{\mathrm{red}}
\longrightarrow
\mathscr N_{f,v}(G/R)(\Qp),
\qquad
x\longmapsto[{}^{\vphantom{G}}_{x}P_{b_i}^G].
\]
After forgetting the neutral reductive condition, these are the
restrictions of Kantor's relative Malcev Kummer maps.
\end{proposition}

\begin{proof}
By Proposition~\ref{proposition:kantor_kummer}, the relative path
torsors attached to integral points satisfy the relevant global and
local \(p\)-adic Bloch--Kato conditions.

If \(x\) and \(b_i\) lie in the same reductive stratum, their
Galois-equivariant tensor fibre functors on $\langle \mathcal{L}^\et\rangle_\otimes$ are isomorphic. Hence
\[
{}^{\vphantom{R}}_{x}P_{b_i}^R
=
{}^{\vphantom{G}}_{x}P_{b_i}^G\times^G R
\]
is isomorphic to the trivial equivariant \(R\)-torsor. The relative
Kummer class therefore lies canonically in the corresponding neutral
fibre.
\end{proof}

\begin{remark}[Fixed-point inertia]
\label{remark:effective_reductive_inertia}
The fixed-point group \(R^\Gamma\) need not be trivial. In the
Kodaira--Parshin setting, it generally contains the simultaneous
central element
\[
(-1,\ldots,-1)\in R^\circ.
\]
The neutral-fibre targets are therefore naturally stacky. Their
finite-level algebraicity and effective inertia are studied in
\S\ref{section:neutral_selmer_representability}.
\end{remark}

\subsection{Finite Descent and Kummer Stratification}
\label{section:finite_descent_kummer}

Fix
\[
S_p\subseteq T'\subseteq T,
\]
where \(S_p\) is the set of places above \(p\), and retain the
finite-type motivic quotient
\[
1\longrightarrow U
\longrightarrow G
\longrightarrow R
\longrightarrow1
\]
from \S\ref{section:neutral_fibre_selmer}. For each reductive
stratum, we regard this quotient as rebased at the chosen point of
that stratum.

We construct, at every place
\[
w\in T\setminus T',
\]
a local partition on whose cells the relative \(G\)-valued path
torsors are unramified.

\begin{lemma}[Finite descent to pro-\(p\) monodromy]
\label{lemma:pro_p_finite_descent}
Let \(w\nmid p\). There is a finite collection of finite \'{e}tale
covers
\[
\pi_{w,j}:X_{w,j}\longrightarrow X_{K_w}
\]
such that:

\begin{enumerate}
\item the images of the sets \(X_{w,j}(K_w)\) cover \(X(K_w)\);

\item for every \(j\), the closure of the image of
\[
\pi_1^{\et}(\overline{X}_{w,j})
\longrightarrow
G(\Qp)
\]
is a pro-\(p\) group.
\end{enumerate}
\end{lemma}

\begin{proof}
Let
\[
\widetilde\rho:
\pi_1^{\et}(X_{\overline{K}},\overline b)
\longrightarrow G(\Qp)
\]
be the finite-level relative-monodromy homomorphism, and let \(H\) be
the closure of its image. The group \(H\) is a compact \(p\)-adic
analytic group, stable under \(G_w\), and hence contains a
\(G_w\)-stable open normal pro-\(p\) subgroup \(H_p\). Indeed, one may
start with any open normal pro-\(p\) subgroup and intersect its finite
orbit under \(G_w\).

Set
\[
N\coloneqq \widetilde\rho^{-1}(H_p).
\]
Then \(N\) is an open normal \(G_w\)-stable subgroup of
\(\pi_1^{\et}(X_{\overline{K}},\overline b)\). The corresponding pointed
geometric cover descends to a finite \'{e}tale torsor over \(K_w\)
under the finite \'{e}tale group scheme associated with
\[
Q\coloneqq 
\pi_1^{\et}(X_{\overline{K}},\overline b)/N.
\]

The set \(H^1(K_w,Q)\) is finite. Let
\[
\{X_{w,j}\longrightarrow X_{K_w}\}_j
\]
be the resulting finite collection of twists. Every
\(x\in X(K_w)\) lifts to the twist determined by the class of the
fibre over \(x\), so their images cover \(X(K_w)\). After geometric
base change, every twist has geometric fundamental group \(N\), up to
inner automorphism, and its relative monodromy therefore has image
contained in \(H_p\).
\end{proof}

\begin{lemma}
[{Local invariant paths;
cf.~\cite[Theorem~1.3.1 and the proof of
Theorem~3.1.6]{BETTS_DOGRA}}]
\label{lemma:local_invariant_path}
Let
\[
\pi:C\longrightarrow X_{K_w}
\]
be one of the covers in
Lemma~\ref{lemma:pro_p_finite_descent}, and let
\(x,y\in C(K_w)\). Suppose that their reductions lie on the same
irreducible component of the geometric special fibre of a regular
semistable model of \(C\) over the ring of integers of a finite
extension of \(K_w\). Then
\[
{}^{\vphantom{G}}_{\pi(x)}P_{\pi(y)}^G
\]
has an \(I_w\)-fixed point. Equivalently,
\[
[{}^{\vphantom{G}}_{\pi(x)}P_{\pi(y)}^G]|_{I_w}=1
\qquad\text{in }\mathscr H^1(I_w,G).
\]
\end{lemma}

\begin{proof}
Let \(\ell_w\neq p\) be the residue characteristic of \(K_w\).
Betts--Dogra's relative Oda--Tamagawa criterion gives an
inertia-invariant path between \(x\) and \(y\). More precisely, the
proof of their nonabelian Picard--Lefschetz theorem establishes the
corresponding statement for the maximal
pro-prime-to-\(\ell_w\) \'{e}tale fundamental groupoid.

The maximal pro-\(p\) groupoid is a quotient of the
pro-prime-to-\(\ell_w\) groupoid. By
Lemma~\ref{lemma:pro_p_finite_descent}, the relative monodromy on
\(C\) factors through this pro-\(p\) quotient. Pushing forward the
invariant path therefore gives an \(I_w\)-fixed point of the
\(G\)-valued path torsor.
\end{proof}

For the remainder of this subsection, we work in the genus $\ge 2$
case of Convention~\ref{convention:fixed_kp_family}; thus
\(g(\widehat X)\geq2\), and the relative completion is formed on
\(\widehat X\) and evaluated at points of
\(X\subseteq\widehat X\). This is the only setting in which we shall
use the finite Kummer stratification below.

Accordingly, the covers in
Lemma~\ref{lemma:pro_p_finite_descent} are the restrictions
\[
X_{w,j}
=
\widehat X_{w,j}\times_{\widehat X}X
\]
of finite \'{e}tale covers
\[
\widehat X_{w,j}\longrightarrow\widehat X_{K_w}.
\]

\begin{proposition}[Finite local partitions]
\label{proposition:local_partition}
For every \(w\in T\setminus T'\), there exist a finite set \(A_w\)
and a partition
\[
X(K_w)=\bigsqcup_{a\in A_w}C_{w,a}
\]
such that, whenever \(x,y\in C_{w,a}\), the relative path torsor
\[
{}^{\vphantom{G}}_{x}P_y^G
\]
is unramified at \(w\).
\end{proposition}

\begin{proof}
For each \(j\), choose, after a finite extension of \(K_w\), a
regular semistable model of the proper curve \(\widehat X_{w,j}\).
Partition \(X_{w,j}(K_w)\subseteq\widehat X_{w,j}(K_w)\) according
to the irreducible component of the geometric special fibre met by
reduction. Since \(\widehat X_{w,j}\) is proper and its special fibre
has finitely many irreducible components, this is a finite partition.

If two points lie in the same component cell,
Lemma~\ref{lemma:local_invariant_path} shows that the relative path
torsor between their images in \(X(K_w)\) is unramified. By
Lemma~\ref{lemma:pro_p_finite_descent}, the images of the finitely
many component cells cover \(X(K_w)\). Taking the nonempty Boolean atoms of these finitely many subsets gives
the required partition. We denote the finite set of such atoms by
\(A_w\), and the corresponding atoms by \(C_{w,a}\), \(a\in A_w\).
\end{proof}

For each reductive stratum \(\Sigma_i^{\mathrm{red}}\), apply
Proposition~\ref{proposition:local_partition} to the finite-level
quotient rebased at \(b_i\). Denote the resulting finite indexing set
by \(A_{i,w}\), and write
\[
X(K_w)
=
\bigsqcup_{a\in A_{i,w}}C_{i,w,a},
\qquad
w\in T\setminus T',
\]
for the resulting local partitions.

\begin{definition}[Kummer stratification]
\label{definition:kummer_stratification}
The \emph{Kummer stratification} associated with \(G\) and \(T'\) is
the collection of nonempty subsets
\[
\Sigma_{i,\mathbf a}^{\mathrm{Kummer}}
\coloneqq 
\Sigma_i^{\mathrm{red}}
\cap
\bigcap_{w\in T\setminus T'}
\operatorname{loc}_w^{-1}(C_{i,w,a_w}),
\]
where
\[
1\leq i\leq s,
\qquad
\mathbf a=(a_w)_{w\in T\setminus T'}.
\]
These subsets form a finite partition
\[
\mathcal X(\mathcal O_{K,S})
=
\bigsqcup_{i,\mathbf a}
\Sigma_{i,\mathbf a}^{\mathrm{Kummer}}.
\]
\end{definition}

\begin{theorem}[Kummer stratification and factorization]
\label{thm:kummer_stratification}
Let \(\Sigma\) be a Kummer stratum and choose \(b\in\Sigma\).
After rebasing the finite-type motivic quotient at \(b\), the
relative Malcev Kummer map restricts canonically to
\[
\kappa_{\Sigma}^{\mathrm{neu}}:
\Sigma
\longrightarrow
\mathscr N_{f,T'}(G/R)(\Qp),
\qquad
x\longmapsto[{}^{\vphantom{G}}_{x}P_b^G].
\]
After forgetting the neutral reductive and local conditions, this is
the restriction of Kantor's relative Malcev Kummer map.
\end{theorem}

\begin{proof}
Since \(\Sigma\) is contained in a single reductive stratum, one has
\[
x\sim_R b
\]
for every \(x\in\Sigma\). Hence the reductive path torsor
\[
{}^{\vphantom{R}}_{x}P_b^R
\]
is neutral as a \(G_T\)-equivariant \(R\)-torsor.

By Proposition~\ref{proposition:kantor_kummer}, the relative path
torsor
\[
{}^{\vphantom{G}}_{x}P_b^G
\]
satisfies the global \(p\)-adic Bloch--Kato conditions. Moreover, for
every \(w\in T\setminus T'\), the points
\[
\operatorname{loc}_w(x),
\qquad
\operatorname{loc}_w(b)
\]
belong to the same local cell. Therefore
Proposition~\ref{proposition:local_partition} gives
\[
[{}^{\vphantom{G}}_{x}P_b^G]|_{I_w}=1
\qquad\text{in}\qquad
\mathscr H^1(I_w,G).
\]
These are precisely the defining conditions of
\(\mathscr N_{f,T'}(G/R)\).
\end{proof}

\subsection{Algebraicity of Finite-Level Neutral Selmer Stacks}
\label{section:neutral_selmer_representability}

We now work on the category of \(\Qp\)-schemes with the \'{e}tale
topology. Retain the finite-type \(G_T\)-equivariant motivic quotient
\[
1\longrightarrow U
\longrightarrow G
\longrightarrow R
\longrightarrow1
\]
of \eqref{eq:finite_level_relative_extension}. We first establish the
algebraicity of the unipotent cohomology and Bloch--Kato functors that
provide atlases for the neutral-fibre stacks.

\paragraph{Unipotent cohomology and Bloch--Kato schemes.}

Let \(\Gamma\) be a profinite group, and let \(U\) be a unipotent
algebraic group over \(\Qp\), equipped with a continuous
\(\Gamma\)-action and a finite \(\Gamma\)-stable central series
\[
U=U^1\supseteq U^2\supseteq\cdots\supseteq U^{r+1}=1.
\]
Set
\[
V_i\coloneqq \Lie(U^i/U^{i+1}).
\]

\begin{theorem}
[Kim;
{\cite[Proposition~2]{KIM2};
cf.~\cite[Theorem~2.2.1 and Proposition~2.2.3]{BETTS_WEIGHT}}]
\label{theorem:kim_unrestricted_representability}
Assume that
\[
H^0(\Gamma,V_i)=0
\]
and that \(H^j(\Gamma,V_i)\) is finite-dimensional for every \(i\)
and \(j\). Then the functor of continuous cohomology classes
\[
A\longmapsto H^1\bigl(\Gamma,U(A)\bigr)
\]
is represented by an affine \(\Qp\)-scheme of finite type, denoted
\(H^1(\Gamma,U)\). Moreover, the cohomology stack
\[
\mathscr H^1(\Gamma,U)
\]
is represented by this scheme.
\end{theorem}

\begin{proof}
The representability assertion is
\cite[Proposition~2]{KIM}. If \(c\) is a \(U\)-valued cocycle, its
automorphism group in the cohomology stack is
\[
H^0(\Gamma,{}^cU).
\]
Twisting by \(c\) acts by inner automorphisms and therefore does not
alter the action on the central graded pieces \(V_i\). Induction
along the central series gives
\[
H^0(\Gamma,{}^cU)=1.
\]
Thus the stack has trivial inertia and is represented by its sheaf of
isomorphism classes.
\end{proof}

We use the following unfiltered form of Betts's local result.
\begin{proposition}
[{\cite[Propositions~3.1.3 and~3.1.6,
Remark~3.1.7]{BETTS_WEIGHT}}]
\label{proposition:betts_local_bk_representability}
Let \(w\mid p\), and suppose that \(H^1(G_w,U)\) is represented by
an affine scheme. Assume that each \(V_i\) is crystalline and that
\(1\) is not an eigenvalue of Frobenius on \(D_{\cris}(V_i)\). Then
\[
\mathscr H^1_f(G_w,U)
\]
is represented by a closed affine subscheme
\[
H^1_f(G_w,U)
\hookrightarrow
H^1(G_w,U).
\]
\end{proposition}

\begin{proof}
Betts proves that the functor of Bloch--Kato cohomology classes is
represented by a closed affine subscheme under a purity hypothesis on
the graded pieces. In his proof, purity is used only to ensure that
\(\varphi-1\) is invertible on their crystalline realizations; the
same induction therefore applies under the stated hypothesis.

It remains to pass from isomorphism classes to the cohomology stack.
The Frobenius hypothesis gives
\[
H^0(G_w,V_i)
\hookrightarrow
D_{\cris}(V_i)^{\varphi=1}
=0.
\]
For a \(U\)-valued cocycle \(c\), its automorphism group in the
cohomology stack is \(H^0(G_w,{}^cU)\). Inner twisting does not alter
the \(G_w\)-action on the central graded pieces \(V_i\), so induction
along the central series gives
\[
H^0(G_w,{}^cU)=1.
\]
Thus the Bloch--Kato stack has trivial inertia and is represented by
its sheaf of isomorphism classes, namely Betts's closed affine
subscheme.
\end{proof}

\begin{proposition}[Kodaira--Parshin unipotent Selmer schemes]
\label{proposition:kp_unipotent_selmer_algebraicity}
For the unipotent radical \(U\) of every finite-type
Kodaira--Parshin motivic quotient, the stacks
\[
\mathscr H^1(G_T,U),
\qquad
\mathscr H^1(G_w,U),
\qquad
\mathscr H^1_f(G_w,U)
\quad(w\mid p),
\]
and
\[
\mathscr H^1_f(G_T,U)
\]
are represented by affine \(\Qp\)-schemes of finite type. Moreover,
the Bloch--Kato schemes are closed in the corresponding unrestricted
cohomology schemes.
\end{proposition}

\begin{proof}
Equip \(U\) with its DCS filtration and write
\[
V_n\coloneqq \Gr_n^{\DCS}\Lie(U).
\]
For every \(w\mid p\), the \'{e}tale--crystalline comparison shows
that the \(V_n|_{G_w}\) are crystalline. By
\cite[Proposition~6.2.4]{KANTORTHESIS}, their crystalline
realizations have strictly negative Weil weights. In particular,
they contain no copy of the unit \(F\)-isocrystal, and hence
\[
H^0(G_w,V_n)=0.
\]
It follows that
\[
H^0(G_T,V_n)=0.
\]

The cohomology groups of the finite-dimensional representations
\(V_n\) are finite-dimensional, both locally and globally.
Theorem~\ref{theorem:kim_unrestricted_representability} therefore
applies to \(G_T\) and to every \(G_w\), \(w\mid p\).
Proposition~\ref{proposition:betts_local_bk_representability} gives
the local Bloch--Kato schemes. Finally,
\[
H^1_f(G_T,U)
=
H^1(G_T,U)
\times_{\prod_{w\mid p}H^1(G_w,U)}
\prod_{w\mid p}H^1_f(G_w,U)
\]
is a closed affine finite-type subscheme of \(H^1(G_T,U)\).
\end{proof}

\paragraph{The \(p\)-adic neutral fibres.}

Combining
Proposition~\ref{proposition:kp_unipotent_selmer_algebraicity}
with
Proposition~\ref{proposition:bk_neutral_fibre_comparison}
gives the required algebraic models.

\begin{corollary}
\label{corollary:neutral_bk_algebraicity}
For every \(w\mid p\), there are equivalences of algebraic stacks
\[
\mathscr N_{f,w}(G/R)
\simeq
\left[
H^1_f(G_w,U)/R^{G_w}
\right]
\]
and
\[
\mathscr N_{f,T}(G/R)
\simeq
\left[
H^1_f(G_T,U)/R^{G_T}
\right].
\]
In particular, these are Artin stacks of finite type over \(\Qp\)
with affine inertia.
\end{corollary}

\begin{proof}
The groups \(R^{G_w}\) and \(R^{G_T}\) are closed affine algebraic
subgroups of \(R\). The result follows from
Proposition~\ref{proposition:kp_unipotent_selmer_algebraicity} and the
quotient presentations of Proposition~\ref{proposition:bk_neutral_fibre_comparison}.
\end{proof}

\paragraph{Conditions away from \(p\).}

Recall from
Definition~\ref{definition:relative_neutral_inertia_locus}
that, for \(w\nmid p\),
\[
\mathscr Z_w(G/U)
=
\mathscr H^1(I_w,U)
\times_{\mathscr H^1(I_w,G)}
\mathscr B_{I_w}(G),
\]
and from
Proposition~\ref{proposition:restricted_neutral_fibre_presentation}
that
\[
\mathscr Z_w(G/U)
\simeq
\left[R^{I_w}/G^{I_w}\right].
\]

\begin{proposition}[Algebraicity of the relative neutral inertia locus]
\label{proposition:local_neutral_locus}
The stacks
\[
\mathscr H^1(I_w,U)
\qquad\text{and}\qquad
\mathscr Z_w(G/U)
\]
are Artin stacks of finite type over \(\Qp\). Moreover, the natural
morphism
\[
\mathscr Z_w(G/U)
\longrightarrow
\mathscr H^1(I_w,U)
\]
is representable and of finite type.
\end{proposition}

\begin{proof}
Let
\[
U
=
\mathcal F_{\mathrm{DCS}}^1U
\supseteq
\mathcal F_{\mathrm{DCS}}^2U
\supseteq
\cdots
\supseteq
\mathcal F_{\mathrm{DCS}}^{c+1}U
=
1
\]
be the descending central series of \(U\). It is \(I_w\)-stable,
since every \(\mathcal F_{\mathrm{DCS}}^nU\) is characteristic in
\(U\), and each quotient
\[
V_n
\coloneqq \Gr_n^{\DCS} U = 
\mathcal F_{\mathrm{DCS}}^nU/
\mathcal F_{\mathrm{DCS}}^{n+1}U
\]
is a finite-dimensional vector group over \(\Qp\).

We first treat vector groups uniformly in families. Let
\[
J_w
\coloneqq 
\ker\left(
I_w\longrightarrow I_w^{(p)}
\right),
\qquad
I_w^{(p)}
\simeq
\mathbb Z_p(1),
\]
where \(I_w^{(p)}\) is the maximal pro-\(p\) quotient, and choose a
topological generator \(\gamma\) of \(I_w^{(p)}\). Since \(J_w\) has
pro-order prime to \(p\), its invariants functor is exact. Thus, for
every finite-dimensional continuous \(I_w\)-representation \(V\), if
\[
W\coloneqq V^{J_w},
\]
then for every \(\Qp\)-algebra \(A\) there is a functorial
quasi-isomorphism
\[
R\Gamma_{\mathrm{cts}}
\left(
I_w,V\otimes_{\Qp}A
\right)
\simeq
\left[
W\otimes_{\Qp}A
\xrightarrow{\gamma-1}
W\otimes_{\Qp}A
\right].
\]
Consequently,
\[
H^i(I_w,V)\otimes_{\Qp}A
\xrightarrow{\sim}
H^i
\left(
I_w,V\otimes_{\Qp}A
\right)
\]
for every \(i\), and both sides vanish for \(i\geq2\). Comparing the
exact sequences
\[
0
\longrightarrow
V^{I_w}
\longrightarrow
V
\longrightarrow
Z^1(I_w,V)
\longrightarrow
H^1(I_w,V)
\longrightarrow
0
\]
before and after scalar extension gives
\[
Z^1(I_w,V)\otimes_{\Qp}A
\xrightarrow{\sim}
Z^1
\left(
I_w,V\otimes_{\Qp}A
\right).
\]
Since \(Z^1(I_w,V)\) is finite-dimensional, the corresponding
cocycle functor is therefore represented by a finite-dimensional
vector group.

Now consider an \(I_w\)-equivariant central extension
\[
1\longrightarrow V
\longrightarrow U_1
\longrightarrow U_0
\longrightarrow1.
\]
For every \(\Qp\)-algebra \(A\), the obstruction to lifting a cocycle
in \(Z^1(I_w,U_0(A))\) lies in
\[
H^2
\left(
I_w,V\otimes_{\Qp}A
\right)
=
0.
\]
Since \(V\) is central, twisting by a \(U_0\)-cocycle does not alter
its \(I_w\)-action. Hence
\[
Z^1(I_w,U_1)
\longrightarrow
Z^1(I_w,U_0)
\]
is a torsor under the vector group representing \(Z^1(I_w,V)\).
Induction along the central series shows that \(Z^1(I_w,U)\) is
represented by an affine scheme of finite type.
Hence
\[
\mathscr H^1(I_w,U)
\simeq
[Z^1(I_w,U)/U]
\]
is an Artin stack of finite type.

By
Proposition~\ref{proposition:restricted_neutral_fibre_presentation},
there is a canonical equivalence
\[
\mathscr Z_w(G/U)
\simeq
[R^{I_w}/G^{I_w}].
\]
The fixed-point group schemes \(G^{I_w}\) and \(R^{I_w}\) are closed
affine subgroup schemes of the finite-type groups \(G\) and \(R\).
Thus \(\mathscr Z_w(G/U)\) is an Artin stack of finite type.

Finally, by
Definition~\ref{definition:relative_neutral_inertia_locus}, the square
\[
\begin{tikzcd}
\mathscr Z_w(G/U)
  \arrow[r]
  \arrow[d]
&
\mathscr B_{I_w}(G)
  \arrow[d]
\\
\mathscr H^1(I_w,U)
  \arrow[r]
&
\mathscr H^1(I_w,G)
\end{tikzcd}
\]
is \(2\)-cartesian. The right vertical morphism is fully faithful,
since \(\mathscr B_{I_w}(G)\) is a strictly full substack.
Consequently, its base change
\[
\mathscr Z_w(G/U)
\longrightarrow
\mathscr H^1(I_w,U)
\]
is a monomorphism of algebraic stacks and is therefore representable
by algebraic spaces. Since its source and target are of finite type
over \(\Qp\), it is of finite type.
\end{proof}

\begin{theorem}[Algebraicity of finite-level neutral Selmer stacks]
\label{theorem:neutral_selmer_algebraicity}
For every
\[
S_p\subseteq T'\subseteq T,
\]
the trivialized neutral Selmer stack
\[
\widetilde{\mathscr N}_{f,T'}(G/R)
\]
is represented by an algebraic space of finite type over \(\Qp\).
Consequently,
\[
\mathscr N_{f,T'}(G/R)
\simeq
\left[
\widetilde{\mathscr N}_{f,T'}(G/R)/R^{G_T}
\right]
\]
is an Artin stack of finite type with affine inertia.
\end{theorem}

\begin{proof}
By
Proposition~\ref{proposition:restricted_neutral_fibre_presentation},
there is a canonical equivalence
\[
\widetilde{\mathscr N}_{f,T'}(G/R)
\simeq
H^1_f(G_T,U)
\times_{
\prod_{w\in T\setminus T'}
\mathscr H^1(I_w,U)
}
\prod_{w\in T\setminus T'}
\mathscr Z_w(G/U).
\]

The scheme \(H^1_f(G_T,U)\) is affine of finite type by
Proposition~\ref{proposition:kp_unipotent_selmer_algebraicity}.
Moreover, by
Proposition~\ref{proposition:local_neutral_locus}, each morphism
\[
\mathscr Z_w(G/U)
\longrightarrow
\mathscr H^1(I_w,U)
\]
is representable and of finite type. The displayed fibre product is
therefore an algebraic space of finite type.

The quotient presentation follows from
Definition~\ref{definition:neutral_fibre_selmer}. Since \(R^{G_T}\)
is a closed affine algebraic subgroup of \(R\), the resulting quotient
is an Artin stack of finite type with affine inertia.
\end{proof}
﻿\section{\texorpdfstring{The \(v\)-adic Period Map and the
Bloch--Kato Logarithm}{The v-adic Period Map and the Bloch--Kato
Logarithm}}
\label{section:dr}
This section develops the Hodge-theoretic foundations of our method.
We construct the finite-level \(v\)-adic period map and neutral Bloch--Kato logarithm,
prove their compatibility with the relative Kummer map, and combine
this compatibility with the density of the period map to obtain the
dimension-cutter criterion.

\begin{itemize}
\item
\S\ref{section:dr_moduli} defines the Frobenius-fixed group, the
Hodge filtration space, and the de Rham moduli stack, and identifies
the latter with the moduli stack of admissible de Rham torsors.

\item
\S\ref{section:analyticity} constructs the finite-level lifted
\(v\)-adic period map on the residue tube from the universal
crystalline and de Rham path torsors and their family-level
comparison. The required Hodge filtration on the universal relative
de Rham path torsor is constructed in
Appendix~\ref{appendix:filtered_de_rham_path_groupoid}.

\item
\S\ref{section:density} compares the \(v\)-adic and complex local
period maps and proves that the lifted period map is algebraically
Zariski dense in the connected component containing the identity
coset.

\item
\S\ref{chapter:bk_log} constructs the local and global neutral
Bloch--Kato logarithms by relative \(p\)-adic Hodge realization and
proves their compatibility with the neutral Kummer maps and the
stack-valued period map.

\item
\S\ref{section:selmer_fibred_prod} forms the Selmer pullback, bounds
the dimension of its image in the Hodge filtration space, and proves
the resulting dimension-cutter criteria.
\end{itemize}

\begin{notation}
\label{notation:finite_level_period_data}
Throughout this section,
\[
1\longrightarrow
\mathcal U_n^\bullet
\longrightarrow
\mathcal G_n^\bullet
\longrightarrow
\mathcal R^\bullet
\longrightarrow1,
\qquad
\bullet\in\{\et,\cri,\dR\},
\]
denotes a compatible finite-type motivic quotient of the relative
completion in the three realizations. Compatibility is understood
with respect to the realization comparison isomorphisms. All
constructions below are natural under morphisms of compatible
finite-type motivic quotients.
\end{notation}

\subsection{The de Rham Moduli Stack}
\label{section:dr_moduli}

We give a finite-level formulation of Kantor's de Rham double quotient and its
interpretation as the moduli stack of torsors satisfying his basic
admissibility condition; cf.~\cite[Definition~6.4.4 and
Proposition~6.4.6]{KANTORTHESIS}.

The crystalline--de Rham comparison of
Theorem~\ref{theorem:crystalline_de_rham_comparison} gives a canonical
isomorphism
\[
\iota_n:
\mathcal G_n^{\cri}\otimes_{K_0}K_v
\xrightarrow{\sim}
\mathcal G_n^{\dR}.
\]
Write
\[
\varphi:
\sigma^*\mathcal G_n^{\cri}
\xrightarrow{\sim}
\mathcal G_n^{\cri}
\]
for crystalline Frobenius, where \(\sigma\) is Frobenius on \(K_0\),
and let
\[
F^0\mathcal G_n^{\dR}
\subseteq
\mathcal G_n^{\dR}
\]
denote the Hodge subgroup of
Proposition~\ref{proposition:de_rham_torsor_filtration}.

\begin{definition}
[Frobenius-fixed group;
cf.~{\cite[Definition~6.4.3]{KANTORTHESIS}}]
\label{definition:de_rham_frobenius_fixed_group}
For every \(\Qp\)-algebra \(A\), set
\[
\mathcal G_n^{\varphi=1}(A)
\coloneqq 
\left\{
g\in
\mathcal G_n^{\cri}
\bigl(A\otimes_{\Qp}K_0\bigr)
\ \middle|\
\varphi(g)=g
\right\}.
\]
This functor is represented by the affine \(\Qp\)-group
\[
\mathcal G_n^{\varphi=1}
=
\left(
\operatorname{Res}_{K_0/\Qp}
\mathcal G_n^{\cri}
\right)^{\varphi=1}.
\]
We similarly write
\[
\mathcal R^{\varphi=1}
\coloneqq 
\left(
\operatorname{Res}_{K_0/\Qp}
\mathcal R^{\cri}
\right)^{\varphi=1}
\]
for the Frobenius-fixed group of the reductive quotient.

Set
\[
\mathcal G_{n,K_v}^{\varphi=1}
\coloneqq 
\mathcal G_n^{\varphi=1}\otimes_{\Qp}K_v.
\]
The natural evaluation morphism, followed by crystalline--de Rham
comparison, gives a canonical closed immersion of \(K_v\)-group
schemes
\[
\mathcal G_{n,K_v}^{\varphi=1}
\hookrightarrow
\mathcal G_n^{\cri}\otimes_{K_0}K_v
\xrightarrow[\sim]{\iota_n}
\mathcal G_n^{\dR}.
\]
\end{definition}

Indeed, if \(f=[K_0:\Qp]\), then after base change to \(K_0\) the
evaluation morphism identifies
\(\mathcal G_n^{\varphi=1}\otimes_{\Qp}K_0\) with the closed
equalizer
\[
\operatorname{Eq}
\left(
\varphi^f,\id:
\mathcal G_n^{\cri}
\rightrightarrows
\mathcal G_n^{\cri}
\right).
\]

\begin{definition}
[{Admissible de Rham torsors;
cf.~\cite[Definition~6.4.4]{KANTORTHESIS}}]
\label{definition:dr_torsors}
Let \(T\) be a \(K_v\)-scheme. An \emph{admissible
\(\mathcal G_n^{\dR}\)-torsor} over \(T\) is a right
\(\mathcal G_{n,T}^{\dR}\)-torsor \(P\), together with reductions of
structure group
\[
P^H
\quad\text{to}\quad
F^0\mathcal G_{n,T}^{\dR},
\qquad
P^\varphi
\quad\text{to}\quad
\mathcal G_{n,K_v,T}^{\varphi=1}.
\]
We call \(P^H\) and \(P^\varphi\) the \emph{Hodge} and
\emph{Frobenius reductions}, respectively. No trivialization of
either reduction is part of the data.
\end{definition}

\begin{definition}[Hodge filtration space and de Rham moduli stack]
\label{definition:de_rham_period_stack}
Define the \emph{Hodge filtration space}
\[
\operatorname{Fil}_n^{\dR}
\coloneqq 
\mathcal G_n^{\dR}/F^0\mathcal G_n^{\dR},
\]
where the quotient is taken in the fppf topology. The preceding
closed immersion induces an action of
\(\mathcal G_{n,K_v}^{\varphi=1}\) on
\(\operatorname{Fil}_n^{\dR}\). The \emph{de Rham moduli stack} is
\[
\mathscr M_n^{\dR}
\coloneqq 
\left[
\mathcal G_{n,K_v}^{\varphi=1}
\backslash
\operatorname{Fil}_n^{\dR}
\right].
\]
Equivalently,
\[
\mathscr M_n^{\dR}
=
\left[
\mathcal G_{n,K_v}^{\varphi=1}
\backslash
\mathcal G_n^{\dR}/F^0\mathcal G_n^{\dR}
\right].
\]
We write
\[
\pi_n:
\operatorname{Fil}_n^{\dR}
\longrightarrow
\mathscr M_n^{\dR}
\]
for the quotient morphism.
\end{definition}

Since \(\mathcal G_n^{\dR}\) is an affine algebraic group over
\(K_v\) and \(F^0\mathcal G_n^{\dR}\) is a closed subgroup, the fppf
quotient
\[
\operatorname{Fil}_n^{\dR}
=
\mathcal G_n^{\dR}/F^0\mathcal G_n^{\dR}
\]
is represented by a quasi-projective \(K_v\)-scheme; see
\cite[Example~1.3]{CONRAD_QUOTIENT_FORMALISM}. The quotient
morphism
\[
\mathcal G_n^{\dR}
\longrightarrow
\operatorname{Fil}_n^{\dR}
\]
is an \(F^0\mathcal G_n^{\dR}\)-torsor. Since
\(F^0\mathcal G_n^{\dR}\) is smooth in characteristic zero, this
torsor is \'{e}tale-locally trivial, and the fppf and \'{e}tale
sheaf quotients agree.

The group \(\mathcal G_{n,K_v}^{\varphi=1}\) is affine, smooth, and
of finite type. Hence
\[
\operatorname{Fil}_n^{\dR}
\longrightarrow
\mathscr M_n^{\dR}
=
\left[
\mathcal G_{n,K_v}^{\varphi=1}
\backslash
\operatorname{Fil}_n^{\dR}
\right]
\]
is a smooth surjective atlas. Thus \(\mathscr M_n^{\dR}\) is an
Artin stack of finite type over \(K_v\). Its inertia is affine,
since after pullback to \(\operatorname{Fil}_n^{\dR}\) it is the
stabilizer, a closed subgroup scheme of
\[
\mathcal G_{n,K_v}^{\varphi=1}
\times
\operatorname{Fil}_n^{\dR}.
\]

\begin{remark}[Comparison with Kantor's notation]
Our notation \(\mathscr M_n^{\dR}\) denotes the full double quotient
classifying torsors satisfying Kantor's basic admissibility
condition. Kantor later writes \(M^{\dR}\) for the further locus of
doubly admissible torsors; cf.~\cite[Definition~6.4.23]{KANTORTHESIS}---that refinement will not be used here.
\end{remark}

\begin{proposition}
[{Moduli interpretation;
cf.~\cite[Proposition~6.4.6]{KANTORTHESIS}}]
\label{proposition:de_rham_period_moduli}
The stack \(\mathscr M_n^{\dR}\) is canonically equivalent to the
moduli stack of admissible \(\mathcal G_n^{\dR}\)-torsors.
\end{proposition}

\begin{proof}
There is a canonical equivalence
\[
\mathscr M_n^{\dR}
\simeq
B\mathcal G_{n,K_v}^{\varphi=1}
\times_{B\mathcal G_n^{\dR}}
B F^0\mathcal G_n^{\dR},
\]
where \(BH\) denotes the classifying stack of \(H\)-torsors. The
right-hand side classifies a \(\mathcal G_n^{\dR}\)-torsor equipped
with reductions of structure group to
\(\mathcal G_{n,K_v}^{\varphi=1}\) and
\(F^0\mathcal G_n^{\dR}\), which is precisely an admissible torsor.
\end{proof}

\begin{remark}[Local description and lifts]
\label{remark:moduli_map}
Let \(P\) be an admissible \(\mathcal G_n^{\dR}\)-torsor over a
\(K_v\)-scheme \(T\). After passing to an fppf cover \(T_1\to T\),
choose sections
\[
p^\varphi\in P^\varphi(T_1),
\qquad
p^H\in P^H(T_1).
\]
There is a unique
\[
g\in\mathcal G_n^{\dR}(T_1)
\]
such that
\[
p^H=p^\varphi\cdot g.
\]
Replacing the chosen sections by
\[
p^\varphi h_\varphi,
\qquad
p^Hh_H,
\]
where
\[
h_\varphi\in
\mathcal G_{n,K_v}^{\varphi=1}(T_1),
\qquad
h_H\in
F^0\mathcal G_n^{\dR}(T_1),
\]
replaces \(g\) by
\[
h_\varphi^{-1}gh_H.
\]
These local representatives, together with their descent data,
therefore define the corresponding object of
\(\mathscr M_n^{\dR}(T)\).

The successive morphisms
\[
\mathcal G_n^{\dR}
\longrightarrow
\operatorname{Fil}_n^{\dR}
\xrightarrow{\pi_n}
\mathscr M_n^{\dR}
\]
have the following moduli interpretation. A lift of the object
represented by \(P\) along \(\pi_n\) is equivalent to a
trivialization of its Frobenius reduction \(P^\varphi\). Relative to
such a lift, a further lift to \(\mathcal G_n^{\dR}\) is equivalent
to a trivialization of its Hodge reduction \(P^H\). Thus the lifted
period map constructed below requires only a canonical
trivialization of the Frobenius reduction.
\end{remark}

\begin{proposition}[Path torsors on a residue disk]
\label{proposition:dr_path_admissible}
Let
\[
x\in]\underline b[(K_v).
\]
Then the finite-level de Rham path torsor
\[
{}^{\vphantom{\dR}}_{x}P_{b,n}^{\dR}
\]
is admissible. Its Frobenius reduction has a canonical \(K_v\)-point,
and hence its class in \(\mathscr M_n^{\dR}(K_v)\) admits a canonical
lift to \(\operatorname{Fil}_n^{\dR}(K_v)\).
\end{proposition}

\begin{proof}
The Hodge filtration on the de Rham path torsor gives a reduction
\[
F^0{}^{\vphantom{\dR}}_{x}P_{b,n}^{\dR}
\subseteq
{}^{\vphantom{\dR}}_{x}P_{b,n}^{\dR}
\]
to \(F^0\mathcal G_n^{\dR}\); see
Proposition~\ref{proposition:de_rham_torsor_filtration}.

Since \(x\) and \(b\) lie in the same residue disk,
\[
\underline x=\underline b.
\]
The corresponding crystalline path torsor is therefore canonically
\[
{}^{\vphantom{\cri}}_{\underline x}P_{\underline b,n}^{\cri}
=
{}^{\vphantom{\cri}}_{\underline b}P_{\underline b,n}^{\cri}
=
\mathcal G_n^{\cri}.
\]
Its identity element is fixed by crystalline Frobenius. Under the
comparison isomorphism
\[
{}^{\vphantom{\cri}}_{\underline x}P_{\underline b,n}^{\cri}
\otimes_{K_0}K_v
\xrightarrow{\sim}
{}^{\vphantom{\dR}}_{x}P_{b,n}^{\dR},
\]
it determines a canonical point
\[
p_x^\varphi
\in
{}^{\vphantom{\dR}}_{x}P_{b,n}^{\dR}(K_v).
\]
The orbit
\[
p_x^\varphi\mathcal G_{n,K_v}^{\varphi=1}
\subseteq
{}^{\vphantom{\dR}}_{x}P_{b,n}^{\dR}
\]
is the Frobenius reduction, and \(p_x^\varphi\) gives its canonical
trivialization. Together with the Hodge reduction, this proves
admissibility and, by
Remark~\ref{remark:moduli_map}, determines the asserted lift to
\(\operatorname{Fil}_n^{\dR}(K_v)\).
\end{proof}

These constructions are natural under morphisms of compatible
finite-type motivic quotients. Hence, for every compatible
cofiltered system \(\{\mathcal G_n^\bullet\}_n\), the Hodge
filtration spaces, de Rham moduli stacks, and canonical pointwise
lifts form compatible inverse systems. The next subsection proves
that the resulting maps
\[
]\underline b[(K_v)
\longrightarrow
\operatorname{Fil}_n^{\dR}(K_v)
\]
extend to rigid-analytic maps on the residue tube.

\subsection{\texorpdfstring{Analyticity of the \(v\)-adic Period Map}{Analyticity of the v-adic Period Map}}
\label{section:analyticity}

This subsection shows that the canonical pointwise lifts of
Proposition~\ref{proposition:dr_path_admissible} assemble into
compatible rigid-analytic maps on the residue tube
\[
]\underline b[
\subseteq
X_v^{\an}.
\]
We first construct the universal de Rham and crystalline path torsors
separately, then compare their realizations over the residue tube.
Horizontal continuation of the crystalline identity gives the
canonical Frobenius trivialization used, together with the universal
Hodge reduction, to define the lifted period map.

\begin{itemize}
\item
\S\ref{subsubsection:universal_de_rham_path_torsor} records the
finite-level \(K_v\)-form of the universal de Rham path torsor,
together with its canonical connection, Hodge filtration, and Hodge
reduction. These structures are obtained by base change from
Appendix~\ref{appendix:filtered_de_rham_path_groupoid}.

\item
\S\ref{subsubsection:universal_crystalline_path_torsor} constructs
the universal crystalline path torsor intrinsically in isocrystals,
together with its intrinsic Frobenius structure, and, after extension
of coefficients to \(K_v\), realizes it on the fixed formal
\(\mathcal O_{K_v}\)-enlargement as a rigid-analytic torsor with
convergent connection.

\item
\S\ref{subsubsection:comparison_residue_tube} establishes the
family-level crystalline--de Rham comparison over the residue tube.
It horizontally extends the identity of the crystalline fibre and
transports the resulting horizontal section to the universal de Rham
torsor, recovering at every \(K_v\)-point the canonical Frobenius
trivialization of
Proposition~\ref{proposition:dr_path_admissible}.

\item
\S\ref{subsubsection:v_adic_period_map} combines this horizontal
trivialization with the universal Hodge reduction to define the
finite-level \(v\)-adic period map. It proves rigid analyticity,
identifies its values with the canonical pointwise lifts, and records
compatibility under morphisms of finite-level motivic quotients.
\end{itemize}

\subsubsection{The Universal de Rham Path Torsor}
\label{subsubsection:universal_de_rham_path_torsor}

The universal relative de Rham path torsor, its canonical connection,
and its family-level Hodge filtration are constructed over \(K\) in
Appendix~\ref{appendix:filtered_de_rham_path_groupoid}. We record
here the finite-level \(K_v\)-forms used in the construction of the
period map; see especially
\S\ref{appendix:finite_level_filtered_path_torsors}.

Let
\[
\mathcal T_n^{\dR}
\subseteq
\langle\mathcal L^{\dR}\rangle_\otimes^{\mathrm{ext}}
\]
denote the full Tannakian subcategory corresponding to the quotient
\[
\mathcal G^{\dR}\twoheadrightarrow\mathcal G_n^{\dR}.
\]
In the genus $\ge 2$ case, the constructions below are formed on
\(\widehat X_v\) and then restricted to \(X_v\), in accordance with
Convention~\ref{convention:fixed_kp_family}.

Consider the tensor functors
\[
\omega_{b,n}^{\dR}\otimes_{K_v}\mathcal O_{X_v},
\quad
\omega_{X_v,n}^{\dR}:
\mathcal T_n^{\dR}
\longrightarrow
\operatorname{VB}(X_v),
\]
where
\[
\omega_{b,n}^{\dR}
\coloneqq 
\left.\omega_b^{\dR}\right|_{\mathcal T_n^{\dR}},
\]
and \(\omega_{X_v,n}^{\dR}\) sends a bundle with integrable
connection to its underlying vector bundle.

\begin{definition}
[{Universal de Rham path torsor;
cf.~\S\ref{appendix:finite_level_filtered_path_torsors}}]
\label{definition:universal_de_rham_path_torsor}
The \emph{universal finite-level de Rham path torsor based at \(b\)}
is
\[
\mathcal P_n^{\dR}
\coloneqq 
\operatorname{Isom}^{\otimes}
\left(
\omega_{b,n}^{\dR}\otimes_{K_v}\mathcal O_{X_v},
\omega_{X_v,n}^{\dR}
\right).
\]
\end{definition}

Let
\[
p_{n,K}:
\mathcal P_{n,b,-,K}^{\dR}
\longrightarrow
X^\dagger
\]
denote the structure morphism of the global finite-level torsor
constructed in
\S\ref{appendix:finite_level_filtered_path_torsors}, and set
\[
\mathscr A_{n,b,K}^{\dR}
\coloneqq 
(p_{n,K})_*
\mathcal O_{\mathcal P_{n,b,-,K}^{\dR}}.
\]
By
Corollary~\ref{corollary:finite_level_filtered_path_groupoid}, there
is a canonical identification
\[
\mathcal P_n^{\dR}
\simeq
\left.
\left(
\mathcal P_{n,b,-,K}^{\dR}
\otimes_KK_v
\right)
\right|_{X_v}.
\]

Write
\[
p_n:
\mathcal P_n^{\dR}
\longrightarrow
X_v
\]
for the structure morphism, and set
\[
\mathscr A_n^{\dR}
\coloneqq 
(p_n)_*
\mathcal O_{\mathcal P_n^{\dR}}.
\]
Equivalently,
\[
\mathscr A_n^{\dR}
\simeq
\left.
\left(
\mathscr A_{n,b,K}^{\dR}
\otimes_KK_v
\right)
\right|_{X_v}.
\]

\begin{proposition}
[{Universal finite-level de Rham path torsor;
see
Proposition~\ref{proposition:tannakian_path_torsor_connection}
and
Corollary~\ref{corollary:finite_level_filtered_path_groupoid}}]
\label{proposition:universal_de_rham_path_torsor}
The scheme \(\mathcal P_n^{\dR}\) is an affine right
\(\mathcal G_n^{\dR}\)-torsor over \(X_v\), equipped with its
canonical integrable connection. Dually,
\(\mathscr A_n^{\dR}\) is a commutative algebra object in
\[
\operatorname{Ind}
\operatorname{VIC}(X_v/K_v).
\]

For every \(K_v\)-scheme \(T\) and every morphism
\[
x:T\longrightarrow X_v,
\]
there is a canonical identification
\[
x^*\mathcal P_n^{\dR}
\xrightarrow{\sim}
{}^{\vphantom{\dR}}_{x}P_{b,n}^{\dR}.
\]
In particular,
\[
\left.
\mathcal P_n^{\dR}
\right|_x
\simeq
{}^{\vphantom{\dR}}_{x}P_{b,n}^{\dR}
\qquad
\bigl(x\in X_v(K_v)\bigr).
\]
\end{proposition}

\begin{proposition}
[{Universal finite-level Hodge filtration;
see
Corollary~\ref{corollary:finite_level_filtered_path_groupoid}}]
\label{proposition:universal_de_rham_hodge_filtration}
The algebra \(\mathscr A_n^{\dR}\) carries a canonical decreasing,
exhaustive, and multiplicative Hodge filtration
\[
\cdots
\supseteq
F^p\mathscr A_n^{\dR}
\supseteq
F^{p+1}\mathscr A_n^{\dR}
\supseteq
\cdots .
\]
It is the base change of the global finite-level filtration:
\[
F^p\mathscr A_n^{\dR}
\simeq
\left.
\left(
F^p\mathscr A_{n,b,K}^{\dR}
\otimes_KK_v
\right)
\right|_{X_v},
\]
and has the following properties.

\begin{enumerate}
\item
The canonical connection satisfies Griffiths transversality:
\[
\nabla
\left(
F^p\mathscr A_n^{\dR}
\right)
\subseteq
F^{p-1}\mathscr A_n^{\dR}
\otimes_{\mathcal O_{X_v}}
\Omega^1_{X_v/K_v}.
\]

\item
The torsor coaction
\[
\Delta_n:
\mathscr A_n^{\dR}
\longrightarrow
\mathscr A_n^{\dR}
\otimes_{K_v}
\mathcal O(\mathcal G_n^{\dR})
\]
is a morphism of filtered algebras, where
\(\mathcal O(\mathcal G_n^{\dR})\) carries the filtration obtained by
specializing at \(b\).

\item
The filtration commutes with base change on \(X_v\). For every
\(K_v\)-scheme \(T\) and every \(x:T\to X_v\), the induced filtration
on
\[
\mathcal O
\left(
{}^{\vphantom{\dR}}_{x}P_{b,n}^{\dR}
\right)
\]
is the finite-level Hodge filtration on the de Rham path torsor.

\item
The construction is natural under morphisms of compatible
finite-type motivic quotients.
\end{enumerate}
\end{proposition}

\begin{definition}
[{Universal Hodge locus;
cf.~Corollary~\ref{corollary:finite_level_universal_hodge_reduction}}]
\label{definition:universal_hodge_locus}
Let
\[
\mathcal P_{n,b,-,K}^{H}
\subseteq
\mathcal P_{n,b,-,K}^{\dR}
\]
be the global finite-level Hodge reduction of
Corollary~\ref{corollary:finite_level_universal_hodge_reduction}.
Define
\[
\mathcal P_n^H
\coloneqq 
\left.
\left(
\mathcal P_{n,b,-,K}^{H}
\otimes_KK_v
\right)
\right|_{X_v}
\subseteq
\mathcal P_n^{\dR}.
\]
We also write
\[
\mathcal P_n^H
=
F^0\mathcal P_n^{\dR}.
\]

Equivalently, for a \(K_v\)-scheme \(T\) and a morphism
\(x:T\to X_v\), the \(T\)-points of \(x^*\mathcal P_n^H\) are the
points of \(x^*\mathcal P_n^{\dR}\) whose corresponding algebra
morphisms
\[
x^*\mathscr A_n^{\dR}
\longrightarrow
\mathcal O_T
\]
preserve the Hodge filtration, where \(\mathcal O_T\) carries the
trivial filtration.
\end{definition}

\begin{proposition}
[{Universal Hodge reduction;
see
Corollary~\ref{corollary:finite_level_universal_hodge_reduction}}]
\label{proposition:universal_hodge_reduction}
The closed subscheme
\[
\mathcal P_n^H
\subseteq
\mathcal P_n^{\dR}
\]
is a reduction of structure group to
\(F^0\mathcal G_n^{\dR}\). Thus the restricted action induces an
isomorphism
\[
\mathcal P_n^H
\times_{K_v}
F^0\mathcal G_n^{\dR}
\xrightarrow{\sim}
\mathcal P_n^H
\times_{X_v}
\mathcal P_n^H,
\]
and
\[
\mathcal P_n^H
\longrightarrow
X_v
\]
is an fppf torsor under \(F^0\mathcal G_n^{\dR}\).

Its formation commutes with base change on \(X_v\) and with
morphisms of compatible finite-type motivic quotients. For every
\(K_v\)-scheme \(T\) and every morphism \(x:T\to X_v\), there is a
canonical identification
\[
x^*\mathcal P_n^H
\xrightarrow{\sim}
F^0{}^{\vphantom{\dR}}_{x}P_{b,n}^{\dR}.
\]
\end{proposition}

\subsubsection{The Universal Crystalline Path Torsor}
\label{subsubsection:universal_crystalline_path_torsor}
This construction is the crystalline analogue of
\S\ref{subsubsection:universal_de_rham_path_torsor}.
The universal object is first defined intrinsically as an affine
scheme object in isocrystals, together with its intrinsic Frobenius
structure. After extension of coefficients to \(K_v\), evaluation on
the fixed formal \(\mathcal O_{K_v}\)-enlargement gives a
rigid-analytic torsor with convergent connection. Its compatibility
with the universal de Rham torsor is established in
Proposition~\ref{proposition:universal_crystalline_de_rham_comparison},
using the comparison results of
\S\ref{section:crystalline_de_rham_comparison_appendix}.

Let
\[
\mathcal T_n^{\cri}
\subseteq
\langle\mathcal L^{\cri}\rangle_\otimes^{\mathrm{ext}}
\]
denote the full Tannakian subcategory corresponding to the quotient
\[
\mathcal G^{\cri}\twoheadrightarrow\mathcal G_n^{\cri}.
\]
Write
\[
\omega_{\underline b,n}^{\cri}
\coloneqq 
\left.
\omega_{\underline b}^{\cri}
\right|_{\mathcal T_n^{\cri}},
\]
and let
\[
\omega_{X_k,n}^{\cri}:
\mathcal T_n^{\cri}
\longrightarrow
\operatorname{Isoc}(X_k/K_0)
\]
be the tautological tensor functor. The fibre functor at
\(\underline b\) determines the constant-isocrystal functor
\[
\omega_{\underline b,n}^{\cri}
\otimes_{K_0}\mathbf 1_{X_k}.
\]

\begin{definition}
[{Universal crystalline path torsor;
cf.~\cite[\S5]{DELIGNE_MONOGRAPH}}]
\label{definition:universal_crystalline_path_torsor}
The \emph{universal finite-level crystalline path torsor based at
\(\underline b\)} is the affine scheme object
\[
\mathcal P_n^{\cri}
\coloneqq 
\operatorname{Isom}^{\otimes}
\left(
\omega_{\underline b,n}^{\cri}
\otimes_{K_0}\mathbf 1_{X_k},
\omega_{X_k,n}^{\cri}
\right)
\]
in \(\operatorname{Isoc}(X_k/K_0)\). Its coordinate
algebra is the corresponding commutative algebra object in
\[
\operatorname{Ind}\operatorname{Isoc}(X_k/K_0).
\]
\end{definition}

Set
\[
\underline{\mathcal G}_{n,X_k}^{\cri}
\coloneqq 
\mathcal G_n^{\cri}
\otimes_{K_0}\mathbf 1_{X_k},
\]
viewed as the constant affine group scheme object in
\(\operatorname{Isoc}(X_k/K_0)\).

\begin{proposition}[The intrinsic crystalline torsor]
\label{proposition:universal_crystalline_path_torsor}
The object \(\mathcal P_n^{\cri}\) is a right
\(\underline{\mathcal G}_{n,X_k}^{\cri}\)-torsor in affine scheme
objects in \(\operatorname{Isoc}(X_k/K_0)\). For every
\(\underline x\in X_k(k)\), there is a canonical identification of
affine \(K_0\)-schemes
\[
\underline x^*\mathcal P_n^{\cri}
\xrightarrow{\sim}
{}^{\vphantom{\cri}}_{\underline x}P_{\underline b,n}^{\cri}.
\]
In particular,
\[
\underline b^*\mathcal P_n^{\cri}
\simeq
\mathcal G_n^{\cri}.
\]

The Frobenius pullback on \(\mathcal T_n^{\cri}\), together with its
compatibility with the two tensor functors, induces a Frobenius
isomorphism
\[
\varphi_{\mathcal P_n}:
\sigma^*\mathcal P_n^{\cri}
\xrightarrow{\sim}
\mathcal P_n^{\cri},
\]
equivariant with respect to crystalline Frobenius on
\(\underline{\mathcal G}_{n,X_k}^{\cri}\).
\end{proposition}

\begin{proof}
Relative Tannakian duality gives the torsor structure. Formation of
the tensor-isomorphism torsor commutes with pullback, yielding the
fibre identifications. Applying
\(\operatorname{Isom}^{\otimes}\) to the Frobenius-transformed tensor
functors gives the final assertion.
\end{proof}

\paragraph{Rigid realization on the fixed formal enlargement.}

Choose an affine open neighborhood
\[
\underline b\in U_k\subseteq X_k.
\]
Let \(\mathfrak X_v\) denote the \(p\)-adic formal completion of the
fixed smooth \(\mathcal O_{K_v}\)-model along its special fibre, and
let
\[
\mathfrak U_v
\coloneqq 
\left.\mathfrak X_v\right|_{U_k}
\]
be the open formal subscheme with special fibre \(U_k\). Write
\[
\mathfrak U_{v,\eta}
\]
for its rigid generic fibre.

After extension of coefficients from \(K_0\) to \(K_v\), evaluation
on the enlargement \(\mathfrak U_v\) gives a \(K_v\)-linear tensor
functor
\[
\omega_{\mathfrak U_v,n}^{\cri,\nabla}:
\mathcal T_n^{\cri}\otimes_{K_0}K_v
\longrightarrow
\operatorname{VIC}^{\mathrm{conv}}
\bigl(
\mathfrak U_{v,\eta}/K_v
\bigr).
\]
Let
\[
\omega_{\mathfrak U_v,n}^{\cri}:
\mathcal T_n^{\cri}\otimes_{K_0}K_v
\longrightarrow
\operatorname{VB}(\mathfrak U_{v,\eta})
\]
denote the underlying vector-bundle functor.

\begin{definition}[Rigid realization]
\label{definition:rigid_universal_crystalline_torsor}
For a rigid-analytic space
\[
f:Y\longrightarrow\mathfrak U_{v,\eta},
\]
set
\[
\mathcal P_{n,\mathfrak U_v}^{\cri,\rig}(Y)
\coloneqq 
\operatorname{Isom}^{\otimes}
\left(
f^*
\left(
\left(
\omega_{\underline b,n}^{\cri}
\otimes_{K_0}K_v
\right)
\otimes_{K_v}
\mathcal O_{\mathfrak U_{v,\eta}}
\right),
f^*\omega_{\mathfrak U_v,n}^{\cri}
\right),
\]
where the right-hand side denotes tensor-natural isomorphisms of
functors
\[
\mathcal T_n^{\cri}\otimes_{K_0}K_v
\longrightarrow
\operatorname{VB}(Y).
\]
This presheaf on the big rigid-analytic \'{e}tale site of
\(\mathfrak U_{v,\eta}\) is the \emph{rigid realization}, after
extension of coefficients to \(K_v\), of
\(\mathcal P_n^{\cri}|_{U_k}\) on the fixed formal enlargement
\(\mathfrak U_v\).
\end{definition}

\begin{proposition}[Representability of the rigid realization]
\label{proposition:rigid_universal_crystalline_torsor}
The presheaf
\(\mathcal P_{n,\mathfrak U_v}^{\cri,\rig}\) is an \'{e}tale sheaf
represented by a rigid-analytic space locally of finite type over
\(\mathfrak U_{v,\eta}\). It is a right torsor on the \'{e}tale site
under the constant rigid-analytic group
\[
\mathcal G_{n,\mathfrak U_{v,\eta}}^{\cri,\an}
\coloneqq 
\left(
\mathcal G_n^{\cri}\otimes_{K_0}K_v
\right)^{\an}
\times_{\operatorname{Sp}K_v}
\mathfrak U_{v,\eta}.
\]

More precisely, if
\[
V=\operatorname{Sp}(A)
\subseteq
\mathfrak U_{v,\eta}
\]
is an affinoid subdomain, then the restriction of the torsor to \(V\)
is represented by the relative analytification of the affine
\(A\)-scheme
\[
\operatorname{Isom}^{\otimes}
\left(
\left(
\omega_{\underline b,n}^{\cri}
\otimes_{K_0}K_v
\right)
\otimes_{K_v}A,
\Gamma
\left(
V,
\omega_{\mathfrak U_v,n}^{\cri}
\right)
\right)
\]
provided by relative Tannakian duality.

If \(\widetilde x\) is a rigid point of
\(\mathfrak U_{v,\eta}\), with residue field
\(K(\widetilde x)\), specializing to
\(\underline x\in U_k(k)\), then there is a canonical identification
\[
\left(
\mathcal P_{n,\mathfrak U_v}^{\cri,\rig}
\right)_{\widetilde x}
\xrightarrow{\sim}
\left(
{}^{\vphantom{\cri}}_{\underline x}P_{\underline b,n}^{\cri}
\otimes_{K_0}K(\widetilde x)
\right)^{\an}.
\]
\end{proposition}

\begin{proof}
Tensor isomorphisms of rigid-analytic vector bundles satisfy
\'{e}tale descent, and tensor compatibility is an \'{e}tale-local
condition. On an affinoid subdomain, relative Tannakian duality
represents the tensor-isomorphism functor by an affine scheme of
finite presentation; see
\cite[\S3]{DelMil_TC22}. Its relative analytification gives the
stated local representing space, and these spaces glue because they
represent the restrictions of the same sheaf.

Relative Tannakian duality also gives the torsor structure. Since
\(\mathcal G_n^{\cri}\otimes_{K_0}K_v\) is smooth, the torsor is
\'{e}tale-locally trivial. Finally, tensor-isomorphism torsors commute
with base change, while the crystal property identifies the fibre at
\(\widetilde x\) with
\[
\omega_{\underline x,n}^{\cri}
\otimes_{K_0}K(\widetilde x),
\]
giving the asserted fibre identification.
\end{proof}

By the crystal property, this realization is compatible with
restriction to smaller open neighborhoods of \(\underline b\).

\begin{proposition}
[Connection;
cf.~Proposition~\ref{proposition:tannakian_path_torsor_connection}]
\label{proposition:rigid_universal_crystalline_connection}
The torsor
\(\mathcal P_{n,\mathfrak U_v}^{\cri,\rig}\) carries a canonical
convergent integrable connection for which the torsor action is
horizontal.
\end{proposition}

\begin{proof}
The tensor functor
\(\omega_{\mathfrak U_v,n}^{\cri,\nabla}\) takes values in
rigid-analytic vector bundles with convergent integrable connection.
The rigid-analytic form of the tensor-isomorphism construction in
Proposition~\ref{proposition:tannakian_path_torsor_connection}
therefore gives the canonical connection and the horizontality of the
torsor action.
\end{proof}

\subsubsection{Comparison and Horizontal Trivialization on the Residue Tube}
\label{subsubsection:comparison_residue_tube}

Retain \(\mathfrak X_v\) and \(\mathfrak U_v\) from above. The
specialization map identifies
\[
]\underline b[
=
\operatorname{sp}^{-1}(\underline b)
\subseteq
\mathfrak U_{v,\eta}
\]
with the residue tube in \(X_v^{\an}\).

\begin{notation}
For a finite-type \(X_v\)-scheme \(Z\), write
\[
Z_{]\underline b[}^{\an}
\coloneqq 
]\underline b[
\times_{X_v^{\an}}
Z^{\an}
\]
for the restriction of its rigid analytification to
\(]\underline b[\).
\end{notation}

Set
\[
\mathcal P_{n,]\underline b[}^{\dR,\an}
\coloneqq 
\left(\mathcal P_n^{\dR}\right)_{]\underline b[}^{\an},
\qquad
\mathcal P_{n,]\underline b[}^{H,\an}
\coloneqq 
\left(\mathcal P_n^H\right)_{]\underline b[}^{\an}.
\]
Set
\[
\mathcal P_{n,]\underline b[}^{\cri,\rig}
\coloneqq 
\left.
\mathcal P_{n,\mathfrak U_v}^{\cri,\rig}
\right|_{]\underline b[}.
\]

\begin{proposition}
[{Family-level crystalline--de Rham comparison;
cf.~\cite[Chapter~4]{OLSSON_NA_P_HODGE}}]
\label{proposition:universal_crystalline_de_rham_comparison}
There is a canonical isomorphism of rigid-analytic torsors over
\(]\underline b[\),
\[
\iota_{\mathcal P_n}:
\mathcal P_{n,]\underline b[}^{\cri,\rig}
\xrightarrow{\sim}
\mathcal P_{n,]\underline b[}^{\dR,\an}.
\]
Its structure-group isomorphism is the analytification of
\[
\mathcal G_n^{\cri}\otimes_{K_0}K_v
\xrightarrow[\sim]{\iota_n}
\mathcal G_n^{\dR}.
\]
The comparison is compatible with the integrable connections and,
for every \(x\in]\underline b[(K_v)\), restricts to the
analytification of the pointwise comparison
\[
{}^{\vphantom{\cri}}_{\underline x}P_{\underline b,n}^{\cri}
\otimes_{K_0}K_v
\xrightarrow{\sim}
{}^{\vphantom{\dR}}_{x}P_{b,n}^{\dR}.
\]
In particular, it identifies the identity elements in the fibres at
\(b\).
\end{proposition}

\begin{proof}
Set
\[
T\coloneqq ]\underline b[.
\]

Because \(\mathcal G_n^\bullet\) is a compatible finite-level
motivic quotient, the crystalline--de Rham realization equivalence
of
Theorem~\ref{theorem:olsson_crystalline_de_rham_comparison}
restricts to a tensor equivalence
\[
\Psi_n:
\mathcal T_n^{\cri}\otimes_{K_0}K_v
\xrightarrow{\sim}
\mathcal T_n^{\dR}.
\]

The realization functor in
Theorem~\ref{theorem:olsson_crystalline_de_rham_comparison}
is defined by evaluation on the fixed formal
\(\mathcal O_{K_v}\)-enlargement. Consequently, after restriction to
\(T\), it gives a tensor-natural horizontal isomorphism
\[
c_{\mathfrak U_v,n}:
\left.
\omega_{\mathfrak U_v,n}^{\cri,\nabla}
\right|_T
\xrightarrow{\sim}
\left(
\omega_{X_v,n}^{\dR,\nabla}
\right)^{\an}\big|_T
\circ\Psi_n,
\]
where
\[
\omega_{X_v,n}^{\dR,\nabla}:
\mathcal T_n^{\dR}
\longrightarrow
\operatorname{VIC}(X_v/K_v)
\]
is the tautological connection-valued de Rham functor.

Pullback along the lift \(b\in T(K_v)\) gives a compatible tensor
isomorphism of fibre functors
\[
c_{b,n}:
\omega_{\underline b,n}^{\cri}
\otimes_{K_0}K_v
\xrightarrow{\sim}
\omega_{b,n}^{\dR}\circ\Psi_n.
\]

After identifying the finite-level categories through \(\Psi_n\),
let
\[
f:Y\longrightarrow T
\]
be a rigid-analytic space over \(T\). For a tensor isomorphism
\[
\alpha:
f^*
\left(
\left(
\omega_{\underline b,n}^{\cri}
\otimes_{K_0}K_v
\right)
\otimes_{K_v}\mathcal O_T
\right)
\xrightarrow{\sim}
f^*
\left(
\left.
\omega_{\mathfrak U_v,n}^{\cri}
\right|_T
\right),
\]
define
\[
\alpha^{\dR}
\coloneqq 
f^*c_{\mathfrak U_v,n}
\circ
\alpha
\circ
\left(
f^*c_{b,n}
\right)^{-1}.
\]
This construction is functorial in \(Y\), compatible with morphisms,
tensor products, duals, and the tensor unit, and invertible by the
same formula using the inverse comparison isomorphisms. It therefore
defines an isomorphism of tensor-isomorphism sheaves on the big
rigid-analytic \'{e}tale site of \(T\):
\[
\mathcal P_{n,]\underline b[}^{\cri,\rig}
\xrightarrow{\sim}
\mathcal P_{n,]\underline b[}^{\dR,\an}.
\]
Since both sheaves are represented, this is the asserted
isomorphism \(\iota_{\mathcal P_n}\) of rigid-analytic torsors.

Applying the same construction to tensor automorphisms of the
constant fibre functor gives the structure-group isomorphism
\[
\left(
\mathcal G_n^{\cri}
\otimes_{K_0}K_v
\right)^{\an}
\xrightarrow{\sim}
\left(
\mathcal G_n^{\dR}
\right)^{\an}.
\]
By Tannakian duality, this is the analytification of
\(\iota_n\), and the preceding conjugation formula is equivariant
for these group actions.

The comparison maps
\(c_{\mathfrak U_v,n}\) and \(c_{b,n}\) are horizontal. Functoriality
of the induced Hom connections therefore shows that
\(\iota_{\mathcal P_n}\) is horizontal.

Finally, let
\[
x:\operatorname{Sp}K_v\longrightarrow T
\]
be a \(K_v\)-point, specializing to
\(\underline x\in X_k(k)\). The fibre identifications of
Propositions~\ref{proposition:rigid_universal_crystalline_torsor}
and~\ref{proposition:universal_de_rham_path_torsor} give
\[
\left(
\mathcal P_{n,]\underline b[}^{\cri,\rig}
\right)_x
\simeq
\left(
{}^{\vphantom{\cri}}_{\underline x}P_{\underline b,n}^{\cri}
\otimes_{K_0}K_v
\right)^{\an}
\]
and
\[
\left(
\mathcal P_{n,]\underline b[}^{\dR,\an}
\right)_x
\simeq
\left(
{}^{\vphantom{\dR}}_{x}P_{b,n}^{\dR}
\right)^{\an}.
\]
Compatibility of the realization comparison on the fixed formal
enlargement with pullback to \(x\) identifies the restriction of
\(\iota_{\mathcal P_n}\) with the analytification of the canonical
pointwise comparison
\[
\iota_{x,b,n}:
{}^{\vphantom{\cri}}_{\underline x}P_{\underline b,n}^{\cri}
\otimes_{K_0}K_v
\xrightarrow{\sim}
{}^{\vphantom{\dR}}_{x}P_{b,n}^{\dR}.
\]

At \(x=b\), the variable-endpoint comparison \(c_{x,n}\) is
\(c_{b,n}\). Hence the formula
\[
\alpha\longmapsto
c_{x,n}\circ\alpha\circ c_{b,n}^{-1}
\]
sends the identity tensor isomorphism to the identity.
\end{proof}

\begin{lemma}[Canonical horizontal trivialization]
\label{lemma:frobenius_horizontal_trivialization}
The identity element
\[
1\in
\left(
\mathcal P_{n,]\underline b[}^{\cri,\rig}
\right)_b(K_v)
\simeq
\left(
\mathcal G_n^{\cri}\otimes_{K_0}K_v
\right)(K_v)
\]
extends uniquely to a horizontal section
\[
s_n^{\cri}:
]\underline b[
\longrightarrow
\mathcal P_{n,]\underline b[}^{\cri,\rig}.
\]

Transport through \(\iota_{\mathcal P_n}\) gives a horizontal section
\[
s_n^\varphi:
]\underline b[
\longrightarrow
\mathcal P_{n,]\underline b[}^{\dR,\an}.
\]
For every \(x\in]\underline b[(K_v)\), let
\[
\iota_{x,b,n}:
{}^{\vphantom{\cri}}_{\underline x}P_{\underline b,n}^{\cri}
\otimes_{K_0}K_v
\xrightarrow{\sim}
{}^{\vphantom{\dR}}_{x}P_{b,n}^{\dR}
\]
denote the pointwise crystalline--de Rham comparison. Since
\(\underline x=\underline b\), there is a canonical identity element
\[
1\in
{}^{\vphantom{\cri}}_{\underline x}P_{\underline b,n}^{\cri}
=
\mathcal G_n^{\cri},
\]
and
\[
s_n^\varphi(x)
=
\iota_{x,b,n}(1).
\]
In particular, \(s_n^\varphi(b)=1\), and
\(s_n^\varphi(x)\) is the canonical trivializing point of the
Frobenius reduction of \({}^{\vphantom{\dR}}_{x}P_{b,n}^{\dR}\) constructed in
Proposition~\ref{proposition:dr_path_admissible}.
\end{lemma}

\begin{proof}
For \(E\in\mathcal T_n^{\cri}\otimes_{K_0}K_v\), write
\[
E_{\mathfrak U_v}
\coloneqq 
\omega_{\mathfrak U_v,n}^{\cri,\nabla}(E).
\]
Its convergent connection determines a convergent stratification on
the tube of the diagonal. Pullback along
\[
(b,\operatorname{id}_{]\underline{b}[}):
]\underline{b}[\longrightarrow]\Delta_{U_k}[
\]
therefore gives a horizontal Taylor isomorphism
\[
\operatorname{Tay}_{b,-,E}:
(E_{\mathfrak U_v})_b\otimes_{K_v}\mathcal O_{]\underline{b}[}
\xrightarrow{\sim}
E_{\mathfrak U_v}|_{]\underline{b}[},
\]
normalized to be the identity at \(b\).

Let
\[
\lambda_{b,E}:
\omega_{\underline b,n}^{\cri}(E)\otimes_{K_0}K_v
\xrightarrow{\sim}
(E_{\mathfrak U_v})_b
\]
be the canonical crystal-evaluation isomorphism. The tensor-natural
family
\[
\operatorname{Tay}_{b,-,E}\circ\lambda_{b,E}
\]
defines \(s_n^{\cri}\). Its uniqueness follows from uniqueness of
horizontal continuation from \(b\).

Now let \(x\in ]\underline{b}[(K_v)\), and let
\[
\lambda_{x,E}:
\omega_{\underline x,n}^{\cri}(E)\otimes_{K_0}K_v
\xrightarrow{\sim}
(E_{\mathfrak U_v})_x
\]
be the crystal-evaluation isomorphism used in the canonical
identification of the fibre at \(x\). Since
\(\underline x=\underline b\), the cocycle identity for the
stratification gives
\[
\operatorname{Tay}_{b,x,E}\circ\lambda_{b,E}
=
\lambda_{x,E}.
\]
The fibre identification sends a tensor isomorphism \(\alpha\) to
\(\lambda_{x,E}^{-1}\circ\alpha_E\). It therefore sends
\(s_n^{\cri}(x)\) to
\[
\lambda_{x,E}^{-1}
\circ
\operatorname{Tay}_{b,x,E}
\circ
\lambda_{b,E}
=
\operatorname{id}.
\]
Thus \(s_n^{\cri}(x)=1\) in
\({}^{\vphantom{\cri}}_{\underline x}P_{\underline b,n}^{\cri}\).

Transport through the horizontal comparison
\(\iota_{\mathcal P_n}\) defines \(s_n^\varphi\). Compatibility of
the family-level comparison with pullback to \(x\) gives
\[
s_n^\varphi(x)=\iota_{x,b,n}(1).
\]
At \(x=b\), the comparison preserves the identity; for general \(x\),
the right-hand side is precisely the point \(p_x^\varphi\) constructed
in Proposition~\ref{proposition:dr_path_admissible}.
\end{proof}

\subsubsection{\texorpdfstring{The \(v\)-adic Period Map}{The v-adic Period Map}}
\label{subsubsection:v_adic_period_map}

The section \(s_n^\varphi\) gives a rigid-analytic trivialization
\[
\begin{aligned}
\tau_n:
]\underline b[
\times
\left(\mathcal G_n^{\dR}\right)^{\an}
&\xrightarrow{\sim}
\mathcal P_{n,]\underline b[}^{\dR,\an},
\\
(x,g)
&\longmapsto
s_n^\varphi(x)\cdot g.
\end{aligned}
\]
Under this trivialization, the Hodge reduction
\[
\mathcal P_{n,]\underline b[}^{H,\an}
\subseteq
\mathcal P_{n,]\underline b[}^{\dR,\an}
\]
determines a section of the associated homogeneous-space bundle.

\begin{definition}[Finite-level \(v\)-adic period map]
\label{definition:v_adic_period_map}
The map
\[
\Phi_{n,v}:
]\underline b[
\longrightarrow
\left(\operatorname{Fil}_n^{\dR}\right)^{\an}
\]
determined by the trivialization \(\tau_n\) and the Hodge reduction
is the \emph{finite-level \(v\)-adic period map}. The corresponding
stack-valued period map is
\[
[\Phi_{n,v}]
\coloneqq 
\pi_n^{\an}\circ\Phi_{n,v}:
]\underline b[
\longrightarrow
\left(\mathscr M_n^{\dR}\right)^{\an}.
\]
\end{definition}

\begin{proposition}[Analyticity and pointwise interpretation]
\label{proposition:an_rep}
The map \(\Phi_{n,v}\) is rigid analytic. For every
\(x\in]\underline b[(K_v)\),
\[
[\Phi_{n,v}](x)
=
[{}^{\vphantom{\dR}}_{x}P_{b,n}^{\dR}]
\qquad\text{in}\qquad
\mathscr M_n^{\dR}(K_v),
\]
and \(\Phi_{n,v}(x)\) is the canonical lift of this class to
\(\operatorname{Fil}_n^{\dR}(K_v)\) furnished by
Proposition~\ref{proposition:dr_path_admissible}.
\end{proposition}

\begin{proof}
The algebraic Hodge reduction
\(\mathcal P_n^H\subseteq\mathcal P_n^{\dR}\) restricts to a
rigid-analytic reduction
\[
\mathcal P_{n,]\underline b[}^{H,\an}
\subseteq
\mathcal P_{n,]\underline b[}^{\dR,\an}
\]
to \((F^0\mathcal G_n^{\dR})^{\an}\). Under \(\tau_n\), such a
reduction is equivalently a rigid-analytic section of
\[
]\underline b[
\times
\left(
\mathcal G_n^{\dR}/F^0\mathcal G_n^{\dR}
\right)^{\an}.
\]
This proves that \(\Phi_{n,v}\) is rigid analytic.

At \(x\), the point \(s_n^\varphi(x)\) trivializes the Frobenius
reduction of \({}^{\vphantom{\dR}}_{x}P_{b,n}^{\dR}\). Relative to this trivialization,
its Hodge reduction determines \(\Phi_{n,v}(x)\). The moduli
interpretation of
Proposition~\ref{proposition:de_rham_period_moduli} therefore gives
\[
\pi_n\bigl(\Phi_{n,v}(x)\bigr)
=
[{}^{\vphantom{\dR}}_{x}P_{b,n}^{\dR}],
\]
and the pointwise assertion follows.
\end{proof}

\paragraph{Compatibility with finite-level quotients.}

Let
\[
q^\bullet:
\mathcal G_{n'}^\bullet
\twoheadrightarrow
\mathcal G_n^\bullet
\]
be a morphism of compatible finite-type motivic quotients. The
induced morphisms of the universal crystalline and de Rham torsors
preserve the connections, Frobenius structures, Hodge reductions,
and comparison isomorphisms. By uniqueness of horizontal
continuation, they also carry \(s_{n'}^\varphi\) to
\(s_n^\varphi\). Hence the diagram
\[
\begin{tikzcd}
]\underline b[
  \arrow[r, "\Phi_{n',v}"]
  \arrow[dr, "\Phi_{n,v}"']
&
\left(\operatorname{Fil}_{n'}^{\dR}\right)^{\an}
  \arrow[d]
\\
&
\left(\operatorname{Fil}_n^{\dR}\right)^{\an}
\end{tikzcd}
\]
commutes.

For a compatible cofiltered system
\(\{\mathcal G_n^\bullet\}_n\), the maps
\[
\Phi_v
\coloneqq 
(\Phi_{n,v})_n
\]
form a compatible family
\[
]\underline b[
\longrightarrow
\left\{
\left(\operatorname{Fil}_n^{\dR}\right)^{\an}
\right\}_n.
\]
We call this family the \emph{pro-analytic \(v\)-adic period map};
equivalently, every finite-level projection \(\Phi_{n,v}\) is rigid
analytic. The compatible family
\[
\bigl([\Phi_{n,v}]\bigr)_n
\]
is the corresponding stack-valued pro-period map.

\subsection{\texorpdfstring{Density of the \(v\)-adic Period Map}{Density of the v-adic Period Map}}
\label{section:density}
We prove that the finite-level \(v\)-adic period map is algebraically
Zariski dense in the connected component of the Hodge filtration space
containing the identity coset, by comparison with the corresponding
complex period map.

Fix a compatible finite-type motivic quotient
\(\mathcal G_n^\bullet\). Retain the global \(K\)-models
\[
\mathcal G_{n,K}^{\dR},
\qquad
F^0\mathcal G_{n,K}^{\dR},
\qquad
\mathcal P_{n,b,-,K}^{\dR},
\qquad
\mathcal P_{n,b,-,K}^{H}
\]
constructed in
\S\ref{appendix:finite_level_filtered_path_torsors}; see
Corollaries~\ref{corollary:finite_level_filtered_path_groupoid}
and~\ref{corollary:finite_level_universal_hodge_reduction}.

Set
\[
\operatorname{Fil}_{n,K}^{\dR}
\coloneqq 
\mathcal G_{n,K}^{\dR}/F^0\mathcal G_{n,K}^{\dR}.
\]
Its base change to \(K_v\) is
\(\operatorname{Fil}_n^{\dR}\). Let
\[
\operatorname{Fil}_{n,K}^{\dR,\circ}
\coloneqq 
\mathcal G_{n,K}^{\dR,\circ}\cdot[1]
\subseteq
\operatorname{Fil}_{n,K}^{\dR}
\]
be the connected component containing the identity coset, and set
\[
\operatorname{Fil}_n^{\dR,\circ}
\coloneqq 
\operatorname{Fil}_{n,K}^{\dR,\circ}\otimes_KK_v.
\]
This is the connected component of
\(\operatorname{Fil}_n^{\dR}\) containing the identity coset. Since
\(]\underline b[\) is connected and
\[
\Phi_{n,v}(b)=[1],
\]
the period map factors through
\[
\Phi_{n,v}:
]\underline b[
\longrightarrow
\left(
\operatorname{Fil}_n^{\dR,\circ}
\right)^{\an}.
\]

A rigid-analytic map
\[
f:Y\longrightarrow Z^{\an},
\]
with \(Z\) an algebraic space of finite type, is
\emph{algebraically Zariski dense} if its image is not contained in
the analytification of any proper closed algebraic subspace of \(Z\).

\paragraph{Comparison with the complex period map.}

Fix an embedding
\[
\iota_\infty:
K\hookrightarrow\mathbb C,
\]
and write
\[
X^\dagger_{\mathbb C}
\coloneqq 
X^\dagger\otimes_{K,\iota_\infty}\mathbb C,
\qquad
\mathcal G_{n,\mathbb C}^{\dR}
\coloneqq 
\mathcal G_{n,K}^{\dR}
\otimes_{K,\iota_\infty}\mathbb C.
\]
Let
\[
D_{b,\mathbb C}
\subseteq
\left(
X^\dagger_{\mathbb C}
\right)^{\an}
\]
be a simply connected analytic neighborhood of \(b\). The unique
horizontal trivialization of
\(\mathcal P_{n,b,-,K}^{\dR}\otimes_{K,\iota_\infty}\mathbb C\)
whose value at \(b\) is the identity, together with the base change
of its Hodge reduction, defines a holomorphic period map
\[
\Phi_{n,\mathbb C}:
D_{b,\mathbb C}
\longrightarrow
\left(
\operatorname{Fil}_{n,K}^{\dR,\circ}
\otimes_{K,\iota_\infty}\mathbb C
\right)^{\an}.
\]

\begin{lemma}
[{Comparison of local Zariski closures;
cf.~\cite[\S\S3.3--3.4, especially
Lemmas~3.1--3.2]{LV}}]
\label{lemma:period_closure_comparison}
After shrinking to neighborhoods
\[
D_v\subseteq]\underline b[,
\qquad
D_{\mathbb C}\subseteq D_{b,\mathbb C}
\]
of \(b\), there is a closed \(K\)-subscheme
\[
Z_{n,K}
\subseteq
\operatorname{Fil}_{n,K}^{\dR,\circ}
\]
such that
\[
\overline{\Phi_{n,v}(D_v)}^{\,\mathrm{Zar}}
=
Z_{n,K}\otimes_KK_v
\]
and
\[
\overline{\Phi_{n,\mathbb C}(D_{\mathbb C})}^{\,\mathrm{Zar}}
=
Z_{n,K}\otimes_{K,\iota_\infty}\mathbb C.
\]
\end{lemma}

\begin{proof}
Choose an \'{e}tale local parameter \(t\) at \(b\) and a locally
closed immersion
\[
\operatorname{Fil}_{n,K}^{\dR,\circ}
\hookrightarrow
\mathbb P_K^N.
\]
The universal connection and Hodge reduction are defined over \(K\)
by
Corollaries~\ref{corollary:finite_level_filtered_path_groupoid}
and~\ref{corollary:finite_level_universal_hodge_reduction}.
The \(v\)-adic section of
Lemma~\ref{lemma:frobenius_horizontal_trivialization} and the complex
horizontal trivialization are the unique horizontal sections
normalized by the identity at \(b\). Their formal Taylor expansions,
and hence those of the associated period maps, are therefore obtained
from the same \(K\)-defined connection and initial condition.

After shrinking \(D_v\) and \(D_{\mathbb C}\), the two maps to
projective space are consequently represented by the same tuple
\[
[B_0(t):\cdots:B_N(t)],
\qquad
B_i(t)\in K[[t]],
\]
converging on the respective neighborhoods. The assertion follows
from \cite[Lemma~3.2]{LV}, after intersecting the common projective
Zariski closure with
\(\operatorname{Fil}_{n,K}^{\dR,\circ}\).
\end{proof}

\paragraph{Complex density.}

Let
\[
\widetilde X^\dagger_{\mathbb C}
\longrightarrow
\left(
X^\dagger_{\mathbb C}
\right)^{\an}
\]
be the universal cover, and choose a lift \(\widetilde b\) of \(b\).
Horizontal continuation extends the local complex period map to
\[
\widetilde\Phi_{n,\mathbb C}:
\widetilde X^\dagger_{\mathbb C}
\longrightarrow
\left(
\operatorname{Fil}_{n,K}^{\dR,\circ}
\otimes_{K,\iota_\infty}\mathbb C
\right)^{\an}.
\]

Let
\[
\widetilde\rho_n:
\pi_1^{\Betti}(X^\dagger_{\mathbb C},b)
\longrightarrow
\mathcal G_{n,\mathbb C}^{\dR}(\mathbb C)
\]
be the finite-level relative-completion homomorphism. Hain's relative
de Rham theorem and the equivariance of the universal connection give
\[
\widetilde\Phi_{n,\mathbb C}(\gamma z)
=
\widetilde\rho_n(\gamma)\,
\widetilde\Phi_{n,\mathbb C}(z);
\]
see
\cite[Theorem~10.1 and Proposition~14.4]{HAIN_HDR_RMC}.

\begin{lemma}[Density of the complex period map]
\label{lemma:complex_period_density}
For every analytic neighborhood
\[
D\subseteq D_{b,\mathbb C}
\]
of \(b\), the image of
\[
\Phi_{n,\mathbb C}|_D:
D
\longrightarrow
\left(
\operatorname{Fil}_{n,K}^{\dR,\circ}
\otimes_{K,\iota_\infty}\mathbb C
\right)^{\an}
\]
is algebraically Zariski dense.
\end{lemma}

\begin{proof}
Let
\[
Z
\subseteq
\operatorname{Fil}_{n,K}^{\dR,\circ}
\otimes_{K,\iota_\infty}\mathbb C
\]
be the Zariski closure of
\(\Phi_{n,\mathbb C}(D)\). Then
\[
\widetilde\Phi_{n,\mathbb C}^{-1}(Z^{\an})
\]
is a closed complex-analytic subset of the connected universal cover
containing a neighborhood of \(\widetilde b\). It is therefore the
entire universal cover.

Set
\[
\Gamma_n^\circ
\coloneqq 
\widetilde\rho_n^{-1}
\left(
\mathcal G_{n,\mathbb C}^{\dR,\circ}(\mathbb C)
\right).
\]
For every \(\gamma\in\Gamma_n^\circ\), equivariance gives
\[
\widetilde\Phi_{n,\mathbb C}
(\gamma\widetilde b)
=
\widetilde\rho_n(\gamma)\cdot[1]
\in
Z(\mathbb C).
\]
Since the relative-completion homomorphism has Zariski-dense image,
\(\widetilde\rho_n(\Gamma_n^\circ)\) is Zariski dense in
\(\mathcal G_{n,\mathbb C}^{\dR,\circ}\). Hence its orbit through
\([1]\) is dense in
\[
\mathcal G_{n,\mathbb C}^{\dR,\circ}\cdot[1]
=
\operatorname{Fil}_{n,K}^{\dR,\circ}
\otimes_{K,\iota_\infty}\mathbb C.
\]
Since \(Z\) contains this orbit, it is the entire filtration space.
\end{proof}

\begin{proposition}[Density of the \(v\)-adic period map]
\label{proposition:p-adic-dense}
For every compatible finite-type motivic quotient
\(\mathcal G_n^\bullet\), the period map
\[
\Phi_{n,v}:
]\underline b[
\longrightarrow
\left(
\operatorname{Fil}_{n,K}^{\dR,\circ}
\otimes_KK_v
\right)^{\an}
\]
is algebraically Zariski dense. In fact, its restriction to every
sufficiently small neighborhood of \(b\) is already algebraically
Zariski dense.
\end{proposition}

\begin{proof}
Let \(U\subseteq]\underline b[\) be a sufficiently small neighborhood
of \(b\). Apply
Lemma~\ref{lemma:period_closure_comparison} with
\(D_v\subseteq U\). By
Lemma~\ref{lemma:complex_period_density},
\[
\overline{\Phi_{n,\mathbb C}(D_{\mathbb C})}^{\,\mathrm{Zar}}
=
\operatorname{Fil}_{n,K}^{\dR,\circ}
\otimes_{K,\iota_\infty}\mathbb C.
\]
It follows that
\[
Z_{n,K}
=
\operatorname{Fil}_{n,K}^{\dR,\circ},
\]
and hence
\[
\overline{\Phi_{n,v}(D_v)}^{\,\mathrm{Zar}}
=
\operatorname{Fil}_{n,K}^{\dR,\circ}\otimes_KK_v.
\]
Thus \(\Phi_{n,v}|_U\), and therefore \(\Phi_{n,v}\), is
algebraically Zariski dense.
\end{proof}

\begin{remark}
If \(\mathcal G_{n,K}^{\dR}\) is connected, then
\[
\operatorname{Fil}_{n,K}^{\dR,\circ}
=
\operatorname{Fil}_{n,K}^{\dR}.
\]
In general, the connected residue tube meets only the component
containing the identity coset. All components are translates of one
another and have the same dimension.
\end{remark}

Consequently, for every compatible cofiltered system
\(\{\mathcal G_n^\bullet\}_n\) in the Kodaira--Parshin motivic tower,
the pro-analytic period map
\[
\Phi_v:
]\underline b[
\longrightarrow
\left\{
\left(
\operatorname{Fil}_n^{\dR}
\right)^{\an}
\right\}_n
\]
is \emph{levelwise algebraically Zariski dense}: each finite-level
projection is dense in the component of its target containing the
identity coset.

\subsection{The Neutral Bloch--Kato Logarithm}
\label{chapter:bk_log}
This subsection constructs the local and global neutral
Bloch--Kato logarithms by relative \(p\)-adic Hodge realization. A
local Bloch--Kato finite torsor is first realized as a de Rham torsor
with canonical Hodge and Frobenius pseudo-reductions. On the neutral
reductive fibre these loci are genuine reductions, so the realization
defines an object of the de Rham moduli stack.

\begin{itemize}
\item
\S\ref{subsubsection:relative_padic_hodge_realization} constructs
the relative crystalline realization of a Bloch--Kato finite
\(\mathcal G_n^{\et}\)-torsor, together with its Frobenius and
\(B_{\cris}\)-comparison. After base change to \(K_v\), it constructs the corresponding de Rham
torsor, its selected \(B_{\dR}\)-comparison, and its canonical Hodge
and Frobenius pseudo-reductions.

\item
\S\ref{subsubsection:neutral_realization_admissibility} restricts
to the neutral reductive fibre and proves that both
pseudo-reductions are fppf-locally nonempty. The Hodge reduction is
obtained from the unipotent \(B_{\dR}^+\)-trivialization, while the
Frobenius reduction follows from the unipotent Frobenius--Lang
isomorphism and the negative-weight hypothesis. The realized torsor
is therefore admissible.

\item
\S\ref{subsubsection:neutral_bk_logarithms} stackifies this local
realization to define
\[
\BK_{n,v}:
\mathscr S_{n,v}
\longrightarrow
\mathscr M_n^{\dR},
\]
where for $S_p\subseteq T'\subseteq T$ and, for $\lambda\in\{v,T'\}$, we set
\[
\mathscr S_{n,\lambda}
\coloneqq 
\mathscr N_{f,\lambda}
\bigl(
\mathcal G_n^{\et}/\mathcal R^{\et}
\bigr)_{K_v}.
\]
It then defines the global logarithm by localization at \(v\), as
well as its pullback to the trivialized neutral Selmer space.

\item
\S\ref{subsubsection:bk_compatibility} proves compatibility of the
neutral Bloch--Kato logarithm with the neutral Kummer map and the
stack-valued period map. It also establishes functoriality under
morphisms of compatible finite-type motivic quotients and hence for
the resulting pro-systems.
\end{itemize}

\subsubsection{\texorpdfstring{Relative \(p\)-adic Hodge Realization}{Relative p-adic Hodge Realization}}
\label{subsubsection:relative_padic_hodge_realization}
Let
\[
N=\Spec(A)
\]
be an affine \(K_v\)-scheme, and set
\[
A_0\coloneqq A\otimes_{\Qp}K_0,
\qquad
A_{\cris}\coloneqq A\otimes_{\Qp}B_{\cris},
\qquad
N_0\coloneqq \Spec(A_0),
\]
where \(A_{\cris}\) is regarded as an \(A_0\)-algebra through the
canonical map induced by \(K_0\hookrightarrow B_{\cris}\).

We use the finite-level comparison induced by
Theorem~\ref{theorem:etale_crystalline_comparison}:
\[
\iota_{n,\cris}:
\mathcal G_n^{\et}\otimes_{\Qp}B_{\cris}
\xrightarrow{\sim}
\mathcal G_n^{\cri}\otimes_{K_0}B_{\cris}.
\]
It is compatible with the \(G_v\)-actions and with crystalline
Frobenius.

\begin{definition}[Relative crystalline realization]
\label{definition:relative_crystalline_realization}
Let \(P\) be a Bloch--Kato finite \(G_v\)-equivariant
\(\mathcal G_n^{\et}\)-torsor over \(N\). Define a presheaf on
\(A_0\)-algebras by
\[
\mathbf D_{\cris,N}(P)(B)
\coloneqq 
P\bigl(B\otimes_{K_0}B_{\cris}\bigr)^{G_v}.
\]
Here \(B\otimes_{K_0}B_{\cris}\) is regarded as an \(A\)-algebra
through
\[
A\longrightarrow A_0\longrightarrow B
\longrightarrow B\otimes_{K_0}B_{\cris},
\]
and \(G_v\) acts diagonally, trivially on \(B\).
\end{definition}

\begin{proposition}[Relative crystalline descent]
\label{proposition:relative_crystalline_descent}
The presheaf \(\mathbf D_{\cris,N}(P)\) is an fppf sheaf represented
by an affine \(N_0\)-scheme of finite type
\[
\mathcal D_{\cris,N}(P).
\]
It is a right \(\mathcal G_{n,N_0}^{\cri}\)-torsor and admits a
canonical equivariant comparison isomorphism
\[
\epsilon_{P,\cris}:
\mathcal D_{\cris,N}(P)
\times_{N_0}\Spec(A_{\cris})
\xrightarrow{\sim}
P\times_N\Spec(A_{\cris}),
\]
where the structure groups are identified through
\(\iota_{n,\cris}\). Its coordinate algebra is canonically
\[
\mathcal O\bigl(\mathcal D_{\cris,N}(P)\bigr)
\simeq
\left(
\mathcal O(P)\otimes_AA_{\cris}
\right)^{G_v}.
\]

Let
\[
\sigma_{N_0}:
N_0\longrightarrow N_0
\]
be induced by
\[
\id_A\otimes\sigma:
A_0\longrightarrow A_0.
\]
Frobenius on \(B_{\cris}\) induces a semilinear isomorphism
\[
\varphi_P:
\sigma_{N_0}^*
\mathcal D_{\cris,N}(P)
\xrightarrow{\sim}
\mathcal D_{\cris,N}(P),
\]
equivariant with respect to crystalline Frobenius
\[
\varphi_n:
\sigma^*\mathcal G_n^{\cri}
\xrightarrow{\sim}
\mathcal G_n^{\cri}.
\]

These constructions commute with pullback in \(N\), with
isomorphisms of torsors, and with extension of structure group along
morphisms of compatible finite-type motivic quotients.
\end{proposition}

\begin{proof}
For every \(A_0\)-algebra \(B\), the comparison
\(\iota_{n,\cris}\) and the identity
\[
\left(
B\otimes_{K_0}B_{\cris}
\right)^{G_v}
=
B
\]
give a natural isomorphism
\[
\mathcal G_n^{\et}
\bigl(
B\otimes_{K_0}B_{\cris}
\bigr)^{G_v}
\xrightarrow{\sim}
\mathcal G_n^{\cri}(B).
\]
Thus \(\mathbf D_{\cris,N}(P)\) is a
\(\mathcal G_{n,N_0}^{\cri}\)-pseudo-torsor.

After an fppf base change
\[
N_i=\Spec(A_i)\longrightarrow N,
\]
the Bloch--Kato condition supplies a fixed point
\[
p_i\in
P\bigl(A_i\otimes_{\Qp}B_{\cris}\bigr)^{G_v}.
\]
On an overlap \(N_{ij}\), write
\[
p_j=p_i\cdot g_{ij}.
\]
Since \(p_i\) and \(p_j\) are \(G_v\)-fixed,
\[
g_{ij}
\in
\mathcal G_n^{\et}
\bigl(
A_{ij}\otimes_{\Qp}B_{\cris}
\bigr)^{G_v}
\simeq
\mathcal G_n^{\cri}
\bigl(
A_{ij}\otimes_{\Qp}K_0
\bigr).
\]
The cocycle \(\{g_{ij}\}\) therefore descends the trivial
crystalline torsors over \(N_{i,0}\) to
\(\mathcal D_{\cris,N}(P)\), and the same local trivializations give
\(\epsilon_{P,\cris}\). Taking coordinate rings and \(G_v\)-invariants
gives the displayed formula.

Frobenius on \(B_{\cris}\) preserves the fixed-point sheaf and is
semilinear over \(\id_A\otimes\sigma\), yielding \(\varphi_P\).
Compatibility of \(\iota_{n,\cris}\) with Frobenius and
\(G_v\)-equivariance of the original torsor action give
\[
\varphi_P(p\cdot g)
=
\varphi_P(p)\cdot\varphi_n(g).
\]
The remaining functorialities follow directly from the fixed-point
construction.
\end{proof}

\begin{definition}[Relative de Rham realization]
\label{definition:relative_de_rham_realization}
Let
\[
m_A:
A_0\longrightarrow A,
\qquad
a\otimes c\longmapsto ac,
\]
be induced by \(K_0\hookrightarrow K_v\). Define
\[
\mathcal D_{\dR,N}(P)
\coloneqq 
\mathcal D_{\cris,N}(P)
\times_{N_0,m_A}N.
\]
\end{definition}

\begin{proposition}[Relative de Rham descent]
\label{proposition:relative_padic_hodge_realization}
Through the crystalline--de Rham comparison
\[
\mathcal G_n^{\cri}\otimes_{K_0}K_v
\xrightarrow{\sim}
\mathcal G_n^{\dR},
\]
the scheme \(\mathcal D_{\dR,N}(P)\) is a right
\(\mathcal G_{n,N}^{\dR}\)-torsor.

Set
\[
A_{\dR}
\coloneqq 
A\otimes_{K_v}B_{\dR}.
\]
The inclusion \(B_{\cris}\hookrightarrow B_{\dR}\) induces an
\(A_0\)-algebra homomorphism
\[
\jmath_{A,\dR}:
A_{\cris}
=
A\otimes_{\Qp}B_{\cris}
\longrightarrow
A\otimes_{K_v}B_{\dR}
=
A_{\dR},
\]
where the \(A_0\)-algebra structure on \(A_{\dR}\) factors through
\[
m_A:A_0\longrightarrow A.
\]
Base change of \(\epsilon_{P,\cris}\) along
\(\jmath_{A,\dR}\) induces a canonical \(G_v\)-equivariant,
structure-group-equivariant isomorphism
\[
\epsilon_{P,\dR}:
\mathcal D_{\dR,N}(P)
\times_N\Spec(A_{\dR})
\xrightarrow{\sim}
P\times_N\Spec(A_{\dR}),
\]
where \(G_v\) acts on the left through its action on \(A_{\dR}\).

The construction commutes with pullback in \(N\), with isomorphisms
of torsors, and with extension of structure group along morphisms of
compatible finite-type motivic quotients.
\end{proposition}

\begin{proof}
The torsor assertion follows by base change from
Proposition~\ref{proposition:relative_crystalline_descent} and the
crystalline--de Rham comparison of the structure groups. Since the
\(A_0\)-algebra structure on \(A_{\dR}\) factors through \(m_A\),
base change of \(\epsilon_{P,\cris}\) along
\[
\jmath_{A,\dR}:A_{\cris}\longrightarrow A_{\dR}
\]
gives \(\epsilon_{P,\dR}\). Its equivariance and the asserted
functorialities follow from the corresponding properties of the
crystalline comparison.
\end{proof}

\begin{definition}[Pseudo-reductions]
\label{definition:pseudo_reduction}
Let \(H\hookrightarrow G\) be a closed immersion of affine group
schemes over a scheme \(N\), and let \(P\) be a right \(G\)-torsor.
An \emph{\(H\)-pseudo-reduction} of \(P\) is an \(H\)-stable
subsheaf \(Q\subseteq P\) such that the action induces an isomorphism
\[
Q\times_NH
\xrightarrow{\sim}
Q\times_NQ.
\]
It is a \emph{reduction of structure group} if \(Q\) is
fppf-locally nonempty, equivalently, if \(Q\to N\) is an
\(H\)-torsor.
\end{definition}

For an \(A\)-algebra \(B\), write
\[
B_{\dR,B}
\coloneqq 
B\otimes_{K_v}B_{\dR},
\qquad
B_{\dR,B}^+
\coloneqq 
B\otimes_{K_v}B_{\dR}^+.
\]
By
Proposition~\ref{proposition:selected_finite_level_filtered_comparison},
the selected-component finite-level filtered comparison identifies,
functorially in \(B\),
\begin{equation}
\label{eq:selected_filtered_comparison}
\begin{aligned}
\mathcal G_{n,N}^{\dR}(B)
&\xrightarrow{\sim}
\mathcal G_n^{\et}(B_{\dR,B})^{G_v},
\\
F^0\mathcal G_{n,N}^{\dR}(B)
&\xrightarrow{\sim}
\mathcal G_n^{\et}(B_{\dR,B}^+)^{G_v},
\end{aligned}
\end{equation}
where the second identification is compatible with the inclusions
into the first.

\begin{proposition}[Canonical Hodge pseudo-reduction]
\label{proposition:relative_hodge_pseudoreduction}
For every \(A\)-algebra \(B\), the comparison isomorphism and the
identity
\[
(B_{\dR,B})^{G_v}=B
\]
induce a natural isomorphism
\[
\theta_{P,B}:
\mathcal D_{\dR,N}(P)(B)
\xrightarrow{\sim}
P(B_{\dR,B})^{G_v},
\qquad
d\longmapsto
\epsilon_{P,\dR}(d_{B_{\dR,B}}).
\]
Define
\[
\mathcal D_{\dR,N}(P)^H
\subseteq
\mathcal D_{\dR,N}(P)
\]
to be the fppf subfunctor consisting of those \(d\) for which
\(\theta_{P,B}(d)\) lifts, fppf-locally on \(\Spec(B)\), to a point
of
\[
P(B_{\dR,B}^+)^{G_v}.
\]
Then \(\mathcal D_{\dR,N}(P)^H\) is an
\(F^0\mathcal G_{n,N}^{\dR}\)-pseudo-reduction.

The construction commutes with pullback in \(N\) and is functorial
under isomorphisms of torsors. A morphism of compatible finite-type
motivic quotients induces a natural equivariant morphism between the
corresponding Hodge pseudo-reductions.
\end{proposition}

\begin{proof}
The assertion that \(\theta_{P,B}\) is an isomorphism is fppf-local
on \(\Spec(B)\). After choosing a section of
\(\mathcal D_{\dR,N}(P)\), it becomes the first identification in
\eqref{eq:selected_filtered_comparison}, and hence follows by
descent.

Let
\[
d,d'\in\mathcal D_{\dR,N}(P)^H(B),
\]
and write uniquely
\[
d'=d\cdot g,
\qquad
g\in\mathcal G_{n,N}^{\dR}(B).
\]
After passing to an fppf \(A\)-algebra \(B'\) over \(B\), choose
lifts
\[
q,q'\in P(B_{\dR,B'}^+)^{G_v}
\]
of \(\theta_{P,B}(d)\) and \(\theta_{P,B}(d')\), respectively.
There is then a unique
\[
h\in\mathcal G_n^{\et}(B_{\dR,B'}^+)^{G_v}
\]
such that
\[
q'=q\cdot h.
\]
Compatibility of \(\epsilon_{P,\dR}\) with the torsor actions
identifies the image of \(h\) under
\eqref{eq:selected_filtered_comparison} with the restriction of
\(g\) to \(B'\). Hence this restriction belongs to
\[
F^0\mathcal G_{n,N}^{\dR}(B').
\]
Since
\[
F^0\mathcal G_{n,N}^{\dR}
\hookrightarrow
\mathcal G_{n,N}^{\dR}
\]
is a closed immersion, membership descends faithfully flatly, and
therefore
\[
g\in F^0\mathcal G_{n,N}^{\dR}(B).
\]

Conversely, the action of
\(F^0\mathcal G_{n,N}^{\dR}\) preserves the selected
\(B_{\dR}^+\)-lifting condition. Thus the torsor isomorphism
restricts to an isomorphism of fppf sheaves
\[
\mathcal D_{\dR,N}(P)^H
\times_NF^0\mathcal G_{n,N}^{\dR}
\xrightarrow{\sim}
\mathcal D_{\dR,N}(P)^H
\times_N
\mathcal D_{\dR,N}(P)^H.
\]
This proves the pseudo-reduction assertion.

Finally, the selected \(B_{\dR}^+\)-lifting condition commutes with
base change in \(N\) and is preserved by isomorphisms of torsors.
A morphism of compatible finite-type motivic quotients carries
selected points forward and therefore induces the asserted
equivariant morphism.
\end{proof}

\begin{proposition}[Canonical Frobenius pseudo-reduction]
\label{proposition:relative_frobenius_pseudoreduction}
For an \(A\)-algebra \(B\), set
\[
B_0\coloneqq B\otimes_{\Qp}K_0
\]
and let
\[
m_B:B_0\longrightarrow B,
\qquad
b\otimes c\longmapsto bc,
\]
be the multiplication map. Define
\[
\mathcal D_{\dR,N}(P)^\varphi(B)
\coloneqq 
\left\{
p\in
\mathcal D_{\cris,N}(P)(B_0)
\ \middle|\
\varphi_P(p)=p
\right\},
\]
viewed in \(\mathcal D_{\dR,N}(P)(B)\) through \(m_B\).

The functor \(\mathcal D_{\dR,N}(P)^\varphi\) is represented by a
closed affine \(N\)-subscheme of \(\mathcal D_{\dR,N}(P)\), and is a
\(\mathcal G_{n,K_v,N}^{\varphi=1}\)-pseudo-reduction.

The construction commutes with pullback in \(N\) and is functorial
under isomorphisms of torsors. A morphism of compatible finite-type
motivic quotients induces a natural equivariant morphism between the
corresponding Frobenius pseudo-reductions.
\end{proposition}

\begin{proof}
Set
\[
\mathscr R_P
\coloneqq 
\operatorname{Res}_{N_0/N}
\mathcal D_{\cris,N}(P).
\]
Since \(N_0\to N\) is finite étale, \(\mathscr R_P\) is an affine
\(N\)-scheme of finite type. The semilinear Frobenius
\[
\varphi_P:
\sigma_{N_0}^*\mathcal D_{\cris,N}(P)
\xrightarrow{\sim}
\mathcal D_{\cris,N}(P)
\]
induces an \(N\)-automorphism
\[
\Phi_P:
\mathscr R_P
\xrightarrow{\sim}
\mathscr R_P.
\]
Its fixed locus
\[
\mathscr R_P^{\Phi_P=1}
\coloneqq 
\operatorname{Eq}
\left(
\Phi_P,\id_{\mathscr R_P}
\right)
\]
is a closed affine \(N\)-subscheme, and for every \(A\)-algebra \(B\)
one has
\[
\mathscr R_P^{\Phi_P=1}(B)
=
\left\{
p\in
\mathcal D_{\cris,N}(P)(B_0)
\ \middle|\
\varphi_P(p)=p
\right\}.
\]

Evaluation along the section
\[
m_A:
N\longrightarrow N_0
\]
induces a morphism
\[
\mathscr R_P^{\Phi_P=1}
\longrightarrow
\mathcal D_{\dR,N}(P).
\]
This is a closed immersion. Indeed, if
\[
f=[K_0:\Qp],
\]
then for every \(A\)-algebra \(B\) the canonical decomposition
\[
B\otimes_{\Qp}K_0
\simeq
\prod_{i=0}^{f-1}B,
\]
indexed by the powers of \(\sigma\), shows that a Frobenius-fixed
point is uniquely determined by the component selected by \(m_B\).
That component is cut out by the equality between the \(f\)-fold
Frobenius and the identity. The image of the displayed closed
immersion is therefore precisely
\[
\mathcal D_{\dR,N}(P)^\varphi.
\]

Applying the same construction to
\(\mathcal G_{n,N_0}^{\cri}\) identifies its fixed locus with
\[
\mathcal G_{n,K_v,N}^{\varphi=1}.
\]
The crystalline torsor isomorphism
\[
\mathcal D_{\cris,N}(P)
\times_{N_0}
\mathcal G_{n,N_0}^{\cri}
\xrightarrow{\sim}
\mathcal D_{\cris,N}(P)
\times_{N_0}
\mathcal D_{\cris,N}(P)
\]
is Frobenius equivariant. Weil restriction and passage to fixed
equalizers therefore give an isomorphism
\[
\mathcal D_{\dR,N}(P)^\varphi
\times_N
\mathcal G_{n,K_v,N}^{\varphi=1}
\xrightarrow{\sim}
\mathcal D_{\dR,N}(P)^\varphi
\times_N
\mathcal D_{\dR,N}(P)^\varphi.
\]
Thus \(\mathcal D_{\dR,N}(P)^\varphi\) is a
\(\mathcal G_{n,K_v,N}^{\varphi=1}\)-pseudo-reduction.

Finally, Weil restriction along \(N_0\to N\), fixed equalizers, and
evaluation along \(m_A\) commute with base change in \(N\).
Isomorphisms of torsors preserve these constructions, while a
morphism of compatible finite-type motivic quotients carries
Frobenius-fixed points forward and hence induces the asserted
equivariant morphism.
\end{proof}

\subsubsection{Admissibility on the Neutral Fibre}
\label{subsubsection:neutral_realization_admissibility}

We now restrict to the neutral reductive fibre and show that the two
canonical pseudo-reductions are locally nonempty.

\begin{proposition}[The Hodge reduction]
\label{proposition:neutral_hodge_reduction}
Let \(P\) belong to the local neutral Bloch--Kato prestack over an
affine \(K_v\)-scheme \(N\). Then
\[
\mathcal D_{\dR,N}(P)^H
\]
is fppf-locally nonempty and hence is a reduction of structure group
to \(F^0\mathcal G_{n,N}^{\dR}\).
\end{proposition}

\begin{proof}
Choose an fppf cover
\[
N'\longrightarrow N
\]
over which the reductive pushout of \(P\) is trivial. By
Proposition~\ref{proposition:bk_neutral_fibre_comparison}, there is a
Bloch--Kato finite \(\mathcal U_n^{\et}\)-torsor \(Q\) over \(N'\)
and an isomorphism
\[
P_{N'}
\simeq
Q\times^{\mathcal U_n^{\et}}\mathcal G_n^{\et}.
\]

After a further fppf cover
\[
N''=\Spec(A'')
\longrightarrow
N',
\]
choose a trivialization of \(Q_{N''}\), and let
\[
\xi:
G_v\longrightarrow
\mathcal U_n^{\et}(A'')
\]
be the resulting cocycle. Since \(Q\) is Bloch--Kato finite, the
class of \(\xi\) lies in
\[
H^1_f
\bigl(
G_v,\mathcal U_n^{\et}(A'')
\bigr).
\]
The weight hypotheses for \(\mathcal U_n^{\et}\) and
\cite[Lemma~3.1.12]{BETTS_WEIGHT} give a
\(B_{\dR}^+\)-trivialization of \(\xi\), equivalently a fixed point
\[
q_{\dR}
\in
Q_{N''}
\left(
A''\otimes_{\Qp}B_{\dR}^+
\right)^{G_v}.
\]
Let
\[
q_{\dR,v}
\in
Q_{N''}
\left(
A''\otimes_{K_v}B_{\dR}^+
\right)^{G_v}
\]
be its image under the natural \(G_v\)-equivariant quotient map
\[
A''\otimes_{\Qp}B_{\dR}^+
\longrightarrow
A''\otimes_{K_v}B_{\dR}^+.
\]
Applying
Proposition~\ref{proposition:relative_padic_hodge_realization}
with \(\mathcal G_n^\bullet\) replaced by its compatible unipotent
radical \(\mathcal U_n^\bullet\), write
\[
\epsilon_{Q,\dR}:
\mathcal D_{\dR,N''}(Q_{N''})
\times_{N''}
\Spec\bigl(A''\otimes_{K_v}B_{\dR}\bigr)
\xrightarrow{\sim}
Q_{N''}
\times_{N''}
\Spec\bigl(A''\otimes_{K_v}B_{\dR}\bigr)
\]
for the resulting comparison isomorphism.

The image of \(q_{\dR,v}\) in
\[
Q_{N''}
\left(
A''\otimes_{K_v}B_{\dR}
\right)^{G_v}
\]
corresponds under \(\epsilon_{Q,\dR}\) to a point
\[
\gamma_Q^H
\in
\mathcal D_{\dR,N''}(Q_{N''})(N'').
\]
By definition of the selected \(B_{\dR}^+\)-locus,
\[
\gamma_Q^H
\in
\mathcal D_{\dR,N''}(Q_{N''})^H(N'');
\]
cf.~\cite[Remark~3.1.13]{BETTS_WEIGHT}.

Relative de Rham realization commutes with extension of structure
group, giving
\[
\mathcal D_{\dR,N'}(P_{N'})
\simeq
\mathcal D_{\dR,N'}(Q)
\times^{\mathcal U_n^{\dR}}
\mathcal G_n^{\dR}.
\]
The image of \(\gamma_Q^H\) under this pushout belongs to
\[
\mathcal D_{\dR,N''}(P_{N''})^H(N'').
\]
Thus \(\mathcal D_{\dR,N}(P)^H\) is fppf-locally nonempty. Since it is an \(F^0\mathcal G_{n,N}^{\dR}\)-pseudo-reduction by
Proposition~\ref{proposition:relative_hodge_pseudoreduction},
it is a reduction of structure group.
\end{proof}

\begin{lemma}
[{The Frobenius--Lang isomorphism;
cf.~\cite[Theorem~3.1 and Corollary~3.2]{BESSER_COLEMAN}}]
\label{lemma:unipotent_frobenius_lang}
The morphism
\[
L_\varphi:
\operatorname{Res}_{K_0/\Qp}\mathcal U_n^{\cri}
\longrightarrow
\operatorname{Res}_{K_0/\Qp}\mathcal U_n^{\cri},
\qquad
u\longmapsto u^{-1}\varphi(u),
\]
is an isomorphism. Consequently, every
\(\mathcal U_n^{\cri}\)-torsor equipped with a compatible
semilinear Frobenius has a unique Frobenius-fixed section.
\end{lemma}

\begin{proof}
Set
\[
V_i
\coloneqq 
\Gr_i^{\DCS}\Lie(\mathcal U_n^{\cri}),
\qquad
f\coloneqq [K_0:\Qp].
\]
By
\cite[Proposition~6.2.4]{KANTORTHESIS},
the \(K_0\)-linear operator \(\varphi^f\) has strictly negative Weil
weights on every \(V_i\). In particular,
\[
V_i^{\varphi^f}=0.
\]

Let
\[
W_i\coloneqq \operatorname{Res}_{K_0/\Qp}V_i.
\]
After base change to \(K_0\), there is a canonical decomposition
\[
(W_i)_{K_0}
\simeq
\prod_{j=0}^{f-1}\sigma^{j*}V_i,
\]
under which \(\varphi\) cyclically permutes the factors. Projection
onto any one factor therefore identifies the \(\varphi\)-fixed locus
with the \(\varphi^f\)-fixed locus of \(V_i\). Hence
\[
W_i^{\varphi}=0.
\]
Thus the Lang map on \(W_i\) is the injective \(\Qp\)-linear map
\[
\varphi-1:W_i\longrightarrow W_i.
\]
Since \(W_i\) is finite-dimensional, this map is an isomorphism.

The DCS filtration gives successive central
extensions whose kernels are the vector groups \(W_i\). Inductively,
one solves the Lang equation first in the quotient and then uniquely
corrects a lift in the central kernel. Since the Lang map is an
isomorphism on every \(W_i\), this proves that
\[
L_\varphi:
\operatorname{Res}_{K_0/\Qp}\mathcal U_n^{\cri}
\longrightarrow
\operatorname{Res}_{K_0/\Qp}\mathcal U_n^{\cri}
\]
is an isomorphism.

Now let \(Q\) be a Frobenius torsor. Fppf-locally choose \(q\in Q\)
and write
\[
\varphi(q)=q\cdot c.
\]
There is a unique \(u\) satisfying
\[
u^{-1}\varphi(u)=c.
\]
Then
\[
\varphi(q\cdot u^{-1})=q\cdot u^{-1}.
\]
The resulting local fixed sections are unique and therefore agree on
overlaps, so they descend to a unique Frobenius-fixed section of
\(Q\).
\end{proof}

\begin{proposition}[Admissibility of the neutral realization]
\label{proposition:neutral_realization_admissible}
Let \(P\) belong to the local neutral Bloch--Kato prestack over an
affine \(K_v\)-scheme \(N\). Then
\[
\left(
\mathcal D_{\dR,N}(P),
\mathcal D_{\dR,N}(P)^H,
\mathcal D_{\dR,N}(P)^\varphi
\right)
\]
is an admissible \(\mathcal G_n^{\dR}\)-torsor over \(N\) in the
sense of Definition~\ref{definition:dr_torsors}. This construction is
functorial in \(P\) and commutes with pullback in \(N\).
\end{proposition}

\begin{proof}
We first show that the Frobenius pseudo-reduction is fppf-locally nonempty.
After an fppf base change, write
\[
P\simeq
Q\times^{\mathcal U_n^{\et}}\mathcal G_n^{\et}
\]
for a Bloch--Kato finite \(\mathcal U_n^{\et}\)-torsor \(Q\). Its
crystalline realization is a Frobenius torsor under
\(\mathcal U_n^{\cri}\), and hence has a unique Frobenius-fixed
section by
Lemma~\ref{lemma:unipotent_frobenius_lang}. Its image in
\(\mathcal D_{\cris,N}(P)\) gives a local point of
\(\mathcal D_{\dR,N}(P)^\varphi\).

Thus the Frobenius pseudo-reduction is a genuine reduction. Together
with Proposition~\ref{proposition:neutral_hodge_reduction}, this
proves admissibility. Functoriality follows from the intrinsic Hodge
and fixed-point constructions.
\end{proof}

\subsubsection{The Local and Global Logarithms}
\label{subsubsection:neutral_bk_logarithms}

\begin{proposition}[Local neutral Bloch--Kato logarithm]
\label{proposition:local_neutral_bk_logarithm}
The assignment
\[
P
\longmapsto
\left(
\mathcal D_{\dR,N}(P),
\mathcal D_{\dR,N}(P)^H,
\mathcal D_{\dR,N}(P)^\varphi
\right)
\]
on the local neutral Bloch--Kato prestack induces a canonical
morphism of algebraic stacks
\[
\BK_{n,v}:
\mathscr S_{n,v}
\longrightarrow
\mathscr M_n^{\dR}.
\]
\end{proposition}

\begin{proof}
Propositions~\ref{proposition:relative_padic_hodge_realization},
\ref{proposition:relative_hodge_pseudoreduction},
and~\ref{proposition:relative_frobenius_pseudoreduction}
construct the realized torsor and its canonical pseudo-reductions,
functorially under pullback and isomorphisms. On the neutral fibre,
Propositions~\ref{proposition:neutral_hodge_reduction}
and~\ref{proposition:neutral_realization_admissible} show that these
are genuine reductions.

We therefore obtain a pseudonatural functor from the local neutral
Bloch--Kato prestack to the stack of admissible
\(\mathcal G_n^{\dR}\)-torsors. It extends uniquely across
stackification, and the moduli interpretation of
Proposition~\ref{proposition:de_rham_period_moduli} gives the
displayed morphism.
\end{proof}

Let
\[
\widetilde{\mathscr S}_{n,T'}
\coloneqq 
\widetilde{\mathscr N}_{f,T'}
\bigl(
\mathcal G_n^{\et}/\mathcal R^{\et}
\bigr)_{K_v},
\]
and let
\[
q_{n,T'}:
\widetilde{\mathscr S}_{n,T'}
\longrightarrow
\mathscr S_{n,T'}
\]
be the morphism forgetting the chosen trivialization of the
reductive pushout.

\begin{definition}[Neutral Bloch--Kato logarithm]
\label{definition:neutral_bk_logarithm}
Localization at \(v\) induces
\[
\operatorname{loc}_v:
\mathscr S_{n,T'}
\longrightarrow
\mathscr S_{n,v}.
\]
The \emph{finite-level neutral Bloch--Kato logarithm} is
\[
\BK_{n,T'}
\coloneqq 
\BK_{n,v}\circ\operatorname{loc}_v:
\mathscr S_{n,T'}
\longrightarrow
\mathscr M_n^{\dR}.
\]
When \(T'\) is fixed, we write
\[
\BK_n\coloneqq \BK_{n,T'}.
\]

The corresponding morphism on the trivialized neutral fibre is
\[
\widetilde{\BK}_{n,T'}
\coloneqq 
\BK_{n,T'}\circ q_{n,T'}:
\widetilde{\mathscr S}_{n,T'}
\longrightarrow
\mathscr M_n^{\dR}.
\]
When \(T'\) is fixed, we write
\[
\widetilde{\BK}_n\coloneqq \widetilde{\BK}_{n,T'}.
\]
\end{definition}

\subsubsection{Kummer Compatibility and Functoriality}
\label{subsubsection:bk_compatibility}

Let \(\Sigma^{\mathrm{glob}}\subseteq \mathcal{X}(\mathcal{O}_{K,S})\) be a global Kummer stratum containing \(b\), and set
\[
\Sigma_{b}^{\mathrm{glob}}
\coloneqq 
\left\{
x\in\Sigma^{\mathrm{glob}}
\ \middle|\
\operatorname{loc}_v(x)\in]\underline b[(K_v)
\right\}.
\]
The finite-level neutral Kummer map of
Theorem~\ref{thm:kummer_stratification} restricts to
\[
\kappa_{\Sigma,n}^{\mathrm{neu}}:
\Sigma_{b}^{\mathrm{glob}}
\longrightarrow
\mathscr S_{n,T'}(K_v).
\]

\begin{theorem}[Compatibility with the period map]
\label{theorem:bk_period_compatibility}
The diagram
\[
\begin{tikzcd}[column sep=large,row sep=large]
\Sigma_{b}^{\mathrm{glob}}
  \arrow[d, "\kappa_{\Sigma,n}^{\mathrm{neu}}"]
  \arrow[r, "\operatorname{loc}_v"]
&
  ]\underline b[(K_v)
  \arrow[d, "{[\Phi_{n,v}]}"']
  \\
  \mathscr S_{n,T'}(K_v)
  \arrow[r, "\BK_n"]&
\mathscr M_n^{\dR}(K_v)
\end{tikzcd}
\]
is \(2\)-commutative.
\end{theorem}

\begin{proof}
Let \(x\in\Sigma_{b}^{\mathrm{glob}}\), and write
\[
x_v\coloneqq \operatorname{loc}_v(x).
\]
The localization of the neutral Kummer class is represented by the
finite-level étale path torsor
\[
{}^{\vphantom{\et}}_{x_v}P_{b,n}^{\et}.
\]
The étale--crystalline and crystalline--de Rham comparison
isomorphisms induce a canonical isomorphism
\[
\mathcal D_{\dR,K_v}
\left(
{}^{\vphantom{\et}}_{x_v}P_{b,n}^{\et}
\right)
\xrightarrow{\sim}
{}^{\vphantom{\dR}}_{x_v}P_{b,n}^{\dR}.
\]
Under this isomorphism,
Proposition~\ref{proposition:selected_finite_level_filtered_comparison}
identifies the \(B_{\dR}^+\)-condition defining
\[
\mathcal D_{\dR,K_v}
\left(
{}^{\vphantom{\et}}_{x_v}P_{b,n}^{\et}
\right)^H
\]
with the canonical Hodge reduction
\[
F^0{}^{\vphantom{\dR}}_{x_v}P_{b,n}^{\dR},
\]
while Frobenius compatibility identifies the corresponding
Frobenius reductions. Hence there is a canonical isomorphism in
\(\mathscr M_n^{\dR}(K_v)\),
\[
\BK_n
\left(
\kappa_{\Sigma,n}^{\mathrm{neu}}(x)
\right)
\simeq
[{}^{\vphantom{\dR}}_{x_v}P_{b,n}^{\dR}].
\]
By Proposition~\ref{proposition:an_rep}, the right-hand side is
canonically isomorphic to
\[
[\Phi_{n,v}](x_v),
\]
which proves the assertion.
\end{proof}

\begin{proposition}[Functoriality]
\label{proposition:neutral_bk_functoriality}
Let
\[
\pi_{n',n}^\bullet:
\mathcal G_{n'}^\bullet
\twoheadrightarrow
\mathcal G_n^\bullet
\]
be a morphism of compatible finite-type motivic quotients. For
\[
\lambda\in\{v,T'\},
\]
the diagram
\[
\begin{tikzcd}[column sep=large,row sep=large]
\mathscr S_{n',\lambda}
  \arrow[r, "\BK_{n',\lambda}"]
  \arrow[d]
&
\mathscr M_{n'}^{\dR}
  \arrow[d]
\\
\mathscr S_{n,\lambda}
  \arrow[r, "\BK_{n,\lambda}"']
&
\mathscr M_n^{\dR}
\end{tikzcd}
\]
is \(2\)-commutative, where
\(\BK_{n,v}\) is the local logarithm and
\(\BK_{n,T'}\) is the global logarithm.

After precomposition with
\[
q_{n',T'}:
\widetilde{\mathscr S}_{n',T'}
\longrightarrow
\mathscr S_{n',T'},
\]
the global square gives the corresponding \(2\)-commutative diagram
for the trivialized logarithms.

Consequently, for every compatible cofiltered system
\(\{\mathcal G_n^\bullet\}_n\), the family
\[
\BK
\coloneqq 
(\BK_{n,T'})_n
\]
defines a morphism of pro-stacks
\[
\left\{
\mathscr S_{n,T'}
\right\}_n
\longrightarrow
\left\{
\mathscr M_n^{\dR}
\right\}_n.
\]
\end{proposition}

\begin{proof}
We first prove the local assertion. Let \(N\) be an affine
\(K_v\)-scheme and let \(P_{n'}\) be an object of the local neutral
Bloch--Kato prestack at level \(n'\). Set
\[
P_n
\coloneqq 
P_{n'}\times^{\mathcal G_{n'}^{\et}}\mathcal G_n^{\et}.
\]

Functoriality of the finite-level étale--crystalline comparison gives
a canonical isomorphism
\[
\mathcal D_{\cris,N}(P_n)
\xrightarrow{\sim}
\mathcal D_{\cris,N}(P_{n'})
\times^{\mathcal G_{n'}^{\cri}}
\mathcal G_n^{\cri}.
\]
After extension to \(K_v\), compatibility of the
crystalline--de Rham comparison with
\(\pi_{n',n}^\bullet\) gives
\[
\mathcal D_{\dR,N}(P_n)
\xrightarrow{\sim}
\mathcal D_{\dR,N}(P_{n'})
\times^{\mathcal G_{n'}^{\dR}}
\mathcal G_n^{\dR}.
\]
Functoriality of the selected \(B_{\dR}^+\)-condition and of the
Frobenius fixed-point construction gives equivariant morphisms from
the Hodge and Frobenius reductions obtained on the right by extension
of structure group to the corresponding reductions on the left. On
the neutral locus all these pseudo-reductions are genuine torsors.
The induced morphisms after contracted product are therefore
isomorphisms. Thus the two resulting admissible
\(\mathcal G_n^{\dR}\)-torsors are canonically isomorphic.

These isomorphisms are pseudonatural in \(N\) and in \(P_{n'}\), and
therefore descend under stackification to the asserted
\(2\)-commutative local square.

For the global assertion, extension of structure group commutes with
restriction from \(G_{T}\) to \(G_v\). Hence the localization square
\[
\begin{tikzcd}[column sep=large]
\mathscr S_{n',T'}
  \arrow[r, "\operatorname{loc}_v"]
  \arrow[d]
&
\mathscr S_{n',v}
  \arrow[d]
\\
\mathscr S_{n,T'}
  \arrow[r, "\operatorname{loc}_v"']
&
\mathscr S_{n,v}
\end{tikzcd}
\]
is \(2\)-commutative. Since
\[
\BK_{n,T'}
=
\BK_{n,v}\circ\operatorname{loc}_v,
\]
the global square follows by composing this square with the local
one.

The forgetful maps from the trivialized neutral fibres are natural
under \(\pi_{n',n}^{\et}\), giving the corresponding trivialized
diagram.
\end{proof}

\subsection{The Dimension Cutter}
\label{section:selmer_fibred_prod}
We form the Selmer pullback and combine its dimension bound with the
density of the \(v\)-adic period map to obtain the analytic and
numerical dimension-cutter criteria.

Retain the trivialized neutral Selmer space and logarithm
\[
\widetilde{\mathscr S}_{n,T'},
\qquad
\widetilde{\BK}_n:
\widetilde{\mathscr S}_{n,T'}
\longrightarrow
\mathscr M_n^{\dR}
\]
from
\S\ref{chapter:bk_log}. By
Theorem~\ref{theorem:neutral_selmer_algebraicity},
\(\widetilde{\mathscr S}_{n,T'}\) is an algebraic space of finite
type over \(K_v\).

\paragraph{The Selmer image.}

\begin{definition}[The Selmer pullback]
\label{definition:selmer_pullback}
Define the \(2\)-fibre product
\[
\widetilde{\mathscr C}_n
\coloneqq 
\widetilde{\mathscr S}_{n,T'}
\times_{\mathscr M_n^{\dR}}
\operatorname{Fil}_n^{\dR}.
\]
Thus \(\widetilde{\mathscr C}_n\) parametrizes a trivialized neutral
Selmer class together with a trivialization of the Frobenius
reduction of its de Rham realization, equivalently a lift of that
realization from
\(\mathscr M_n^{\dR}\) to
\(\operatorname{Fil}_n^{\dR}\).
\end{definition}

\begin{proposition}[Dimension of the Selmer pullback]
\label{proposition:selmer_pullback_dimension}
The space \(\widetilde{\mathscr C}_n\) is an algebraic space of finite
type over \(K_v\), and the projection
\[
\widetilde{\mathscr C}_n
\longrightarrow
\widetilde{\mathscr S}_{n,T'}
\]
is a
\(\mathcal G_{n,K_v}^{\varphi=1}\)-torsor. In particular,
\[
\dim\widetilde{\mathscr C}_n
=
\dim\widetilde{\mathscr S}_{n,T'}
+
\dim\mathcal G_{n,K_v}^{\varphi=1}.
\]
\end{proposition}

\begin{proof}
The quotient morphism
\[
\pi_n:
\operatorname{Fil}_n^{\dR}
\longrightarrow
\left[
\mathcal G_{n,K_v}^{\varphi=1}
\backslash
\operatorname{Fil}_n^{\dR}
\right]
=
\mathscr M_n^{\dR}
\]
is a
\(\mathcal G_{n,K_v}^{\varphi=1}\)-torsor. Its base change along
\(\widetilde{\BK}_n\) is precisely
\[
\widetilde{\mathscr C}_n
\longrightarrow
\widetilde{\mathscr S}_{n,T'}.
\]
The representability, finite-type assertion, and dimension formula
follow.
\end{proof}

Retain
\(\operatorname{Fil}_n^{\dR,\circ}\)
from \S\ref{section:density}, and set
\[
\widetilde{\mathscr C}_n^\circ
\coloneqq 
\widetilde{\mathscr C}_n
\times_{\operatorname{Fil}_n^{\dR}}
\operatorname{Fil}_n^{\dR,\circ}.
\]

Let
\[
\mathscr W_n
\hookrightarrow
\operatorname{Fil}_n^{\dR,\circ}
\]
be the scheme-theoretic image of the finite-type morphism
\[
\widetilde{\mathscr C}_n^\circ
\longrightarrow
\operatorname{Fil}_n^{\dR,\circ};
\]
thus \(\mathscr W_n\) is the smallest closed subscheme through which
this morphism factors. Its underlying topological space is the
Zariski closure of the image of
\(\widetilde{\mathscr C}_n^\circ\), and hence
\[
\dim\mathscr W_n
\leq
\dim\widetilde{\mathscr C}_n^\circ
\leq
\dim\widetilde{\mathscr C}_n.
\]
By
Proposition~\ref{proposition:selmer_pullback_dimension}, this gives
\begin{equation}
\label{eq:selmer_image_dimension_bound}
\dim\mathscr W_n
\leq
\dim\widetilde{\mathscr S}_{n,T'}
+
\dim\mathcal G_{n,K_v}^{\varphi=1}.
\end{equation}

Fix a Kummer stratum \(\Sigma^{\mathrm{glob}}\) containing \(b\), and retain the
notation
\[
\Sigma_{b}^{\mathrm{glob}}
\coloneqq 
\left\{
x\in\Sigma^{\mathrm{glob}}
\ \middle|\
\operatorname{loc}_v(x)\in]\underline b[(K_v)
\right\}.
\]

\begin{lemma}[Containment of the integral period image]
\label{lemma:integral_period_image_containment}
One has
\[
\Phi_{n,v}
\left(
\operatorname{loc}_v(\Sigma_{b}^{\mathrm{glob}})
\right)
\subseteq
\mathscr W_n(K_v).
\]
\end{lemma}

\begin{proof}
Let \(x\in\Sigma_{b}^{\mathrm{glob}}\), and set
\[
x_v\coloneqq \operatorname{loc}_v(x).
\]
By
Theorem~\ref{theorem:bk_period_compatibility}, there is a canonical
isomorphism in \(\mathscr M_n^{\dR}(K_v)\),
\[
\BK_n
\left(
\kappa_{\Sigma,n}^{\mathrm{neu}}(x)
\right)
\simeq
\pi_n
\left(
\Phi_{n,v}(x_v)
\right).
\]

After an fppf extension \(A/K_v\), lift the neutral Kummer class to
\[
\widetilde\kappa_x
\in
\widetilde{\mathscr S}_{n,T'}(A).
\]
Together with the preceding isomorphism, this gives an \(A\)-point of
\(\widetilde{\mathscr C}_n\) mapping to
\(\Phi_{n,v}(x_v)_A\). Since the period map takes values in
\(\operatorname{Fil}_n^{\dR,\circ}\), this point lies in
\(\widetilde{\mathscr C}_n^\circ(A)\), and hence
\(\Phi_{n,v}(x_v)_A\) factors through \(\mathscr W_{n,A}\).

Factorization through a closed subscheme descends under faithfully
flat base change. Therefore
\[
\Phi_{n,v}(x_v)\in\mathscr W_n(K_v).
\]
\end{proof}

\paragraph{The analytic cutter.}

\begin{theorem}[The dimension cutter]
\label{theorem:dim_ineq}
Suppose that
\[
\dim\mathscr W_n
<
\dim\operatorname{Fil}_n^{\dR,\circ}.
\]
Then \(\Sigma_{b}^{\mathrm{glob}}\) is finite.

More precisely, there are a closed affinoid disk
\[
D_b^+
\subseteq
]\underline b[
\]
containing \(]\underline b[(K_v)\) and a nonzero analytic function
\[
F_{n,b}\in\mathcal O(D_b^+)
\]
such that
\[
F_{n,b}
\left(
\operatorname{loc}_v(x)
\right)
=
0
\qquad
(x\in\Sigma_{b}^{\mathrm{glob}}).
\]
\end{theorem}

\begin{proof}
Choose an \'{e}tale local parameter \(t\) at \(b\), identifying
\(]\underline b[\) with the open unit disk, and fix a uniformizer
\(\varpi_v\) of \(K_v\). Since \(K_v\) is discretely valued, every
\(K_v\)-point of the tube satisfies
\[
|t|\leq|\varpi_v|.
\]
Thus
\[
D_b^+
\coloneqq 
\left\{
x\in]\underline b[
\ \middle|\
|t(x)|\leq|\varpi_v|
\right\}
\]
is a closed affinoid disk containing \(]\underline b[(K_v)\).

Consider the closed analytic subspace
\[
\mathcal Z_{n,b}
\coloneqq 
\left(
\Phi_{n,v}|_{D_b^+}
\right)^{-1}
\left(
\mathscr W_n^{\an}
\right)
\subseteq
D_b^+.
\]
It is proper. Indeed, if
\(\mathcal Z_{n,b}=D_b^+\), then every sufficiently small
neighborhood \(D_v\subseteq D_b^+\) of \(b\) would satisfy
\[
\Phi_{n,v}(D_v)
\subseteq
\mathscr W_n^{\an},
\]
contradicting the algebraic Zariski density of
\(\Phi_{n,v}|_{D_v}\) from
Proposition~\ref{proposition:p-adic-dense}, since
\(\mathscr W_n\) is a proper closed subscheme of
\(\operatorname{Fil}_n^{\dR,\circ}\).

The defining ideal of \(\mathcal Z_{n,b}\) in
\(\mathcal O(D_b^+)\) is therefore nonzero. Choose
\[
0\neq F_{n,b}
\in
\mathcal O(D_b^+)
\]
in this ideal. By
Lemma~\ref{lemma:integral_period_image_containment},
\(F_{n,b}\) vanishes at
\(\operatorname{loc}_v(x)\) for every \(x\in\Sigma_{b}^{\mathrm{glob}}\).

A nonzero analytic function on a closed one-dimensional affinoid disk
has only finitely many zeros, by Weierstrass preparation. Hence
\[
\operatorname{loc}_v(\Sigma_{b}^{\mathrm{glob}})
\]
is finite. Finally,
\[
X(K)\longrightarrow X(K_v)
\]
is injective by fpqc descent, so \(\Sigma_{b}^{\mathrm{glob}}\) is finite.
\end{proof}

\begin{corollary}[Numerical dimension criterion]
\label{corollary:dim_ineq}
If
\begin{equation}
\label{eq:neutral_dimension_inequality}
\dim\widetilde{\mathscr S}_{n,T'}
+
\dim\mathcal G_{n,K_v}^{\varphi=1}
<
\dim\operatorname{Fil}_n^{\dR,\circ},
\end{equation}
then
\[
\#\Sigma_{b}^{\mathrm{glob}}<\infty.
\]

Equivalently, using invariance of dimension under scalar extension,
it suffices that
\[
\dim
\widetilde{\mathscr N}_{f,T'}
\left(
\mathcal G_n^{\et}/\mathcal R^{\et}
\right)
+
\dim\mathcal G_n^{\varphi=1}
+
\dim F^0\mathcal G_n^{\dR}
<
\dim\mathcal G_n^{\dR}.
\]
\end{corollary}

\begin{proof}
Condition~\eqref{eq:neutral_dimension_inequality} and
\eqref{eq:selmer_image_dimension_bound} give
\[
\dim\mathscr W_n
<
\dim\operatorname{Fil}_n^{\dR,\circ},
\]
so Theorem~\ref{theorem:dim_ineq} applies. The equivalent form uses
\[
\dim\operatorname{Fil}_n^{\dR,\circ}
=
\dim\mathcal G_n^{\dR}
-
\dim F^0\mathcal G_n^{\dR}.
\]
\end{proof}

\begin{lemma}[Reduction to the unipotent radical]
\label{lemma:unipotent_dimension_reduction}
With notation as above, one has:
\begin{enumerate}
\item $
\dim\mathcal G_n^{\varphi=1}
=
\dim\mathcal R^{\varphi=1}$, and

\item $\dim F^0\mathcal G_n^{\dR}
\leq
\dim F^0\mathcal U_n^{\dR}
+
\dim F^0\mathcal R^{\dR}$.
\end{enumerate}
\end{lemma}

\begin{proof}
The first statement follows from~\cite[Proposition~6.4.7]{KANTORTHESIS}.

For the second statement, observe that
\[
\mathrm{Ker}\,\left(F^0\mathcal G_n^{\dR}
\longrightarrow
F^0\mathcal R^{\dR}\right)\subseteq F^0\mathcal U_n^{\dR}.
\]
\end{proof}

The following is an immediate consequence of
Lemma~\ref{lemma:unipotent_dimension_reduction} and
Corollary~\ref{corollary:dim_ineq}.

\begin{corollary}[Unipotent numerical criterion]
\label{corollary:unipotent_dimension_criterion}
If
\begin{equation}
\label{eq:unipotent_dimension_criterion}
\dim\widetilde{\mathscr S}_{n,T'}
+
\dim F^0\mathcal U_n^{\dR}
+
\dim F^0\mathcal R^{\dR}
<
\dim\mathcal U_n^{\dR},
\end{equation}
then
\[
\#\Sigma_{b}^{\mathrm{glob}}<\infty.
\]
\qed
\end{corollary}
﻿\section{Applications to Fundamental Groups}
\label{chapter:applications}

This section applies the unipotent numerical criterion of
Corollary~\ref{corollary:unipotent_dimension_criterion} to two
classes of curves over \(\mathbb Q\). Assuming the relevant
Bloch--Kato dimension formula, we construct finite-type motivic
quotients satisfying the dimension inequality and deduce the
finiteness of integral points.

\begin{itemize}
\item
\S\ref{section:eis} treats modular curves with enough Eisenstein
classes. A suitable Eisenstein quotient of the abelianized
unipotent radical satisfies the dimension inequality; no higher
unipotent truncation is required.

\item
\S\ref{section:hyperbolic} treats curves of genus at least \(2\).
Using the finite Kummer stratification and sufficiently deep
truncations of the adjoint-supported motivic tower constructed in
\S\ref{sec:param_families}, we obtain the dimension inequality at
every relevant Kummer stratum.
\end{itemize}

Throughout this section,
\[
K=\mathbb Q,
\qquad
p\notin S,
\qquad
T=S\cup\{p\}.
\]

\subsection{Modular Curves with Enough Eisenstein Classes}
\label{section:eis}

The use of Eisenstein quotients in relative Chabauty goes back to
\cite[\S7.3]{KANTORTHESIS}. In the modular setting, the required
dimension inequality is already obtained from a quotient of the
abelianized unipotent radical. We take \(T'=T\), so that the
trivialized neutral fibre is an ordinary unipotent Selmer space.

\subsubsection{Eisenstein Quotients}
\label{subsubsection:eisenstein_quotients}

Let
\[
f:
\mathcal E
\longrightarrow
X
\]
be the universal elliptic curve over the modular curve
\(X/\mathbb Q\). For each realization
\(\bullet\in\{\et,\cri,\dR\}\), set
\[
\mathbb V^\bullet
\coloneqq 
R^1_\bullet f_*\mathbf 1,
\qquad
V^\bullet
\coloneqq 
\omega_b^\bullet(\mathbb V^\bullet),
\qquad
\mathbb V_m^\bullet
\coloneqq 
\Sym^m\mathbb V^\bullet,
\qquad
V_m^\bullet
\coloneqq 
\Sym^mV^\bullet.
\]
The geometric monodromy group is
\(\SL_2^\bullet\), and \(V^\bullet\) is its standard representation.

\begin{definition}[Motivic Eisenstein class]
\label{definition:motivic_eisenstein_class}
A \emph{motivic Eisenstein class of degree \(m\)} is a compatible
system of rank-one direct summands
\[
\mathfrak e_m^\bullet
\hookrightarrow
H^1_\bullet
\left(
X,\mathbb V_m^\bullet
\right),
\qquad
\bullet\in\{\et,\cri,\dR\},
\]
under the realization comparison isomorphisms, such that its
geometric étale realization is the pure Tate representation
\[
\mathfrak e_m^{\et}
\simeq
\mathbb Q_p(-m-1)
\]
of \(G_{\mathbb Q}\). We call \(m\) a \emph{motivic Eisenstein
value}.
\end{definition}

Let \(\mathcal U^\bullet\) be the unipotent radical of the relative
completion. By
Fact~\ref{fact:u_ab_revised}, the
\(\mathbb V_m^\bullet\)-isotypic quotient of its abelianized Lie
algebra is
\[
\Lie(\mathcal U^\bullet)^{\ab}
\twoheadrightarrow
V_m^\bullet
\otimes
H^1_\bullet
\left(
X,\mathbb V_m^\bullet
\right)^\vee.
\]
A motivic Eisenstein class therefore induces the compatible quotient
\begin{equation}
\label{eq:eis_quotient}
\mathcal Q_m^\bullet
\coloneqq 
V_m^\bullet(m+1)
\simeq
\left(
V_m^\bullet
\right)^\vee(1)
\end{equation}
of \((\mathcal U^\bullet)^{\ab}\), where we identify a representation
with its associated vector group. Denote by $\mathcal{G}_m^\bullet$ the pushout of $\mathcal{G}^\bullet$ 
along the quotient $\mathcal{U}^\bullet \twoheadrightarrow \mathcal{Q}_m^\bullet$.

\begin{remark}[Betti interpretation]
Fix an embedding \(\mathbb Q\hookrightarrow\mathbb C\), write
\[
X(\mathbb C)
=
X_\Gamma
=
\Gamma\backslash\mathbb H,
\qquad
D
=
\Gamma\backslash\mathbb P^1(\mathbb Q),
\]
and let
\[
\mathbb V_\Gamma
=
R^1_{\Betti}(f_\Gamma)_*\mathbb Q.
\]
Eichler--Shimura identifies the weight-\(m+1\) part of
\[
H^1_{\Betti}
\left(
X_\Gamma,\Sym^m\mathbb V_\Gamma
\right)
\]
with the classes of cusp forms of weight \(m+2\) and their complex
conjugates. For \(m>0\), its mixed Hodge structure fits into the exact
sequence
\begin{equation}
\label{equation:zucker}
0
\longrightarrow
W_{m+1}
H^1_{\Betti}
\left(
X_\Gamma,\Sym^m\mathbb V_\Gamma
\right)
\longrightarrow
H^1_{\Betti}
\left(
X_\Gamma,\Sym^m\mathbb V_\Gamma
\right)
\longrightarrow
\bigoplus_{P\in D}\mathbb Q(-m-1)
\longrightarrow
0;
\end{equation}
see
\cite[Theorems~11.4--11.5]{Hain_HodgeDeRham_Chapter2016}
and~\cite{ZUCKER}.
The right-hand side is the Eisenstein quotient: a normalized
Eisenstein series determines a class mapping nontrivially to it, and
a Betti Eisenstein class is equivalently a splitting of one of
its one-dimensional factors.
\end{remark}

\begin{definition}
\label{definition:eis_quotient}
The modular curve \(X\) has \emph{enough Eisenstein classes} if there
are infinitely many integers \(m\geq0\) for which
\[
H^1_\bullet
\left(
X,\mathbb V_m^\bullet
\right)
\]
contains a motivic Eisenstein class.
\end{definition}

\begin{fact}[{\cite[\S3.1]{LOEFFLER}}]
\label{fact:what_has_enough_eis}
For every \(N\geq4\), the modular curve \(Y_1(N)\) has enough
Eisenstein classes. In particular,
\[
Y_1(4)_{\mathbb Z[1/2]}
\simeq
\mathbb P^1_{\mathbb Z[1/2]}
\setminus\{0,1,\infty\}
\]
has enough Eisenstein classes.
\end{fact}

\subsubsection{Selmer Estimates}
\label{subsubsection:eisenstein_ramification}

Let \(W\) be a finite-dimensional \(p\)-adic representation of
\(G_T\).

\begin{definition}
[{Bloch--Kato Selmer groups;
cf.~\cite[(3.7)]{BLOCH_KATO}
and Definition~\ref{definition:restricted_ramification_stack}}]
\label{definition:linear_bloch_kato_selmer_groups}
The \emph{\(T\)-ramified Bloch--Kato Selmer group} is
\[
H^1_f(G_T,W)
\coloneqq 
\ker
\left(
H^1(G_T,W)
\longrightarrow
\frac{H^1(G_p,W)}{H^1_f(G_p,W)}
\right),
\]
where
\[
H^1_f(G_p,W)
\coloneqq 
\ker
\left(
H^1(G_p,W)
\longrightarrow
H^1
\left(
G_p,W\otimes_{\Qp}B_{\cris}
\right)
\right).
\]

For
\[
S_p\subseteq T'\subseteq T,
\]
define
\[
H^1_{f,T'}(G_T,W)
\coloneqq 
\ker
\left(
H^1_f(G_T,W)
\longrightarrow
\bigoplus_{\ell\in T\setminus T'}
H^1(I_\ell,W)
\right).
\]
Thus there are natural inclusions
\[
H^1_{f,S_p}(G_T,W)
\hookrightarrow
H^1_{f,T'}(G_T,W)
\hookrightarrow
H^1_{f,T}(G_T,W)
=
H^1_f(G_T,W).
\]
We write
\[
H^1_f(\mathbb Q,W)
\coloneqq 
H^1_{f,S_p}(G_T,W)
\]
for the usual Bloch--Kato Selmer group, with the finite local
condition imposed at every finite place.
\end{definition}

We use the following dimension consequence of the Bloch--Kato
conjecture.

\begin{conjecture}
[{Bloch--Kato dimension formula;
cf.~\cite[Lemma~2.10]{BETTS_CORWIN}}]
\label{conj:BCL_BK}
Let \(W\) be a pure geometric \(p\)-adic representation of
\(G_{\mathbb Q}\) of weight at most \(-3\). Then
\[
\dim H^1_f(\mathbb Q,W)
=
\dim H^1_f(\mathbb Q_p,W)
-
\dim W^+,
\]
where \(W^+\) is the subspace fixed by complex conjugation.
\end{conjecture}

We shall also use the following local dimension formula.

\begin{proposition}
[{\cite[Proposition~2.8]{BELLAICHE}}]
\label{proposition:2.8_ballaiche}
Let \(W\) be a de Rham representation of \(G_p\). Then
\[
\dim H^1_f(\mathbb Q_p,W)
-
\dim H^0(\mathbb Q_p,W)
=
\dim
\frac{D_{\dR}(W)}{F^0D_{\dR}(W)}.
\]
In particular, if
\[
H^0(\mathbb Q_p,W)=0,
\]
then
\[
\dim H^1_f(\mathbb Q_p,W)
=
\dim W^{\dR}
-
\dim F^0W^{\dR}.
\]
\end{proposition}

\begin{lemma}[Uniform Eisenstein ramification bound]
\label{lemma:eisenstein_ramification_bound}
For every \(\ell\in S\) and every motivic Eisenstein value of \(m\),
\[
\dim
H^1
\left(
I_\ell,\mathcal Q_m^{\et}
\right)^{G_\ell/I_\ell}
\leq1.
\]
Consequently,
\begin{equation}
\label{eq:eisenstein_relaxed_selmer_bound}
\dim H^1_f(G_T,\mathcal Q_m^{\et})
\leq
\dim H^1_f(\mathbb Q,\mathcal Q_m^{\et})
+
\#S.
\end{equation}
\end{lemma}

\begin{proof}
By
Definition~\ref{definition:linear_bloch_kato_selmer_groups}
and the Hochschild--Serre sequence for
\[
1\longrightarrow I_\ell
\longrightarrow G_\ell
\longrightarrow G_\ell/I_\ell
\longrightarrow1,
\]
restriction to inertia factors through
\[
H^1(G_\ell,W)
\longrightarrow
H^1(I_\ell,W)^{G_\ell/I_\ell}.
\]
Consequently, there is an exact sequence
\[
0
\longrightarrow
H^1_f(\mathbb Q,\mathcal Q_m^{\et})
\longrightarrow
H^1_f(G_T,\mathcal Q_m^{\et})
\longrightarrow
\bigoplus_{\ell\in S}
H^1
\left(
I_\ell,\mathcal Q_m^{\et}
\right)^{G_\ell/I_\ell}.
\]
It therefore suffices to bound each local term by \(1\).

Fix \(\ell\in S\), and set
\[
W\coloneqq \mathcal Q_m^{\et}.
\]
Let
\[
\rho_W:
I_\ell\longrightarrow\operatorname{GL}(W)
\]
denote the inertia action on \(W\). By the monodromy theorem, there
is an open normal subgroup
\[
I'_\ell\subseteq I_\ell
\]
such that
\[
\rho_W(\sigma)
=
\exp\bigl(t_\ell(\sigma)N\bigr),
\qquad
\sigma\in I'_\ell,
\]
where
\[
t_\ell:I'_\ell\longrightarrow\mathbb Z_p(1)
\]
is the \(p\)-primary tame character and \(N\) is nilpotent. The kernel of \(t_\ell\) has
vanishing continuous cohomology with \(\Qp\)-coefficients. Thus, if
\(\tau\in I'_\ell\) maps to a topological generator of an open
subgroup of \(\mathbb Z_p(1)\), then
\[
H^1(I'_\ell,W)
\simeq
W/(\tau-1)W,
\]
a cocycle is determined by its value at \(\tau\), while the
coboundaries are the elements of \((\tau-1)W\).

If $m=0$, then
\[
W=\Qp(1).
\]
Since $\ell\neq p$, inertia acts trivially on $W$, and hence
\[
H^1(I'_\ell,W)
\simeq
W/(\tau-1)W
=
W
\]
is one-dimensional. Inflation--restriction and passage to
$G_\ell/I_\ell$-invariants therefore give
\[
\dim
H^1(I_\ell,W)^{G_\ell/I_\ell}
\leq1.
\]
We may henceforth assume that $m>0$.

Let $N_V$ denote the monodromy operator on $V^{\et}$, and let
$N_W$ denote the induced operator on
\[
W
\simeq
\left(\Sym^mV^{\et}\right)^\vee(1).
\]
Thus the operator $N$ above is $N_W$.
Moreover, as endomorphisms of \(W\),
\[
\tau-1
=
\rho_W(\tau)-1
=
\exp\bigl(t_\ell(\tau)N_W\bigr)-1
=
N_Wu_\tau(N_W),
\]
where
\[
u_\tau(N_W)
\coloneqq 
t_\ell(\tau)\id_W
+
\frac{t_\ell(\tau)^2}{2!}N_W
+
\frac{t_\ell(\tau)^3}{3!}N_W^2
+\cdots .
\]
This is a finite sum because \(N_W\) is nilpotent, and it is
invertible because its constant term \(t_\ell(\tau)\) is nonzero.
Consequently,
\[
(\tau-1)W=N_WW
\qquad\text{and}\qquad
\operatorname{coker}(\tau-1)
\simeq
\operatorname{coker}(N_W).
\]

\smallskip
\noindent
\textbf{Case $N_W\neq0$.}
Then $N_V\neq0$. Choose a basis $e_0,e_1$ of $V^{\et}$ such that
\[
N_V(e_0)=0,
\qquad
N_V(e_1)=e_0.
\]
The induced monodromy operator on $\Sym^mV^{\et}$ satisfies
\[
N_{\Sym^m V}
\left(e_0^{m-i}e_1^i\right)
=
i\,e_0^{m-i+1}e_1^{i-1},
\]
and therefore acts through a single Jordan block. Dualizing and Tate twisting preserve the Jordan-block sizes, so
\(N_W\), and hence \(\tau-1\), has one-dimensional cokernel. The
identification
\[
H^1(I'_\ell,W)\simeq W/(\tau-1)W
\]
and inflation--restriction therefore give
\[
\dim
H^1(I_\ell,W)^{G_\ell/I_\ell}
\leq1.
\]

\smallskip
\noindent
\textbf{Case $N_W=0$.}
We claim that $N_V=0$. Indeed, if $N_V\neq0$, then after choosing
$e_0,e_1$ as above one has
\[
N_{\Sym^m V}(e_1^m)
=
m\,e_0e_1^{m-1}
\neq0,
\]
since $m>0$. Dualizing and Tate twisting do not change whether the
monodromy operator vanishes, contradicting $N_W=0$.

Thus \(N_V=0\). Choose a finite Galois extension
\(L/\mathbb Q_\ell\) such that
\[
I_L\subseteq I'_\ell
\]
and \(V^{\et}|_{G_L}\) is unramified. By the
Néron--Ogg--Shafarevich criterion, the elliptic curve has good
reduction over \(L\).

Since \(I_\ell/I_L\) is finite and the coefficients are in
\(\mathbb Q_p\), Hochschild--Serre gives
\[
H^1(I_\ell,W)
\xrightarrow{\sim}
H^1(I_L,W)^{I_\ell/I_L}.
\]
Consequently, restriction induces an injection
\[
H^1(I_\ell,W)^{G_\ell/I_\ell}
\hookrightarrow
H^1(I_L,W)^{G_L/I_L}.
\]
Since \(W\) is unramified over \(L\), there is a
\(G_L/I_L\)-equivariant identification
\[
H^1(I_L,W)
\simeq
W(-1).
\]
Moreover,
\[
W(-1)
\simeq
\left(\Sym^mV^{\et}\right)^\vee
\]
is pure of weight \(-m<0\), because the elliptic curve has good
reduction over \(L\). Hence
\[
H^1(I_L,W)^{G_L/I_L}
\simeq
W(-1)^{G_L/I_L}
=
0.
\]
It follows that
\[
H^1(I_\ell,W)^{G_\ell/I_\ell}=0,
\]
which proves the required bound in this case.

Thus every local summand has dimension at most \(1\). Summing over
\(\ell\in S\) gives
\[
\dim H^1_f(G_T,\mathcal Q_m^{\et})
\leq
\dim H^1_f(\mathbb Q,\mathcal Q_m^{\et})
+
\#S,
\]
as required.
\end{proof}

\subsubsection{The Eisenstein Dimension Estimate}
\label{section:modular2}

The Hodge subgroup
\[
F^0\SL_2^{\dR}
\]
is the stabilizer of the Hodge line in
\(H^1_{\dR}(\mathcal E_b)\). It is therefore a Borel subgroup, and
\begin{equation}
\label{eq:hodge_sl2_dimension}
\dim F^0\SL_2^{\dR}=2.
\end{equation}

\begin{proposition}[Eisenstein dimension estimate]
\label{proposition:modular_main}
Let \(m\geq1\) be a motivic Eisenstein value satisfying
\begin{equation}
\label{eq:eisenstein_m_bound}
\left\lfloor
\frac{m+1}{2}
\right\rfloor
>
\#S+2.
\end{equation}
Assume Conjecture~\ref{conj:BCL_BK}. Then
\begin{equation}
\label{eq:modular_unipotent_inequality}
\dim H^1_f(G_T,\mathcal Q_m^{\et})
+
\dim F^0\mathcal Q_m^{\dR}
+
\dim F^0\SL_2^{\dR}
<
\dim\mathcal Q_m^{\dR}.
\end{equation}
\end{proposition}

\begin{proof}
The representation \(\mathcal Q_m^{\et}\) is geometric and pure of
weight \(-m-2\), so
\[
H^0(\mathbb Q_p,\mathcal Q_m^{\et})=0.
\]
Proposition~\ref{proposition:2.8_ballaiche} therefore gives
\[
\dim H^1_f(\mathbb Q_p,\mathcal Q_m^{\et})
=
\dim\mathcal Q_m^{\dR}
-
\dim F^0\mathcal Q_m^{\dR}.
\]
Together with
Conjecture~\ref{conj:BCL_BK} and
Lemma~\ref{lemma:eisenstein_ramification_bound}, this yields
\[
\dim H^1_f(G_T,\mathcal Q_m^{\et})
\leq
\dim\mathcal Q_m^{\dR}
-
\dim F^0\mathcal Q_m^{\dR}
-
\dim(\mathcal Q_m^{\et})^+
+
\#S.
\]

Let \(\sigma\) denote complex conjugation. Since the Weil pairing is
anti-invariant under \(\sigma\), the standard representation
\(V^{\et}\) has one-dimensional \(+1\)- and \(-1\)-eigenspaces.
Choosing corresponding eigenvectors \(e_+\) and \(e_-\), one finds
that
\[
\sigma
\left(
e_+^{\,m-i}e_-^{\,i}
\right)
=
(-1)^i e_+^{\,m-i}e_-^{\,i}.
\]
Hence
\[
\dim
\left(
\Sym^mV^{\et}
\right)^-
=
\left\lfloor
\frac{m+1}{2}
\right\rfloor.
\]
Since
\[
\mathcal Q_m^{\et}
\simeq
\left(
\Sym^mV^{\et}
\right)^\vee(1),
\]
the Tate twist interchanges the two eigenspaces, and therefore
\[
\dim(\mathcal Q_m^{\et})^+
=
\left\lfloor
\frac{m+1}{2}
\right\rfloor.
\]
The hypothesis~\eqref{eq:eisenstein_m_bound}, together with
\eqref{eq:hodge_sl2_dimension}, now gives
\eqref{eq:modular_unipotent_inequality}.
\end{proof}

\begin{theorem}[The modular dimension inequality]
\label{theorem:modular_main}
Let \(X/\mathbb Q\) be a modular curve with enough Eisenstein
classes, and assume Conjecture~\ref{conj:BCL_BK}. Then, for
some motivic Eisenstein value \(m\),
\begin{equation}
\label{ineq:main_modular}
\dim
\widetilde{\mathscr N}_{f,T}
\left(
\mathcal G_m^{\et}/\SL_2^{\et}
\right)
+
\dim F^0\mathcal Q_m^{\dR}
+
\dim F^0\SL_2^{\dR}
<
\dim\mathcal Q_m^{\dR}.
\end{equation}
\end{theorem}

\begin{proof}
The set of motivic Eisenstein values is infinite and hence unbounded.
Choose one satisfying
\eqref{eq:eisenstein_m_bound}. Since \(T'=T\), the neutral-fibre
comparison gives
\[
\widetilde{\mathscr N}_{f,T}
\left(
\mathcal G_m^{\et}/\SL_2^{\et}
\right)
\simeq
H^1_f(G_T,\mathcal Q_m^{\et}).
\]
The result follows from
Proposition~\ref{proposition:modular_main}.
\end{proof}

\begin{corollary}[Finiteness for modular curves]
\label{corollary:modular_finiteness}
Under the hypotheses of
Theorem~\ref{theorem:modular_main},
\[
\#\mathcal X(\mathbb Z[1/S])<\infty.
\]
\end{corollary}

\begin{proof}
The reductive Galois stratification of
Definition~\ref{definition:global_partition} is finite, as is the
set of residue tubes at \(p\). For each nonempty intersection of a
stratum with a residue tube, rebase at an integral point and apply
Theorem~\ref{theorem:modular_main} and
Corollary~\ref{corollary:unipotent_dimension_criterion}. Each such
intersection contains only finitely many integral points, and the
result follows by taking their finite union.
\end{proof}

\subsection{\texorpdfstring{The Genus \(\geq2\) Case}
{The Genus at least two Case}}
\label{section:hyperbolic}

We now work in the genus \(\geq2\) case. The finite
Kummer stratification of
Theorem~\ref{thm:kummer_stratification} allows us to impose all
away-from-\(p\) conditions, and hence to take
\[
T'=S_p.
\]

For \(n\geq1\), write
\[
\mathcal G_n^\bullet
\coloneqq 
\mathcal G_{\mathrm{ad},n}^\bullet
\]
for the \(n\)-th quotient in the adjoint-supported motivic tower of
Theorem~\ref{theorem:special_motivic_exists}, and let
\[
1
\longrightarrow
\mathcal U_n^\bullet
\longrightarrow
\mathcal G_n^\bullet
\longrightarrow
\mathcal R^\bullet
\longrightarrow
1
\]
be its defining extension. For \(k\leq n\), set
\[
W_k^\bullet
\coloneqq 
\Gr_k^{\DCS}\Lie(\mathcal U_n^\bullet).
\]
By Theorem~\ref{theorem:special_motivic_exists} and
Lemma~\ref{lemma:purity_adjoint_supported_gradeds}, these graded
pieces are nonzero, adjoint-supported, independent of \(n\geq k\),
and pure of weight \(-k\).

\subsubsection{The Neutral Selmer Estimate}
\label{subsubsection:hyperbolic_neutral_selmer}

We bound the dimension of the restricted neutral Selmer space by
separating the finitely many neutral inertia parameters from the
strict unipotent Selmer directions.

Recall that the trivialized restricted neutral Selmer space is the
rigidified neutral fibre
\[
\widetilde{\mathscr N}_{f,S_p}
\left(
\mathcal G_n^{\et}/\mathcal R^{\et}
\right)
\coloneqq 
\mathscr H^1_{f,S_p}
\left(
G_T,\mathcal G_n^{\et}
\right)
\times_{\mathscr H^1(G_T,\mathcal R^{\et})}
\{\mathbf1_{\mathcal R^{\et}}\}.
\]
Thus it classifies restricted Selmer
\(\mathcal G_n^{\et}\)-torsors together with a trivialization of
their reductive pushout. By
Proposition~\ref{proposition:kp_unipotent_selmer_algebraicity},
the unipotent Bloch--Kato stack
\[
\mathscr H^1_f
\left(
G_T,\mathcal U_n^{\et}
\right)
\]
is represented by an affine \(\Qp\)-scheme of finite type, which we
denote by
\[
H^1_f
\left(
G_T,\mathcal U_n^{\et}
\right).
\]
Since
\[
T\setminus S_p=S,
\]
Proposition~\ref{proposition:restricted_neutral_fibre_presentation}
gives a canonical equivalence
\[
\widetilde{\mathscr N}_{f,S_p}
\left(
\mathcal G_n^{\et}/\mathcal R^{\et}
\right)
\simeq
H^1_f
\left(
G_T,\mathcal U_n^{\et}
\right)
\times_{
\prod_{\ell\in S}
\mathscr H^1(I_\ell,\mathcal U_n^{\et})
}
\prod_{\ell\in S}
\mathscr Z_\ell
\left(
\mathcal G_n^{\et}/\mathcal U_n^{\et}
\right).
\]

Set
\[
\mathcal R_S^I
\coloneqq 
\prod_{\ell\in S}
\left(
\mathcal R^{\et}
\right)^{I_\ell},
\qquad
\mathcal G_{n,S}^I
\coloneqq 
\prod_{\ell\in S}
\left(
\mathcal G_n^{\et}
\right)^{I_\ell}.
\]
Write
\[
\mathscr Z_{n,S}
\coloneqq 
\prod_{\ell\in S}
\mathscr Z_\ell
\left(
\mathcal G_n^{\et}/\mathcal U_n^{\et}
\right).
\]
The same proposition identifies
\[
\mathscr Z_{n,S}
\simeq
\left[
\mathcal R_S^I/
\mathcal G_{n,S}^I
\right].
\]

Define
\[
\mathscr Y_n
\coloneqq 
\widetilde{\mathscr N}_{f,S_p}
\left(
\mathcal G_n^{\et}/\mathcal R^{\et}
\right)
\times_{\mathscr Z_{n,S}}
\mathcal R_S^I.
\]
Then
\[
\mathscr Y_n
\longrightarrow
\widetilde{\mathscr N}_{f,S_p}
\left(
\mathcal G_n^{\et}/\mathcal R^{\et}
\right)
\]
is a \(\mathcal G_{n,S}^I\)-torsor.

\begin{definition}
[Serre twisting;
cf.~{\cite[Definition~3.1.5 and Proposition~3.1.6]
{KANTORTHESIS}}]
\label{definition:serre_twisting}
Let \(\Gamma\) be a profinite group, let \(H/\Qp\) be an affine
algebraic group with continuous \(\Gamma\)-action, and let \(S\) be
a test object in the big \'{e}tale site of \(\Spec(\Qp)\), on which
\(\Gamma\) acts trivially. For a \(\Gamma\)-equivariant right
\(H_S\)-torsor \(Q\), its \emph{Serre twist} is the affine \(S\)-group
\[
{}^QH_S
\coloneqq 
\underline{\operatorname{Aut}}_{H_S}(Q)
\simeq
Q\times^{H_S}H_S,
\]
where \(H_S\) acts on itself by conjugation. It carries the induced
\(\Gamma\)-action. If \(Q\) represents
\(\alpha\in\mathscr H^1(\Gamma,H)(S)\), we also write
\({}^\alpha H_S\), or simply \({}^\alpha H\) when \(S\) is understood.

If \(Q'\) is another \(\Gamma\)-equivariant right \(H_S\)-torsor,
then
\[
\underline{\operatorname{Isom}}_{H_S}(Q,Q')
\]
is naturally a left \({}^{Q'}H_S\)-torsor and a right
\({}^{Q}H_S\)-torsor, via postcomposition and precomposition,
respectively.
\end{definition}

\begin{lemma}[Twisting and the fibres]
\label{lemma:neutral_selmer_twisting_fibres}
Let \(\Omega/\Qp\) be algebraically closed, and let
\[
r:\Spec(\Omega)\longrightarrow\mathcal R_S^I
\]
be a geometric point such that
\[
\mathscr Y_{n,r}
\coloneqq 
\mathscr Y_n\times_{\mathcal R_S^I}\Spec(\Omega)
\]
is nonempty. Choose \(y\in\mathscr Y_{n,r}(\Omega)\), let
\(Q_\alpha\) be its underlying
\(\mathcal U_{n,\Omega}^{\et}\)-torsor, and let
\[
\alpha\in
H^1_f(G_T,\mathcal U_n^{\et})(\Omega)
\]
be its class. Serre twisting relative to \(Q_\alpha\), followed by
forgetting the chosen inertial identifications, defines a morphism
\[
\mathscr Y_{n,r}
\longrightarrow
\mathscr H^1_{f,S_p}
\left(
G_T,{}^\alpha\mathcal U_n^{\et}
\right)
\]
whose geometric fibres have dimension at most
\[
\dim\mathcal G_{n,S}^I.
\]
Moreover,
\[
\dim
\mathscr H^1_{f,S_p}
\left(
G_T,{}^\alpha\mathcal U_n^{\et}
\right)
\leq
\sum_{k=1}^n
\dim H^1_f(\mathbb Q,W_k^{\et}).
\]
Consequently,
\begin{equation}
\label{eq:neutral_selmer_fibre_bound}
\dim\mathscr Y_{n,r}
\leq
\dim\mathcal G_{n,S}^I
+
\sum_{k=1}^n
\dim H^1_f(\mathbb Q,W_k^{\et}).
\end{equation}
\end{lemma}

\begin{proof}
Let \(Z/\Omega\) be a scheme and let
\(y'\in\mathscr Y_{n,r}(Z)\), with underlying
\(\mathcal U_{n,Z}^{\et}\)-torsor \(Q_\beta\). By
Definition~\ref{definition:serre_twisting},
\[
D_{\alpha,\beta}
\coloneqq 
\underline{\operatorname{Isom}}_
{\mathcal U_{n,Z}^{\et}}
(Q_{\alpha,Z},Q_\beta)
\]
is a right
\(({}^\alpha\mathcal U_n^{\et})_Z\)-torsor. This construction
commutes with base change and descent.

At \(p\), crystalline trivializations of \(Q_{\alpha,Z}\) and
\(Q_\beta\) induce a crystalline trivialization of
\(D_{\alpha,\beta}\). Writing
\[
r=(r_\ell)_{\ell\in S},
\]
the chosen identifications with the local object \(r_\ell\) induce
an \(I_\ell\)-fixed point of
\(D_{\alpha,\beta}|_{I_\ell}\). Thus \(D_{\alpha,\beta}\) satisfies
the strict Selmer conditions and defines the asserted morphism.

Evaluation gives a canonical isomorphism
\[
D_{\alpha,\beta}
\times^{({}^\alpha\mathcal U_n^{\et})_Z}
Q_{\alpha,Z}
\xrightarrow{\sim}
Q_\beta.
\]
Hence, after fixing the resulting strict Selmer object, the only
forgotten data are the inertial identifications. Their choices are
controlled by
\[
K_r
\coloneqq 
\prod_{\ell\in S}
\operatorname{Stab}_{(\mathcal G_n^{\et})^{I_\ell}}(r_\ell),
\]
which is a closed subgroup of \(\mathcal G_{n,S}^I\). Every nonempty
geometric fibre is therefore a \(K_r\)-torsor and has
dimension at most \(\dim\mathcal G_{n,S}^I\).

The descending central series is characteristic and therefore
descends to the Serre twist. Since inner automorphisms act trivially
on its graded pieces, there are canonical isomorphisms
\[
\Gr_k^{\DCS}
\Lie\left({}^\alpha\mathcal U_n^{\et}\right)
\simeq
W_k^{\et}\otimes_{\Qp}\Omega
\]
of \(G_T\)-representations. Set
\[
H\coloneqq {}^{\alpha}\mathcal U_n^{\et}.
\]
The proof of
Proposition~\ref{proposition:kp_unipotent_selmer_algebraicity}
therefore applies over \(\Omega\) and shows that
\(\mathscr H^1_f(G_T,H)\) is represented by an affine scheme of
finite type. Moreover, the inertia argument in
Proposition~\ref{proposition:local_neutral_locus} shows that
\[
\mathscr B_{I_\ell}(H)
\longrightarrow
\mathscr H^1(I_\ell,H)
\]
is representable and of finite type for every \(\ell\in S\). Since
\[
\mathscr H^1_{f,S_p}(G_T,H)
\simeq
\mathscr H^1_f(G_T,H)
\times_{
\prod_{\ell\in S}\mathscr H^1(I_\ell,H)
}
\prod_{\ell\in S}\mathscr B_{I_\ell}(H),
\]
the restricted Selmer stack is an algebraic space of finite type over
\(\Omega\).

At the \(k\)-th stage of the DCS tower, the
obstruction to lifting lies in
\(H^2(G_T,W_k^{\et})\otimes_{\Qp}\Omega\), while every nonempty
geometric lifting fibre is a torsor under a subquotient of
\[
H^1_{f,S_p}(G_T,W_k^{\et})\otimes_{\Qp}\Omega
=
H^1_f(\mathbb Q,W_k^{\et})\otimes_{\Qp}\Omega.
\]
The obstruction can only cut down the lifting locus. Induction
therefore gives
\[
\dim
\mathscr H^1_{f,S_p}
\left(
G_T,{}^\alpha\mathcal U_n^{\et}
\right)_\Omega
\leq
\sum_{k=1}^n
\dim H^1_f(\mathbb Q,W_k^{\et}).
\]
Combining this with the preceding fibre-dimension estimate proves the
result.
\end{proof}

\begin{proposition}[Uniform neutral Selmer bound]
\label{proposition:neutral_selmer_dimension_bound}
For every \(n\geq1\),
\begin{equation}
\label{eq:neutral_selmer_dimension_bound}
\dim
\widetilde{\mathscr N}_{f,S_p}
\left(
\mathcal G_n^{\et}/\mathcal R^{\et}
\right)
\leq
\dim\mathcal R_{n,S}^I
+
\sum_{k=1}^n
\dim H^1_f
\left(
\mathbb Q,W_k^{\et}
\right).
\end{equation}
\end{proposition}

\begin{proof}
If
\[
\widetilde{\mathscr N}_{f,S_p}
\left(
\mathcal G_n^{\et}/\mathcal R^{\et}
\right)
\]
is empty, there is nothing to prove. Otherwise, by the construction
of \(\mathscr Y_n\), the morphism
\[
\mathscr Y_n
\longrightarrow
\widetilde{\mathscr N}_{f,S_p}
\left(
\mathcal G_n^{\et}/\mathcal R^{\et}
\right)
\]
is a \(\mathcal G_{n,S}^I\)-torsor. Hence
\[
\dim
\widetilde{\mathscr N}_{f,S_p}
\left(
\mathcal G_n^{\et}/\mathcal R^{\et}
\right)
=
\dim\mathscr Y_n
-
\dim\mathcal G_{n,S}^I.
\]

By the dimension-of-fibres inequality and
Lemma~\ref{lemma:neutral_selmer_twisting_fibres},
\[
\dim\mathscr Y_n
\leq
\dim\mathcal R_S^I
+
\dim\mathcal G_{n,S}^I
+
\sum_{k=1}^n
\dim H^1_f
\left(
\mathbb Q,W_k^{\et}
\right).
\]
Subtracting
\(\dim\mathcal G_{n,S}^I\) gives
\eqref{eq:neutral_selmer_dimension_bound}.
\end{proof}

\subsubsection{Complex Conjugation and Dimension Growth}
\label{subsubsection:hyperbolic_dimension_growth}

Let
\[
c\in G_{\mathbb Q}
\]
be complex conjugation, and, for a $G_{\mathbb Q}$-representation $W$, write
\[
W^+\coloneqq W^{c=1}.
\]

\begin{lemma}[Complex conjugation on adjoint-supported representations]
\label{lemma:adjoint_supported_complex_conjugation}
Let \(W\) be a nonzero finite-dimensional
\(G_{\mathbb Q}\)-equivariant, adjoint-supported
\(\mathcal R^{\et}\)-representation. Then
\[
W^+\neq0.
\]
\end{lemma}

\begin{proof}
It suffices to prove the assertion after extension of scalars to
\(\overline{\mathbb Q}_p\), since taking invariants under the
involution \(c\) commutes with scalar extension. Full monodromy gives
\[
\left(
\mathcal R^{\et,\circ}
\right)_{\overline{\mathbb Q}_p}
\simeq
\prod_{i\in I}\Sp(V_i).
\]
Set
\[
\mathfrak a_i\coloneqq \mathfrak{sp}(V_i).
\]
Since \(W\) is adjoint-supported, semisimplicity gives the canonical
isotypic decomposition
\[
W_{\overline{\mathbb Q}_p}
=
\bigoplus_{i\in I}W_i,
\qquad
W_i
\simeq
\mathfrak a_i\otimes M_i,
\]
where
\[
M_i
\coloneqq 
\Hom_{\mathcal R^{\et,\circ}}
\left(
\mathfrak a_i,
W_{\overline{\mathbb Q}_p}
\right).
\]
Complex conjugation permutes the summands \(W_i\).

Choose \(i\) with \(W_i\neq0\). If \(c(i)\neq i\), then \(c\)
exchanges the distinct summands \(W_i\) and \(W_{c(i)}\). Thus, for
any \(0\neq w\in W_i\), the vector
\[
w+c(w)
\]
is nonzero and \(c\)-fixed.

Suppose instead that \(c(i)=i\). The polarization on \(V_i\) takes
values in \(\mathbb Q_p(-1)\), so \(c\) acts anti-symplectically on
\(V_i\). Hence its \(+1\)- and \(-1\)-eigenspaces in \(V_i\) are
nonzero Lagrangian subspaces. The canonical identification
\[
\mathfrak a_i
\simeq
\Sym^2(V_i)(1)
\]
then shows that
\[
\mathfrak a_i^+\neq0,
\qquad
\mathfrak a_i^-\neq0.
\]

The induced action of \(c\) on \(M_i\) makes the evaluation
isomorphism
\[
\mathfrak a_i\otimes M_i
\xrightarrow{\sim}
W_i
\]
\(c\)-equivariant. Therefore
\[
W_i^+
=
\left(
\mathfrak a_i^+\otimes M_i^+
\right)
\oplus
\left(
\mathfrak a_i^-\otimes M_i^-
\right).
\]
Since \(M_i\neq0\), at least one of \(M_i^+\) and \(M_i^-\) is
nonzero; since both eigenspaces of \(\mathfrak a_i\) are nonzero, it
follows that \(W_i^+\neq0\).

Thus
\[
\left(
W_{\overline{\mathbb Q}_p}
\right)^+\neq0,
\]
and consequently \(W^+\neq0\).
\end{proof}

Since the \(W_k^{\et}\) are nonzero, adjoint-supported, 
Lemma~\ref{lemma:adjoint_supported_complex_conjugation} gives
\begin{equation}
\label{eq:adjoint_graded_positive_invariants}
\dim
\left(
W_k^{\et}
\right)^+
\geq1
\qquad
(k\geq1).
\end{equation}

\begin{theorem}[The genus \(\geq2\) dimension inequality]
\label{theorem:main_theorem}
Assume Conjecture~\ref{conj:BCL_BK}. Then, for all sufficiently
large \(n\),
\begin{equation}
\label{ineq:main}
\dim
\widetilde{\mathscr N}_{f,S_p}
\left(
\mathcal G_n^{\et}/\mathcal R^{\et}
\right)
+
\dim F^0\mathcal U_n^{\dR}
+
\dim F^0\mathcal R^{\dR}
<
\dim\mathcal U_n^{\dR}.
\end{equation}
\end{theorem}

\begin{proof}
By the exactness and strictness assertion in
Corollary~\ref{corollary:finite_level_filtered_path_groupoid}, the
DCS extensions of
\(\mathcal U_n^{\dR}\) are strict for their canonical filtered de
Rham realizations. Additivity of dimensions therefore
gives
\[
\dim\mathcal U_n^{\dR}
-
\dim F^0\mathcal U_n^{\dR}
=
\sum_{k=1}^n
\left(
\dim W_k^{\dR}
-
\dim F^0W_k^{\dR}
\right).
\]

Since \(W_k^{\et}\) is pure of weight \(-k<0\),
\[
H^0
\left(
\mathbb Q_p,W_k^{\et}
\right)
=
0.
\]
Proposition~\ref{proposition:2.8_ballaiche} therefore gives
\begin{equation}
\label{eq:hyperbolic_local_dimension}
\dim W_k^{\dR}
-
\dim F^0W_k^{\dR}
=
\dim H^1_f
\left(
\mathbb Q_p,W_k^{\et}
\right).
\end{equation}
For \(k\geq3\), Conjecture~\ref{conj:BCL_BK} gives
\begin{equation}
\label{eq:hyperbolic_bk_dimension}
\dim H^1_f
\left(
\mathbb Q,W_k^{\et}
\right)
=
\dim H^1_f
\left(
\mathbb Q_p,W_k^{\et}
\right)
-
\dim
\left(
W_k^{\et}
\right)^+.
\end{equation}

Combining
Proposition~\ref{proposition:neutral_selmer_dimension_bound},
\eqref{eq:hyperbolic_local_dimension}, and
\eqref{eq:hyperbolic_bk_dimension}, and absorbing the fixed
contributions from \(k=1,2\), the local neutral loci, and
\(F^0\mathcal R^{\dR}\) into a constant \(C_0\), we obtain
\[
\begin{aligned}
&
\dim\mathcal U_n^{\dR}
-
\dim F^0\mathcal U_n^{\dR}
-
\dim
\widetilde{\mathscr N}_{f,S_p}
\left(
\mathcal G_n^{\et}/\mathcal R^{\et}
\right)
-
\dim F^0\mathcal R^{\dR}
\\
&\qquad\geq
\sum_{k=3}^n
\dim
\left(
W_k^{\et}
\right)^+
-
C_0.
\end{aligned}
\]
By
\eqref{eq:adjoint_graded_positive_invariants}, the right-hand side is
at least
\[
n-2-C_0,
\]
and is therefore positive for all sufficiently large \(n\). This is
precisely
\eqref{ineq:main}.
\end{proof}

\begin{corollary}[Finiteness in genus at least \(2\)]
\label{corollary:hyperbolic_main}
Assume Conjecture~\ref{conj:BCL_BK}. Then
\[
\#\mathcal X(\mathbb Z[1/S])<\infty.
\]
\end{corollary}

\begin{proof}
The estimate in
Theorem~\ref{theorem:main_theorem} is uniform under rebasing at an
integral point. Indeed, rebasing changes the \(G_T\)-action on the
unipotent radical by an inner twist, which leaves the graded
representations \(W_k^{\et}\) and the corresponding Hodge dimensions
unchanged, while
Lemma~\ref{lemma:neutral_selmer_twisting_fibres} gives the neutral
Selmer estimate uniformly under such twists. Since only finitely many
reductive types occur, we may therefore choose a single sufficiently
large \(n\) for which
\eqref{ineq:main} holds after every relevant rebasing.

Fix such an \(n\), and apply
Theorem~\ref{thm:kummer_stratification} to the quotient
\(\mathcal G_n^\bullet\) with \(T'=S_p\). This gives a finite Kummer
stratification of \(\mathcal X(\mathbb Z[1/S])\). There are also only
finitely many residue tubes at \(p\).

For every nonempty intersection of a level-\(n\) Kummer stratum with
a residue tube, choose an integral point \(b\) in that intersection
and rebase \(\mathcal G_n^\bullet\) at \(b\). By the choice of \(n\),
Theorem~\ref{theorem:main_theorem} applies, and
Corollary~\ref{corollary:unipotent_dimension_criterion} shows that
this intersection is finite. Since only finitely many such
intersections occur, their union
\(\mathcal X(\mathbb Z[1/S])\) is finite.
\end{proof}
﻿\appendix

\section{Comparison and Semisimplicity of Kodaira--Parshin Realizations}
\label{appendix:rmc}
This chapter is dedicated to the comparison and semisimplicity
properties of the Kodaira--Parshin realizations introduced in
\S\ref{section:rmc_realizations}.
\begin{itemize}
    \item \S\ref{section:nonabelian_comparisons} recalls the
    non-abelian comparison theorems of Olsson and verifies their
    hypotheses for the Kodaira--Parshin Gauss--Manin system.
    \item \S\ref{section:semisimplicity} establishes the
    semisimplicity of these realizations and proves
    Proposition~\ref{prop:semisimplicity_kp}.
\end{itemize}

\subsection{Non-Abelian Comparison Isomorphisms}
\label{section:nonabelian_comparisons}
This subsection is dedicated to the non-abelian comparison
isomorphisms for the relative completions and path torsors associated
with the Kodaira--Parshin realizations.
\begin{itemize}
    \item \S\ref{section:etale_crystalline_comparison_appendix}
    treats the \'{e}tale--crystalline comparison.
    \item \S\ref{section:crystalline_de_rham_comparison_appendix}
    treats the crystalline--de Rham comparison.
\end{itemize}
In each case, we recall the relevant theorem of Olsson, verify its
hypotheses for the Kodaira--Parshin Gauss--Manin system, and deduce
the comparison statements recorded in
\S\ref{section:rmc_realizations}.

\subsubsection{The \'{e}tale--Crystalline Comparison}
\label{section:etale_crystalline_comparison_appendix}

We deduce Theorem~\ref{theorem:etale_crystalline_comparison} from
Olsson's general non-abelian comparison theorem.

Let \(v\mid p\), let \(k\) be the residue field of \(K_v\), and set
\[
K_0\coloneqq \operatorname{Frac}W(k).
\]
Let \(\widehat{\mathcal X}\) be a smooth proper
\(\mathcal O_{K_v}\)-scheme, let
\(\mathcal D\subseteq\widehat{\mathcal X}\) be a relative
normal-crossings divisor, and set
\[
\mathcal X^\circ\coloneqq \widehat{\mathcal X}\setminus\mathcal D.
\]
Write
\[
X_v\coloneqq \mathcal X^\circ_{K_v}
=\widehat X_v\setminus D_v,
\qquad
X_k\coloneqq \mathcal X^\circ_k
=\widehat X_k\setminus D_k.
\]

Let
\[
(\mathfrak X,\mathfrak D)
\]
be the formal completion of the pair
\((\widehat{\mathcal X},\mathcal D)\) along its special fibre
\((\widehat X_k,D_k)\). Thus
\[
\mathfrak X
=
\widehat{\mathcal X}^{\,\wedge}_{\widehat X_k},
\qquad
\mathfrak D
=
\mathcal D\times_{\widehat{\mathcal X}}\mathfrak X.
\]
In particular, \((\mathfrak X,\mathfrak D)\) is a smooth formal
\(\mathcal O_{K_v}\)-lift of \((\widehat X_k,D_k)\). Let
\[
\mathfrak X_\eta^\circ
\coloneqq 
\mathfrak X_\eta\setminus\mathfrak D_\eta
\]
be its open rigid generic fibre. Denote by
\[
a:
\bigl(\mathfrak X_\eta^\circ\bigr)_{\et}
\longrightarrow
(X_v)_{\et}
\]
the natural morphism of \'{e}tale sites.

Let \(\mathcal V^{\et}\) be a smooth \(\mathbb Q_p\)-local system on
\(X_v\), and write
\[
\mathcal V^{\et,\mathrm{an}}\coloneqq a^*\mathcal V^{\et}
\]
for its pullback to the rigid-analytic \'{e}tale site. Let
\(\mathcal V^{\cri}\) be a filtered logarithmic convergent
\(F\)-isocrystal on the log special fibre
\((\widehat X_k,D_k)\), and write
\[
(\mathcal E,\nabla,\Fil^\bullet)
\]
for its realization on
\((\mathfrak X_\eta,\mathfrak D_\eta)\) as a vector bundle with
logarithmic integrable connection and Griffiths-transverse
filtration.

Fix a basis \(\mathscr B\) of the \'{e}tale site of
\(\widehat{\mathcal X}\) consisting of connected affine \'{e}tale charts
\[
U=\Spec(R)\longrightarrow\widehat{\mathcal X}
\]
that are small in the sense of
\cite[\S6.1]{OLSSON_NA_P_HODGE}, namely, that admit local
\'{e}tale coordinates compatible with the normal-crossings divisor
\(\mathcal D\). Such charts form an \'{e}tale basis by loc.~cit.

For \(U\in\mathscr B\), set
\[
U^\circ
\coloneqq 
U\times_{\widehat{\mathcal X}}\mathcal X^\circ,
\qquad
\mathfrak U
\coloneqq 
U^\wedge_{U_k}.
\]
Let
\[
\bar\eta\longrightarrow U^\circ_{K_v}
\]
be a geometric generic point. Olsson's relative Fontaine
construction associates with the pair \((\mathfrak U,\bar\eta)\)
filtered Frobenius period rings
\[
A_{\cris}(\mathfrak U,\bar\eta),
\qquad
B_{\cris}(\mathfrak U,\bar\eta)
=
A_{\cris}(\mathfrak U,\bar\eta)[1/p,1/t],
\]
equipped with an action of the corresponding local Galois group; see
\cite[\S5.2 and \S\S6.10.1,~6.13]
{OLSSON_NA_P_HODGE}.
We suppress \(\bar\eta\) from the notation when no ambiguity can
arise.

The logarithmic divided-power construction defining
\(A_{\cris}(\mathfrak U,\bar\eta)\) determines a logarithmic
convergent enlargement of the special fibre. Since
\(\mathcal V^{\cri}\) is a convergent isocrystal, it may be evaluated
on the rationalized enlargement. After extension of scalars to
\(B_{\cris}(\mathfrak U,\bar\eta)\), this gives a filtered $\varphi$-module denoted by
\[
\mathcal V^{\cri}
\bigl(B_{\cris}(\mathfrak U,\bar\eta)\bigr).
\]

\begin{definition}
[Association;
cf.~{\cite[\S6.13]{OLSSON_NA_P_HODGE}}
and~{\cite[Definition~7.7]{PRIDHAM_GALOIS_ACTIONS}}]
\label{definition:associated_etale_crystalline}
We say that \(\mathcal V^{\et}\) and \(\mathcal V^{\cri}\) are
\emph{associated} if they are equipped, for every pair
\((U,\bar\eta)\) as above, with an isomorphism
\[
\iota_{U,\bar\eta}:
\mathcal V_{\bar\eta}^{\et}
\otimes_{\mathbb Q_p}
B_{\cris}(\mathfrak U,\bar\eta)
\xrightarrow{\sim}
\mathcal V^{\cri}
\bigl(B_{\cris}(\mathfrak U,\bar\eta)\bigr)
\]
of filtered Frobenius
\(B_{\cris}(\mathfrak U,\bar\eta)\)-modules, equivariant for the
corresponding local Galois action. These isomorphisms are required to
be compatible with morphisms of small \'{e}tale charts and with
changes of geometric generic point.
\end{definition}

For a \(K_0\)-linear Tannakian category \(\mathcal C\), write
\[
\mathcal C_{K_v}
\]
for its Tannakian scalar extension to \(K_v\), and, for
\(E\in\mathcal C\), write \(E_{K_v}\) for its image in
\(\mathcal C_{K_v}\). We use the same convention for tensor
subcategories and tensor functors. Thus every exact
\(K_0\)-linear tensor functor from \(\mathcal C\) to a
\(K_v\)-linear Tannakian category factors canonically through
\(\mathcal C_{K_v}\).

\begin{lemma}
[Restriction to the open crystalline category;
see~{\cite[\S4.18 and \S4.19(ii),
especially~(4.19.2)]{OLSSON_NA_P_HODGE}}]
\label{lemma:logarithmic_crystalline_restriction}
Let
\[
j_k:
X_k\hookrightarrow\widehat X_k
\]
be the complementary open immersion. After forgetting Frobenius and
the filtration, restriction to the trivial-log open defines an exact
fully faithful \(K_0\)-linear tensor functor
\[
j_k^*:
\operatorname{Isoc}_{\mathrm{nilp}}
\bigl((\widehat X_k,D_k)/K_0\bigr)
\longrightarrow
\operatorname{Isoc}(X_k/K_0).
\]
Its essential image is closed under finite direct sums, tensor
products, duals, subquotients, and extensions. The functor \(j_k^*\)
is also compatible with evaluation at points of \(X_k\).

Consequently, if \(\mathcal E\) is a logarithmic convergent
isocrystal with nilpotent residues and
\[
\mathcal E_\circ\coloneqq j_k^*\mathcal E,
\]
then \(j_k^*\) identifies the thin and thick Tannakian categories
generated by \(\mathcal E\) with the corresponding categories
generated by \(\mathcal E_\circ\) in
\(\operatorname{Isoc}(X_k/K_0)\).
\end{lemma}

For the remainder of this subsection, set
\[
\mathcal V_\circ^{\cri}
\coloneqq 
j_k^*\mathcal V^{\cri},
\]
where the \(F\)-structure and filtration on
\(\mathcal V^{\cri}\) are forgotten when forming its Tannakian
category.

Assume that the following conditions hold:
\begin{enumerate}
\item \(\mathcal V^{\et}\) and \(\mathcal V^{\cri}\) are associated
in the sense of Definition~\ref{definition:associated_etale_crystalline};

\item the Tannakian monodromy groups of
\[
\left.\mathcal V^{\et}\right|_{X_{\overline K_v}}
\qquad\text{and}\qquad
\mathcal V_\circ^{\cri}
\]
are reductive;

\item \(\mathcal V^{\cri}\) has unipotent local monodromy along
\(D_k\), i.e.~its logarithmic residues are nilpotent; see
\cite[\S4.14]{OLSSON_NA_P_HODGE}.
\end{enumerate}

Let
\[
b,x\in\mathcal X^\circ(\mathcal O_{K_v}),
\]
and denote their geometric generic fibres by
\(\overline b,\overline x\) and their reductions by
\(\underline b,\underline x\).

For \(\epsilon\in\{0,\mathrm{ext}\}\), set
\[
\begin{aligned}
\mathcal C_{\et}^{0}
&\coloneqq 
\left\langle
\left.\mathcal V^{\et}\right|_{X_{\overline K_v}}
\right\rangle_\otimes,
&
\mathcal C_{\et}^{\mathrm{ext}}
&\coloneqq 
\left\langle
\left.\mathcal V^{\et}\right|_{X_{\overline K_v}}
\right\rangle_\otimes^{\mathrm{ext}},
\\
\mathcal C_{\cri}^{0}
&\coloneqq 
\langle\mathcal V_\circ^{\cri}\rangle_\otimes,
&
\mathcal C_{\cri}^{\mathrm{ext}}
&\coloneqq 
\langle\mathcal V_\circ^{\cri}\rangle_\otimes^{\mathrm{ext}}.
\end{aligned}
\]
We write
\[
{}^{\vphantom{\et,\epsilon}}_{\overline x}
P_{\overline b}^{\et,\epsilon}
\coloneqq 
\operatorname{Isom}^{\otimes}
\left(
\left.\omega_{\overline b}^{\et}\right|_{\mathcal C_{\et}^{\epsilon}},
\left.\omega_{\overline x}^{\et}\right|_{\mathcal C_{\et}^{\epsilon}}
\right)
\]
and
\[
{}^{\vphantom{\cri,\epsilon}}_{\underline x}
P_{\underline b}^{\cri,\epsilon}
\coloneqq 
\operatorname{Isom}^{\otimes}
\left(
\left.\omega_{\underline b}^{\cri}\right|_{\mathcal C_{\cri}^{\epsilon}},
\left.\omega_{\underline x}^{\cri}\right|_{\mathcal C_{\cri}^{\epsilon}}
\right).
\]
Thus \(\epsilon=0\) gives the path torsors for the thin Tannakian
categories, while \(\epsilon=\mathrm{ext}\) gives the path
torsors for the thick Tannakian categories.

For fixed \(b\) and
\(\epsilon\in\{0,\mathrm{ext}\}\), let
\[
\mathscr A_{b,v}^{\et,\epsilon}
\qquad\text{and}\qquad
\bigl(
\mathscr A_{\underline b}^{\cri,\epsilon},
\Fil^\bullet
\bigr)
\]
denote, respectively, the regular coordinate ind-sheaf and the
filtered logarithmic convergent Frobenius coordinate ind-isocrystal
obtained by applying Olsson's constructions
\(\mathbb V(\mathcal O(G_{\et}))\) and
\(\mathbb L(\mathcal O(G_{\dR}))\) to the corresponding thin or
thick Tannakian category; see
\cite[\S\S7.5--7.6, especially (7.6.1)--(7.6.4), and
Remark~7.8]{OLSSON_NA_P_HODGE}.

\begin{theorem}
[{Olsson's relative filtered \'{e}tale--crystalline association;
\cite[\S\S7.5--7.7, especially Proposition~7.7]
{OLSSON_NA_P_HODGE}}]
\label{theorem:olsson_etale_crystalline_comparison}
Under the preceding hypotheses, \(\Fil^\bullet\) is Olsson's
filtration defined by the filtered de Rham matrix coefficients, and
\[
\mathscr A_{b,v}^{\et,\epsilon}
\qquad\text{and}\qquad
\bigl(
\mathscr A_{\underline b}^{\cri,\epsilon},
\Fil^\bullet
\bigr)
\]
are associated, as algebra objects, in the sense of
Definition~\ref{definition:associated_etale_crystalline}.
Association for these ind-objects is understood stagewise.

For every integral section \(y\), writing \(\overline y\) and
\(\underline y\) for its geometric generic fibre and reduction,
evaluation gives canonical identifications
\[
\overline y^*
\mathscr A_{b,v}^{\et,\epsilon}
\xrightarrow{\sim}
\mathcal O
\left(
{}^{\vphantom{\et,\epsilon}}_{\overline y}
P_{\overline b}^{\et,\epsilon}
\right),
\qquad
\underline y^*
\mathscr A_{\underline b}^{\cri,\epsilon}
\xrightarrow{\sim}
\mathcal O
\left(
{}^{\vphantom{\cri,\epsilon}}_{\underline y}
P_{\underline b}^{\cri,\epsilon}
\right),
\]
compatible with the associated comparison.
\end{theorem}

\begin{corollary}
[{Olsson's pointwise \'{e}tale--crystalline comparison;
\cite[Theorems~1.8 and~1.11]{OLSSON_NA_P_HODGE};
cf.~\cite[Theorem~3.3.12]{KANTORTHESIS}}]
\label{corollary:olsson_pointwise_etale_crystalline_comparison}
For every integral section \(x\) and
\(\epsilon\in\{0,\mathrm{ext}\}\), there is a canonical isomorphism
of affine \(B_{\cris}\)-schemes
\[
\iota_{x,b}^{\epsilon}:
{}^{\vphantom{\et,\epsilon}}_{\overline x}
P_{\overline b}^{\et,\epsilon}
\otimes_{\mathbb Q_p}B_{\cris}
\xrightarrow{\sim}
{}^{\vphantom{\cri,\epsilon}}_{\underline x}
P_{\underline b}^{\cri,\epsilon}
\otimes_{K_0}B_{\cris}.
\]
For \(z\in\{b,x\}\), the corresponding equal-endpoint comparison
\(\iota_{z,z}^{\epsilon}\) is an isomorphism of affine
\(B_{\cris}\)-group schemes. With respect to these structure-group
comparisons, \(\iota_{x,b}^{\epsilon}\) is an isomorphism of
bitorsors.

These comparisons are \(G_v\)-equivariant, where \(G_v\) acts
diagonally on the \'{e}tale side and through \(B_{\cris}\) on the
crystalline side, and are compatible with Frobenius, composition of
paths, and the natural morphisms from
\(\epsilon=\mathrm{ext}\) to \(\epsilon=0\).
\end{corollary}

\begin{proof}
Specializing the relative associated ind-objects along the integral
sections and using Olsson's specialization from the relative period
rings to \(B_{\cris}\) gives the coordinate-algebra comparison of
\cite[Theorems~1.8 and~1.11]{OLSSON_NA_P_HODGE}. We define
\(\iota_{x,b}^{\epsilon}\) as the inverse of its spectrum. Applying
the same construction with fixed basepoint \(b\), and again with
fixed basepoint \(x\), gives the right and left structure-group
comparisons. Compatibility with the torsor actions, path
composition, Frobenius, and the maps
\(\mathrm{ext}\to0\) follows from the tensorial and basepoint-free
functoriality of Olsson's construction; see
\cite[\S\S8.29--8.31]{OLSSON_NA_P_HODGE}.
\end{proof}

For an ind-crystalline \(G_v\)-representation
\[
V=\varinjlim_\lambda V_\lambda,
\]
written as the filtered union of its finite-dimensional crystalline
\(G_v\)-subrepresentations, set
\[
D_{\cris}(V)
\coloneqq 
\varinjlim_\lambda D_{\cris}(V_\lambda)
=
\left(
V\otimes_{\mathbb Q_p}B_{\cris}
\right)^{G_v}.
\]

\begin{corollary}[Recovery by \(D_{\cris}\)]
For every integral section \(x\) and
\(\epsilon\in\{0,\mathrm{ext}\}\),
\[
\mathcal O
\left(
{}^{\vphantom{\et,\epsilon}}_{\overline x}
P_{\overline b}^{\et,\epsilon}
\right)
\]
is an ind-crystalline \(G_v\)-algebra, and there is a canonical
isomorphism of Frobenius algebras
\[
D_{\cris}
\left(
\mathcal O
\left(
{}^{\vphantom{\et,\epsilon}}_{\overline x}
P_{\overline b}^{\et,\epsilon}
\right)
\right)
\xrightarrow{\sim}
\mathcal O
\left(
{}^{\vphantom{\cri,\epsilon}}_{\underline x}
P_{\underline b}^{\cri,\epsilon}
\right).
\]
\end{corollary}

\begin{proof}
By \cite[Theorem~1.11]{OLSSON_NA_P_HODGE}, the \'{e}tale coordinate
algebra is ind-crystalline. Taking \(G_v\)-invariants stagewise in
the coordinate-algebra comparison, and using
\(B_{\cris}^{G_v}=K_0\), gives the assertion.
\end{proof}

\begin{proof}[Proof of
Theorem~\ref{theorem:etale_crystalline_comparison}]
Let
\[
\fp:
\mathcal Y^\circ\longrightarrow\mathcal X^\circ
\]
be the model fixed in
Convention~\ref{convention:fixed_kp_family}, and denote its generic,
geometric generic, and special fibres by
\[
f_v:Y_v\longrightarrow X_v,
\qquad
\overline f:Y_{\overline K_v}\longrightarrow X_{\overline K_v},
\qquad
\underline f:Y_k\longrightarrow X_k.
\]
Set
\[
\mathbb L_v^{\et}
\coloneqq 
R^1_{\et}(f_v)_*\mathbb Q_p,
\qquad
\mathcal L^{\et}
\coloneqq 
R^1_{\et}\overline f_*\mathbb Q_p,
\qquad
\mathcal L^{\cri}
\coloneqq 
R^1_{\cri}\underline f_*\mathbf1.
\]
Smooth proper base change identifies
\[
\left.\mathbb L_v^{\et}\right|_{X_{\overline K_v}}
\xrightarrow{\sim}
\mathcal L^{\et},
\]
and the resulting descent datum is the natural \(G_v\)-action on
\(\mathcal L^{\et}\).

Let \(\mathcal L_{\log}^{\cri}\) be the logarithmic Gauss--Manin
\(F\)-isocrystal associated with the fixed proper log-smooth
compactification. Its restriction to the trivial-log open is
canonically identified with \(\mathcal L^{\cri}\):
\[
j_k^*\mathcal L_{\log}^{\cri}
\xrightarrow{\sim}
\mathcal L^{\cri}.
\]
For every \(U\in\mathscr B\) and every geometric generic point
\(\bar\eta\to U_{K_v}^\circ\), Faltings' relative crystalline
comparison gives an isomorphism
\[
(\mathbb L_v^{\et})_{\bar\eta}
\otimes_{\mathbb Q_p}
B_{\cris}(\mathfrak U,\bar\eta)
\xrightarrow{\sim}
\mathcal L_{\log}^{\cri}
\bigl(B_{\cris}(\mathfrak U,\bar\eta)\bigr).
\]
These isomorphisms are compatible with morphisms of charts, changes
of geometric generic point, the local Galois actions, the
filtrations, and the semilinear Frobenius operators; see
\cite[\S\S6.10.1 and~6.13]{OLSSON_NA_P_HODGE}.
Thus \(\mathbb L_v^{\et}\) and
\(\mathcal L_{\log}^{\cri}\) are associated in the sense of
Definition~\ref{definition:associated_etale_crystalline}, with
restrictions \(\mathcal L^{\et}\) and \(\mathcal L^{\cri}\).

In the modular case, the fixed semistable compactification gives
nilpotent residues for \(\mathcal L_{\log}^{\cri}\); in the
genus $\ge 2$ case, the boundary is empty. Moreover,
Proposition~\ref{prop:semisimplicity_kp} gives reductivity of the
corresponding Tannakian monodromy groups, and the chosen endpoints are
integral sections of \(\mathcal X^\circ\). Hence
Theorem~\ref{theorem:olsson_etale_crystalline_comparison} applies, and
Corollary~\ref{corollary:olsson_pointwise_etale_crystalline_comparison},
with \(\epsilon=\mathrm{ext}\), gives precisely the path-torsor
comparison asserted in the theorem.
\end{proof}

\subsubsection{The Crystalline--de Rham Comparison}
\label{section:crystalline_de_rham_comparison_appendix}

We now deduce
Theorem~\ref{theorem:crystalline_de_rham_comparison}
from Olsson's logarithmic crystalline realization together with
Shiho's proper log-smooth direct-image comparison.

Notations are as in
\S\ref{section:etale_crystalline_comparison_appendix}.
Retain the fixed smooth proper logarithmic
$\mathcal O_{K_v}$-model
\[
(\widehat{\mathcal X},\mathcal D)
\]
of $(\widehat X_v,D_v)$, and let
\[
(\mathfrak X,\mathfrak D)
\]
be its $p$-adic formal completion.

Let
\[
\operatorname{MIC}_{\mathrm{nilp}}(\widehat X_v,D_v)
\]
denote the category of vector bundles on $\widehat X_v$ equipped
with logarithmic integrable connections having nilpotent residues
along $D_v$. Restriction to $X_v$ identifies this with a full
Tannakian subcategory of $\operatorname{VIC}(X_v)$, closed under
subquotients and extensions; see
\cite[4.17.2]
{OLSSON_NA_P_HODGE}.

Retain the restriction functor \(j_k^*\) of
Lemma~\ref{lemma:logarithmic_crystalline_restriction}. Write
\[
j_v:
X_v\hookrightarrow\widehat X_v
\]
for the complementary open immersion.

Pull back the logarithmic convergent site along
\[
W(k)\longrightarrow\mathcal O_{K_v}.
\]
The closed immersion of the reduced special fibre
\[
(\widehat X_k,D_k)
\hookrightarrow
(\mathfrak X,\mathfrak D)
\]
makes \((\mathfrak X,\mathfrak D)\) an enlargement in the resulting
logarithmic convergent site. Evaluation on this enlargement, followed
by proper rigid GAGA, gives a \(K_0\)-linear tensor functor from
\(\operatorname{Isoc}_{\mathrm{nilp}}
\bigl((\widehat X_k,D_k)/K_0\bigr)\) to logarithmic connections over \(K_v\).
By the universal property of Tannakian scalar extension, it extends,
uniquely up to unique tensor isomorphism, to a \(K_v\)-linear tensor
functor
\[
\operatorname{ev}_{\log,v}:
\left(\operatorname{Isoc}_{\mathrm{nilp}}
\bigl((\widehat X_k,D_k)/K_0\bigr)\right)_{K_v}
\longrightarrow
\operatorname{MIC}_{\mathrm{nilp}}(\widehat X_v,D_v).
\]
Define its restriction to the trivial-log open by
\[
\operatorname{ev}_v
\coloneqq 
j_v^*\circ\operatorname{ev}_{\log,v}:
\left(\operatorname{Isoc}_{\mathrm{nilp}}
\bigl((\widehat X_k,D_k)/K_0\bigr)\right)_{K_v}
\longrightarrow
\operatorname{VIC}(X_v).
\]

\begin{theorem}
[Olsson's logarithmic crystalline realization and its Tannakian
consequences;
see~{\cite[Lemmas~4.17.2 and 4.24]{OLSSON_NA_P_HODGE}}]
\label{theorem:olsson_crystalline_de_rham_comparison}
The functors
\[
\operatorname{ev}_{\log,v}
\qquad\text{and}\qquad
\operatorname{ev}_v
\]
are exact and fully faithful, and their essential images are closed
under finite direct sums, tensor products, duals, subquotients, and
extensions.

Consequently, let
\[
\mathcal V_{\log}^{\cri}
\in
\operatorname{Isoc}_{\mathrm{nilp}}
\bigl((\widehat X_k,D_k)/K_0\bigr),
\qquad
\mathcal V_\circ^{\cri}
\coloneqq 
j_k^*\mathcal V_{\log}^{\cri},
\]
and set
\[
\mathcal V_{\log}^{\dR}
\coloneqq 
\operatorname{ev}_{\log,v}
\bigl((\mathcal V_{\log}^{\cri})_{K_v}\bigr),
\qquad
\mathcal V_v^{\dR}
\coloneqq 
j_v^*\mathcal V_{\log}^{\dR}
=
\operatorname{ev}_v
\bigl((\mathcal V_{\log}^{\cri})_{K_v}\bigr).
\]
Combining
Lemma~\ref{lemma:logarithmic_crystalline_restriction}
with \(\operatorname{ev}_v\) gives tensor equivalences
\[
\left(
\langle\mathcal V_\circ^{\cri}\rangle_\otimes
\right)_{K_v}
\xrightarrow{\sim}
\langle\mathcal V_v^{\dR}\rangle_\otimes
\]
and
\[
\left(
\langle\mathcal V_\circ^{\cri}\rangle_\otimes^{\mathrm{ext}}
\right)_{K_v}
\xrightarrow{\sim}
\langle\mathcal V_v^{\dR}\rangle_\otimes^{\mathrm{ext}}.
\]

For an \(\mathcal O_{K_v}\)-valued point lifting a point of the
special fibre, these equivalences are compatible with the
corresponding fibre functors. Hence they identify the associated
Tannakian fundamental groups and path torsors after extension of
scalars to \(K_v\).

If \(\mathcal V_{\log}^{\cri}\) carries an \(F\)-structure and its
de Rham realization carries a filtration, the \(K_0\)-realization
with Frobenius together with its filtered \(K_v\)-realization forms
the corresponding filtered Frobenius realization.
\end{theorem}

\begin{proof}[Proof of
Theorem~\ref{theorem:crystalline_de_rham_comparison}]
Let \(v\mid p\), and let
\[
\fp:
\mathcal Y^\circ\longrightarrow\mathcal X^\circ
\]
be the smooth proper \(\mathcal O_{K_v}\)-model supplied by
Convention~\ref{convention:fixed_kp_family}, with generic and special
fibres
\[
f_v:
Y_v\longrightarrow X_v,
\qquad
\underline f:
Y_k\longrightarrow X_k.
\]

In the modular case, write
\[
\widehat f_v:
(\widehat Y_v,E_v)
\longrightarrow
(\widehat X_v,D_v)
\]
for the generic fibre of the fixed proper log-smooth semistable
compactification of \(\fp\); its restriction to the trivial-log open
is \(f_v\). In the genus $\ge 2$ case, put
\[
\widehat f_v\coloneqq f_v,
\qquad
E_v=D_v=\varnothing.
\]

Retain the logarithmic crystalline Gauss--Manin object
\(\mathcal L_{\log}^{\cri}\) introduced above and its canonical
identification
\[
j_k^*\mathcal L_{\log}^{\cri}
\xrightarrow{\sim}
\mathcal L^{\cri}
=
R^1_{\cri}\underline f_*\mathbf1,
\]
and let $\mathcal L_{\log}^{\dR}$ denote the logarithmic de Rham Gauss--Manin bundle, equipped with its
logarithmic Gauss--Manin connection. Restriction to the trivial-log
open gives a canonical identification
\[
j_v^*\mathcal L_{\log}^{\dR}
\xrightarrow{\sim}
\mathcal L^{\dR} = R^1_{\dR}(f_v)_*\mathbf1.
\]

The fixed compactification is proper log smooth and integral and has
log-smooth parameter. After base change along
\[
W(k)\longrightarrow\mathcal O_{K_v},
\]
Shiho's crystalline--convergent comparison identifies the
logarithmic crystalline Gauss--Manin object with the convergent
direct-image isocrystal; see
\cite[Theorem~2.36 and Theorem~4.8]
{SHIHO_RELATIVE_LOG_CONVERGENT}.
By
\cite[Definition~4.1 and Corollary~4.10]
{SHIHO_RELATIVE_LOG_CONVERGENT},
the evaluation of this isocrystal on the fixed formally log-smooth
enlargement is the relative logarithmic analytic cohomology, namely
the higher direct image of the logarithmic de Rham complex on the
rigid generic fibre. Since the fixed enlargement is obtained from
the proper algebraic log-smooth compactification, proper rigid GAGA
identifies this analytic Gauss--Manin bundle with the analytification
of
\[
\mathcal L_{\log}^{\dR}
=
R^1_{\log\text{-}\dR}
(\widehat f_v^{\log})_*\mathbf1.
\]
Compatibility with the base change
\(W(k)\to\mathcal O_{K_v}\) follows from
\cite[Corollary~3.10 and Lemma~3.11(1)]
{SHIHO_RELATIVE_LOG_CONVERGENT}.
Thus there is a canonical horizontal isomorphism
\[
\operatorname{ev}_{\log,v}
\bigl((\mathcal L_{\log}^{\cri})_{K_v}\bigr)
\xrightarrow{\sim}
\mathcal L_{\log}^{\dR}.
\]
Restricting along \(j_v\) gives
\[
\operatorname{ev}_v
\bigl((\mathcal L_{\log}^{\cri})_{K_v}\bigr)
\xrightarrow{\sim}
\mathcal L^{\dR}.
\]

Theorem~\ref{theorem:olsson_crystalline_de_rham_comparison}
therefore induces tensor equivalences
\[
\left(
\langle\mathcal L^{\cri}\rangle_\otimes
\right)_{K_v}
\xrightarrow{\sim}
\langle\mathcal L^{\dR}\rangle_\otimes
\]
and
\[
\left(
\langle\mathcal L^{\cri}\rangle_\otimes^{\mathrm{ext}}
\right)_{K_v}
\xrightarrow{\sim}
\langle\mathcal L^{\dR}\rangle_\otimes^{\mathrm{ext}}.
\]

For every
\[
z\in\mathcal X^\circ(\mathcal O_{K_v})
\]
with reduction \(\underline z\), functoriality under pullback gives
\[
\omega_{\underline z}^{\cri}\otimes_{K_0}K_v
\xrightarrow{\sim}
\omega_z^{\dR}.
\]
The preceding equivalences are therefore compatible with the fibre
functors. Tannakian duality gives
\[
\mathcal R^{\cri}\otimes_{K_0}K_v
\xrightarrow{\sim}
\mathcal R^{\dR},
\qquad
\mathcal G^{\cri}\otimes_{K_0}K_v
\xrightarrow{\sim}
\mathcal G^{\dR},
\]
and, for
\[
b,x\in\mathcal X^\circ(\mathcal O_{K_v}),
\]
gives
\[
{}^{\vphantom{\cri}}_{\underline x}P_{\underline b}^{\cri}
\otimes_{K_0}K_v
\xrightarrow{\sim}
{}^{\vphantom{\dR}}_{x}P_b^{\dR}.
\]

The Frobenius is carried by the crystalline \(K_0\)-realization.
The canonical Hodge filtration on the de Rham group and path torsors
is constructed independently in
Theorem~\ref{theorem:relative_filtered_path_groupoid}; its
finite-level counterpart is given by
Corollary~\ref{corollary:finite_level_filtered_path_groupoid}.
Together with the underlying crystalline--de Rham comparison above,
these structures give the asserted filtered Frobenius realizations.
\end{proof}

\subsection{Semisimplicity of Kodaira--Parshin Realizations}
\label{section:semisimplicity}

\begin{proof}[Proof of Proposition~\ref{prop:semisimplicity_kp}]
Retain the standing abelian-by-finite factorization
\[
Y^\dagger
\xrightarrow{g^\dagger}
X'^\dagger
\xrightarrow{\pi^\dagger}
X^\dagger,
\qquad
f^\dagger=\pi^\dagger\circ g^\dagger.
\]
Thus, in the modular case,
\[
X'^\dagger=X^\dagger,
\qquad
\pi^\dagger=\id_{X^\dagger},
\]
while in the genus $\ge 2$ case, this is the fixed compactified
Kodaira--Parshin family.

Let
\[
\mathcal L^{\mathrm B}
\coloneqq 
R^1_{\mathrm B}(f^\dagger)_*\mathbf1
\]
be the Betti realization, and let \(\mathcal R^{\mathrm B}\) denote
its algebraic monodromy group. For a complex point
\[
x_0\in X^\dagger(\mathbb C),
\]
the abelian-by-finite decomposition gives
\[
\mathcal L^{\mathrm B}_{x_0}
\simeq
\bigoplus_{x'\in(\pi^\dagger)^{-1}(x_0)}
H^1_{\mathrm B}(Y^\dagger_{x'},\mathbb Q).
\]
The monodromy permutes these summands and preserves their
polarizations. Full monodromy therefore gives
\[
\prod_{x'\in(\pi^\dagger)^{-1}(x_0)}
\operatorname{Sp}
\bigl(
H^1_{\mathrm B}(Y^\dagger_{x'},\mathbb Q),
\omega_{x'}
\bigr)
\subseteq
\mathcal R^{\mathrm B}
\subseteq
\left(
\prod_{x'\in(\pi^\dagger)^{-1}(x_0)}
\operatorname{Sp}
\bigl(
H^1_{\mathrm B}(Y^\dagger_{x'},\mathbb Q),
\omega_{x'}
\bigr)
\right)\rtimes\Pi
\]
for a finite permutation group \(\Pi\). Hence
\[
(\mathcal R^{\mathrm B})^\circ
=
\prod_{x'\in(\pi^\dagger)^{-1}(x_0)}
\operatorname{Sp}
\bigl(
H^1_{\mathrm B}(Y^\dagger_{x'},\mathbb Q),
\omega_{x'}
\bigr),
\]
so \(\mathcal R^{\mathrm B}\) is reductive and
\(\mathcal L^{\mathrm B}\) is semisimple.

Artin comparison identifies the geometric \'{e}tale monodromy of
\(\mathcal L^{\et}\), after extension of scalars, with the Betti
monodromy of \(\mathcal L^{\mathrm B}\). Thus
\(\mathcal R^{\et}\) is reductive and \(\mathcal L^{\et}\) is
semisimple.

Likewise, applying de Rham--Betti comparison to the global \(K\)-form
of the Gauss--Manin connection and then base-changing from \(K\) to
\(K_v\) shows that \(\mathcal R^{\dR}\) is reductive and
\(\mathcal L^{\dR}\) is semisimple.

Finally, Theorem~\ref{theorem:crystalline_de_rham_comparison} gives
\[
\mathcal R^{\cri}\otimes_{K_0}K_v
\simeq
\mathcal R^{\dR}.
\]
Since reductivity descends under field extension,
\(\mathcal R^{\cri}\) is reductive. Hence
\(\mathcal L^{\cri}\) is semisimple as well.
\end{proof}

\section{The Hodge Filtration on the Universal Relative de Rham Path Torsor}
\label{appendix:filtered_de_rham_path_groupoid}
This chapter constructs over \(K\) the Hodge filtration on the
universal relative de Rham path torsor based at \(b\), and on the
compatible finite-type motivic quotients used in
\S\ref{section:analyticity}. After recording the canonical
connection on Tannakian families of path torsors, we proceed as
follows.

\begin{itemize}
    \item \S\ref{appendix:filtered_reductive_path_torsor} constructs
    the filtered reductive path-coordinate algebra through its matrix
    coefficients.
    \item \S\ref{appendix:relative_de_rham_bar_model} constructs the
    relative two-sided bar model and identifies its zeroth cohomology
    with the universal path-coordinate algebra.
    \item \S\ref{appendix:universal_relative_hodge_filtration}
    transports the bar filtration to the universal path torsor and
    compares it with Hain's Hodge filtration and Olsson's filtered
    realization.
    \item \S\ref{appendix:finite_level_filtered_path_torsors} passes
    to compatible finite-type motivic quotients and constructs their
    Hodge reductions and filtered period comparison.
\end{itemize}

Set
\[
X^\dagger
\coloneqq 
\begin{cases}
X, & \text{in the modular case},\\
\widehat X, & \text{in the genus $\ge 2$ case},
\end{cases}
\]
as in Convention~\ref{convention:fixed_kp_family}, and fix
\[
b\in X^\dagger(K).
\]

We first record the canonical connection on a Tannakian family of
path torsors. Let
\[
\mathcal T
\subseteq
\operatorname{VIC}(X^\dagger/K)
\]
be a full Tannakian subcategory, and write
\[
\omega_b
\coloneqq 
\left.b^*\right|_{\mathcal T},
\qquad
\omega_{X^\dagger}
:
\mathcal T
\longrightarrow
\operatorname{VB}(X^\dagger)
\]
for the fibre functor at \(b\) and the functor obtained by forgetting
the connection. Set
\[
\mathcal G_{\mathcal T}
\coloneqq 
\operatorname{Aut}^{\otimes}(\mathcal T,\omega_b)
\]
and
\[
\mathcal P_{\mathcal T,b,-}
\coloneqq 
\operatorname{Isom}^{\otimes}
\left(
\omega_b\otimes_K\mathcal O_{X^\dagger},
\omega_{X^\dagger}
\right).
\]
Writing
\[
p_{\mathcal T}:
\mathcal P_{\mathcal T,b,-}
\longrightarrow
X^\dagger
\]
for the structure morphism, set
\[
\mathscr A_{\mathcal T,b}
\coloneqq 
(p_{\mathcal T})_*
\mathcal O_{\mathcal P_{\mathcal T,b,-}}.
\]

\begin{proposition}[Canonical connection on Tannakian path torsors]
\label{proposition:tannakian_path_torsor_connection}
The affine right \(\mathcal G_{\mathcal T}\)-torsor
\[
\mathcal P_{\mathcal T,b,-}
\longrightarrow
X^\dagger
\]
carries a canonical integrable connection. Dually,
\(\mathscr A_{\mathcal T,b}\) is a commutative algebra object in
\[
\operatorname{Ind}(\mathcal T)
\subseteq
\operatorname{Ind}\operatorname{VIC}(X^\dagger/K),
\]
whose canonical connection
\[
\nabla_{\mathcal T}:
\mathscr A_{\mathcal T,b}
\longrightarrow
\mathscr A_{\mathcal T,b}
\otimes_{\mathcal O_{X^\dagger}}
\Omega^1_{X^\dagger/K}
\]
is an integrable derivation.

The torsor action is horizontal. Equivalently, the unit,
multiplication, and coaction
\[
\Delta_{\mathcal T}:
\mathscr A_{\mathcal T,b}
\longrightarrow
\mathscr A_{\mathcal T,b}
\otimes_K
\mathcal O(\mathcal G_{\mathcal T})
\]
are horizontal, where
\(\mathcal O(\mathcal G_{\mathcal T})\) carries the trivial
connection.

The construction commutes with scalar extension and pullback on
\(X^\dagger\). Moreover, if
\[
\mathcal T'\subseteq\mathcal T
\]
is a full Tannakian subcategory, then the canonical pushout
isomorphism
\[
\mathcal P_{\mathcal T,b,-}
\times^{\mathcal G_{\mathcal T}}
\mathcal G_{\mathcal T'}
\xrightarrow{\sim}
\mathcal P_{\mathcal T',b,-}
\]
is horizontal.
\end{proposition}

\begin{proof}
Equip
\[
\omega_b\otimes_K\mathcal O_{X^\dagger}
\]
with its trivial connection, and let
\[
\omega_{X^\dagger}^{\nabla}:
\mathcal T
\longrightarrow
\operatorname{VIC}(X^\dagger/K)
\]
be the tautological tensor functor.

We first suppose that \(\mathcal T\) admits a tensor generator \(E\).
The tensor-isomorphism torsor
\(\mathcal P_{\mathcal T,b,-}\) is then a closed subscheme of
\[
\underline{\operatorname{Isom}}
\left(
\omega_b(E)\otimes_K\mathcal O_{X^\dagger},
\omega_{X^\dagger}(E)
\right),
\]
cut out by the equations expressing compatibility with all morphisms
in \(\mathcal T\), tensor products, duals, and the tensor unit.

The vector bundle
\[
\underline{\operatorname{Hom}}
\left(
\omega_b(E)\otimes_K\mathcal O_{X^\dagger},
\omega_{X^\dagger}(E)
\right)
\]
carries the induced connection
\[
\nabla(\alpha)
=
\nabla_E\circ\alpha
-
(\alpha\otimes1)\circ\nabla_{\mathrm{triv}}.
\]
It is integrable because \(\nabla_E\) and
\(\nabla_{\mathrm{triv}}\) are integrable. The induced connection on
the symmetric algebra of the dual Hom bundle is a derivation and
extends to the open isomorphism locus.

Every morphism in \(\mathcal T\) is horizontal, and the tensor,
duality, and unit constraints in
\(\operatorname{VIC}(X^\dagger/K)\) are horizontal. Hence the ideal
cutting out the tensor-isomorphism torsor is preserved by the
connection. The connection therefore descends to
\(\mathcal P_{\mathcal T,b,-}\), or equivalently to an integrable
derivation \(\nabla_{\mathcal T}\) on its coordinate algebra.

For a general \(\mathcal T\), let \(\mathcal T_0\) range over its
full finitely generated Tannakian subcategories. Then
\[
\mathcal P_{\mathcal T,b,-}
\simeq
\varprojlim_{\mathcal T_0}
\mathcal P_{\mathcal T_0,b,-},
\qquad
\mathscr A_{\mathcal T,b}
\simeq
\varinjlim_{\mathcal T_0}
\mathscr A_{\mathcal T_0,b}.
\]
The preceding connections are compatible under restriction of tensor
isomorphisms and therefore induce the asserted connection on the
limit torsor and its coordinate algebra. Relative Tannakian duality
identifies \(\mathscr A_{\mathcal T,b}\) with the realization of the
regular ind-representation of \(\mathcal G_{\mathcal T}\); hence it
is a commutative algebra object of
\(\operatorname{Ind}(\mathcal T)\).

Composition of tensor isomorphisms is horizontal for the induced Hom
connections. It follows that the torsor action is horizontal and,
dually, that the coaction \(\Delta_{\mathcal T}\) is horizontal. The
unit and multiplication are horizontal because
\(\nabla_{\mathcal T}\) is a derivation.

The entire construction is defined from tensor-isomorphism functors
and induced connections, and hence commutes with scalar extension and
pullback. If \(\mathcal T'\subseteq\mathcal T\), restriction of tensor
isomorphisms gives the pushout identification
\[
\mathcal P_{\mathcal T,b,-}
\times^{\mathcal G_{\mathcal T}}
\mathcal G_{\mathcal T'}
\xrightarrow{\sim}
\mathcal P_{\mathcal T',b,-},
\]
and the construction above shows that this identification is
horizontal.
\end{proof}

We now apply the preceding construction to the fixed
Kodaira--Parshin family. Let
\[
\mathcal L_K^{\dR}
\in
\operatorname{VIC}(X^\dagger/K)
\]
be its Gauss--Manin object, and set
\[
\mathcal D_K^{\dR}
\coloneqq 
\langle\mathcal L_K^{\dR}\rangle_\otimes,
\qquad
\mathcal C_K^{\dR}
\coloneqq 
\langle\mathcal L_K^{\dR}\rangle_\otimes^{\mathrm{ext}}.
\]
Using the notation introduced above, write
\[
\mathcal R_K^{\dR}
\coloneqq 
\mathcal G_{\mathcal D_K^{\dR}},
\qquad
\mathcal G_K^{\dR}
\coloneqq 
\mathcal G_{\mathcal C_K^{\dR}},
\]
and
\[
\mathcal P_{\mathcal R,b,-,K}^{\dR}
\coloneqq 
\mathcal P_{\mathcal D_K^{\dR},b,-},
\qquad
\mathcal P_{b,-,K}^{\dR}
\coloneqq 
\mathcal P_{\mathcal C_K^{\dR},b,-}.
\]
Their coordinate algebras are denoted
\[
\mathscr A_{\mathcal R,b,K}^{\dR}
\coloneqq 
\mathscr A_{\mathcal D_K^{\dR},b},
\qquad
\mathscr A_{b,K}^{\dR}
\coloneqq 
\mathscr A_{\mathcal C_K^{\dR},b}.
\]

For every field extension \(K'/K\) and every
\(x\in X^\dagger(K')\), there are canonical identifications
\[
\left.
\mathcal P_{\mathcal R,b,-,K}^{\dR}
\right|_x
\simeq
{}^{\vphantom{\dR}}_{x}P_{b,\mathcal R}^{\dR},
\qquad
\left.
\mathcal P_{b,-,K}^{\dR}
\right|_x
\simeq
{}^{\vphantom{\dR}}_{x}P_b^{\dR}.
\]
Since
\[
\mathcal D_K^{\dR}
\subseteq
\mathcal C_K^{\dR},
\]
Proposition~\ref{proposition:tannakian_path_torsor_connection} gives a
canonical horizontal pushout isomorphism
\[
\mathcal P_{b,-,K}^{\dR}
\times^{\mathcal G_K^{\dR}}
\mathcal R_K^{\dR}
\xrightarrow{\sim}
\mathcal P_{\mathcal R,b,-,K}^{\dR}.
\]

We write
\[
\nabla_{\mathcal R}
\quad\text{and}\quad
\nabla
\]
for the canonical integrable connections on
\(\mathscr A_{\mathcal R,b,K}^{\dR}\) and
\(\mathscr A_{b,K}^{\dR}\), respectively. The natural inclusion
\[
\mathscr A_{\mathcal R,b,K}^{\dR}
\hookrightarrow
\mathscr A_{b,K}^{\dR}
\]
dual to the reductive pushout is horizontal.

\subsection{The Filtered Reductive Path Torsor}
\label{appendix:filtered_reductive_path_torsor}

We construct the Hodge filtration directly on the reductive
path-coordinate algebra
\[
\mathscr A_{\mathcal R,b,K}^{\dR}
\]
through its matrix coefficients. Throughout, the polarizations of
the Gauss--Manin factors retain their Tate-valued targets.

Choose a finite Galois extension
\[
L/K
\]
such that
\[
H
\coloneqq 
\left(
\mathcal R_L^{\dR}
\right)^\circ
\]
is split, the component group
\[
\Pi
\coloneqq 
\pi_0
\left(
\mathcal R_L^{\dR}
\right)
\]
is constant, and every simple object of
\[
\mathcal D_L^{\dR}
\coloneqq 
\mathcal D_K^{\dR}\otimes_KL
\]
is absolutely simple. Let
\[
\left\{
\mathbb V_\lambda
\right\}_{\lambda\in\Lambda}
\]
be representatives of their isomorphism classes.

Tannakian duality identifies
\(\mathscr A_{\mathcal R,b,K}^{\dR}\otimes_KL\) with the regular
coend of the two endpoint fibre functors. Since
\(\mathcal D_L^{\dR}\) is semisimple, the Tannakian Peter--Weyl
decomposition gives
\begin{equation}
\label{eq:K_rational_matrix_coefficient_decomposition}
\mathscr A_{\mathcal R,b,K}^{\dR}
\otimes_KL
\simeq
\bigoplus_{\lambda\in\Lambda}
\omega_{b,L}^{\dR}
\left(
\mathbb V_\lambda
\right)^\vee
\otimes_L
\mathbb V_\lambda.
\end{equation}

For each \(\lambda\), choose a filtered algebraic lift
\[
\mathbb V_\lambda^F
\]
constructed from the filtered Gauss--Manin factors, their
Tate-valued polarizations, and the finite component-group data.
Equip the \(\lambda\)-th matrix-coefficient summand with the
tensor-product filtration on
\begin{equation}
\label{eq:filtered_matrix_coefficient_summand}
\omega_{b,L}^{\dR}
\left(
\mathbb V_\lambda^F
\right)^\vee
\otimes_L
\mathbb V_\lambda^F,
\end{equation}
and denote the resulting direct-sum filtration by
\[
F_L^\bullet
\left(
\mathscr A_{\mathcal R,b,K}^{\dR}\otimes_KL
\right).
\]
The opposite Tate twists on the fixed-endpoint dual factor and the
variable-endpoint factor make this filtration independent of the
choice of Hodge lifts, as the following lemma records.

\begin{lemma}
[\(K\)-rationality of the matrix-coefficient filtration;
cf.~{\cite[
Lemma~13.5,
Corollaries~13.6--13.7,
and Theorem~13.10]{HAIN_HDR_RMC}}]
\label{lemma:K_rational_reductive_matrix_coefficients}
The filtration \(F_L^\bullet\) is independent of the chosen filtered
lifts and invariant under
\(\operatorname{Gal}(L/K)\). It therefore descends uniquely to a
decreasing filtration
\[
F^\bullet
\mathscr A_{\mathcal R,b,K}^{\dR}.
\]
The descended filtration is independent of \(L\), of the
representatives \(\mathbb V_\lambda\), and of all auxiliary choices.

The unit and multiplication, the reductive torsor coaction, the
identity evaluation
\[
\epsilon_{\mathcal R}:
\mathcal O
\left(
\mathcal R_K^{\dR}
\right)
\longrightarrow
K,
\]
and reductive path composition are filtered. The construction
commutes with extension of the ground field and with specialization
to either endpoint. For every embedding
\[
\iota_\infty:
K\hookrightarrow\mathbb C,
\]
its base change through \(\iota_\infty\) is Hain's canonical
weight-zero Hodge filtration on the reductive path-coordinate
algebra.
\end{lemma}

\begin{proof}
By the full-monodromy description of the Kodaira--Parshin system,
after enlarging \(L\) if necessary one has
\[
H
\simeq
\prod_{i\in I}
\Sp(\mathbb W_i),
\]
where the \(\mathbb W_i\) are the standard Gauss--Manin factors. The
component group \(\Pi\) permutes these factors and their
Tate-valued polarizations.

Let \(\mathbb W_\xi\) be an irreducible representation of \(H\).
Weyl's construction realizes \(\mathbb W_\xi\) from the standard
representations \(\mathbb W_i\) by tensor products, Schur functors,
and intersections of kernels of contraction morphisms. We retain the
Tate-valued targets of the polarizations in these contractions.
Consequently, Weyl's construction takes place in filtered algebraic
connections and gives a lift
\[
\mathbb W_\xi^F
\]
of \(\mathbb W_\xi\). After extension along any embedding
\(L\hookrightarrow\mathbb C\), this is an admissible polarized
variation of Hodge structure.

We normalize this lift by the weight determined by the chosen Weyl
construction. After extension to \(\mathbb C\), Hain's uniqueness
theorem implies that any two Hodge lifts of
\(\mathbb W_\xi\) differ by a Tate twist. The weight normalization
therefore makes
\(\mathbb W_\xi^F\) independent, up to filtered isomorphism, of the
chosen Weyl presentation.

The component group \(\Pi\) acts on the set of irreducible
\(H\)-representations. Since it permutes the filtered standard
factors and their Tate-valued polarizations, it carries the
normalized filtered lift of \(\mathbb W_\xi\) to the normalized
filtered lift of its conjugate
\[
{}^\pi\mathbb W_\xi
\qquad
(\pi\in\Pi).
\]
Thus the normalized Weyl lifts are compatible with the
\(\Pi\)-action.

Now fix a simple object
\(\mathbb V_\lambda\in\mathcal D_L^{\dR}\), choose an irreducible
constituent \(\mathbb W_\xi\) of
\(\left.\mathbb V_\lambda\right|_H\), and let
\(\Pi_\xi\subseteq\Pi\) be its stabilizer. Clifford theory reconstructs
\(\mathbb V_\lambda\) from the \(\Pi\)-orbit of
\(\mathbb W_\xi\) together with finite projective
\(\Pi_\xi\)-data.

For \(\pi\in\Pi_\xi\), the normalized lifts
\(\mathbb W_\xi^F\) and \({}^{\pi}\mathbb W_\xi^F\) have the same
weight and underlie isomorphic simple local systems. Hence
\cite[Lemma~13.5]{HAIN_HDR_RMC} shows, after extension to
\(\mathbb C\), that their stabilizer intertwiners are filtered: the
possible Tate twist is trivial by the weight normalization. Faithful
flatness gives the same conclusion over \(L\). The resulting Clifford
cocycle is scalar and therefore has type \((0,0)\). Equipping the
finite projective stabilizer data with its type-\((0,0)\) filtration,
Clifford induction produces a filtered algebraic lift
\[
\mathbb V_\lambda^F
\]
of \(\mathbb V_\lambda\), whose complexification is an admissible
variation of Hodge structure.

Equivalently, writing
\[
M_{\lambda,\xi}
\coloneqq 
\Hom_H
\left(
\mathbb W_\xi,
\mathbb V_\lambda
\right),
\]
the \(H\)-isotypic decomposition of the corresponding
matrix-coefficient object is
\[
\begin{aligned}
&
\omega_{b,L}^{\dR}
\left(
\mathbb V_\lambda
\right)^\vee
\otimes_L
\mathbb V_\lambda
\\
&\qquad\simeq
\bigoplus_{\xi,\eta}
\omega_{b,L}^{\dR}
\left(
\mathbb W_\xi
\right)^\vee
\otimes_L
\mathbb W_\eta
\otimes_L
M_{\lambda,\xi}^\vee
\otimes_L
M_{\lambda,\eta},
\end{aligned}
\]
where \(\xi\) and \(\eta\) range over the common
\(\Pi\)-orbit of constituents. The multiplicity spaces carry the
trivial filtration, while the first two factors carry the
tensor-product filtration induced by the normalized Weyl lifts. The
component group permutes these summands through filtered
isomorphisms, so this filtration is precisely the filtration on
\eqref{eq:filtered_matrix_coefficient_summand}.

We next verify independence and descent. After extension along any
embedding \(L\hookrightarrow\mathbb C\), the chosen lift
\(\mathbb V_\lambda^F\) and any other compatible Hodge lift of
\(\mathbb V_\lambda\) are admissible variations on the same simple
local system. By
\cite[Lemma~13.5]{HAIN_HDR_RMC}, they differ by a Tate twist.
Such a twist cancels in the matrix-coefficient object:
\[
\omega_{b,L}^{\dR}
\left(
\mathbb V_\lambda^F(r)
\right)^\vee
\otimes_L
\mathbb V_\lambda^F(r)
\simeq
\omega_{b,L}^{\dR}
\left(
\mathbb V_\lambda^F
\right)^\vee
\otimes_L
\mathbb V_\lambda^F
\]
as filtered objects. Faithful flatness therefore shows that the
filtration on each matrix-coefficient summand is independent of the
chosen lift.

Let \(\Gamma\coloneqq \Gal(L/K)\). For \(\sigma\in\Gamma\), the same argument
applied to the conjugate simple object shows that \(\sigma\) carries
the filtration on the \(\lambda\)-summand to that on the
\(\sigma\lambda\)-summand. Hence the direct-sum filtration on
\[
\mathscr A_{\mathcal R,b,K}^{\dR}\otimes_KL
\]
is \(\Gamma\)-stable and descends uniquely to a filtration on
\(\mathscr A_{\mathcal R,b,K}^{\dR}\). Passing to a common finite
Galois extension proves independence of \(L\), of the representatives
\(\mathbb V_\lambda\), and of all auxiliary choices.

Multiplication of matrix coefficients is induced by tensor products
of representations and their decomposition into simple summands.
After extension to \(\mathbb C\), the construction above is Hain's
weight-zero matrix-coefficient construction, for which
multiplication, evaluation, coevaluation, the reductive torsor
coaction, and reductive path composition are filtered. Faithful
flatness therefore gives the same assertions over \(L\) and over
\(K\). In particular, the unit and the counit obtained by evaluation
at the identity are filtered.

The matrix-coefficient description and uniqueness of the descended
filtration show that the construction commutes with scalar extension
and endpoint specialization. Finally, after base change through any
embedding
\[
K\hookrightarrow\mathbb C,
\]
the summands
\[
\omega_b
\left(
\mathbb V_\lambda^F
\right)^\vee
\otimes
\mathbb V_\lambda^F
\]
are exactly Hain's weight-zero matrix-coefficient variations.
Therefore the resulting filtration is Hain's canonical Hodge
filtration.
\end{proof}

\begin{proposition}
[{Filtered reductive coefficient algebra;
cf.~\cite[Corollaries~13.6--13.7]{HAIN_HDR_RMC}}]
\label{proposition:filtered_reductive_path_torsor}
The commutative ind-algebra
\[
\mathscr A_{\mathcal R,b,K}^{\dR},
\]
equipped with its canonical connection
\(\nabla_{\mathcal R}\), carries a canonical decreasing
multiplicative Hodge filtration
\[
\cdots
\supseteq
F^p\mathscr A_{\mathcal R,b,K}^{\dR}
\supseteq
F^{p+1}\mathscr A_{\mathcal R,b,K}^{\dR}
\supseteq
\cdots
\]
with the following properties.

\begin{enumerate}
\item
There is a filtered system of finite-rank
\(\nabla_{\mathcal R}\)-stable subbundles
\[
\mathscr V_\alpha
\subseteq
\mathscr A_{\mathcal R,b,K}^{\dR}
\]
such that
\[
\mathscr A_{\mathcal R,b,K}^{\dR}
=
\varinjlim_\alpha\mathscr V_\alpha.
\]
The system may be chosen to contain the unit and to be
multiplicatively cofinal: for every \(\alpha,\beta\), there is a
\(\gamma\) such that multiplication restricts to
\[
\mathscr V_\alpha\otimes
\mathscr V_\beta
\longrightarrow
\mathscr V_\gamma.
\]
For every \(\alpha\), the subsheaves
\[
F^p\mathscr V_\alpha
\coloneqq 
\mathscr V_\alpha
\cap
F^p\mathscr A_{\mathcal R,b,K}^{\dR}
\]
are subbundles and define a finite, exhaustive, and separated
filtration on \(\mathscr V_\alpha\).

\item
The connection satisfies Griffiths transversality:
\[
\nabla_{\mathcal R}
\left(
F^p\mathscr A_{\mathcal R,b,K}^{\dR}
\right)
\subseteq
F^{p-1}\mathscr A_{\mathcal R,b,K}^{\dR}
\otimes_{\mathcal O_{X^\dagger}}
\Omega^1_{X^\dagger/K}.
\]

\item
The unit and multiplication of
\(\mathscr A_{\mathcal R,b,K}^{\dR}\), the reductive torsor
coaction
\[
\Delta_{\mathcal R}:
\mathscr A_{\mathcal R,b,K}^{\dR}
\longrightarrow
\mathscr A_{\mathcal R,b,K}^{\dR}
\otimes_K
\mathcal O(\mathcal R_K^{\dR}),
\]
the identity evaluation
\[
\epsilon_{\mathcal R}:
\mathcal O(\mathcal R_K^{\dR})
\longrightarrow
K,
\]
and the reductive path-composition morphisms are filtered. Here
\(\mathcal O(\mathcal R_K^{\dR})\) carries the filtration obtained
by specializing
\(\mathscr A_{\mathcal R,b,K}^{\dR}\) at \(b\), and \(K\) carries
the trivial filtration.

\item
The construction commutes with extension of the ground field. For
every field extension \(K'/K\), there is a canonical identification
\[
F^p\mathscr A_{\mathcal R,b,K}^{\dR}
\otimes_KK'
\xrightarrow{\sim}
F^p\mathscr A_{\mathcal R,b,K'}^{\dR}.
\]
For every
\[
x\in X^\dagger(K'),
\]
specialization at \(x\) gives a canonical filtered identification
\[
x^*\mathscr A_{\mathcal R,b,K'}^{\dR}
\xrightarrow{\sim}
\mathcal O
\left(
{}^{\vphantom{\dR}}_{x}P_{b,\mathcal R}^{\dR}
\right).
\]
After base change through any embedding
\[
K\hookrightarrow\mathbb C,
\]
the resulting filtration is Hain's canonical weight-zero Hodge
filtration on the reductive path-coordinate algebra.
\end{enumerate}
\end{proposition}

\begin{proof}
Let
\[
F^\bullet
\mathscr A_{\mathcal R,b,K}^{\dR}
\]
be the filtration of
Lemma~\ref{lemma:K_rational_reductive_matrix_coefficients}.
It remains to establish the finite-rank exhaustion and Griffiths
transversality.

Choose \(L/K\) and simple representatives
\[
\left\{
\mathbb V_\lambda
\right\}_{\lambda\in\Lambda}
\]
as in the matrix-coefficient construction. For every finite
\(\operatorname{Gal}(L/K)\)-stable subset
\(\Lambda_0\subseteq\Lambda\) containing the tensor unit, set
\[
\mathscr V_{\Lambda_0,L}
\coloneqq 
\bigoplus_{\lambda\in\Lambda_0}
\omega_{b,L}^{\dR}
\left(
\mathbb V_\lambda
\right)^\vee
\otimes_L
\mathbb V_\lambda.
\]
Galois invariance of the matrix-coefficient filtration descends this
to a filtered finite-rank subbundle
\[
\mathscr V_{\Lambda_0}
\subseteq
\mathscr A_{\mathcal R,b,K}^{\dR}.
\]
These subbundles form a filtered system and exhaust
\(\mathscr A_{\mathcal R,b,K}^{\dR}\).

If \(\Lambda_0,\Lambda_1\subseteq\Lambda\) are two such subsets, let
\(\Lambda_2\) be a finite Galois-stable set containing every simple
constituent of
\[
\mathbb V_\lambda\otimes\mathbb V_\mu,
\qquad
\lambda\in\Lambda_0,\quad
\mu\in\Lambda_1.
\]
Multiplication of matrix coefficients then factors through
\[
\mathscr V_{\Lambda_0,L}
\otimes_L
\mathscr V_{\Lambda_1,L}
\longrightarrow
\mathscr V_{\Lambda_2,L},
\]
and hence, after descent, through the corresponding map over \(K\).
Thus the exhaustion is multiplicatively cofinal.

After base change through an embedding
\(K\hookrightarrow\mathbb C\), each
\(\mathscr V_{\Lambda_0}\) becomes a finite direct sum of Hain's
weight-zero matrix-coefficient variations. Its filtered pieces are
therefore subbundles, and its filtration is finite, exhaustive, and
separated. These properties descend to \(K\) by faithful flatness.

The canonical connection of
Proposition~\ref{proposition:tannakian_path_torsor_connection}
restricts on each matrix-coefficient summand to its tensor-product
connection. After base change to \(\mathbb C\), it is therefore the
connection on Hain's weight-zero variation and satisfies Griffiths
transversality. The asserted inclusion descends to \(K\) by faithful
flatness.

Filteredness of the unit, identity evaluation, multiplication,
reductive torsor coaction, and reductive path composition, together
with compatibility with scalar extension and endpoint
specialization, follows from
Lemma~\ref{lemma:K_rational_reductive_matrix_coefficients}. The same
lemma identifies the complexification of the filtration with Hain's
canonical weight-zero Hodge filtration.
\end{proof}

\subsection{The Relative Bar Model}
\label{appendix:relative_de_rham_bar_model}

Choose a smooth proper compactification
\[
j:
X^\dagger
\hookrightarrow
\overline X
\]
whose boundary
\[
D
\coloneqq 
\overline X\setminus X^\dagger
\]
is a normal-crossings divisor. In the projective case, \(D\) is
empty.

Retain a multiplicatively cofinal filtered system
\[
\left\{
\mathscr V_\alpha
\right\}_\alpha
\]
as in
Proposition~\ref{proposition:filtered_reductive_path_torsor}. For each \(\alpha\), let
\[
\overline{\mathscr V}_\alpha
\]
be the Deligne canonical logarithmic extension of
\(\mathscr V_\alpha\) to \(\overline X\), and set
\[
F^p\overline{\mathscr V}_\alpha
\coloneqq 
\overline{\mathscr V}_\alpha
\cap
j_*F^p\mathscr V_\alpha
\subseteq
j_*\mathscr V_\alpha.
\]
After extension along any embedding \(K\hookrightarrow\mathbb C\),
the \(\mathscr V_\alpha\) are finite direct sums of admissible
polarizable variations of pure weight zero, and the transition and
multiplication maps are morphisms of such variations; see
\cite[Lemma~13.5 and Corollary~13.6]{HAIN_HDR_RMC}.
Their local monodromy is unipotent by the standing semistability
hypothesis and the stability of unipotence under tensor operations
and subquotients.  The canonical-extension theorem for admissible
variations therefore makes each
\(F^p\overline{\mathscr V}_\alpha\) a subbundle of
\(\overline{\mathscr V}_\alpha\); cf.~the discussion preceding
\cite[Theorem~13.4]{HAIN_HDR_RMC}.  Finally, morphisms of polarizable
pure variations are strict and their category is semisimple.
Decomposing each transition or multiplication map into its kernel,
image, and complementary summands, and using functoriality of the
Deligne canonical extension, shows that these maps extend strictly
as morphisms of filtered logarithmic connections.
Since the residues are nilpotent, Deligne canonical extension is
compatible with tensor products. These assertions descend to \(K\)
by faithful flatness.
Consequently,
\[
\overline{\mathscr A}_{\mathcal R,b,K}^{\dR}
\coloneqq 
\varinjlim_\alpha
\overline{\mathscr V}_\alpha
\]
is a filtered commutative algebra in logarithmic ind-connections on
\((\overline X,D)\). Write
\[
\overline\nabla_{\mathcal R}:
\overline{\mathscr A}_{\mathcal R,b,K}^{\dR}
\longrightarrow
\overline{\mathscr A}_{\mathcal R,b,K}^{\dR}
\otimes_{\mathcal O_{\overline X}}
\Omega^1_{\overline X/K}(\log D)
\]
for its logarithmic connection.

Let
\[
\pr_1:
\overline X\times_KX^\dagger
\longrightarrow
\overline X,
\qquad
\pr_2:
\overline X\times_KX^\dagger
\longrightarrow
X^\dagger
\]
be the projections.

\begin{definition}
[{Coefficient logarithmic de Rham dga;
cf.~\cite[\S\S12--13]{HAIN_HDR_RMC}}]
\label{definition:coefficient_logarithmic_de_rham_dga}
For every \(\alpha\) and \(m\geq0\), set
\[
\mathscr K_{\alpha,b,K}^m
\coloneqq 
\Omega^m_{
\overline X\times_KX^\dagger/X^\dagger
}
\bigl(
\log(D\times_KX^\dagger)
\bigr)
\otimes
\pr_1^*
\overline{\mathscr V}_\alpha.
\]
For a local section
\[
\omega\otimes a
\in
\mathscr K_{\alpha,b,K}^m,
\]
define
\[
d_{\mathscr K}(\omega\otimes a)
\coloneqq 
d_{\mathrm{rel}}\omega\otimes a
+
(-1)^m
\omega\wedge
\pr_1^*\overline\nabla_{\mathcal R}(a),
\]
and equip \(\mathscr K_{\alpha,b,K}^\bullet\) with the filtration
\[
F^p\mathscr K_{\alpha,b,K}^m
\coloneqq 
\Omega^m_{
\overline X\times_KX^\dagger/X^\dagger
}
\bigl(
\log(D\times_KX^\dagger)
\bigr)
\otimes
F^{p-m}
\pr_1^*
\overline{\mathscr V}_\alpha.
\]

The transition and multiplication maps on the
\(\overline{\mathscr V}_\alpha\) induce compatible maps
\[
\mathscr K_{\alpha,b,K}^\bullet
\otimes
\mathscr K_{\beta,b,K}^\bullet
\longrightarrow
\mathscr K_{\gamma,b,K}^\bullet,
\qquad
(\omega\otimes a)(\eta\otimes c)
=
(\omega\wedge\eta)\otimes ac.
\]
Set
\[
\bigl(
\mathscr K_{b,K}^\bullet,F^\bullet
\bigr)
\coloneqq 
\varinjlim_\alpha
\bigl(
\mathscr K_{\alpha,b,K}^\bullet,F^\bullet
\bigr).
\]
Equivalently,
\[
\mathscr K_{b,K}^m
=
\Omega^m_{
\overline X\times_KX^\dagger/X^\dagger
}
\bigl(
\log(D\times_KX^\dagger)
\bigr)
\otimes
\pr_1^*
\overline{\mathscr A}_{\mathcal R,b,K}^{\dR}.
\]

Integrability of \(\overline\nabla_{\mathcal R}\) gives
\[
d_{\mathscr K}^2=0,
\]
while its compatibility with multiplication gives the graded Leibniz
rule. Griffiths transversality and multiplicativity of the
coefficient filtration therefore make
\[
\bigl(
\mathscr K_{b,K}^\bullet,F^\bullet
\bigr)
\]
a filtered commutative dga in
\(\pr_2^{-1}\mathcal O_{X^\dagger}\)-modules.
\end{definition}

Fix Navarro--Aznar's filtered Godement--Thom--Whitney model
\[
R_{\mathrm{TW}}\pr_{2,*}.
\]
We apply its filtered-module version at each finite coefficient stage.
The Thom--Whitney simple and the Godement derived direct-image
construction, including their filtered and multiplicative forms, are
developed in
\cite[\S\S2--4 and~6.6--6.15]{NAVARRO_HODGE_DELIGNE}.
The relative logarithmic \v{C}ech model and its functoriality for
sections are treated in
\cite[\S\S6.1,~6.3, and~6.8]
{NAVARRO_GAUSS_MANIN_HOMOTOPY}.
Thus transition maps, multiplication maps, and section-induced
augmentations are carried functorially to the Thom--Whitney models.
All ind-objects below are formed by applying this fixed model at each
finite coefficient stage and then taking the filtered colimit
degreewise.

For every \(\alpha\), set
\[
\bigl(
\mathscr E_{\alpha,b,K}^\bullet,F^\bullet
\bigr)
\coloneqq 
R_{\mathrm{TW}}\pr_{2,*}
\bigl(
\mathscr K_{\alpha,b,K}^\bullet,F^\bullet
\bigr),
\]
and define
\[
\bigl(
\mathscr E_{b,K}^\bullet,F^\bullet
\bigr)
\coloneqq 
\varinjlim_\alpha
\bigl(
\mathscr E_{\alpha,b,K}^\bullet,F^\bullet
\bigr).
\]
The multiplicative structure of the filtered Thom--Whitney functor
and the multiplicative cofinality of
\(\{\mathscr V_\alpha\}_\alpha\) make
\[
\bigl(
\mathscr E_{b,K}^\bullet,F^\bullet
\bigr)
\]
a filtered commutative differential graded
\(\mathcal O_{X^\dagger}\)-algebra. Its underlying complex is
naturally quasi-isomorphic to
\[
R\pr_{2,*}\mathscr K_{b,K}^\bullet.
\]

\begin{proposition}[Endpoint augmentations and scalar extension]
\label{proposition:filtered_relative_de_rham_dga}
Let
\[
s_b:
X^\dagger
\longrightarrow
\overline X\times_KX^\dagger,
\qquad
x\longmapsto(b,x),
\]
and let
\[
\Delta:
X^\dagger
\longrightarrow
\overline X\times_KX^\dagger,
\qquad
x\longmapsto(x,x),
\]
be the diagonal. These sections induce augmentations of filtered
commutative dgas
\[
\epsilon_b:
\bigl(
\mathscr E_{b,K}^\bullet,F^\bullet
\bigr)
\longrightarrow
\bigl(
\mathcal O_{X^\dagger},F^\bullet_{\mathrm{triv}}
\bigr)
\]
and
\[
\epsilon_-:
\bigl(
\mathscr E_{b,K}^\bullet,F^\bullet
\bigr)
\longrightarrow
\bigl(
\mathscr A_{\mathcal R,b,K}^{\dR},F^\bullet
\bigr),
\]
where the targets are regarded as dgas concentrated in degree zero
and
\[
F^p_{\mathrm{triv}}\mathcal O_{X^\dagger}
=
\begin{cases}
\mathcal O_{X^\dagger}, & p\leq0,\\
0, & p>0.
\end{cases}
\]

The first augmentation is induced by evaluation along \(s_b\),
followed by the identity evaluation
\[
\epsilon_{\mathcal R}:
\mathcal O(\mathcal R_K^{\dR})
\longrightarrow
K.
\]
The second is induced by evaluation along \(\Delta\).

For every field extension \(K'/K\), there is a canonical filtered
multiplicative quasi-isomorphism
\[
\bigl(
\mathscr E_{b,K}^\bullet,F^\bullet
\bigr)
\otimes_KK'
\xrightarrow{\simeq}
\bigl(
\mathscr E_{b,K'}^\bullet,F^\bullet
\bigr),
\]
compatible with both augmentations after scalar extension.
\end{proposition}

\begin{proof}
At each coefficient stage, a section of \(\pr_2\) gives the usual
augmentation of the relative de Rham complex: it evaluates functions
and coefficients in degree zero and vanishes in positive
relative-form degrees.

Along \(s_b\), the coefficient algebra specializes to
\[
\mathcal O
\left(
{}^{\vphantom{\dR}}_{b}P_{b,\mathcal R}^{\dR}
\right)
=
\mathcal O(\mathcal R_K^{\dR}),
\]
and composition with \(\epsilon_{\mathcal R}\) gives
\(\epsilon_b\). Along \(\Delta\), it specializes to
\(\mathscr A_{\mathcal R,b,K}^{\dR}\), giving \(\epsilon_-\).
These maps are filtered by
Proposition~\ref{proposition:filtered_reductive_path_torsor}.
Functoriality of the filtered Thom--Whitney direct image and passage
to the filtered colimit give the asserted augmentations.

For every \(\alpha\), formation of the filtered logarithmic de Rham
complex of \(\overline{\mathscr V}_\alpha\) commutes with scalar
extension. Comparison of the fixed Thom--Whitney model with the
ordinary derived direct image, together with flat proper base change
for \(\pr_2\), therefore gives a natural filtered multiplicative
quasi-isomorphism
\[
\bigl(
\mathscr E_{\alpha,b,K}^\bullet,F^\bullet
\bigr)
\otimes_KK'
\xrightarrow{\simeq}
\bigl(
\mathscr E_{\alpha,b,K'}^\bullet,F^\bullet
\bigr).
\]
Functoriality with respect to the sections \(s_b\) and \(\Delta\)
gives compatibility with both augmentations. Since filtered colimits
of quasi-coherent modules are exact, passage to the stagewise
filtered colimit proves the assertion.
\end{proof}

\begin{definition}
[{Relative two-sided bar complex;
cf.~\cite[\S7]{HAIN_HDR_RMC}}]
\label{definition:relative_de_rham_bar_complex}
We use Hain's reduced two-sided bar construction
\[
B(M^\bullet,A^\bullet,N^\bullet),
\]
where \(A^\bullet\) is a commutative dga,
\(M^\bullet\) is a right dg \(A^\bullet\)-module, and
\(N^\bullet\) is a left dg \(A^\bullet\)-module. We use the
differential, balancing relations, and shuffle product of
\cite[\S7]{HAIN_HDR_RMC}. Bar elements are denoted
\[
m[a_1|\cdots|a_r]n.
\]

Regard \(\mathcal O_{X^\dagger}\) as a right
\(\mathscr E_{b,K}^\bullet\)-module through \(\epsilon_b\), and
regard \(\mathscr A_{\mathcal R,b,K}^{\dR}\) as a left
\(\mathscr E_{b,K}^\bullet\)-module through \(\epsilon_-\). Define
the \emph{relative two-sided bar complex} by
\[
\mathscr B_{b,K}^\bullet
\coloneqq 
B
\left(
\mathcal O_{X^\dagger},
\mathscr E_{b,K}^\bullet,
\mathscr A_{\mathcal R,b,K}^{\dR}
\right),
\]
where all tensor products are taken over
\(\mathcal O_{X^\dagger}\). The shuffle product makes
\(\mathscr B_{b,K}^\bullet\) a commutative dga in
ind-quasi-coherent \(\mathcal O_{X^\dagger}\)-modules.

Equip each bar-length summand with the tensor-product filtration
induced by the trivial filtration on \(\mathcal O_{X^\dagger}\),
the filtration on \(\mathscr E_{b,K}^\bullet\), and the Hodge
filtration on \(\mathscr A_{\mathcal R,b,K}^{\dR}\). Explicitly,
\(F^p\mathscr B_{b,K}^\bullet\) is generated by bar elements
\[
m[e_1|\cdots|e_r]\phi
\]
for which
\[
m\in F_{\mathrm{triv}}^{p_0}\mathcal O_{X^\dagger},
\qquad
e_i\in F^{p_i}\mathscr E_{b,K}^\bullet,
\qquad
\phi\in
F^{p_{r+1}}\mathscr A_{\mathcal R,b,K}^{\dR},
\]
and
\[
p_0+p_1+\cdots+p_{r+1}\geq p.
\]
This is the \emph{Hodge convolution filtration}. Since the
augmentations and module structures are filtered, the bar
differential and shuffle product preserve it.
\end{definition}

Let
\[
\mathscr B_{b,b,K}^\bullet
\coloneqq 
b^*\mathscr B_{b,K}^\bullet
\]
denote the specialization at the endpoint \(b\).

\begin{proposition}
[{Composition coproduct;
cf.~\cite[\S12, especially~(14)]{HAIN_HDR_RMC}}]
\label{proposition:relative_bar_composition_coproduct}
Composition of reductive paths induces a canonical morphism of
filtered commutative dgas
\[
\Delta_b^{\mathrm{bar}}:
\mathscr B_{b,K}^\bullet
\longrightarrow
\mathscr B_{b,K}^\bullet
\otimes_K
\mathscr B_{b,b,K}^\bullet.
\]
It is the chain-level coproduct corresponding to the right action of
loops based at \(b\) on paths from \(b\) to the variable endpoint.

Consequently, setting
\[
\mathscr H_{b,K}
\coloneqq 
\mathcal H^0(\mathscr B_{b,K}^\bullet),
\qquad
\mathscr H_{b,b,K}
\coloneqq 
\mathcal H^0(\mathscr B_{b,b,K}^\bullet),
\]
there is an induced morphism of filtered algebras
\[
\Delta_b^{\mathrm{bar}}:
\mathscr H_{b,K}
\longrightarrow
\mathscr H_{b,K}
\otimes_K
\mathscr H_{b,b,K}.
\]
\end{proposition}

\begin{proof}
The right action
\[
\mathcal P_{\mathcal R,b,-,K}^{\dR}
\times_K
\mathcal R_K^{\dR}
\longrightarrow
\mathcal P_{\mathcal R,b,-,K}^{\dR}
\]
induces compatible coactions on the reductive coefficient algebra
and on
\(\mathscr E_{b,K}^\bullet\).

Combining these coactions with deconcatenation of bar words gives
Hain's three-endpoint coproduct
\cite[(14)]{HAIN_HDR_RMC}, specialized to the case in which the
intermediate endpoint is \(b\). After applying the symmetry
isomorphism to match our right-torsor convention, this gives
\(\Delta_b^{\mathrm{bar}}\).

Hain's formula shows that the coproduct commutes with the internal
and bar differentials and with the shuffle product. Each constituent
map is filtered, so it also preserves the Hodge convolution
filtration. Passing to zeroth cohomology gives the final assertion.
\end{proof}

\paragraph{The algebraic comparison morphism.}

Let
\[
q_{\mathcal Q}:
\mathcal Q_{b,-,K}
\longrightarrow
X^\dagger
\]
be the affine \(X^\dagger\)-scheme, generally not of finite type,
with coordinate algebra
\[
(q_{\mathcal Q})_*
\mathcal O_{\mathcal Q_{b,-,K}}
=
\mathscr H_{b,K}.
\]

Let
\[
V\in\mathcal C_K^{\dR}.
\]
Choose a finite filtration of \(V\) whose graded pieces belong to
\(\mathcal D_K^{\dR}\), and split it locally on \(X^\dagger\).
Relative to the universal reductive transport
\(\tau_V^{\mathcal R}\), the extension part of the connection is then
represented by a nilpotent Maurer--Cartan matrix
\[
\Omega_V
\in
\mathscr E_{b,K}^1\otimes_K\mathfrak n_V,
\qquad
d\Omega_V+\Omega_V\wedge\Omega_V=0,
\]
where \(\mathfrak n_V\) is the algebra of strictly
filtration-lowering endomorphisms of
\(\omega_{b,K}^{\dR}(V)\). Set
\[
T_V
\coloneqq 
\left(
\sum_{r\geq0}
[
\underbrace{
\Omega_V|\cdots|\Omega_V
}_{r}
]
\right)
\tau_V^{\mathcal R}.
\]
The sum is finite because \(\mathfrak n_V\) is nilpotent.

Locally, \(T_V\) is a degree-zero section of the matrix-valued bar
complex
\[
\mathscr B_{b,K}^\bullet
\otimes_{\mathcal O_{X^\dagger}}
\mathscr Hom_{\mathcal O_{X^\dagger}}
\left(
\omega_{b,K}^{\dR}(V)\otimes_K\mathcal O_{X^\dagger},
\omega_{X^\dagger,K}^{\dR}(V)
\right),
\]
where the second factor is placed in degree zero. Thus the assertion
that \(T_V\) is a cocycle means that each of its matrix coefficients
is closed for the full two-sided bar differential.

\begin{lemma}
[Relative bar holonomy;
cf.~{\cite[\S\S8--10 and~12]{HAIN_HDR_RMC}}]
\label{lemma:relative_bar_holonomy}
The bar series \(T_V\) is a degree-zero cocycle and defines a
canonical algebraic isomorphism
\[
\mathbf T_V:
\omega_{b,K}^{\dR}(V)
\otimes_K
\mathcal O_{\mathcal Q_{b,-,K}}
\xrightarrow{\sim}
q_{\mathcal Q}^*
\omega_{X^\dagger,K}^{\dR}(V).
\]
It is independent of the chosen filtration and local splitting,
functorial in \(V\), compatible with tensor products, duals, and the
tensor unit, and commutes with extension of the ground field.
Under deconcatenation, the isomorphisms \(\mathbf T_V\) are compatible
with composition of paths.
\end{lemma}

\begin{proof}
Each occurrence of \(\Omega_V\) has cohomological degree one, while
bar suspension lowers degree by one. Consequently every term
\[
[\Omega_V|\cdots|\Omega_V]\tau_V^{\mathcal R}
\]
has total degree zero.

We now apply the full two-sided bar differential, consisting of the
internal differential, the contractions of adjacent bar entries, and
the two endpoint terms. The endpoint terms vanish because both
augmentations have targets concentrated in degree zero and hence
annihilate \(\mathscr E_{b,K}^1\). The internal terms insert
\(d\Omega_V\), while the adjacent contractions insert
\(\Omega_V\wedge\Omega_V\), with the signs of the reduced bar
differential. Therefore
\[
d_{\mathrm{bar}}T_V
=
\sum
[\cdots|
d\Omega_V+\Omega_V\wedge\Omega_V
|\cdots]\tau_V^{\mathcal R}
=
0
\]
by the Maurer--Cartan equation. The reductive equivariance of
\(\Omega_V\), inherited from its construction relative to the
universal reductive transport, ensures that this series respects the
balancing relations in the two-sided bar construction. This is the
matrix-valued form of
\cite[Proposition~9.1]{HAIN_HDR_RMC}.

It follows that \(T_V\) determines a local section
\[
[T_V]\in
\mathscr H_{b,K}
\otimes_{\mathcal O_{X^\dagger}}
\mathscr Hom_{\mathcal O_{X^\dagger}}
\left(
\omega_{b,K}^{\dR}(V)\otimes_K\mathcal O_{X^\dagger},
\omega_{X^\dagger,K}^{\dR}(V)
\right).
\]
If two local splittings differ by a unipotent gauge transformation
\(g\), then
\[
\Omega_V'
=
g^{-1}\Omega_Vg+g^{-1}dg,
\qquad
T_V'
=
g(-)^{-1}T_Vg(b).
\]
These are precisely the changes of frames in the target and source.
Hence the local sections glue.

Since
\[
(q_{\mathcal Q})_*
\mathcal O_{\mathcal Q_{b,-,K}}
=
\mathscr H_{b,K},
\]
affine adjunction and the projection formula identify the resulting
global section with an algebraic morphism
\[
\mathbf T_V:
\omega_{b,K}^{\dR}(V)\otimes_K
\mathcal O_{\mathcal Q_{b,-,K}}
\longrightarrow
q_{\mathcal Q}^*
\omega_{X^\dagger,K}^{\dR}(V).
\]
Its bar-length-zero term is the invertible reductive transport
\(\tau_V^{\mathcal R}\), while the remaining factor lies in
\(1+\mathfrak n_V\). Since \(\mathfrak n_V\) is nilpotent, this
unipotent factor has a finite inverse. Thus \(\mathbf T_V\) is an
isomorphism.

The preceding gauge-transformation formula proves independence of the
local splitting. Moreover, the resulting bar class induces the full
monodromy representation of \(V\). By~\cite[Corollaries~9.2--9.3 and Theorem~10.1]{HAIN_HDR_RMC},
this representation factors uniquely through the relative
completion, so the class is also independent of the chosen extension
filtration. After extension along \(K\hookrightarrow\mathbb C\) this
is Hain's uniqueness statement, and the conclusion descends to \(K\)
by faithful scalar extension.

Naturality of the connection matrices gives functoriality in \(V\),
while the shuffle identity gives compatibility with tensor products.
Compatibility with the tensor unit and duals then follows from the
unit, evaluation, and coevaluation morphisms. Finally, Hain's
three-endpoint coproduct
\cite[\S12, especially Theorem~12.2 and~(14)]{HAIN_HDR_RMC}
identifies the transport along a composite path with the composite of
the two transports. All constructions are algebraic over \(K\) and
commute with extension of the ground field.
\end{proof}
\begin{corollary}[Algebraic bar transport]
\label{proposition:algebraic_bar_transport}
The isomorphisms of
Lemma~\ref{lemma:relative_bar_holonomy} assemble into a canonical
tensor isomorphism
\[
\mathbf T_{b,-}:
\omega_{b,K}^{\dR}
\otimes_K
\mathcal O_{\mathcal Q_{b,-,K}}
\xrightarrow{\sim}
q_{\mathcal Q}^*
\omega_{X^\dagger,K}^{\dR}
\]
on \(\mathcal C_K^{\dR}\). Equivalently, they define a canonical
morphism
\[
\eta_{b,-}:
\mathcal Q_{b,-,K}
\longrightarrow
\mathcal P_{b,-,K}^{\dR},
\]
or, dually, a morphism of \(\mathcal O_{X^\dagger}\)-algebras
\[
\vartheta_{b,-}:
\mathscr A_{b,K}^{\dR}
\longrightarrow
\mathscr H_{b,K}.
\]
These morphisms are natural in the variable endpoint and commute with
extension of the ground field.

Let
\[
\vartheta_{b,b}:
\mathcal O(\mathcal G_K^{\dR})
\longrightarrow
\mathscr H_{b,b,K}
\]
denote the specialization at \(b\). Then the diagram
\[
\begin{tikzcd}[column sep=large]
\mathscr A_{b,K}^{\dR}
\arrow[r, "\vartheta_{b,-}"]
\arrow[d, "\Delta_{\mathcal P}"']
&
\mathscr H_{b,K}
\arrow[d, "\Delta_b^{\mathrm{bar}}"]
\\
\mathscr A_{b,K}^{\dR}
\otimes_K
\mathcal O(\mathcal G_K^{\dR})
\arrow[r, "\vartheta_{b,-}\otimes\vartheta_{b,b}"']
&
\mathscr H_{b,K}
\otimes_K
\mathscr H_{b,b,K}
\end{tikzcd}
\]
commutes.
\end{corollary}

\begin{proof}
The first assertion is
Lemma~\ref{lemma:relative_bar_holonomy}. Since
\[
\mathcal P_{b,-,K}^{\dR}
=
\operatorname{Isom}^{\otimes}
\left(
\omega_{b,K}^{\dR}\otimes_K\mathcal O_{X^\dagger},
\omega_{X^\dagger,K}^{\dR}
\right),
\]
its representing property gives \(\eta_{b,-}\), and affineness gives
the dual morphism \(\vartheta_{b,-}\).

Compatibility of \(\mathbf T_{b,-}\) with path composition says that
\(\eta_{b,-}\) is equivariant with respect to its equal-endpoint
specialization \(\eta_{b,b}\). Dualizing this equivariance gives the
displayed coaction diagram. Naturality and compatibility with scalar
extension follow from the same lemma.
\end{proof}

\begin{lemma}[Complex comparison with Hain's filtered bar model]
\label{lemma:comparison_with_hain_bar_complex}
Fix an embedding
\[
\iota_\infty:K\hookrightarrow\mathbb C,
\qquad
Y_\infty
\coloneqq 
\left(
X^\dagger\otimes_{K,\iota_\infty}\mathbb C
\right)^{\an},
\]
here \((-)^{\an}\) denotes complex analytification with respect to
\(\iota_\infty\), rather than rigid analytification. Let \(\mathcal O_{Y_\infty}\) and
\(\mathcal C_{Y_\infty}^\infty\) denote, respectively, the sheaves of
holomorphic and smooth complex-valued functions on \(Y_\infty\).
Let
\[
\left(
\mathscr B_{\mathrm H,b,-}^\bullet,
F_{\mathrm H}^\bullet
\right)
\]
denote Hain's filtered two-sided bar family formed from the complex
part of his multiplicative mixed-Hodge complex and the two endpoint
augmentations.

There is a natural isomorphism of filtered
\(\mathcal C_{Y_\infty}^\infty\)-algebras
\[
\begin{aligned}
&
\left(
\mathcal H^0(\mathscr B_{b,K}^\bullet)
\otimes_{K,\iota_\infty}\mathbb C
\right)^{\an}
\otimes_{\mathcal O_{Y_\infty}}
\mathcal C_{Y_\infty}^\infty
\\
&\hspace{3cm}\xrightarrow{\sim}
\mathcal H^0
\left(
\mathscr B_{\mathrm H,b,-}^\bullet
\right).
\end{aligned}
\]
It is horizontal and compatible with endpoint specialization and
the path-groupoid structure maps. Under this comparison, the
complexification of \(\vartheta_{b,-}\) is Hain's bar-transport
morphism. In particular, the filtration on the left identifies with
Hain's Hodge filtration on the relative Malcev path-coordinate
family.
\end{lemma}

\begin{proof}
For every finite diagram of coefficient objects, analytification and
the fine smooth logarithmic resolution give a multiplicative
filtered comparison between the Thom--Whitney model and the complex
part of Hain's mixed-Hodge complex. This comparison is horizontal,
functorial for the endpoint augmentations and the three-endpoint
maps, and a quasi-isomorphism on every associated graded; see
\cite[
Theorem~13.4,
Corollary~13.6,
Proposition~13.9,
and the proof of Theorem~13.3
]{HAIN_HDR_RMC}.

The middle dgas are concentrated in nonnegative cohomological
degrees, and the endpoint modules are concentrated in degree zero.
Moreover, the middle dgas are connected. Indeed,
\[
\begin{aligned}
&H^0\!\left(
E_{\mathrm{fin}}^\bullet
\left(
X^\dagger(\mathbb C),
\mathcal O(P_{\mathcal R,b,-})
\right)
\right)
\\
&\qquad =
\mathcal O(\mathcal R_{\mathbb C}^{\dR})^{
\rho_{\mathcal R}(\pi_1(X^\dagger(\mathbb C),b))}
=
\mathcal O(\mathcal R_{\mathbb C}^{\dR})^{
\mathcal R_{\mathbb C}^{\dR}}
=
\mathbb C,
\end{aligned}
\]
where the second equality follows from Zariski density. After
extension to \(\mathcal C_{Y_\infty}^\infty\), and hence for every
middle dga \(A^\bullet\) in the comparison, this gives
\[
\mathcal H^0(A^\bullet)
\simeq
\mathcal C_{Y_\infty}^\infty.
\]

At each finite coefficient diagram, strictness of the Hodge
filtration gives
\[
\mathcal H^n\!\left(
\Gr_F^\bullet A^\bullet
\right)
\xrightarrow{\sim}
\Gr_F^\bullet\mathcal H^n(A^\bullet).
\]
Since \(\mathcal C_{Y_\infty}^\infty\) has the trivial filtration,
the full associated-graded middle dga is also connected:
\[
\mathcal H^0\!\left(
\Gr_F^\bullet A^\bullet
\right)
\simeq
\Gr_F^\bullet\mathcal C_{Y_\infty}^\infty
\simeq
\mathcal C_{Y_\infty}^\infty.
\]

Hain's bar-invariance theorem
\cite[Proposition~7.1]{HAIN_HDR_RMC},
together with flat scalar extension to
\(\mathcal C_{Y_\infty}^\infty\), therefore applies both to the
underlying comparison and to its full associated-graded comparison.
For the Hodge convolution filtration,
\[
\Gr_F^\bullet B(M^\bullet,A^\bullet,N^\bullet)
\simeq
B\left(
\Gr_F^\bullet M^\bullet,
\Gr_F^\bullet A^\bullet,
\Gr_F^\bullet N^\bullet
\right).
\]
Thus the induced comparison of bar complexes is a quasi-isomorphism
both on the underlying complexes and on every associated graded.
Exactness of filtered colimits and their compatibility with tensor
products permit passage through the multiplicatively cofinal
coefficient system, yielding the displayed filtered isomorphism on
zeroth cohomology.

After forgetting the filtrations, the comparison following
\cite[Corollary~13.8]{HAIN_HDR_RMC}
identifies this construction with Hain's
\(E_{\mathrm{fin}}^\bullet\)-bar model. Finally,
\cite[Theorem~12.2]{HAIN_HDR_RMC}
identifies the induced morphism with Hain's bar transport and proves
its compatibility with path composition. Functoriality for the
endpoint augmentations and the three-endpoint construction gives the
remaining path-groupoid compatibilities.
\end{proof}

\begin{corollary}[Bar model for the universal relative path torsor]
\label{proposition:relative_de_rham_bar_model}
The morphism
\[
\vartheta_{b,-}:
\mathscr A_{b,K}^{\dR}
\longrightarrow
\mathcal H^0(\mathscr B_{b,K}^\bullet)
\]
of
Corollary~\ref{proposition:algebraic_bar_transport}
is an isomorphism. Its inverse defines a canonical isomorphism
\[
\theta_{b,-}:
\mathcal H^0(\mathscr B_{b,K}^\bullet)
\xrightarrow{\sim}
\mathscr A_{b,K}^{\dR}.
\]
Its specialization at \(b\) identifies
\[
\mathscr H_{b,b,K}
\xrightarrow{\sim}
\mathcal O(\mathcal G_K^{\dR}),
\]
and, under these identifications, the bar coproduct is the universal
right-torsor coaction.
\end{corollary}

\begin{proof}
Choose an embedding \(K\hookrightarrow\mathbb C\). By
Lemma~\ref{lemma:comparison_with_hain_bar_complex}, after
analytification and extension to
\(\mathcal C_{Y_\infty}^\infty\), the morphism
\(\vartheta_{b,-}\) identifies with Hain's bar-transport morphism.
The latter is an isomorphism by
\cite[Theorem~12.2]{HAIN_HDR_RMC}.

Faithful flatness of
\[
\mathcal O_{Y_\infty}
\longrightarrow
\mathcal C_{Y_\infty}^\infty,
\]
of algebraic-to-analytic local rings, and of
\(K\hookrightarrow\mathbb C\) therefore shows that
\(\vartheta_{b,-}\) is an isomorphism over \(K\). The same argument
after specialization at \(b\) gives the based-loop identification.
Compatibility with the coproduct follows from
Corollary~\ref{proposition:algebraic_bar_transport}.
\end{proof}

\subsection{The Universal Hodge Filtration}
\label{appendix:universal_relative_hodge_filtration}

Equip
\[
\mathcal H^0
\left(
\mathscr B_{b,K}^\bullet
\right)
\]
with the filtration
\[
F^p
\mathcal H^0
\left(
\mathscr B_{b,K}^\bullet
\right)
\coloneqq 
\operatorname{im}
\left(
\mathcal H^0
\left(
F^p\mathscr B_{b,K}^\bullet
\right)
\longrightarrow
\mathcal H^0
\left(
\mathscr B_{b,K}^\bullet
\right)
\right).
\]
Transporting this filtration through
Corollary~\ref{proposition:relative_de_rham_bar_model} gives a
decreasing filtration
\[
\cdots
\supseteq
F^p\mathscr A_{b,K}^{\dR}
\supseteq
F^{p+1}\mathscr A_{b,K}^{\dR}
\supseteq
\cdots
\]
on the coordinate algebra of the universal relative path torsor.

\begin{lemma}[Strict endpoint specialization]
\label{lemma:strict_endpoint_specialization}
For every \(p\), the quotient
\[
\mathscr A_{b,K}^{\dR}/F^p\mathscr A_{b,K}^{\dR}
\]
is flat over \(X^\dagger\). Consequently, for every field extension
\(K'/K\) and every \(x\in X^\dagger(K')\), the canonical
identification
\[
x^*\mathscr A_{b,K'}^{\dR}
\xrightarrow{\sim}
\mathcal O\left({}^{\vphantom{\dR}}_{x}P_b^{\dR}\right)
\]
is strict:
\[
x^*F^p\mathscr A_{b,K'}^{\dR}
\xrightarrow{\sim}
F^p\mathcal O\left({}^{\vphantom{\dR}}_{x}P_b^{\dR}\right).
\]
\end{lemma}

\begin{proof}
Fix an embedding \(K\hookrightarrow\mathbb C\). By
Lemma~\ref{lemma:comparison_with_hain_bar_complex}, after
analytification and extension to
\(\mathcal C_{Y_\infty}^\infty\), the filtration on
\(\mathscr A_{b,K}^{\dR}\) becomes Hain's Hodge filtration on the
relative path-coordinate family. Its filtered pieces are holomorphic
subbundles by \cite[Theorem~13.10]{HAIN_HDR_RMC}. Hence the displayed
quotient becomes flat after these extensions.

Faithful flatness of algebraic-to-analytic local rings and of
\[
\mathcal O_{Y_\infty}
\longrightarrow
\mathcal C_{Y_\infty}^\infty,
\]
followed by faithful flatness of \(K\hookrightarrow\mathbb C\),
shows that the quotient is flat over \(X^\dagger\). Endpoint pullback
therefore preserves the filtered inclusion. The underlying
identification follows by base change of the universal affine path
torsor, while the filtered comparison identifies the resulting
filtration with Hain's fibrewise filtration. Faithful-flat descent
and scalar extension give the assertion over \(K'\).
\end{proof}

By Lemma~\ref{lemma:strict_endpoint_specialization}, specialization
at the endpoint \(b\) defines
\[
F^p\mathcal O(\mathcal G_K^{\dR})
\coloneqq 
b^*F^p\mathscr A_{b,K}^{\dR}.
\]
Under \(\vartheta_{b,b}\), this agrees with the filtration induced by
the based-loop bar complex.

\begin{proposition}[Hodge algebra and torsor compatibility]
\label{proposition:universal_hodge_algebra_structure}
The filtration on \(\mathscr A_{b,K}^{\dR}\) is multiplicative:
\[
F^p\mathscr A_{b,K}^{\dR}
\cdot
F^q\mathscr A_{b,K}^{\dR}
\subseteq
F^{p+q}\mathscr A_{b,K}^{\dR}.
\]
Moreover, the torsor coaction
\[
\mathscr A_{b,K}^{\dR}
\longrightarrow
\mathscr A_{b,K}^{\dR}
\otimes_K
\mathcal O(\mathcal G_K^{\dR})
\]
is a morphism of filtered algebras. The filtration on
\(\mathcal O(\mathcal G_K^{\dR})\) defined above is a Hopf-algebra
filtration: its coproduct, counit, and antipode are filtered.
\end{proposition}

\begin{proof}
The internal and bar differentials and the shuffle product preserve
the Hodge convolution filtration. Hence the induced filtration on
zeroth cohomology is multiplicative. Likewise, the composition
coproduct of
Proposition~\ref{proposition:relative_bar_composition_coproduct}
is filtered. Under
Corollary~\ref{proposition:relative_de_rham_bar_model}, it
identifies with the universal path-torsor coaction.

After specialization at \(b\), the composition coproduct is the Hopf
coproduct on \(\mathcal O(\mathcal G_K^{\dR})\), while identity
evaluation gives the filtered counit. After extension through an
embedding \(K\hookrightarrow\mathbb C\), the filtered comparison
identifies this filtration with Hain's Hodge filtration. By
\cite[Theorem~13.1]{HAIN_HDR_RMC}, the antipode is a morphism of
mixed Hodge structures, and hence is filtered. Faithful-flat descent
therefore shows that the antipode is filtered over \(K\).
\end{proof}

\begin{theorem}
[{Universal Hodge filtration;
cf.~\cite[
Theorem~12.2,
Theorem~13.3,
Proposition~13.9,
and Theorem~13.10
]{HAIN_HDR_RMC}}]
\label{theorem:relative_filtered_path_groupoid}
The canonical connection satisfies Griffiths transversality:
\[
\nabla
\left(
F^p\mathscr A_{b,K}^{\dR}
\right)
\subseteq
F^{p-1}\mathscr A_{b,K}^{\dR}
\otimes_{\mathcal O_{X^\dagger}}
\Omega^1_{X^\dagger/K}.
\]

The filtered algebra commutes with extension of the ground field.
For every field extension \(K'/K\), there is a canonical
identification
\[
F^p\mathscr A_{b,K}^{\dR}
\otimes_KK'
\xrightarrow{\sim}
F^p\mathscr A_{b,K'}^{\dR}.
\]
For every
\[
x\in X^\dagger(K'),
\]
specialization at \(x\) gives a canonical filtered identification
\[
x^*
\mathscr A_{b,K'}^{\dR}
\xrightarrow{\sim}
\mathcal O
\left(
{}^{\vphantom{\dR}}_{x}P_b^{\dR}
\right),
\]
and hence
\[
x^*
F^p\mathscr A_{b,K'}^{\dR}
\xrightarrow{\sim}
F^p
\mathcal O
\left(
{}^{\vphantom{\dR}}_{x}P_b^{\dR}
\right).
\]

After base change through every embedding
\[
K\hookrightarrow\mathbb C,
\]
the resulting filtration is Hain's Hodge filtration on the
coordinate algebra of the relative Malcev path torsor. After base change to \(K_v\), it is the canonical family-level Hodge
filtration on the universal relative de Rham path torsor.
\end{theorem}

\begin{proof}
Fix an embedding
\[
\iota_\infty:
K\hookrightarrow\mathbb C,
\]
and retain the notation \(Y_\infty\) of
Lemma~\ref{lemma:comparison_with_hain_bar_complex}.
After analytification and extension to
\(\mathcal C_{Y_\infty}^\infty\), that lemma identifies the
filtration on \(\mathscr A_{b,K}^{\dR}\) with Hain's Hodge
filtration and the canonical connection with the flat connection on
his relative path local system. By
\cite[Theorem~13.10]{HAIN_HDR_RMC}, this filtration satisfies
Griffiths transversality.

Equivalently, the composite
\[
F^p\mathscr A_{b,K}^{\dR}
\xrightarrow{\nabla}
\mathscr A_{b,K}^{\dR}
\otimes_{\mathcal O_{X^\dagger}}
\Omega^1_{X^\dagger/K}
\longrightarrow
\left(
\mathscr A_{b,K}^{\dR}/
F^{p-1}\mathscr A_{b,K}^{\dR}
\right)
\otimes_{\mathcal O_{X^\dagger}}
\Omega^1_{X^\dagger/K}
\]
therefore vanishes after these faithfully flat extensions. Hence it
vanishes over \(K\), proving the asserted inclusion.

Proposition~\ref{proposition:filtered_relative_de_rham_dga},
together with base-change compatibility of the endpoint modules and
functoriality of the reduced bar construction, gives a canonical
filtered scalar-extension quasi-isomorphism
\[
\left(
\mathscr B_{b,K}^\bullet,
F^\bullet
\right)
\otimes_KK'
\xrightarrow{\sim}
\left(
\mathscr B_{b,K'}^\bullet,
F^\bullet
\right).
\]
Since \(K'/K\) is flat, scalar extension commutes with zeroth
cohomology and with images. Hence
\[
\begin{aligned}
&
\operatorname{im}
\left(
\mathcal H^0(F^p\mathscr B_{b,K}^\bullet)
\longrightarrow
\mathcal H^0(\mathscr B_{b,K}^\bullet)
\right)
\otimes_KK'
\\
&\hspace{2cm}\xrightarrow{\sim}
\operatorname{im}
\left(
\mathcal H^0(F^p\mathscr B_{b,K'}^\bullet)
\longrightarrow
\mathcal H^0(\mathscr B_{b,K'}^\bullet)
\right).
\end{aligned}
\]
Transport through
Corollary~\ref{proposition:relative_de_rham_bar_model} gives
\[
F^p\mathscr A_{b,K}^{\dR}\otimes_KK'
\xrightarrow{\sim}
F^p\mathscr A_{b,K'}^{\dR}.
\]

The filtered endpoint identifications follow from
Lemma~\ref{lemma:strict_endpoint_specialization}.

The comparison with Hain is
Lemma~\ref{lemma:comparison_with_hain_bar_complex}. Taking
\(K'=K_v\) gives the asserted family-level \(K_v\)-filtration.
\end{proof}

\begin{lemma}
[{Compatibility with Olsson's filtered realization;
cf.~\cite[
\S3,
Corollary~4.6,
diagram~(8.12.2),
and Corollary~8.13
]{OLSSON_BAR_AFFINE_STACKS};
\cite[
Proposition~7.7,
Corollaries~7.15--7.16,
and \S\S8.29--8.31
]{OLSSON_NA_P_HODGE}}]
\label{lemma:olsson_bar_filtration_compatibility}
Set
\[
\mathscr A_{b,K_v}^{\dR}
\coloneqq 
\mathscr A_{b,K}^{\dR}\otimes_KK_v,
\]
and let
\[
F^\bullet_{\mathrm{Ols}}
\mathscr A_{b,K_v}^{\dR}
\]
be the filtration obtained from the relative filtered association of
Theorem~\ref{theorem:olsson_etale_crystalline_comparison}, transported
through the canonical identification of the underlying
ind-connections. Then
\[
F^p_{\mathrm{Ols}}
\mathscr A_{b,K_v}^{\dR}
=
F^p\mathscr A_{b,K}^{\dR}\otimes_KK_v
\]
for every \(p\).

This identification is compatible with multiplication, path
composition, the right torsor coaction, and specialization at \(b\).
In particular, the analogous assertion holds for
\[
\mathcal O(\mathcal G_K^{\dR})\otimes_KK_v.
\]
\end{lemma}

\begin{proof}
Work over a small chart \(U\), and write Olsson's regular
coefficient ind-objects as filtered colimits of associated
finite-rank subobjects. Fix a simplicial degree and a finite
collection of multiplication, unit, action, and endpoint maps.
Choose finite stages containing all source pieces, and then enlarge
the target stages so that every chosen map is represented.
Thus the algebra structure is handled by compatible maps between
finite stages; no individual finite stage is assumed to be closed
under multiplication.

In Olsson's~\cite[
Proposition~7.7,
Remark~7.8,
and Corollary~7.14
]{OLSSON_NA_P_HODGE},
the regular \'{e}tale coefficient ind-object is shown to be
associated, stagewise, with its filtered crystalline counterpart. This association is
constructed from the tensor equivalence and the natural comparison
of fibre functors in
\cite[
Proposition~7.4,
Lemma~7.12,
and \S7.13
]{OLSSON_NA_P_HODGE};
it is therefore natural for the morphisms in the chosen finite
diagram. Its compatibility with the regular group actions is
Corollary~7.15 of loc.~cit.

Applying the functorial equivariant cochain constructions of
\cite[\S\S6.13 and~6.15--6.17]{OLSSON_NA_P_HODGE}
gives the corresponding multiplicative comparison of coefficient
dgas; Corollary~7.16 of loc.~cit.\ records this comparison for the
regular coefficient ind-object. Pullback and the fixed-base-point
augmentation are compatible with the comparison by
\cite[Proposition~6.19 and \S6.20]{OLSSON_NA_P_HODGE}.
The identification of evaluation at the second endpoint with the
regular path-coordinate algebra, and its compatibility with the
period comparison, are given in
\cite[\S\S8.29--8.31]{OLSSON_NA_P_HODGE};
applied chartwise to the diagonal section, these give the
variable-endpoint augmentation. Under our Thom--Whitney
realization, the de Rham member of the resulting augmented diagram
is precisely
\[
\mathcal O_U
\xleftarrow{\epsilon_b}
\mathscr E_{b,K_v}^{\bullet}\big|_U
\xrightarrow{\epsilon_-}
\mathscr A_{\mathcal R,b,K_v}^{\dR}\big|_U.
\]

After extension to \(B_{\dR}\), repeat the construction of
\cite[\S8.12]{OLSSON_BAR_AFFINE_STACKS} with these associated
coefficient objects. At every finite coefficient stage, the maps in
the resulting analogue of the diagram of loc.\ cit.\ are the local
comparison maps of
\cite[\S6.13 and \S\S6.17--6.20]
{OLSSON_NA_P_HODGE}
with finite locally free coefficients. They are compatible with
multiplication and with both endpoint maps, and they are filtered
quasi-isomorphisms. Thus one obtains, over \(B_{\dR}\), a zigzag of
filtered multiplicative quasi-isomorphisms between the corresponding
\'{e}tale and de Rham coefficient spans.

Apply the simplicial two-sided bar construction to these spans.
In every fixed simplicial degree, finite local freeness of the
coefficient pieces and their filtered steps gives a filtered
quasi-isomorphism, over \(B_{\dR}\), between the corresponding
two-sided bar terms. The face and degeneracy
maps are formed from multiplication, the two module structures, and
the unit, so the degreewise comparisons assemble into a comparison
of simplicial filtered differential graded algebras.

Filtered Thom--Sullivan totalization preserves this comparison by
\cite[\S3]{OLSSON_BAR_AFFINE_STACKS}. On the underlying
differential graded algebras, replace the coefficient span by a
cofibrant one. Then
\cite[Corollary~4.6]{OLSSON_BAR_AFFINE_STACKS}
identifies the corresponding simplicial two-sided bar with the
homotopy pushout. Thus the filtered comparison constructed above is
the coefficientwise analogue of
\cite[diagram~(8.12.2)]{OLSSON_BAR_AFFINE_STACKS}.

Under
\cite[\S\S8.29--8.31]{OLSSON_NA_P_HODGE},
the induced map on zeroth cohomology is Olsson's canonical comparison
for the relative path torsor. Under
\[
\theta_{b,-}:
\mathcal H^0(\mathscr B_{b,K}^{\bullet})
\xrightarrow{\sim}
\mathscr A_{b,K}^{\dR},
\]
this is the comparison characterized by normalized tensor transport.
Consequently, the filtration induced by Olsson's period comparison is
the image of the Hodge convolution filtration on
\[
\mathscr B_{b,K}^{\bullet}\otimes_KK_v.
\]
This is the relative-coefficient version of
\cite[Corollary~8.13]{OLSSON_BAR_AFFINE_STACKS}.

Taking zeroth cohomology, passing to the filtered colimit, and gluing
the small charts gives
\[
F^p_{\mathrm{Ols}}\mathscr A_{b,K_v}^{\dR}
=
F^p\mathscr A_{b,K}^{\dR}\otimes_KK_v.
\]
Compatibility with multiplication, path composition, the torsor
coaction, and specialization follows from the functoriality of the
coefficient comparisons, the two-sided bar construction, and
normalized tensor transport.
\end{proof}

\subsection{Finite-Level Quotients}
\label{appendix:finite_level_filtered_path_torsors}

Let
\[
\mathcal G_n^\bullet,
\qquad
\bullet\in\{\et,\cri,\dR\},
\]
be a compatible finite-type motivic quotient in the
Kodaira--Parshin tower, and retain the quotient
\[
\mathcal G_K^{\dR}
\twoheadrightarrow
\mathcal G_{n,K}^{\dR}.
\]
Define the corresponding universal based path torsor by extension of
structure group:
\[
\mathcal P_{n,b,-,K}^{\dR}
\coloneqq 
\mathcal P_{b,-,K}^{\dR}
\times^{\mathcal G_K^{\dR}}
\mathcal G_{n,K}^{\dR}.
\]
Writing
\[
p_n:
\mathcal P_{n,b,-,K}^{\dR}
\longrightarrow
X^\dagger,
\qquad
\mathscr A_{n,b,K}^{\dR}
\coloneqq 
(p_n)_*\mathcal O_{\mathcal P_{n,b,-,K}^{\dR}},
\]
the quotient maps induce inclusions
\[
\mathscr A_{n,b,K}^{\dR}
\hookrightarrow
\mathscr A_{b,K}^{\dR},
\qquad
\mathcal O(\mathcal G_{n,K}^{\dR})
\hookrightarrow
\mathcal O(\mathcal G_K^{\dR}).
\]
We equip the finite-level coordinate algebras with the intersection
filtrations
\[
F^p\mathscr A_{n,b,K}^{\dR}
\coloneqq 
\mathscr A_{n,b,K}^{\dR}
\cap
F^p\mathscr A_{b,K}^{\dR}
\]
and
\[
F^p\mathcal O(\mathcal G_{n,K}^{\dR})
\coloneqq 
\mathcal O(\mathcal G_{n,K}^{\dR})
\cap
F^p\mathcal O(\mathcal G_K^{\dR}).
\]
These filtrations are decreasing, exhaustive, and multiplicative by
construction.
The content of the next result is that they agree with the intrinsic
filtered realizations of the compatible motivic quotient and behave
strictly under base change and specialization.

\begin{corollary}[Finite-level filtered realization]
\label{corollary:finite_level_filtered_path_groupoid}
For the fixed place \(v\mid p\), under the standing comparison
hypotheses, the preceding filtrations become, after extension to
\(K_v\), the filtrations supplied by Olsson's filtered association:
\[
F^p\mathscr A_{n,b,K}^{\dR}\otimes_KK_v
=
F^p_{\mathrm{Ols}}
\mathscr A_{n,b,K_v}^{\dR},
\]
and
\[
F^p\mathcal O(\mathcal G_{n,K}^{\dR})\otimes_KK_v
=
F^p_{\mathrm{Ols}}
\left(
\mathcal O(\mathcal G_{n,K}^{\dR})\otimes_KK_v
\right).
\]
Equivalently, the finite-level coordinate inclusions are strict for
the intrinsic filtered realizations.

The connection on \(\mathscr A_{n,b,K}^{\dR}\) satisfies Griffiths
transversality, the torsor coaction is filtered, and the filtration on
\(\mathcal O(\mathcal G_{n,K}^{\dR})\) is a Hopf-algebra filtration.

For every \(r\in\mathbb Z\), the quotient
\[
\mathscr A_{n,b,K}^{\dR}/F^r\mathscr A_{n,b,K}^{\dR}
\]
is flat over \(X^\dagger\). Consequently, for every field extension
\(K'/K\), every \(K'\)-scheme \(T\), and every morphism
\[
x:T\longrightarrow X^\dagger_{K'},
\]
write
\[
\mathcal P_{n,b,-,K'}^{\dR}
\coloneqq 
\mathcal P_{n,b,-,K}^{\dR}\otimes_KK',
\]
and let
\[
p_{n,x}:
x^*\mathcal P_{n,b,-,K'}^{\dR}
\longrightarrow
T
\]
denote the projection, the canonical affine base-change identification
\[
x^*
\left(
\mathscr A_{n,b,K}^{\dR}\otimes_KK',F^\bullet
\right)
\xrightarrow{\sim}
\left(
(p_{n,x})_*
\mathcal O_{x^*\mathcal P_{n,b,-,K'}^{\dR}},F^\bullet
\right)
\]
is strict, where the right-hand side carries the pulled-back
filtration.
In particular, if \(T=\operatorname{Spec}K'\), then
\[
x^*
\left(
\mathscr A_{n,b,K}^{\dR}\otimes_KK',F^\bullet
\right)
\xrightarrow{\sim}
\left(
\mathcal O({}^{\vphantom{\dR}}_{x}P_{b,n}^{\dR}),F^\bullet
\right)
\]
is a strict filtered isomorphism. Formation of all these filtrations
commutes with extension of the ground field and is natural under
morphisms of compatible finite-type motivic quotients.

Finally, after restriction to \(X_v\),
\[
\left.
\left(
\mathcal P_{n,b,-,K}^{\dR}\otimes_KK_v
\right)
\right|_{X_v}
\]
is canonically the universal finite-level de Rham path torsor
\[
\mathcal P_n^{\dR}\longrightarrow X_v
\]
of Definition~\ref{definition:universal_de_rham_path_torsor}, with
its canonical connection and family-level Hodge filtration.
\end{corollary}

\begin{proof}
By compatibility of \(\mathcal G_n^\bullet\) in the three
realizations, the finite-level coordinate inclusions are morphisms of
the associated path-coordinate ind-objects of
Theorem~\ref{theorem:olsson_etale_crystalline_comparison}. The
reductive regular coefficients underlying that construction are
associated stagewise by
\cite[
Proposition~7.7,
Remark~7.8,
and Corollary~7.14
]{OLSSON_NA_P_HODGE}, compatibly with their actions by
\cite[Corollary~7.15]{OLSSON_NA_P_HODGE}.
The corresponding endpoint identifications are those of
\cite[\S\S8.29--8.31]{OLSSON_NA_P_HODGE}.

Restricting the crystalline objects to the smooth open as in
Lemma~\ref{lemma:logarithmic_crystalline_restriction}, one obtains
associated isocrystals in the sense of Nizio{\l}. Their realization
is exact and tensorial, and exact sequences are strict on formally
smooth widenings; see
\cite[\S2.2.3 and Lemma~3.1]
{NIZIOL_CRYSTALLINE_SMOOTH_SHEAVES}.
Since the finite-rank source objects are finitely presented, each map
factors through a finite target stage. After enlarging that stage, its
kernel vanishes. Thus, after cofinal reindexing, each coordinate
inclusion is represented by strict monomorphisms between finite-rank
associated objects. At every such stage the filtered pieces and their
quotients are locally free.

It follows that the intrinsic finite-level Olsson filtration is the
intersection with the universal Olsson filtration. By
Lemma~\ref{lemma:olsson_bar_filtration_compatibility}, the latter is
the base change of the universal bar filtration. Since intersections
are kernels and \(K_v/K\) is faithfully flat, the scalar extensions of
the intersection filtrations defined above are therefore the
intrinsic finite-level Olsson filtrations.

Filtered colimits of the finite-stage filtered pieces and quotients
are flat, and this flatness descends along \(K_v/K\). The universal
Griffiths-transverse connection and filtered structure maps restrict
to the finite level by strictness; this may be checked after extension
to \(K_v\).

The preceding flatness shows that the defining exact sequence remains
exact after arbitrary pullback on \(X^\dagger\). This proves the
strict base-change assertion; for field-valued endpoints it agrees with
Lemma~\ref{lemma:strict_endpoint_specialization}. Naturality follows
from the functoriality of the associated realizations. The final
identification follows from the compatibility of the finite-level
quotient in the three realizations and the defining
tensor-isomorphism description of the path torsor.
\end{proof}

\begin{definition}[Finite-level Hodge loci]
\label{definition:finite_level_hodge_loci}
Set
\[
\mathcal I_n^G
\coloneqq 
\left\langle
F^1\mathcal O(\mathcal G_{n,K}^{\dR})
\right\rangle
\]
and define
\[
F^0\mathcal G_{n,K}^{\dR}
\coloneqq 
\operatorname{Spec}
\left(
\mathcal O(\mathcal G_{n,K}^{\dR})/\mathcal I_n^G
\right).
\]
Likewise, set
\[
\mathcal I_n^P
\coloneqq 
\left\langle
F^1\mathscr A_{n,b,K}^{\dR}
\right\rangle
\]
and define
\[
\mathcal P_{n,b,-,K}^H
\coloneqq 
F^0\mathcal P_{n,b,-,K}^{\dR}
\coloneqq 
\operatorname{Spec}_{X^\dagger}
\left(
\mathscr A_{n,b,K}^{\dR}/\mathcal I_n^P
\right)
\subseteq
\mathcal P_{n,b,-,K}^{\dR}.
\]

Equivalently, for a commutative \(K\)-algebra \(R\),
\[
F^0\mathcal G_{n,K}^{\dR}(R)
=
\left\{
g\in\mathcal G_{n,K}^{\dR}(R)
\ \middle|\
\operatorname{ev}_g:
\mathcal O(\mathcal G_{n,K}^{\dR})
\longrightarrow R
\text{ is filtered}
\right\},
\]
and, for a \(K\)-scheme \(T\) with a morphism
\(x:T\to X^\dagger\),
each point
\[
q:T\longrightarrow\mathcal P_{n,b,-,K}^{\dR}
\]
over \(X^\dagger\) induces an evaluation map
\[
\operatorname{ev}_q:
x^*\mathscr A_{n,b,K}^{\dR}
\longrightarrow
\mathcal O_T,
\]
and
\[
\operatorname{Hom}_{X^\dagger}
\left(
T,\mathcal P_{n,b,-,K}^H
\right)
=
\left\{
q\in
\operatorname{Hom}_{X^\dagger}
\left(T,\mathcal P_{n,b,-,K}^{\dR}\right)
\ \middle|\
\operatorname{ev}_q\text{ is filtered}
\right\}.
\]
Here \(R\) and \(\mathcal O_T\) carry the trivial filtration.
\end{definition}

\begin{corollary}[Finite-level Hodge reduction]
\label{corollary:finite_level_universal_hodge_reduction}
The ideal \(\mathcal I_n^G\) is a Hopf ideal. Consequently,
\[
F^0\mathcal G_{n,K}^{\dR}
\subseteq
\mathcal G_{n,K}^{\dR}
\]
is a closed subgroup, and
\[
\mathcal P_{n,b,-,K}^H
\longrightarrow
X^\dagger
\]
is an fppf torsor under \(F^0\mathcal G_{n,K}^{\dR}\). Thus it is a
reduction of structure group of
\(\mathcal P_{n,b,-,K}^{\dR}\).

The construction commutes with extension of the ground field. More
generally, for every field extension \(K'/K\), every \(K'\)-scheme
\(S\), and every morphism
\[
x:S\longrightarrow X^\dagger_{K'},
\]
there is a canonical identification
\[
S\times_{X^\dagger_{K'}}
\left(
\mathcal P_{n,b,-,K}^H\otimes_KK'
\right)
\xrightarrow{\sim}
F^0{}^{\vphantom{\dR}}_{x}P_{b,n}^{\dR},
\]
where the right-hand side is the closed subscheme of
\(x^*\mathcal P_{n,b,-,K'}^{\dR}\) defined by the ideal sheaf
generated by the pulled-back \(F^1\).
These identifications are natural under morphisms of compatible
finite-type motivic quotients.
\end{corollary}

\begin{proof}
Filteredness of the Hopf maps gives
\[
\Delta(\mathcal I_n^G)
\subseteq
\mathcal I_n^G\otimes_K
\mathcal O(\mathcal G_{n,K}^{\dR})
+
\mathcal O(\mathcal G_{n,K}^{\dR})\otimes_K
\mathcal I_n^G,
\]
while the counit vanishes on \(\mathcal I_n^G\) and the antipode
preserves it. Thus \(\mathcal I_n^G\) is a Hopf ideal. Likewise, the
filtered torsor coaction restricts the action to
\[
\mathcal P_{n,b,-,K}^H
\times_K
F^0\mathcal G_{n,K}^{\dR}
\longrightarrow
\mathcal P_{n,b,-,K}^H.
\]
Let \(\alpha_n\) denote the resulting morphism
\[
\alpha_n:
\mathcal P_{n,b,-,K}^H
\times_KF^0\mathcal G_{n,K}^{\dR}
\longrightarrow
\mathcal P_{n,b,-,K}^H
\times_{X^\dagger}
\mathcal P_{n,b,-,K}^H,
\qquad
(q,g)\longmapsto(q,q\cdot g).
\]

Fix an embedding \(K\hookrightarrow\mathbb C\). The filtered complex
comparison identifies the schemes in \(\alpha_n\) with the
corresponding Hodge loci in Hain's relative Malcev path groupoid. The
path-composition and bar-reversal formulas are morphisms of mixed
Hodge structures; see
\cite[Theorems~13.1 and~13.3 and their proofs]{HAIN_HDR_RMC}.
Consequently, after analytification and extension to smooth functions,
the ordinary torsor difference map restricts to the inverse of
\(\alpha_n\).

Moreover, Theorem~13.10 and the construction preceding
Proposition~14.4 of \cite{HAIN_HDR_RMC} give local \(C^\infty\)
filtered sections of the relative path bundle: on coordinate algebras
the section is induced by Deligne's multiplicative bigrading and
projection to \(\Gr_0^W\), and the reductive part is supplied by the
Hodge-adapted frames constructed in the proof of
Proposition~\ref{proposition:filtered_reductive_path_torsor}.
Restricting these sections to the strict finite-level coordinate
subalgebra gives local sections of the finite-level Hodge locus.

Thus, on every sufficiently small analytic open \(U\), extension along
\(\mathcal O_U\to\mathcal C_U^\infty\) makes \(\alpha_n\) an
isomorphism and gives
\[
\left.
\left(
\left(
\mathscr A_{n,b,K}^{\dR}/\mathcal I_n^P
\right)
\otimes_K\mathbb C
\right)^{\an}
\right|_U
\otimes_{\mathcal O_U}\mathcal C_U^\infty
\xrightarrow{\sim}
\mathcal C_U^\infty
\otimes_{\mathbb C}
\mathcal O
\left(
F^0\mathcal G_{n,K}^{\dR}\otimes_K\mathbb C
\right).
\]
The right-hand side is faithfully flat over
\(\mathcal C_U^\infty\). Since
\(\mathcal O_U\to\mathcal C_U^\infty\) and the
algebraic-to-analytic local maps are faithfully flat, \(\alpha_n\)
is an isomorphism after extension to \(\mathbb C\), and
\[
\mathcal P_{n,b,-,K}^H\otimes_K\mathbb C
\longrightarrow
X^\dagger_{\mathbb C}
\]
is faithfully flat. Since the finite-level path torsor is of finite
type over the noetherian scheme \(X^\dagger_{\mathbb C}\), this closed
subscheme is of finite presentation. It is therefore an fppf torsor.
Faithfully flat descent along \(\mathbb C/K\) gives both assertions
over \(K\), and hence the claimed reduction of structure group.

Finally, the strict base-change assertion of
Corollary~\ref{corollary:finite_level_filtered_path_groupoid}
identifies the pullback of \(\mathcal I_n^P\) with the ideal generated
by the pulled-back \(F^1\). This proves the base-change formula.
Naturality follows from the strict filtered maps induced by morphisms
of compatible finite-type motivic quotients.
\end{proof}

\begin{proposition}[Finite-level \(B_{\dR}^+\)-comparison]
\label{proposition:selected_finite_level_filtered_comparison}
Let
\[
\mathcal G_n^\bullet,
\qquad
\bullet\in\{\et,\cri,\dR\},
\]
be a compatible finite-type motivic quotient, and put
\[
\mathcal G_{n,K_v}^{\dR}
\coloneqq 
\mathcal G_{n,K}^{\dR}\otimes_KK_v,
\qquad
F^0\mathcal G_{n,K_v}^{\dR}
\coloneqq 
F^0\mathcal G_{n,K}^{\dR}\otimes_KK_v.
\]
For a commutative \(K_v\)-algebra \(B\), write
\[
B_{\dR,B}
\coloneqq 
B\otimes_{K_v}B_{\dR},
\qquad
B_{\dR,B}^+
\coloneqq 
B\otimes_{K_v}B_{\dR}^+,
\]
with trivial \(G_v\)-action on \(B\). All tensor products are
algebraic.

The finite-level period comparison induces functorial
identifications
\[
\mathcal G_{n,K_v}^{\dR}(B)
\xrightarrow{\sim}
\mathcal G_n^{\et}(B_{\dR,B})^{G_v},
\qquad
F^0\mathcal G_{n,K_v}^{\dR}(B)
\xrightarrow{\sim}
\mathcal G_n^{\et}(B_{\dR,B}^+)^{G_v},
\]
compatible with the natural inclusions. Likewise, for
\[
b,x\in\mathcal X^\circ(\mathcal O_{K_v}),
\]
there are functorial identifications
\[
{}^{\vphantom{\dR}}_{x}P_{b,n}^{\dR}(B)
\xrightarrow{\sim}
{}^{\vphantom{\et}}_{\overline x}P_{\overline b,n}^{\et}
(B_{\dR,B})^{G_v},
\qquad
F^0{}^{\vphantom{\dR}}_{x}P_{b,n}^{\dR}(B)
\xrightarrow{\sim}
{}^{\vphantom{\et}}_{\overline x}P_{\overline b,n}^{\et}
(B_{\dR,B}^+)^{G_v}.
\]
The two pairs of identifications are compatible with the inclusions
\(B_{\dR,B}^+\hookrightarrow B_{\dR,B}\), and are natural in \(B\),
in the endpoints, and under morphisms of compatible finite-type
motivic quotients. They are also compatible with group laws, torsor
actions, and composition of paths.
\end{proposition}

\begin{proof}
Write
\[
\mathcal A_n^{\et}
\coloneqq 
\mathcal O(\mathcal G_n^{\et}),
\qquad
\mathcal A_n^{\dR}
\coloneqq 
\mathcal O(\mathcal G_{n,K_v}^{\dR}).
\]
For integral endpoints \(b,x\), write likewise
\[
\mathcal A_{x,b,n}^{\et}
\coloneqq 
\mathcal O
\left(
{}^{\vphantom{\et}}_{\overline x}P_{\overline b,n}^{\et}
\right),
\qquad
\mathcal A_{x,b,n}^{\dR}
\coloneqq 
\mathcal O
\left(
{}^{\vphantom{\dR}}_{x}P_{b,n}^{\dR}
\right).
\]
The filtered comparisons of
Theorems~\ref{theorem:olsson_etale_crystalline_comparison}
and~\ref{theorem:olsson_crystalline_de_rham_comparison}, together with
Corollary~\ref{corollary:finite_level_filtered_path_groupoid}, give
\(G_v\)-equivariant algebra isomorphisms
\[
c_n^G:
\mathcal A_n^{\et}\otimes_{\Qp}B_{\dR}
\xrightarrow{\sim}
\mathcal A_n^{\dR}\otimes_{K_v}B_{\dR}
\]
and
\[
c_{x,b,n}^P:
\mathcal A_{x,b,n}^{\et}\otimes_{\Qp}B_{\dR}
\xrightarrow{\sim}
\mathcal A_{x,b,n}^{\dR}\otimes_{K_v}B_{\dR}.
\]
They are compatible with the algebra and path-groupoid structure
maps. By
Lemma~\ref{lemma:olsson_bar_filtration_compatibility}, they identify
the standard \'{e}tale \(B_{\dR}^+\)-lattices with the Rees lattices
\[
\operatorname{Rees}(\mathcal A)
\coloneqq 
\sum_r
F^r\mathcal A\otimes_{K_v}t^{-r}B_{\dR}^+
\subseteq
\mathcal A\otimes_{K_v}B_{\dR}
\]
for \(\mathcal A=\mathcal A_n^{\dR}\) and
\(\mathcal A=\mathcal A_{x,b,n}^{\dR}\), respectively. This statement
is understood stagewise in the corresponding ind-systems.

Since \(B_{\dR}^{G_v}=K_v\) and \(B\) is a \(K_v\)-vector space,
\[
(B_{\dR,B})^{G_v}=B.
\]
Consequently, extension through \(c_n^G\) and \(c_{x,b,n}^P\)
identifies de Rham \(B\)-points with the corresponding \(G_v\)-fixed
\(B_{\dR,B}\)-valued \'{e}tale points. This gives the two unfiltered
identifications.

It remains to characterize the \(B_{\dR,B}^+\)-valued points. Let
\((\mathcal A,F^\bullet)\) be either of the two de Rham coordinate
algebras above, and let \(f:\mathcal A\to B\) be a \(K_v\)-algebra
map. Since \(B\) has the trivial filtration,
\[
f\text{ is filtered}
\quad\Longleftrightarrow\quad
f(F^1\mathcal A)=0.
\]
The latter condition is equivalent to
\[
(f\otimes1)
\left(
\operatorname{Rees}(\mathcal A)
\right)
\subseteq
B_{\dR,B}^+.
\]
Indeed, the forward implication follows directly from the definition
of the Rees lattice. Conversely, for \(a\in F^1\mathcal A\), the
element \(a\otimes t^{-1}\) belongs to
\(\operatorname{Rees}(\mathcal A)\); if its image lies in
\(B_{\dR,B}^+\), then \(f(a)=0\), since
\[
\left(B\otimes_{K_v}K_v\right)
\cap
\left(B\otimes_{K_v}tB_{\dR}^+\right)
=0
\]
inside \(B_{\dR,B}\).

Through the filtered period comparisons, this says precisely that
the associated \'{e}tale point is defined over \(B_{\dR,B}^+\).
Definition~\ref{definition:finite_level_hodge_loci} now gives the two
\(F^0\)-identifications. Compatibility with inclusions, group laws,
torsor actions, path composition, endpoints, and quotient maps follows
from the corresponding compatibility of the coordinate-algebra
comparisons.
\end{proof}

\section{Auxiliary Tannakian Finiteness Results}
\label{appendix:tannakian_finiteness}

This appendix collects the auxiliary finiteness results used in the
proof of Theorem~\ref{thm:ftg_basic} via abelianizations. We also give
a conceptual alternative proof of that theorem using universal
extensions and a hereditary Tannakian envelope.

\subsection{Auxiliary Finiteness Lemmas}
\label{sec:auxiliary_tannakian_finiteness}

\subsubsection{Pro-unipotent groups}
\label{sec:facts_about_prounip}

Let \(\widehat{\mathbb L}(V)\) denote the free pronilpotent Lie
algebra on a vector space, or pro-vector space, \(V\).

\begin{fact}
[{Pro-unipotent generation;
cf.~\cite[Lemma~2.5]{LubMag82}}]
\label{fact:finite_type_criterion}
Let \(\mathfrak u\) be a pronilpotent Lie algebra over \(F\). Any
continuous linear section
\[
\mathfrak u^{\ab}\longrightarrow\mathfrak u
\]
induces a surjection
\[
\widehat{\mathbb L}(\mathfrak u^{\ab})
\twoheadrightarrow
\mathfrak u.
\]
Equivalently, a pro-unipotent group is generated by any set whose
image generates its abelianization.
\end{fact}

For \(n\geq0\) and \(j\geq1\), set
\[
\Phi_n(j)
\coloneqq 
\dim_F
\Gr_j^{\DCS}\widehat{\mathbb L}(F^n).
\]
Thus \(\Phi_0(j)=0\), while, for \(n\geq1\), Witt's formula gives
\begin{equation}
\label{wittformula}
\Phi_n(j)
=
\frac{1}{j}
\sum_{d\mid j}
\mu(d)n^{j/d},
\end{equation}
where \(\mu\) is the M\"obius function; see
\cite[Satz~3]{Witt37Treue} and
\cite[(1.1)]{Waldinger67Witt}. For \(n\geq0\) and \(k\geq1\), set
\[
d_n(k)\coloneqq \sum_{j=1}^k\Phi_n(j);
\]
in particular, \(d_0(k)=0\).

\begin{corollary}
\label{cor:prounip_dim_bound}
Let \(U\) be a pro-unipotent group over \(F\), and let \(k\geq1\). If
\[
\dim_F\Lie(U)^{\ab}=n
\qquad\text{and}\qquad
\Fc_{\DCS}^{k+1}U=1,
\]
then \(U\) is algebraic and
\[
\dim U\leq d_n(k).
\]
\end{corollary}

\begin{proof}
Set \(\mathfrak u\coloneqq \Lie(U)\). Choose a linear section of
\(\mathfrak u\twoheadrightarrow\mathfrak u^{\ab}\). By
Fact~\ref{fact:finite_type_criterion}, it induces a surjection
\[
\widehat{\mathbb L}(\mathfrak u^{\ab})
\twoheadrightarrow
\mathfrak u.
\]
Since \(\Fc_{\DCS}^{k+1}\mathfrak u=0\), this map factors through
\[
\widehat{\mathbb L}(\mathfrak u^{\ab})/
\Fc_{\DCS}^{k+1}\widehat{\mathbb L}(\mathfrak u^{\ab}),
\]
whose dimension is \(d_n(k)\) by \eqref{wittformula}. Hence
\[
\dim_F\mathfrak u\leq d_n(k).
\]
The equivalence between pronilpotent Lie algebras and pro-unipotent
groups in characteristic zero now gives the result.
\end{proof}

\subsubsection{Algebraic Tannakian categories}
\label{sec:noeth_tann_sub}

\begin{lemma}
\label{lemma:tann_noetherian}
Let \((\Tc,\omega)\) be an algebraic neutral Tannakian category.

\begin{enumerate}
\item If a collection \(A\subseteq\Tc\) tensor-generates \(\Tc\),
then some finite subset \(A_0\subseteq A\) tensor-generates \(\Tc\).

\item Every ascending chain of Tannakian subcategories of \(\Tc\)
stabilizes.
\end{enumerate}
\end{lemma}

\begin{proof}
Since \(\Tc\) is algebraic, it admits a tensor generator \(E\). If
\(A\) generates \(\Tc\), then \(E\) belongs to the directed union of
the subcategories generated by the finite subsets of \(A\). Thus
\[
E\in\langle A_0\rangle_\otimes
\]
for some finite \(A_0\subseteq A\), and hence
\(\Tc=\langle A_0\rangle_\otimes\).

For the second assertion, let
\[
\Tc_0\subseteq\Tc_1\subseteq\Tc_2\subseteq\cdots\subseteq\Tc
\]
be an ascending chain, and set
\[
\Gc\coloneqq \pi_1(\Tc,\omega),
\qquad
K_i\coloneqq 
\ker\bigl(\Gc\to\pi_1(\Tc_i,\omega)\bigr).
\]
If \(I_i\subseteq\mathcal O(\Gc)\) is the defining ideal of \(K_i\),
then the \(I_i\) form an ascending chain. Since \(\Gc\) is of finite
type, \(\mathcal O(\Gc)\) is Noetherian, so the ideals \(I_i\), and
therefore the subcategories \(\Tc_i\), eventually stabilize.
\end{proof}

\begin{lemma}
\label{lemma:unipotent_bounded_dim}
Let \(\{U_i\}_{i\in I}\) be a cofiltered inverse system of unipotent
groups over \(F\), with faithfully flat transition maps. If
\[
\dim U_i\leq d
\qquad(i\in I)
\]
for some \(d\), then the inverse limit
\[
\varprojlim_i U_i
\]
is an algebraic unipotent group.
\end{lemma}

\begin{proof}
Choose \(i_0\in I\) for which \(\dim U_{i_0}\) is maximal. If
\(U_j\to U_{i_0}\) is a transition map, then
\[
\dim U_j\geq\dim U_{i_0}
\]
by faithful flatness, while maximality gives the reverse inequality.
The kernel is therefore a zero-dimensional unipotent group and hence
is trivial in characteristic zero. Thus every term mapping to
\(U_{i_0}\) is isomorphic to \(U_{i_0}\).

The subsystem of such terms is cofinal, so the original inverse
system is essentially constant and
\[
\varprojlim_iU_i\simeq U_{i_0}.
\]
\end{proof}

\subsection{An Alternative Proof via Universal Extensions}
\label{sec:proof_using_univ_ext}

Retain the notation and hypotheses of
Theorem~\ref{thm:ftg_basic}. We give an alternative proof by first
treating hereditary Tannakian categories and then embedding \(\Tc\)
into a hereditary Tannakian category.

Recall that an abelian category \(\Ac\) is \emph{hereditary} if
\[
\Ext_{\Ac}^n(M,N)=0
\qquad
(M,N\in\Ac,\ n\geq2).
\]

By
\cite[Proposition~3.3 and Lemma~3.5]
{AddezioEsnault22UniversalExt},
for \(M,N\in\Tc\) there is a universal extension
\[
0\longrightarrow M
\longrightarrow
\Ec_{\Tc}(N,M)
\longrightarrow
\Ext^1_{\Tc}(N,M)\otimes_F N
\longrightarrow0
\]
in \(\Tc\), characterized by the following property: for every
\(\alpha\in\Ext^1_{\Tc}(N,M)\), pulling back along
\[
\alpha\otimes\id_N:
N\longrightarrow
\Ext^1_{\Tc}(N,M)\otimes_F N
\]
gives an extension representing \(\alpha\).

\begin{lemma}
\label{lemma:universal_embedding}
Let
\[
0\longrightarrow M
\longrightarrow E_\alpha
\longrightarrow N
\longrightarrow0
\]
represent \(\alpha\in\Ext^1_{\Tc}(N,M)\). Then there is a
monomorphism
\[
E_\alpha
\hookrightarrow
\Ec_{\Tc}(N,M)\oplus N.
\]
\end{lemma}

\begin{proof}
This is the special case \(A'=M\) and \(B'=N\) of
\cite[Proposition~3.4]{AddezioEsnault22UniversalExt}.
\end{proof}

Define objects \(\Ec_k\in\Tc\) recursively by
\[
\Ec_0\coloneqq \Vc,
\qquad
\Ec_k
\coloneqq 
\Ec_{\Tc}(\Ec_{k-1},\Vc)\oplus\Ec_{k-1}
\quad(k\geq1).
\]

\begin{proposition}
\label{prop:Ek_embeds_universal}
If \(\Tc\) is hereditary, then every
\[
E_k\in\IExt^k(\Vc)
\]
embeds into \(\Ec_k\).
\end{proposition}

\begin{proof}
We argue by induction on \(k\). The case \(k=0\) is immediate.
Suppose \(k\geq1\), and write
\[
0\longrightarrow\Vc
\longrightarrow E_k
\longrightarrow E'_{k-1}
\longrightarrow0
\]
as in \eqref{ses:e_k'}. By induction, there is an exact sequence
\[
0\longrightarrow E'_{k-1}
\longrightarrow\Ec_{k-1}
\longrightarrow C
\longrightarrow0.
\]
The corresponding long exact sequence contains
\[
\Ext^1_{\Tc}(\Ec_{k-1},\Vc)
\longrightarrow
\Ext^1_{\Tc}(E'_{k-1},\Vc)
\longrightarrow
\Ext^2_{\Tc}(C,\Vc).
\]
Since \(\Tc\) is hereditary, the last group vanishes. Hence the class
of \(E_k\) lifts to a class
\[
\beta\in\Ext^1_{\Tc}(\Ec_{k-1},\Vc).
\]
Let \(E_\beta\) be an extension representing \(\beta\). Pulling it
back along \(E'_{k-1}\hookrightarrow\Ec_{k-1}\) recovers \(E_k\), so
\[
E_k\hookrightarrow E_\beta.
\]
By Lemma~\ref{lemma:universal_embedding},
\[
E_\beta
\hookrightarrow
\Ec_{\Tc}(\Ec_{k-1},\Vc)\oplus\Ec_{k-1}
=
\Ec_k.
\]
\end{proof}

Suppose temporarily that \(\Tc\) is hereditary. Since \(\Rc\) is of
finite type, \(\Tc_0\simeq\operatorname{Rep}_F(\Rc)\) has a tensor
generator \(W_0\). Proposition~\ref{prop:Ek_embeds_universal} gives
\[
\Tc_{\Vc,k}
\subseteq
\left\langle
\Ec_k\oplus W_0
\right\rangle_\otimes.
\]
The category on the right is algebraic, and therefore so is its
Tannakian subcategory \(\Tc_{\Vc,k}\). This proves
Theorem~\ref{thm:ftg_basic} when \(\Tc\) is hereditary.

We now reduce the general case to this one.

\begin{proposition}
\label{prop:embed_hereditary}
There exists a hereditary neutral Tannakian category
\((\Tc^h,\omega^h)\) and a fully faithful exact tensor functor
\[
\iota:\Tc\hookrightarrow\Tc^h
\]
compatible with the fibre functors, such that:

\begin{enumerate}
\item the semisimple subcategory of \(\Tc^h\) is identified with
\(\Tc_0\simeq\operatorname{Rep}_F(\Rc)\);

\item for all \(M,N\in\Tc^h\), the vector space
\[
\Ext^1_{\Tc^h}(M,N)
\]
is finite-dimensional.
\end{enumerate}
\end{proposition}

\begin{proof}
Using the fibre functor, identify
\[
\Tc\simeq\operatorname{Rep}_F(\Gc).
\]
Choose a Levi splitting
\[
\Gc\simeq\Uc\rtimes\Rc,
\qquad
\uf\coloneqq \Lie(\Uc).
\]
Since \(\Rc\) is linearly reductive, the quotient
\[
\uf\twoheadrightarrow\uf^{\ab}
\]
admits a continuous \(\Rc\)-equivariant linear section.

Let
\[
\mathfrak f^h
\coloneqq 
\widehat{\mathbb L}(\uf^{\ab}),
\qquad
\Uc^h\coloneqq \exp(\mathfrak f^h),
\qquad
\Gc^h\coloneqq \Uc^h\rtimes\Rc,
\]
and set
\[
\Tc^h\coloneqq \operatorname{Rep}_F(\Gc^h).
\]
The chosen section induces an \(\Rc\)-equivariant homomorphism
\[
\mathfrak f^h\longrightarrow\uf
\]
which is surjective by
Fact~\ref{fact:finite_type_criterion}. Exponentiating gives a
faithfully flat morphism
\[
\Gc^h\twoheadrightarrow\Gc,
\]
and hence a fully faithful exact tensor functor
\[
\Tc\hookrightarrow\Tc^h
\]
preserving the fibre functors. Both categories have reductive quotient
\(\Rc\), so their semisimple subcategories are identified with
\(\operatorname{Rep}_F(\Rc)\).

We next prove that \(\Tc^h\) is hereditary. For
\(M,N\in\Tc^h\), set
\[
W\coloneqq M^\vee\otimes N.
\]
Then
\[
\Ext^n_{\Tc^h}(M,N)
\simeq
H^n(\Gc^h,W).
\]
The extension
\[
1\longrightarrow\Uc^h
\longrightarrow\Gc^h
\longrightarrow\Rc
\longrightarrow1
\]
gives a Hochschild--Serre spectral sequence
\[
E_2^{a,b}
=
H^a\bigl(\Rc,H^b(\Uc^h,W)\bigr)
\Longrightarrow
H^{a+b}(\Gc^h,W).
\]
Since \(\Rc\) is linearly reductive,
\[
H^a(\Rc,-)=0
\qquad(a>0).
\]
Moreover, a free pro-unipotent group has cohomological dimension one,
so
\[
H^b(\Uc^h,W)=0
\qquad(b>1);
\]
see~\cite[Theorem~2.9]{LubMag82}. Consequently,
\[
\Ext^n_{\Tc^h}(M,N)=0
\qquad(n\geq2),
\]
and \(\Tc^h\) is hereditary.

It remains to prove finite-dimensionality of \(\Ext^1\). If \(M\)
and \(N\) are semisimple, they factor through \(\Rc\), and
Hochschild--Serre gives
\[
\begin{aligned}
\Ext^1_{\Tc^h}(M,N)
&\simeq
\Hom_{\Rc}
\bigl(
(\mathfrak f^h)^{\ab},
\omega(M^\vee\otimes N)
\bigr)
\\
&\simeq
\Hom_{\Rc}
\bigl(
\uf^{\ab},
\omega(M^\vee\otimes N)
\bigr)
\\
&\simeq
\Ext^1_{\Tc}(M,N).
\end{aligned}
\]
This is finite-dimensional by hypothesis. The general case follows
by induction on the lengths of \(M\) and \(N\), using the long exact
sequences for \(\Ext\) and the vanishing of \(\Ext^2\).
\end{proof}

\begin{proof}[Completion of the alternative proof of
Theorem~\ref{thm:ftg_basic}]
Apply the hereditary case to \(\Tc^h\). It follows that
\[
\Tc^h_{\Vc,k}
\]
is algebraic for every \(k\). Since
\(\Tc\hookrightarrow\Tc^h\) is exact and identifies the two
semisimple subcategories, every \(k\)-fold iterated extension of
\(\Vc\) in \(\Tc\) is also one in \(\Tc^h\). Hence
\[
\Tc_{\Vc,k}
\subseteq
\Tc^h_{\Vc,k}.
\]
The Tannaka group of \(\Tc_{\Vc,k}\) is a quotient of the algebraic
Tannaka group of \(\Tc^h_{\Vc,k}\), and is therefore algebraic. This
completes the alternative proof.
\end{proof}

\section{Supported Quotients of Pro-Unipotent Extensions}
\label{appendix:adjoint_supported_quotients}
This appendix develops the formalism of maximal quotients of
pro-unipotent extensions supported on a prescribed finite set
\(\Sigma\) of \(\Rc\)-isotypes; these results are used in
\S\ref{sec:param_families} to construct the adjoint-supported
Kodaira--Parshin tower. In
\S\ref{appendix:support_projectors} we define the support projector
\(Q_\Sigma\); in
\S\S\ref{appendix:supported_lie_quotients}--\ref{appendix:supported_group_quotients}
we construct the corresponding maximal quotients of pronilpotent Lie
algebras and pro-unipotent extensions, obtaining
Proposition~\ref{proposition:maximal_adjoint_supported_quotient};
finally, \S\ref{appendix:supported_free_lie} establishes the
free-Lie and perfect-target results used in the proof of
Theorem~\ref{theorem:special_motivic_exists}.

Let \(F\) be a field of characteristic zero, and let \(\Rc\) be a
reductive algebraic \(F\)-group. Fix a finite
subset
\[
\Sigma\subseteq\operatorname{Irr}(\Rc),
\]
and write
\[
\Sigma^c
\coloneqq 
\operatorname{Irr}(\Rc)\setminus\Sigma
\]
for its complement.

The purpose of this appendix is to construct the maximal quotient of
a pro-unipotent extension whose DCS graded
pieces are supported on \(\Sigma\). In the application of
\S\ref{sec:param_families}, \(\Sigma=\Sigma_{\mathrm{ad}}\) is the
set of simple constituents of
\[
\mathfrak a
\coloneqq 
\Lie\bigl((\Rc^\circ)^{\mathrm{der}}\bigr),
\]

\subsection{Support projectors}
\label{appendix:support_projectors}

\begin{definition}
An object \(W\in\operatorname{Rep}_F(\Rc)\) is
\emph{\(\Sigma\)-supported} if every simple constituent of \(W\)
belongs to \(\Sigma\).
\end{definition}

For \(W\in\operatorname{Rep}_F(\Rc)\), let \(W_\Sigma\subseteq W\)
denote the sum of its simple subobjects belonging to \(\Sigma\), and
let \(W_{\Sigma^c}\subseteq W\) denote the sum of the remaining
isotypic components. Semisimplicity gives a canonical decomposition
\[
W=W_\Sigma\oplus W_{\Sigma^c}.
\]
We write
\[
Q_\Sigma(W)\coloneqq W_\Sigma,
\qquad
C_{\Sigma^c}(W)\coloneqq W_{\Sigma^c},
\]
and regard \(Q_\Sigma(W)\) both as a subobject and as the quotient
\[
W\twoheadrightarrow W/C_{\Sigma^c}(W).
\]

\begin{proposition}
\label{proposition:support_projector}
There is a canonical product decomposition of abelian categories
\[
\operatorname{Rep}_F(\Rc)
\simeq
\operatorname{Rep}_{F,\Sigma}(\Rc)
\times
\operatorname{Rep}_{F,\Sigma^c}(\Rc).
\]
Under this equivalence, \(Q_\Sigma\) and \(C_{\Sigma^c}\) are obtained by
projection to the first and second factors, respectively, followed by
the corresponding inclusions into \(\operatorname{Rep}_F(\Rc)\).

The decomposition extends canonically to
\[
\operatorname{Ind}\operatorname{Rep}_F(\Rc)
\simeq
\operatorname{Ind}\operatorname{Rep}_{F,\Sigma}(\Rc)
\times
\operatorname{Ind}\operatorname{Rep}_{F,\Sigma^c}(\Rc)
\]
and
\[
\operatorname{Pro}\operatorname{Rep}_F(\Rc)
\simeq
\operatorname{Pro}\operatorname{Rep}_{F,\Sigma}(\Rc)
\times
\operatorname{Pro}\operatorname{Rep}_{F,\Sigma^c}(\Rc).
\]
\end{proposition}

\begin{proof}
Since \(\operatorname{Rep}_F(\Rc)\) is semisimple, every object \(W\)
is the direct sum of its simple constituents. Let \(W_\Sigma\) and
\(W_{\Sigma^c}\) be the sums of the constituents belonging to
\(\Sigma\) and \(\Sigma^c\), respectively. Then
\[
W=W_\Sigma\oplus W_{\Sigma^c}.
\]
There are no nonzero morphisms between the two classes of simple
objects, so this decomposition is canonical and functorial in \(W\).
It therefore induces the asserted product decomposition, with
\[
Q_\Sigma(W)=W_\Sigma,
\qquad
C_{\Sigma^c}(W)=W_{\Sigma^c}.
\]

Finally, the \(\operatorname{Ind}\) and \(\operatorname{Pro}\)
constructions commute with finite products of categories, yielding the remaining
equivalences.
\end{proof}

\subsection{Supported quotients of Lie algebras}
\label{appendix:supported_lie_quotients}

Let \(\operatorname{Lie}_{\Rc}\) denote the category of
\(\Rc\)-equivariant Lie algebras whose underlying objects belong to
\(\operatorname{Ind}\operatorname{Rep}_F(\Rc)\).

\begin{definition}
For \(\gf\in\operatorname{Lie}_{\Rc}\), let
\[
I_\Sigma(\gf)
\coloneqq 
\bigl\langle C_{\Sigma^c}(\gf)\bigr\rangle_{\mathrm{Lie}}
\]
be the Lie ideal generated by the non-\(\Sigma\)-supported component
of the underlying representation, and set
\[
\qf_\Sigma(\gf)
\coloneqq 
\gf/I_\Sigma(\gf).
\]
\end{definition}

\begin{proposition}
\label{proposition:maximal_supported_lie_quotient}
Let
\[
\operatorname{Lie}_{\Rc,\Sigma}
\subseteq
\operatorname{Lie}_{\Rc}
\]
denote the full subcategory of Lie algebras whose underlying
\(\Rc\)-representations are \(\Sigma\)-supported. The assignment
\[
\gf\longmapsto\qf_\Sigma(\gf)
\]
defines a functor
\[
\qf_\Sigma:
\operatorname{Lie}_{\Rc}
\longrightarrow
\operatorname{Lie}_{\Rc,\Sigma}
\]
left adjoint to the inclusion.

Equivalently, the quotient map
\[
p_\gf:\gf\twoheadrightarrow\qf_\Sigma(\gf)
\]
exhibits \(\qf_\Sigma(\gf)\) as the maximal
\(\Sigma\)-supported quotient of \(\gf\): every
\(\Sigma\)-supported quotient
\[
\gf\twoheadrightarrow\hf
\]
factors uniquely through a surjection
\[
\qf_\Sigma(\gf)\twoheadrightarrow\hf.
\]
Moreover, the functor \(\qf_\Sigma\) preserves surjections.
\end{proposition}

\begin{proof}
For \(\gf\in\operatorname{Lie}_{\Rc}\), write
\[
p_\gf:
\gf\twoheadrightarrow\qf_\Sigma(\gf)
=
\gf/I_\Sigma(\gf)
\]
for the quotient map.
Let \(f:\gf\to\hf\) be a morphism of \(\Rc\)-equivariant Lie
algebras. By
Proposition~\ref{proposition:support_projector},
\[
f\bigl(C_{\Sigma^c}(\gf)\bigr)
\subseteq
C_{\Sigma^c}(\hf).
\]
Since \(f\) is a Lie algebra homomorphism, it follows that
\[
f\bigl(I_\Sigma(\gf)\bigr)
\subseteq
I_\Sigma(\hf).
\]
Hence \(f\) induces a unique morphism
\[
\qf_\Sigma(f):
\qf_\Sigma(\gf)\longrightarrow\qf_\Sigma(\hf)
\]
satisfying
\[
\qf_\Sigma(f)\circ p_\gf
=
p_\hf\circ f.
\]
Uniqueness of \(p_\gf\) implies that these
assignments preserve identities and compositions. Thus
\(\gf\mapsto\qf_\Sigma(\gf)\) defines a functor.

Since
\[
C_{\Sigma^c}(\gf)\subseteq I_\Sigma(\gf),
\]
the quotient map \(p_\gf\) factors, on the level of underlying
representations, through \(Q_\Sigma(\gf)\). Therefore
\(\qf_\Sigma(\gf)\) is a quotient of the
\(\Sigma\)-supported representation \(Q_\Sigma(\gf)\), and is itself
\(\Sigma\)-supported.

Now let \(\hf\) be \(\Sigma\)-supported. Every
\(\Rc\)-equivariant Lie algebra homomorphism
\[
f:\gf\longrightarrow\hf
\]
vanishes on \(C_{\Sigma^c}(\gf)\). Since \(\ker(f)\) is a Lie ideal, it
contains \(I_\Sigma(\gf)\), so \(f\) factors uniquely through
\(p_\gf\). Thus composition with \(p_\gf\) induces a natural
bijection
\[
\operatorname{Hom}_{\operatorname{Lie}_{\Rc,\Sigma}}
\bigl(\qf_\Sigma(\gf),\hf\bigr)
\xrightarrow{\sim}
\operatorname{Hom}_{\operatorname{Lie}_{\Rc}}
(\gf,\hf).
\]
This proves the adjunction. In particular, if
\[
\gf\twoheadrightarrow\hf
\]
is a \(\Sigma\)-supported quotient, its unique factorization through
\(\qf_\Sigma(\gf)\) is again surjective, proving the stated
maximality.

Finally, if \(f:\gf\twoheadrightarrow\hf\) is any surjection, then
\[
\qf_\Sigma(f)\circ p_\gf
=
p_\hf\circ f
\]
is surjective. Hence \(\qf_\Sigma(f)\) is surjective, so
\(\qf_\Sigma\) preserves surjections.
\end{proof}

We next pass to pronilpotent Lie algebras. By a
\emph{pronilpotent \(\Rc\)-Lie algebra} we mean an inverse limit
\[
\gf=\varprojlim_i\gf_i
\]
of finite-dimensional nilpotent \(\Rc\)-equivariant Lie algebras with
surjective transition maps. We denote their category by
\(\operatorname{PNLie}_{\Rc}\), with morphisms given by morphisms of
pro-objects; equivalently, these are continuous (for the inverse-limit topology) \(\Rc\)-equivariant
Lie algebra homomorphisms.

\begin{definition}
For \(\gf=\varprojlim_i\gf_i\in\operatorname{PNLie}_{\Rc}\), set
\[
\qf^{\mathrm{pro}}_\Sigma(\gf)
\coloneqq 
\varprojlim_i\qf_\Sigma(\gf_i).
\]
Where there is no confusion, we omit the superscript \(\mathrm{pro}\).
\end{definition}

\begin{proposition}
\label{proposition:maximal_supported_pronilpotent_quotient}
The preceding definition is independent of the chosen presentation of
\(\gf\). The functor
\[
\qf_\Sigma:
\operatorname{PNLie}_{\Rc}
\longrightarrow
\operatorname{PNLie}_{\Rc}
\]
is left adjoint to the inclusion of the full subcategory of
pronilpotent Lie algebras whose underlying pro-representations are
\(\Sigma\)-supported.

The construction preserves surjections and commutes with nilpotent
truncations:
\[
\qf_\Sigma(\gf)_n
\xrightarrow{\sim}
\qf_\Sigma
\bigl(\gf_n).
\]
\end{proposition}

\begin{proof}
Let \(\hf=\varprojlim_j\hf_j\) be \(\Sigma\)-supported. Using the
usual formula for morphisms in a pro-category and
Proposition~\ref{proposition:maximal_supported_lie_quotient}, we
obtain
\[
\begin{aligned}
\operatorname{Hom}_{\operatorname{PNLie}_{\Rc,\Sigma}}
\bigl(\qf_\Sigma(\gf),\hf\bigr)
&=
\varprojlim_j\varinjlim_i
\operatorname{Hom}_{\operatorname{Lie}_{\Rc}}
\bigl(\qf_\Sigma(\gf_i),\hf_j\bigr)
\\
&\simeq
\varprojlim_j\varinjlim_i
\operatorname{Hom}_{\operatorname{Lie}_{\Rc}}
(\gf_i,\hf_j)
\\
&=
\operatorname{Hom}_{\operatorname{PNLie}_{\Rc}}
(\gf,\hf).
\end{aligned}
\]
This proves the adjunction and hence independence of the chosen
presentation.

Let \(f:\gf\twoheadrightarrow\hf\) be a surjection. After replacing
the chosen presentations by cofinal ones, \(f\) is represented by
levelwise surjections
\[
f_i:\gf_i\twoheadrightarrow\hf_i.
\]
Proposition~\ref{proposition:maximal_supported_lie_quotient} implies
that each
\[
\qf_\Sigma(f_i):
\qf_\Sigma(\gf_i)\twoheadrightarrow\qf_\Sigma(\hf_i)
\]
is surjective. Hence \(\qf_\Sigma(f)\) is surjective.

Finally, the two sides of the asserted truncation isomorphism
represent the maximal \(\Sigma\)-supported quotient of \(\gf\) having
nilpotence class at most \(n\), so the universal property identifies
them.
\end{proof}

The $\Rc$-action on the Lie algebra $\mathfrak{g}$ of a unipotent extension depends on the choice of a Levi section, while the action on the DCS graded pieces do not. The following elementary observation thus shows the support condition on $\mathfrak{g}$ is intrinsic:

\begin{lemma}
\label{lemma:support_graded_equivalence}
Let \(\gf\) be a finite-dimensional nilpotent \(\Rc\)-Lie algebra.
Then the following are equivalent:
\begin{enumerate}
\item the underlying \(\Rc\)-representation of \(\gf\) is
\(\Sigma\)-supported;
\item every
\[
\Gr_n^{\DCS}\gf
\]
is \(\Sigma\)-supported.
\end{enumerate}
The analogous equivalence holds for pronilpotent
\(\Rc\)-Lie algebras.
\end{lemma}

\begin{proof}
The descending central series is an \(\Rc\)-stable finite filtration.
Since \(\operatorname{Rep}_F(\Rc)\) is semisimple, the filtration
splits as a filtration of \(\Rc\)-representations. Thus the simple
constituents of \(\gf\) are precisely those occurring in its
associated graded.

For a pronilpotent Lie algebra, apply the finite-dimensional statement
to its nilpotent quotients.
\end{proof}

We record two consequences used in the main text.

\begin{proposition}
\label{proposition:supported_abelianization_and_graded}
Let \(\gf\in\operatorname{PNLie}_{\Rc}\), and set
\[
\qf\coloneqq \qf_\Sigma(\gf),
\qquad
A_\Sigma\coloneqq Q_\Sigma(\gf^{\ab}).
\]
Then there is a canonical isomorphism
\[
\qf^{\ab}\xrightarrow{\sim}A_\Sigma.
\]

For every \(n\geq1\), the iterated Lie bracket induces a canonical
surjection
\[
Q_\Sigma\bigl(\mathbb F_n(A_\Sigma)\bigr)
\twoheadrightarrow
\Gr_n^{\DCS}\qf,
\]
where \(\mathbb F_n\) denotes the bracket-length \(n\) component of
the free Lie algebra.

In particular, if \(A_\Sigma\) is finite-dimensional, then every
\[
\Gr_n^{\DCS}\qf
\]
is finite-dimensional, and every nilpotent truncation
\[
\qf_{\leq n}
\coloneqq 
\qf/\Fc_{\DCS}^{n+1}\qf
\]
is finite-dimensional.
\end{proposition}

\begin{proof}
Maps from \(\gf\) to abelian \(\Sigma\)-supported Lie algebras factor
first through \(\gf^{\ab}\), then through
\(Q_\Sigma(\gf^{\ab})\). On the other hand, they factor through the
maximal supported quotient \(\qf\) and then through its
abelianization. The two objects therefore satisfy the same universal
property, giving
\[
\qf^{\ab}\simeq Q_\Sigma(\gf^{\ab}).
\]

The associated graded Lie algebra
\[
\Gr_\bullet^{\DCS}\qf
\]
is generated by its degree-one component
\[
\Gr_1^{\DCS}\qf=\qf^{\ab}=A_\Sigma.
\]
Hence the free graded Lie algebra on \(A_\Sigma\) surjects onto
\(\Gr_\bullet^{\DCS}\qf\). In degree \(n\), this gives
\[
\mathbb F_n(A_\Sigma)
\twoheadrightarrow
\Gr_n^{\DCS}\qf.
\]
The target is \(\Sigma\)-supported, so the map factors through
\(Q_\Sigma(\mathbb F_n(A_\Sigma))\).

If \(A_\Sigma\) is finite-dimensional, then each
\(\mathbb F_n(A_\Sigma)\), and hence each graded piece of \(\qf\), is
finite-dimensional. The same is therefore true of every finite
nilpotent truncation.
\end{proof}

\subsection{Supported quotients of pro-unipotent extensions}
\label{appendix:supported_group_quotients}

Over \(F\), the exponential and logarithm functors identify
pro-unipotent groups with pronilpotent Lie algebras. We use this
equivalence without further comment.

Consider an exact sequence of affine pro-algebraic \(F\)-group schemes
\begin{equation}
\label{eq:appendix_prounipotent_extension}
1\longrightarrow\Uc
\longrightarrow\Gc
\xrightarrow{\pi}\Rc
\longrightarrow1,
\end{equation}
where \(\Rc\) is reductive and \(\Uc\) is pro-unipotent. Since \(F\)
has characteristic zero, the Levi decomposition for pro-algebraic
groups gives a noncanonical section
\[
s:\Rc\longrightarrow\Gc
\]
of \(\pi\). Fixing such a section identifies
\[
\Gc\simeq\Uc\rtimes\Rc
\]
and equips \(\uf\coloneqq \Lie(\Uc)\) with the structure of a pronilpotent
\(\Rc\)-Lie algebra. We shall prove below that $\qf_{\Sigma}(\uf)$ is
independent of the chosen Levi section.

\begin{definition}\label{def:supported_quotient_extension}
A quotient extension
\[
\Gc\twoheadrightarrow\Hc\longrightarrow\Rc
\]
is \emph{\(\Sigma\)-supported} if, writing
\[
\vf
\coloneqq 
\Lie\ker(\Hc\to\Rc),
\]
every
\[
\Gr_n^{\DCS}\vf
\]
is \(\Sigma\)-supported.
\end{definition}

\begin{proposition}[Maximal \(\Sigma\)-supported quotient extension]
\label{proposition:maximal_supported_extension_appendix}
There exists a quotient map
\[
q_\Sigma:
\Gc\twoheadrightarrow\Gc_\Sigma
\]
and an exact sequence
\[
1\longrightarrow\Uc_{\Sigma}
\longrightarrow\Gc_{\Sigma}
\xrightarrow{\pi_\Sigma}\Rc
\longrightarrow1
\]
satisfying
\[
\pi_\Sigma\circ q_\Sigma=\pi.
\]
The group \(\Gc_{\Sigma}\) is \(\Sigma\)-supported, and every
\(\Sigma\)-supported quotient extension
\[
\Gc\twoheadrightarrow\Hc\longrightarrow\Rc
\]
factors uniquely through \(q_\Sigma\):
\[
\Gc
\xrightarrow{q_\Sigma}
\Gc_{\Sigma}
\twoheadrightarrow
\Hc.
\]
This property determines \(\Gc_{\Sigma}\) up to a unique isomorphism
as a quotient extension of \(\Gc\).

For any choice of Levi section, endow
\[
\uf\coloneqq \Lie(\Uc)
\]
with the induced \(\Rc\)-action. Then
\[
\Lie(\Uc_{\Sigma})
\simeq
\qf_{\Sigma}(\uf).
\]

If a group \(\Gamma\) acts on
\eqref{eq:appendix_prounipotent_extension} by automorphisms and the
induced action preserves \(\Sigma\), then the quotient map
\[
q_\Sigma:
\Gc\twoheadrightarrow\Gc_{\Sigma}
\]
is canonically \(\Gamma\)-equivariant.
\end{proposition}

\begin{proof}
Choose a Levi section and identify
\[
\Gc\simeq\Uc\rtimes\Rc.
\]
The induced \(\Rc\)-action makes \(\uf=\Lie(\Uc)\) a pronilpotent
\(\Rc\)-Lie algebra. Set
\[
\uf_{\Sigma}\coloneqq \qf_{\Sigma}(\uf),
\qquad
\Uc_{\Sigma}\coloneqq \exp(\uf_{\Sigma}),
\qquad
\Gc_{\Sigma}\coloneqq \Uc_{\Sigma}\rtimes\Rc.
\]
The quotient
\(\uf\twoheadrightarrow\uf_{\Sigma}\)
induces the quotient map
\[
q_\Sigma:
\Gc\twoheadrightarrow\Gc_{\Sigma}.
\]
Since \(\uf_{\Sigma}\) is \(\Sigma\)-supported, so is
\(\Gc_{\Sigma}\).

Let
\[
\Gc
\xrightarrow{q}
\Hc
\longrightarrow
\Rc
\]
be another \(\Sigma\)-supported quotient extension, and write
\[
\Vc\coloneqq \ker(\Hc\to\Rc),
\qquad
\vf\coloneqq \Lie(\Vc).
\]
The chosen Levi section composes with \(q\), making
\(\vf\) an \(\Rc\)-equivariant quotient of \(\uf\). By assumption,
each
\[
\Gr_n^{\DCS}\vf
\]
is \(\Sigma\)-supported; hence
Lemma~\ref{lemma:support_graded_equivalence} implies that \(\vf\)
itself is \(\Sigma\)-supported as a pro-representation.
Proposition~\ref{proposition:maximal_supported_pronilpotent_quotient}
therefore gives a unique surjection
\[
\uf_{\Sigma}\twoheadrightarrow\vf
\]
compatible with \(\uf\twoheadrightarrow\vf\). Exponentiating and
adjoining \(\Rc\) yields the required factorization
\[
\Gc_{\Sigma}\twoheadrightarrow\Hc.
\]
Its uniqueness follows because
\(\Gc\twoheadrightarrow\Gc_{\Sigma}\) is faithfully flat.

The support condition used above is intrinsic: the conjugation action
of \(\Hc\) on
\[
\Gr_n^{\DCS}\vf
\]
factors canonically through \(\Hc/\Vc\simeq\Rc\), since \(\Vc\) acts
trivially on its DCS graded pieces.
Constructions arising from two different sections therefore satisfy
the same universal property and are uniquely isomorphic as quotient
extensions of \(\Gc\).

Finally, for \(\gamma\in\Gamma\), let
\[
\gamma_\Gc:\Gc\xrightarrow{\sim}\Gc,
\qquad
\gamma_\Rc:\Rc\xrightarrow{\sim}\Rc
\]
denote the induced automorphisms. Let
\({}^{\gamma}\Gc_\Sigma\) denote the group scheme \(\Gc_\Sigma\),
viewed as an extension of \(\Rc\) through the projection
\[
\gamma_\Rc^{-1}\circ\pi_\Sigma:
{}^{\gamma}\Gc_\Sigma\longrightarrow\Rc.
\]
Since \(\gamma_\Rc\) preserves \(\Sigma\), this extension is again
\(\Sigma\)-supported, and
\[
q_\Sigma\circ\gamma_\Gc:
\Gc\twoheadrightarrow{}^{\gamma}\Gc_\Sigma
\]
is a \(\Sigma\)-supported quotient over \(\Rc\). The universal
property therefore gives a unique morphism
\[
\theta_\gamma:
\Gc_\Sigma\longrightarrow{}^{\gamma}\Gc_\Sigma
\]
such that
\[
\theta_\gamma\circ q_\Sigma
=
q_\Sigma\circ\gamma_\Gc.
\]

On the underlying group schemes, \(\theta_\gamma\) defines a morphism
\[
\gamma_\Sigma:\Gc_\Sigma\longrightarrow\Gc_\Sigma
\]
satisfying
\[
\gamma_\Sigma\circ q_\Sigma
=
q_\Sigma\circ\gamma_\Gc,
\qquad
\pi_\Sigma\circ\gamma_\Sigma
=
\gamma_\Rc\circ\pi_\Sigma.
\]
Applying the same construction to \(\gamma^{-1}\) and using uniqueness
shows that \(\gamma_\Sigma\) is an automorphism. The same uniqueness
also gives
\[
(\gamma_1\gamma_2)_\Sigma
=
\gamma_{1,\Sigma}\circ\gamma_{2,\Sigma}.
\]
Thus the automorphisms \(\gamma_\Sigma\) define a compatible
\(\Gamma\)-action on \(\Gc_\Sigma\), and
\(q_\Sigma:\Gc\to\Gc_\Sigma\) is \(\Gamma\)-equivariant.
\end{proof}

\begin{corollary}
\label{corollary:finite_supported_truncations}
With notation as above, set
\[
A_\Sigma\coloneqq Q_\Sigma(\uf^{\ab}).
\]
If \(A_\Sigma\) is finite-dimensional, then
\[
\Gc_{\Sigma,n}
\coloneqq 
\Gc_{\Sigma}/
\Fc_{\DCS}^{n+1}\Uc_{\Sigma}
\]
is an affine algebraic \(F\)-group for every \(n\geq1\).
\end{corollary}

\begin{proof}
By
Proposition~\ref{proposition:supported_abelianization_and_graded},
the Lie algebra of
\[
\Uc_{\Sigma}/\Fc_{\DCS}^{n+1}\Uc_{\Sigma}
\]
is finite-dimensional. The corresponding unipotent group is therefore
algebraic. Since \(\Rc\) is algebraic, the same is true of its
extension \(\Gc_{\Sigma,n}\).
\end{proof}

\begin{remark}
\label{remark:adjoint_support_is_canonical}
Let
\[
\Rc_{\mathrm{ss}}
\coloneqq 
(\Rc^\circ)^{\mathrm{der}},
\]
so that
\[
\mathfrak a=\Lie(\Rc_{\mathrm{ss}}).
\]
Since \(\Rc_{\mathrm{ss}}\) is characteristic in \(\Rc\), every
automorphism of \(\Rc\) preserves \(\Sigma_{\mathrm{ad}}\). Thus the
maximal adjoint-supported quotient is functorial under every
automorphism of the original extension.
\end{remark}

We now establish Proposition~\ref{proposition:maximal_adjoint_supported_quotient}.
\begin{proof}[Proof of
Proposition~\ref{proposition:maximal_adjoint_supported_quotient}]
Apply
Proposition~\ref{proposition:maximal_supported_extension_appendix}
with
\[
\Sigma=\Sigma_{\mathrm{ad}}.
\]
The assertion concerning the abelianization and the quotients of the
free Lie pieces follows from
Proposition~\ref{proposition:supported_abelianization_and_graded}.
If the adjoint-supported abelianization is finite-dimensional,
Corollary~\ref{corollary:finite_supported_truncations} gives the
finite-type truncations. Functoriality follows from
Remark~\ref{remark:adjoint_support_is_canonical} and Proposition~\ref{proposition:maximal_supported_extension_appendix}.
\end{proof}

\subsection{Free Lie algebras and perfect targets}
\label{appendix:supported_free_lie}

Let \(A\in\operatorname{Rep}_F(\Rc)\) be finite-dimensional. Write
\[
\mathbb F(A)
=
\bigoplus_{n\geq1}\mathbb F_n(A)
\]
for the free Lie algebra on \(A\), graded by bracket length, and let
\[
\widehat{\mathbb L}(A)
=
\varprojlim_N
\mathbb F(A)/
\bigoplus_{n>N}\mathbb F_n(A)
\]
be its pronilpotent completion.

Set
\[
D_{\Sigma^c}(A)
\coloneqq 
\bigoplus_{n\geq1}
C_{\Sigma^c}\bigl(\mathbb F_n(A)\bigr)
\subseteq
\mathbb F(A),
\]
and let
\[
J_{\Sigma^c}(A)
\coloneqq 
\bigl\langle D_{\Sigma^c}(A)\bigr\rangle_{\mathrm{Lie}}
\]
be the Lie ideal generated by it. Since \(D_{\Sigma^c}(A)\) is graded,
\(J_{\Sigma^c}(A)\) is a graded ideal. Write
\[
J_{\Sigma^c,n}(A)
\coloneqq 
J_{\Sigma^c}(A)\cap\mathbb F_n(A).
\]

The following proposition bridges the discrete and pronilpotent free Lie algebras by showing that $\mathfrak{q}_{\Sigma}(\mathbb F(A))$ and $\mathfrak{q}_{\Sigma}(\widehat{\mathbb L}(A))$ have the same graded pieces:

\begin{proposition}
\label{proposition:supported_free_lie_description}
There is a canonical isomorphism of pronilpotent
\(\Rc\)-Lie algebras
\[
\qf_\Sigma\bigl(\widehat{\mathbb L}(A)\bigr)
\simeq
\prod_{n\geq1}
\mathbb F_n(A)/J_{\Sigma^c,n}(A),
\]
where the right-hand side carries the completed graded Lie bracket.
It is topologically generated by its degree-one component, and its
DCS filtration agrees with the completed
bracket-degree filtration:
\[
\Fc_{\DCS}^m
\qf_\Sigma\bigl(\widehat{\mathbb L}(A)\bigr)
\simeq
\prod_{n\geq m}
\mathbb F_n(A)/J_{\Sigma^c,n}(A).
\]
Consequently,
\[
\Gr_n^{\DCS}
\qf_\Sigma\bigl(\widehat{\mathbb L}(A)\bigr)
\simeq
\mathbb F_n(A)/J_{\Sigma^c,n}(A).
\]
In particular, this graded piece is a quotient of
\[
Q_\Sigma\bigl(\mathbb F_n(A)\bigr).
\]
\end{proposition}

\begin{proof}
For each \(N\), let
\[
\mathbb F_{\leq N}(A)
\coloneqq 
\mathbb F(A)\Big/
\bigoplus_{n>N}\mathbb F_n(A)
\simeq
\bigoplus_{n=1}^N\mathbb F_n(A)
\]
be the \(N\)-step nilpotent truncation of the free Lie algebra. Under
the displayed vector-space identification, the bracket is the
truncated bracket: a bracket of homogeneous elements is zero whenever
its total degree exceeds \(N\).

The maximal \(\Sigma\)-supported quotient of
\(\mathbb F_{\leq N}(A)\) is obtained by quotienting by the image of
\(J_{\Sigma^c}(A)\). Indeed,
\[
C_{\Sigma^c}\bigl(\mathbb F_{\leq N}(A)\bigr)
=
\bigoplus_{n=1}^N
C_{\Sigma^c}\bigl(\mathbb F_n(A)\bigr),
\]
and the Lie ideal generated by this subrepresentation is precisely the
image of \(J_{\Sigma^c}(A)\).

Thus
\[
\qf_\Sigma\bigl(\mathbb F_{\leq N}(A)\bigr)
\simeq
\bigoplus_{n=1}^N
\mathbb F_n(A)/J_{\Sigma^c,n}(A).
\]

Since $\mathbb F_{\leq N}(A) \cong \widehat{\mathbb L}(A)_N$, we find By Proposition~\ref{proposition:maximal_supported_pronilpotent_quotient} that
\[
\mathfrak{q}_{\Sigma}(\widehat{\mathbb L}(A))_N = \bigoplus_{n=1}^N
\mathbb F_n(A)/J_{\Sigma^c,n}(A).
\]

This graded quotient is generated by its degree-one component, since
it is a quotient of \(\mathbb F_{\leq N}(A)\). Its
DCS filtration therefore agrees with its
bracket-degree filtration:
\[
\Fc_{\DCS}^m
\qf_\Sigma\bigl(\mathbb F_{\leq N}(A)\bigr)
=
\bigoplus_{n=m}^N
\mathbb F_n(A)/J_{\Sigma^c,n}(A).
\]

Since $\mathfrak{q}_{\Sigma}(\widehat{\mathbb L}(A))$ is pronilpotent, we have
\[
\mathfrak{q}_{\Sigma}(\widehat{\mathbb L}(A)) = \varprojlim_N \mathfrak{q}_{\Sigma}(\widehat{\mathbb L}(A))_N
=
\varprojlim_N \bigoplus_{n=1}^N
\mathbb F_n(A)/J_{\Sigma^c,n}(A)
=
\prod_{n \ge 1}
\mathbb F_n(A)/J_{\Sigma^c,n}(A),
\]
as desired.
\end{proof}

The next proposition isolates the special role of $\Sigma=\Sigma_{\mathrm{ad}}$ in \S\ref{sec:param_families} (cf. Remark~\ref{rem:why_ad_supp}). In that application,
\[
\mathfrak a
=
\Lie\bigl((\Rc^\circ)^{\mathrm{der}}\bigr)
\]
is a nonzero perfect \(\Rc\)-Lie algebra. Although
\(\mathfrak a\) has trivial pronilpotent completion, an
\(\Rc\)-equivariant surjection
\[
A\twoheadrightarrow\mathfrak a
\]
forces every bracket-length component of the discrete free Lie
algebra on \(A\) to map onto \(\mathfrak a\). By
Proposition~\ref{proposition:supported_free_lie_description},
\(\mathfrak a\) therefore occurs as an \(\Rc\)-equivariant quotient,
and hence as a direct summand, of every DCS
graded piece of the supported pronilpotent quotient. This is the
nonvanishing mechanism used in the proof of
Theorem~\ref{theorem:special_motivic_exists}.

\begin{proposition}
\label{proposition:perfect_target_all_degrees}
Let \(\mathfrak a\) be a nonzero finite-dimensional
\(\Sigma\)-supported \(\Rc\)-equivariant Lie algebra, and suppose that
\[
p:A\twoheadrightarrow\mathfrak a
\]
is an \(\Rc\)-equivariant surjection. If \(\mathfrak a\) is perfect,
then for every \(n\geq1\) there is an \(\Rc\)-equivariant surjection
\[
\Gr_n^{\DCS}
\qf_\Sigma\bigl(\widehat{\mathbb L}(A)\bigr)
\twoheadrightarrow
\mathfrak a.
\]
In particular,
\(\qf_\Sigma(\widehat{\mathbb L}(A))\) is nonzero in every degree.
\end{proposition}

\begin{proof}
The map \(p\) extends uniquely to an \(\Rc\)-equivariant Lie algebra
homomorphism
\[
\mathbb F(A)\longrightarrow\mathfrak a.
\]
Since \(\mathfrak a\) is \(\Sigma\)-supported, this homomorphism kills
\(D_{\Sigma^c}(A)\), and hence the ideal \(J_{\Sigma^c}(A)\).

The image of \(\mathbb F_n(A)\) is the subspace spanned by the
brackets of length \(n\) in \(\mathfrak a\), namely
\[
\Fc_{\DCS}^n\mathfrak a.
\]
Perfectness gives
\[
\Fc_{\DCS}^n\mathfrak a=\mathfrak a
\qquad(n\geq1).
\]
Therefore the induced map
\[
\mathbb F_n(A)/J_{\Sigma^c,n}(A)
\longrightarrow\mathfrak a
\]
is surjective. Apply
Proposition~\ref{proposition:supported_free_lie_description}.
\end{proof}

The following elementary observation is useful when passing from a
free presentation to a geometrically defined Lie algebra.

\begin{lemma}
\label{lemma:unsupported_relations_vanish}
Let \(B\) be an \(\Rc\)-pro-representation, set
\[
A\coloneqq Q_\Sigma(B),
\]
and consider the canonical composite
\[
\widehat{\mathbb L}(B)
\longrightarrow
\widehat{\mathbb L}(A)
\longrightarrow
\qf_\Sigma\bigl(\widehat{\mathbb L}(A)\bigr).
\]
Let \(Z\subseteq\widehat{\mathbb L}(B)\) be an
\(\Rc\)-subrepresentation such that
\[
Q_\Sigma(Z)=0.
\]
Then the composite kills \(Z\), and hence kills the closed Lie ideal
generated by \(Z\).

Consequently, if
\[
\gf
=
\widehat{\mathbb L}(B)/\overline{\langle Z\rangle}_{\mathrm{Lie}},
\]
then the preceding composite factors through a surjection
\[
\gf\twoheadrightarrow
\qf_\Sigma\bigl(\widehat{\mathbb L}(A)\bigr).
\]
\end{lemma}

\begin{proof}
The target is \(\Sigma\)-supported. Every \(\Rc\)-equivariant map from
\(Z\) to a \(\Sigma\)-supported pro-representation is therefore zero.
The kernel of a Lie algebra homomorphism is a closed Lie ideal, so it
contains the closed ideal generated by \(Z\).
\end{proof}

\printbibliography

\end{document}